\documentclass[twoside,11pt,reqno]{amsart}
\usepackage{amsmath,amssymb,amscd,mathrsfs}

\makeatletter

\newif\ifbook@
\newif\ifappendix@
\book@false\appendix@false

\ifbook@
\else
\fi
\@settopoint\oddsidemargin
\@settopoint\marginparwidth
\@settopoint\evensidemargin
\@settopoint\topmargin
\@addtoreset{equation}{section}

\def\thesubsection{\S\thesection-\alph{subsection}}
\def\Point#1{\addtocounter{subsection}{1}\vspace{2mm}
\noindent\thesubsection. {\bf #1.}\def\@currentlabel{\thesubsection}}

\newtheorem{Proposition}{Proposition}[section]
\newtheorem{Lemma}[Proposition]{Lemma}
\newtheorem{Theorem}[Proposition]{Theorem}
\newtheorem{Corollary}[Proposition]{Corollary}

\newtheorem{Remark}[Proposition]{Remark}

\newtheorem{Example}[Proposition]{Example}

\newbox\squ  
\setbox\squ=\hbox{\vrule width.6pt
   \vbox{\hrule height.6pt width.4em\kern1ex
         \hrule height.6pt}%
   \vrule width.6pt}
\def\sqbox{\copy\squ}
\def\proof{\noindent{\sl Proof.\quad}}
\def\endproof{\quad\sqbox\par\vspace{2mm}}

\def\rep{\operatorname{Rep}}
\def\proj{\operatorname{Proj}}
\def\mod#1{#1\!\operatorname{-mod}}
\def\bimod#1#2{#1\!\operatorname{-mod-}\!#2}
\def\id{\operatorname{id}}
\def\Id{\operatorname{Id}}
\def\C{{\mathbb C}}
\def\Q{{\mathbb Q}}
\def\Z{{\mathbb Z}}
\def\Mtype{\mathtt{M}}
\def\Qtype{\mathtt{Q}}
\def\0{{\bar 0}}
\def\1{{\bar 1}}
\def\pr{{\operatorname{pr}}}

\def\cl{{\operatorname{cl}}}
\def\fn{{\operatorname{fin}}}
\def\Sym{S}
\def\P{{\mathcal P}}
\def\H{{\mathcal H}}
\def\R{{\mathcal R}}
\def\B{{\mathcal B}}
\def\A{{\mathcal A}}
\def\cH{\H^\cl}
\def\fH{\H^\fn}
\def\hom{{\operatorname{Hom}}}
\def\ext{{\operatorname{Ext}}}
\def\End{{\operatorname{End}}}
\def\ind{{\operatorname{ind}}}
\def\Irr{{\operatorname{Irr}\:}}
\def\res{{\operatorname{res}}}
\def\Res{{\operatorname{res}\:}}
\def\im{{\operatorname{im}}}
\def\soc{{\operatorname{soc}\:}}
\def\cosoc{{\operatorname{cosoc}\:}}
\def\ch{{\operatorname{ch}\:}}
\def\wt{{\operatorname{wt}}}
\def\inf{{\operatorname{infl}}}
\def\cont{{\operatorname{cont}}}
\def\underbar{\mathpalette\@underbar}
\def\@underbar#1#2{\settowidth{\@tempdimb}{$#1#2$}\@tempdimb=0.8\@tempdimb
                   \ooalign{$#1#2$\crcr%
                         \hfil\rule[-.5mm]{\@tempdimb}{.4pt}\hfil}}
\def\bi{{\underbar{i}}}
\def\bj{{\underbar{j}}}
\def\bh{{\underbar{h}}}
\def\bil{\hbox{\boldmath{$1_\la$}}}
\def\bid{\hbox{\boldmath{$1$}}}
\def\eps{{\varepsilon}}
\def\phi{{\varphi}}
\def\emptyset{{\varnothing}}
\def\ga{{\gamma}}

\def\la{{\lambda}}

\def\si{{\sigma}}

\def\De{{\Delta}}
\def\al{{\alpha}}

\newdimen\hoogte    \hoogte=12pt    
\newdimen\breedte   \breedte=14pt  
\newdimen\dikte     \dikte=0.5pt 

\newenvironment{Young}{\begingroup
       \def\vr{\vrule height0.89\hoogte width\dikte depth 0.2\hoogte}
       \def\fbox##1{\vbox{\offinterlineskip
                    \hrule height\dikte
                    \hbox to \breedte{\vr\hfill##1\hfill\vr}
                    \hrule height\dikte}}
       \vbox\bgroup \offinterlineskip \tabskip=-\dikte \lineskip=-\dikte
            \halign\bgroup &\fbox{##\unskip}\unskip  \crcr }
       {\egroup\egroup\endgroup}
\def\diagram#1{\relax\ifmmode\vcenter{\,\begin{Young}#1\end{Young}\,}\else%
              $\vcenter{\,\begin{Young}#1\end{Young}\,}$\fi}

\begin{document}
\title[Hecke-Clifford superalgebras]{\boldmath
Hecke-Clifford superalgebras, crystals of type
$A_{2\ell}^{(2)}$ and modular branching rules for $\widehat{S}_n$}
\author{\sc Jonathan Brundan and Alexander Kleshchev}
\address
{Department of Mathematics\\ University of Oregon\\
Eugene\\ OR~97403, USA}
\email{brundan@darkwing.uoregon.edu, klesh@math.uoregon.edu}

\thanks{
{2000 subject classification: 17B67, 20C08, 20C20, 17B10, 17B37.}\\
{\phantom{sp-} Both authors
partially supported by the NSF (grant nos DMS-9801442 and DMS-9900134).}
}

\begin{abstract}
This paper is concerned with the modular representation theory of the
affine Hecke-Clifford superalgebra, the cyclotomic Hecke-Clifford superalgebras,
and projective representations of the symmetric group.
Our approach exploits crystal graphs of affine Kac-Moody algebras.
\end{abstract}

\maketitle

\section{Introduction}

In \cite{LLT}, Lascoux, Leclerc and Thibon made the startling 
combinatorial observation
that the crystal graph of the basic representation
of the affine Kac-Moody
algebra $\mathfrak g = A_{\ell}^{(1)}$, 
determined explicitly by Misra and Miwa \cite{MM},
coincides with the modular branching graph for the symmetric group $S_n$
in characteristic $p=\ell+1$, as in \cite{K1}.
The same observation applies to the modular branching graph for the associated
complex Iwahori-Hecke algebras 
at a primitive $(\ell+1)$th root of unity, see \cite{JWB:branching}.
In this latter case, Lascoux, Leclerc and Thibon 
conjectured moreover that the coefficients of the
canonical basis of the basic representation
coincide with the decomposition numbers 
of the Iwahori-Hecke algebras.

This conjecture was proved by Ariki \cite{A1}\footnote{
Grojnowski has informed us that the result can also be proved 
following his note 
``Representations of affine Hecke algebras (and affine quantum ${\rm GL}\sb n$) at roots of unity'',
{\em Internat. Math. Res. Notices} {\bf 5} (1994), 215--217.}.
More generally, Ariki established a similar result connecting
the canonical basis of an arbitrary integrable highest weight module
of $\mathfrak g$ to the representation theory
of a corresponding cyclotomic Hecke algebra, as defined in \cite{AK}.
Note that this work is concerned with the
cyclotomic Hecke algebras over the ground field $\C$,
but Ariki and Mathas \cite{A2, AM} were later able to 
extend the classification of the irreducible modules (but not the
result on decomposition numbers) to arbitrary fields.
For further developments related to the LLT conjecture,
see \cite{LT2, VV, Sch}.

Subsequently, Grojnowski and Vazirani \cite{G,G2,GV,Vt, V} 
have developed a new approach to (amongst other things)
the classification of the irreducible modules of the cyclotomic
Hecke algebras. The approach is valid over an arbitrary ground field,
and is entirely independent of the ``Specht module theory'' that plays an
important role in Ariki's work.
Branching rules are built in from the outset (like in \cite{OV}),
resulting in an explanation and generalization of the link between modular
branching rules and crystal graphs.
The methods are purely algebraic, exploiting affine Hecke algebras 
in the spirit of \cite{BZ,Z1} and others.
On the other hand, Ariki's result on decomposition numbers does not
follow, since that ultimately depends on the geometric
work of Kazhdan, Lusztig and Ginzburg.

In this article, we use Grojnowski's methods
to develop a parallel theory in the twisted case:
we replace the affine Hecke algebras with the
affine Hecke-Clifford superalgebras of Jones and Nazarov \cite{JNaz},
and the Kac-Moody algebra 
$A_{\ell}^{(1)}$ with the twisted algebra
${\mathfrak g} = A_{2\ell}^{(2)}$.
In particular, we obtain an algebraic construction purely in terms
of the representation theory of Hecke-Clifford superalgebras
of the plus part $U_\Z^+$ of the enveloping algebra of ${\mathfrak g}$, 
as well as of Kashiwara's 
highest weight crystals $B(\infty)$ and $B(\la)$
for each dominant weight $\la$. These emerge as the modular branching
graphs of the Hecke-Clifford superalgebras.
Note there is at present no analogue of the notion 
of Specht module in our theory,
underlining the importance of Grojnowski's methods.
However, we do not obtain an analogue of Ariki's
result on decomposition numbers.

As we work over an arbitrary ground field,
the results of the article have applications to
the modular representation theory of the double covers $\widehat S_n$
of the symmeric groups, as was predicted originally
by Leclerc and Thibon \cite{LT}, see also \cite{BK}.
In particular, the parametrization of irreducibles,
classification of blocks and
analogues of the modular branching rules of the
symmetric group for the double covers over fields
of odd characteristic
follow from the special case $\lambda = \Lambda_0$ of our main results.
These matters are discussed in the final section
of the paper, $\S$\ref{last}.

Let us now describe the main results in more detail.
Let $\H_n$ denote the affine Hecke-Clifford superalgebra of
\cite{JNaz}, over an algebraically closed field $F$ of characteristic
different from $2$ and at defining parameter a primitive $(2 \ell + 1)$-th 
root of unity $q \in F^\times$.
All results also have analogues in the degenerate case $q = 1$, working
instead with the affine Sergeev superalgebra of Nazarov \cite{Naz}, 
when the field $F$ should be taken to be of characteristic
$(2 \ell + 1)$.

We consider
$$
K(\infty) = \bigoplus_{n \geq 0} K(\rep_I \H_n),
$$
the sum of the Grothendieck groups of integral $\Z_2$-graded representations
of $\H_n$ for all $n$ (see \ref{intrep} for the precise definition).
In a familiar way (cf. \cite{Zel}), $K(\infty)$ has a natural structure as a 
commutative graded Hopf algebra over $\Z$, multiplication being induced by 
induction
and comultiplication being induced by restriction.
Hence, the graded dual $K(\infty)^*$ is a cocommutative graded Hopf algebra
over $\Z$.

Next, let ${\mathfrak g}$ denote the twisted affine Kac-Moody algebra of
type $A_{2\ell}^{(2)}$.
We will adopt standard Lie theoretic notation for the root system of
${\mathfrak g}$, summarized in more detail in the main body of the article.
In particular, $U_\Q = U_\Q^- U_\Q^0 U_\Q^+$ 
denotes the $\Q$-form of the universal enveloping algebra
of ${\mathfrak g}$ generated by the Chevalley generators
$e_i, f_i, h_i\:(i \in I)$.
Also $U_{\Z} = U_{\Z}^- U_{\Z}^0 U_{\Z}^+$ denotes the Kostant $\Z$-form
of $U_\Q$.
The first main theorem (Theorem~\ref{thma}) identifies $K(\infty)^*$
with $U_\Z^+$, viewing the latter as a graded Hopf algebra over $\Z$
via the principal grading:

\vspace{2mm}
\noindent
{\bf Theorem A.} {\em
$K(\infty)^*$ and $U_\Z^+$ 
are isomorphic
as graded Hopf algebras.
}

\vspace{2mm}
We also introduce for the first time for each dominant integral
weight $\lambda \in P_+$ a finite dimensional quotient superalgebra
$\H_n^\la$ of $\H_n$.
These superalgebras 
play the role of the cyclotomic Hecke algebras in Grojnowski's
theory.
In the special case $\la = \Lambda_0$, $\H_n^\la$ is the finite
Hecke-Clifford superalgebra introduced by Olshanski \cite{Ol}.
Consider the sum of the Grothendieck groups
$$
K(\la) = \bigoplus_{n \geq 0} K(\rep \H_n^\la)
$$
of finite dimensional $\Z_2$-graded
$\H_n^\la$-modules for all $n$.
In this case, we can identify the graded dual
$K(\la)^*$ with
$$
K(\la)^* = \bigoplus_{n \geq 0} K(\proj \H_n^\la)
$$
where $K(\proj \H_n^\la)$ denotes the Grothendieck
group of finite dimensional 
 $\Z_2$-graded projectives. 
It turns out that the natural Cartan map $K(\la)^* \rightarrow K(\la)$
is injective (Theorem~\ref{injectivity}), so we can view both
$K(\la)^*$ and $K(\la)$ as lattices in 
$K(\la)_\Q = \Q \otimes_{\Z} K(\la) = \Q \otimes_{\Z} K(\la)^*$.

Each $K(\la)$ has a 
natural structure of right $K(\infty)$-comodule,
with $K(\la)^* \subset K(\la)$ being a subcomodule.
In other words, according to Theorem A, $K(\la)^*$ and $K(\la)$
are left $U_\Z^+$-modules. The action of
the Chevalley generator $e_i$ here is essentially a refinement of the
restriction functors from $\H_n^\la$ to $\H_{n-1}^\la$.
We show moreover, by considering refinements of induction, 
that the action of $U_\Z^+$ extends
to an action of all of $U_\Z$ on both $K(\la)$ and $K(\la)^*$.
Hence, on extending scalars, we have an action of $U_\Q$ on
$K(\la)_\Q$, with $K(\la)^* \subset K(\la)$ being two different integral forms.
The second main theorem (Theorem~\ref{thmb}) is the following:

\vspace{2mm}
\noindent
{\bf Theorem B.} {\em 
For each $\la \in P_+$,  
$K(\la)_\Q$ is the integrable highest weight $U_\Q$-module of highest weight
$\la$, with highest weight vector $[\bil]$ corresponding to the irreducible
$\H_0^\la$-module.
Moreover $K(\la)^* \subset K(\la)$
are integral forms for $K(\la)_\Q$ with
$K(\la)^* = U_\Z^- [\bil]$ and $K(\la)$ being the
dual lattice under the Shapovalov form.}

\vspace{2mm}
Now let $B(\infty)$ denote the set of isomorphism classes of irreducible
integral $\H_n$-modules for all $n \geq 0$.
This has a natural crystal graph structure, the action of the crystal operators
$\tilde e_i$ being defined by considering
the {\em socle} of the restriction of an irreducible $\H_n$-module
to $\H_{n-1}$. 
Similarly, writing
$B(\la)$ for the set of isomorphism classes of 
irreducible $\H_n^\la$-modules for all $n \geq 0$, 
each $B(\la)$ has a natural crystal
structure describing branching rules between the
algebras $\H_n^\la$ and $\H_{n-1}^\la$.
We stress the crystal structures on $B(\infty)$ and each $B(\la)$ are 
defined purely in terms of the representation theory
of the Hecke-Clifford superalgebras.
The next main result
(Theorems~\ref{ident1}
and \ref{ident2})
identifies the crystals:

\vspace{2mm}
\noindent
{\bf Theorem C.}
{\em The crystal $B(\infty)$ is isomorphic to Kashiwara's crystal
associated to the crystal base of $U_\Q^-$.
Moreover, for each $\la \in P_+$, the crystal $B(\la)$ is isomorphic to 
Kashiwara's crystal associated to the integrable highest weight
$U_\Q$-module of highest weight $\la$.
}

\vspace{2mm}
\noindent
{\bf Acknowledgements. }
We would like to express our debt to the beautiful ideas of
Ian Grojnowski in \cite{G}. Many of the proofs here, and certainly the overall
strategy adopted in the article, are exactly as in Grojnowski's work.
We would also like to thank Monica Vazirani for explaining \cite{GV}
to us, in discussions which initiated the present work.

\ifbook@\pagebreak\fi

\section{Affine Hecke-Clifford superalgebras}

\Point{Ground field and parameters}
Let $F$ be an algebraically closed field of characteristic different
from $2$, and choose $q \in F^\times$
such that
\begin{itemize}
\item[] {\em either} $q$ is a primitive odd root of unity (including the
possibly $q = 1$),
\item[] {\em or} $q$ is not a root of unity at all.
\end{itemize}
Let $h$ be the ``quantum characteristic'', i.e. the
smallest positive integer such that
$$
q^{1-h} + q^{3-h} + \dots + q^{h-3}+q^{h-1} = 0
$$
or $\infty$ in case no such integer exists.
We refer to the case $q \neq 1$ as the {\em quantum case}
and $q = 1$ as the {\em degenerate case}.
All the results will apply to either situation, but there
are sufficiently many differences that for the purpose of
exposition we
will work in the quantum case in the main body of the text, 
with a summary of modifications
in the degenerate case given at the end of each section
whenever necessary.
We set
$$
\xi = q - q^{-1}
$$
as a convenient shorthand.

\Point{Modules over superalgebras}\label{theory}
We will use freely the basic notions of superalgebra, referring 
the reader to \cite[ch.I]{Leites}, \cite[ch.3, $\S\S$1--2]{Man} 
and 
\cite[$\S$2]{BK}. 
We will denote the {\em parity} 
of a homogeneous vector $v$ of a vector superspace by $\bar{v} \in \Z_2$.
By 
a {\em superalgebra}, we mean a $\Z_2$-graded associative algebra
over the fixed field $F$.
If $A$ and $B$ are two superalgebras, then $A \otimes B = A \otimes_F B$
is again a superalgebra with multiplication satisfying
$$
(a_1 \otimes b_1)(a_2 \otimes b_2) = (-1)^{\bar b_1\bar a_2}
a_1 a_2 \otimes b_1b_2
$$
for $a_i \in A, b_i \in B$. Note this and other such expressions only make
sense for {\em homogeneous} $b_1,a_2$: the intended meaning for
arbitrary elements is to be
obtained by extending linearly from the homogeneous case.

If $A$ is a superalgebra, an {\em $A$-module} means 
a $\Z_2$-graded left $A$-module.
A morphism $f:M \rightarrow N$ of $A$-modules $M$ and $N$ means a (not 
necessarily homogeneous) linear map such that
$f(am) = (-1)^{\bar f \bar a} a f(m)$ for all $a \in A, m \in M$.
The category of all such $A$-modules is denoted
$\mod{A}$. 
By a {\em submodule} of an $A$-module, we always mean a {\em graded}
submodule unless we explicitly say otherwise.
We have the {\em parity change functor}
\begin{equation}\label{pcf}
\Pi:\mod{A} \rightarrow \mod{A}.
\end{equation}
For an object $M$, $\Pi M$ is the same underlying 
vector space but with the opposite $\Z_2$-grading.
The new action of $a \in A$ on $m \in \Pi M$ is defined
in terms of the old action by $a \cdot m := (-1)^{\bar a} am$.

It will occasionally be necessary to consider bimodules over two
superalgebras $A, B$: an $(A,B)$-bimodule is an $A$-module $M$,
as in the previous paragraph,
which is also a $\Z_2$-graded {right} $B$-module such that
$(am)b = a (mb)$ for all $a \in A, b \in B, m \in M$.
Note a morphism $f:M \rightarrow N$ of $(A,B)$-bimodules means a 
morphism of $A$-modules as in the previous paragraph such that
$f(mb) = f(m)b$ for all $m \in M, b \in B$.
This gives us the category $\bimod{A}{B}$ of all $(A,B)$-bimodules.
Also if $M$ is an $(A,B)$-bimodule, $\Pi M$ denotes the $A$-module defined as
in the previous paragraph, with the right $B$-action on $\Pi M$
being the same as the original action on $M$.

If $M$ is a finite dimensional
irreducible $A$-module, Schur's lemma (e.g. \cite[$\S$2]{BK}) says that
$\End_A(M)$ is either one dimensional,
or two dimensional on basis $\id_M, \theta_M$
where $\theta_M$ is an odd involution of $M$,
unique up to a sign.
In the former case, we call $M$ an {\em irreducible of type $\Mtype$},
in the latter case $M$ is an {\em irreducible of type $\Qtype$}.
We will occasionally 
write $M \simeq N$ (as opposed to the usual $M \cong N$)
to indicate that there is actually an {\em even} isomorphism
between $A$-modules $M$ and $N$. 
For instance, if $M$ is an irreducible of type $\Qtype$, then
$M \simeq \Pi M$, while if $M$ is irreducible of type $\Mtype$
then $M \cong \Pi M$ but $M \not \simeq \Pi M$.

Given another superalgebra $B$ 
and  $A$-, $B$-modules $M$, $N$ respectively, 
$M\otimes N = M \otimes_F N$ has a natural 
structure of $A\otimes B$-module with
$$
(a \otimes b) (m \otimes n) = (-1)^{\bar b \bar m} am \otimes bn.
$$
We will call this the {\em outer} tensor product of $M$ and $N$ and denote it 
$M\boxtimes N$. 
If $M$ and $N$ are finite dimensional irreducibles,
$M \boxtimes N$ need not be irreducible (unlike the purely even case).
Indeed,
if both $M$ and $N$ are of type $\Mtype$, 
then $M \boxtimes N$ is also irreducible of type $\Mtype$,
while if one is of type $\Mtype$ and the other of type $\Qtype$,
then $M \boxtimes N$ is irreducible of type $\Qtype$.
In either of these two cases, it will be convenient to write
$M \circledast N$ in place of $M \boxtimes N$.
But if both $M$ and $N$ are of type $\Qtype$,
then $M \boxtimes N$ is decomposable:
let $\theta_M:M \rightarrow M$ be an odd involution of $M$
as an $A$-module and $\theta_N:N \rightarrow N$ be an odd involution of $N$
as a $B$-module. Then
$$
\theta_M \otimes \theta_N:M\boxtimes N \rightarrow M \boxtimes N,\quad
m \otimes n \mapsto (-1)^{\bar m} \theta_M(m) \otimes \theta_N(n)
$$
is an even $A\otimes B$-automorphism of $M \boxtimes N$ whose square is $-1$.
Therefore $M \boxtimes N$ decomposes as a direct sum of two $A\otimes B$-modules,
namely, the $\pm \sqrt{-1}$-eigenspaces of the linear map $\theta_M \otimes \theta_N$.
The map $\theta_M \otimes \id_N$ then gives an odd isomorphism between
the two summands as $A \otimes B$-modules.
In this case, we pick either summand, and denote it
$M \circledast N$: it is an irreducible $A \otimes B$-module 
of type $\Mtype$.
Thus,
\begin{equation*}
M\boxtimes N\simeq
\begin{cases}
M\circledast N\oplus\Pi (M\circledast N), 
& \text{if $M$ and $N$ are both of type $\Qtype$,}\\
M\circledast N, &\text{otherwise.}
\end{cases}
\end{equation*}
We stress that $M \circledast N$ is in general only well-defined up
to isomorphism.

Finally, we make some remarks about antiautomorphisms and duality.
In this paper, all antiautomorphisms $\tau:A \rightarrow A$
of a superalgebra $A$
will be {\em unsigned}, so satisfy $\tau(ab) = \tau(b)\tau(a)$.
If $M$ is a finite dimensional $A$-module,
then we can use $\tau$ to 
make the dual space $M^* = \hom_F(M, F)$ into an $A$-module
by defining $(af)(m) = f(\tau(a)m)$ for all $a \in A, f \in M^*, m \in M$.
We will denote the resulting module by $M^\tau$.
We often use the fact that there is a natural even isomorphism
\begin{equation}\label{switch}
\hom_A(M, N) \rightarrow \hom_A(N^\tau, M^\tau)
\end{equation}
for all finite dimensional $A$-modules $M, N$.
The isomorphism sends $f \in \hom_A(M, N)$ to the
dual map $f^* \in \hom_A(N^\tau, M^\tau)$
defined by $(f^* \theta)(m) = \theta(f m)$ for all $\theta \in N^\tau$.

\Point{Grothendieck groups}\label{ggp}
For a superalgebra $A$, 
the category $\mod{A}$ is 
a {\em superadditive category}: each $\hom_A(M,N)$ is a 
$\Z_2$-graded abelian group in a way that is
 compatible with composition.
Moreover, the {\em underlying even category}, i.e.
the subcategory consisting of the same objects
but only even morphisms, is an abelian category in the usual sense.
This allows us to make use of 
all the basic notions of homological algebra 
on restricting our attention to even morphisms.
For example by a short exact sequence in $\mod{A}$, we mean
a sequence
\begin{equation}\label{ses}
0 \longrightarrow M_1 \longrightarrow M_2 \longrightarrow M_3 \longrightarrow 0
\end{equation}
that is a short exact sequence in
the underlying even category, so in particular 
all the maps are necessarily {even}.
Note all functors between superadditive categories that we shall consider
send even morphisms to even morphisms, so make sense on restriction to
the underlying even subcategory.

Let us write $\rep A$ for the full subcategory of $\mod{A}$
consisting of all finite dimensional $A$-modules.
We define the {\em Grothendieck group} $K(\rep{A})$ to be
the quotient of the free $\Z$-module with generators given by all
finite dimensional $A$-modules by the submodule generated by
\begin{itemize}
\item[(1)]
$M_1 - M_2 + M_3$ for every short exact sequence of the form
(\ref{ses});
\item[(2)]
$M - \Pi M$ for every $A$-module $M$.
\end{itemize}
We will write $[M]$ for the image of the $A$-module $M$ in
$K(\rep{A})$.
In an entirely similar way, we define the Grothendieck group
$K(\proj A)$, where $\proj A$ denotes the full category of 
$\mod{A}$ consisting of finite
dimensional projectives.

Now suppose that $A$ and $B$ are two superalgebras.
Then the Grothendieck groups $K(\rep A)$ and $K(\rep B)$ are free $\Z$-modules
with canonical
bases corresponding to the isomorphism classes of irreducible modules.
Moreover, there is an isomorphism
\begin{equation}\label{stara}
K(\rep A) \otimes_\Z K(\rep B)\rightarrow K(\rep A \otimes B),
\qquad
[L] \otimes [L'] \mapsto [L \circledast L']
\end{equation}
for irreducible modules $L \in \rep A,L' \in \rep B$.
This simple observation explains the importance
of the operation $\circledast$.

\Point{The superalgebras}\label{supalgs}
We now proceed to define the superalgebras we will be interested in.
First, $\P_n$ denotes the algebra
$F[X_1^{\pm 1}, \dots, X_n^{\pm n}]$ of Laurent polynomials, viewed
as a superalgebra concentrated in degree $\0$.
We record the relations:
\begin{align}\label{prel1}
X_i X_i^{-1} =1,\:\:&X_i^{-1} X_i = 1,\\
\label{prel2}
X_k X_j = X_j X_k,
\:
X_k^{-1} X_j^{-1} = X_j^{-1} X_k^{-1},
\:
&X_kX_j^{-1} = X_j^{-1} X_k,
\:
X_k^{-1} X_j = X_j X_k^{-1}
\end{align}
for all $1 \leq i  \leq n, 1 \leq j < k \leq n$.
(Of course, some of these relations are redundant, but the precise form
of the relations is used in the proof of Theorem~\ref{base}.)

Let $\A_n$ denote the {superalgebra} 
with even generators $X_1^{\pm 1}, \dots, X_n^{\pm 1}$ and odd
generators $C_1,\dots, C_n$, where the $X_i^{\pm 1}$
are subject
to the polynomial relations (\ref{prel1}), (\ref{prel2}),
the $C_i$ are subject to the Clifford superalgebra relations
\begin{align}
\label{crel1}
C_i^2 &= 1,\\\label{crel2}
C_k C_j &= - C_j C_k\\\intertext{for all $1 \leq i  \leq n, 1 \leq j < k \leq n$, and there are 
the mixed relations}
\label{cprel}
C_i X_j = X_j C_i,
\quad
C_i X_j^{-1} = X_j^{-1} C_i,
\quad
&C_i X_i = X_i^{-1}C_i,
\quad
C_iX_i^{-1}  = X_iC_i 
\end{align}
for all $1 \leq i,j \leq n$ with $i \neq j$.
For $\alpha = (\alpha_1,\dots,\alpha_n) \in \Z^n$ and 
$\beta = (\beta_1,\dots,\beta_n) \in \Z_2^n$,
we write $X^\alpha$ and $C^\beta$ for the monomials
$X_1^{\alpha_1} X_2^{\alpha_2}\dots X_n^{\alpha_n}$
and $C_1^{\beta_1} C_2^{\beta_2} \dots C_n^{\beta_n}$, respectively.
Then, it is straightforward to show that
the elements $$
\{X^\alpha C^\beta \:|\: \alpha \in \Z^n, \beta \in \Z_2^n\}
$$
form a basis for $\A_n$.
In particular, 
$\A_n \cong \A_1 \otimes \dots \otimes \A_1 \hbox{ ($n$ times)},$
and $\P_n$ can be identified with the subalgebra of $\A_n$
generated by the $X_i^{\pm 1}$.

The symmetric group will be denoted $\Sym_n$, with basic 
transpositions $s_1,\dots,s_{n-1}$ and corresponding
Bruhat ordering denoted $\leq$.
We define a left action of $\Sym_n$ on $\A_n$ by algebra automorphisms
so that 
\begin{equation}
w \cdot X_i = X_{wi}, \qquad w \cdot C_i = C_{wi}
\end{equation}
for each $w \in \Sym_n, i = 1,\dots,n$.

Let $\cH_n$ denote the usual (classical) Hecke algebra of $\Sym_n$,
viewed as a superalgebra concentrated in degree $\0$. 
This can be defined on generators $T_1,\dots,T_{n-1}$ subject to relations
\begin{align}\label{hrel1}
T_i^2 &= \xi T_i + 1,\\\label{hrel2}
T_i T_j = T_j T_i,\quad
&T_i T_{i+1}T_i = T_{i+1} T_i T_{i+1}
\end{align}
for all admissible $i,j$ with $|i-j|>2$.
Recall $\cH_n$ has a basis denoted $\{T_w\:|\:w\in \Sym_n\}$,
where
$T_w = T_{i_1} \dots T_{i_m}$
if $w = s_{i_1} \dots s_{i_m}$ is any reduced expression for $w$.
We also record
\begin{equation}\label{hrel3}
T_i^{-1} = T_i - \xi.
\end{equation}

Now we are ready to define the main object of study:
the {\em affine Hecke-Clifford superalgebra} $\H_n$.
This was first introduced by Jones and Nazarov \cite{JNaz},
being the $q$-analogue of the affine Sergeev superalgebra
of Nazarov \cite{Naz}.
By definition, $\H_n$
has even generators $T_1,\dots,T_{n-1}$, 
$X_1^{\pm 1}, \dots, X_n^{\pm 1}$ and odd generators $C_1,\dots,C_n$,
subject to the same relations as $\A_n$
(\ref{prel1}), (\ref{prel2}), (\ref{crel1}), (\ref{crel2}), (\ref{cprel})
and as $\cH_n$
(\ref{hrel1}), (\ref{hrel2}), together with the new relations:
\begin{align}
\label{wrel0}
T_i C_j &= C_j T_i,\\
\label{wrel1}
T_i C_i &= C_{i+1} T_i,\\
\label{affrel1}
T_i X_j &= X_j T_i,
\qquad 
T_i X_j^{-1} = X_j^{-1} T_i,\\
\label{affrel2}(T_i + \xi C_i C_{i+1}) X_i T_i &= X_{i+1}\\\intertext{for 
all admissible $i,j$ with $j \neq i,i+1$.
The relation (\ref{wrel1}) is equivalent to}
\label{wrel2}
T_i C_{i+1} &= C_i T_i - \xi (C_i -C_{i+1}),\\\intertext{while 
(\ref{affrel2}) is equivalent to any of the following
four statements:}
\label{E290900_2}
T_iX_{i} & = X_{i+1}T_i-\xi(X_{i+1}+C_iC_{i+1}X_i), \\ 
\label{E290800_1}
T_iX_i^{-1} & = X_{i+1}^{-1}T_i+\xi(X_i^{-1}+X_{i+1}^{-1}C_iC_{i+1}),\\
\label{E290900_3}
T_iX_{i+1} & = X_iT_i+\xi(1-C_iC_{i+1})X_{i+1}, \\
\label{E290800_2}
T_iX_{i+1}^{-1} & = X_{i}^{-1}T_i-\xi X_i^{-1}(1-C_iC_{i+1})\\\intertext{for 
all $i = 1,\dots,n-1$.
Using (\ref{E290900_2}) and induction on $j \geq 1$, 
one shows that}
\label{powertx}
(T_i+\xi C_i C_{i+1}) X_i^j &= X_{i+1}^j (T_i-\xi) - \xi \sum_{k=1}^{j-1} 
\left(
X_i^{j-k} X_{i+1}^k + X_i^{-k} X_{i+1}^{j-k} C_i C_{i+1}\right)
\end{align}
for
all $j \geq 1, 1 \leq i < n$. Hence:
\begin{equation}
\label{powertz}
\begin{split}
(T_i+\xi C_i C_{i+1}) X_i^j C_i = &X_{i+1}^j C_{i+1}(T_i+\xi C_iC_{i+1}) \\
&- \xi \sum_{k=1}^{j-1} 
\left(
X_i^{j-k} X_{i+1}^{k} C_i-X_i^{-k} X_{i+1}^{j-k} C_{i+1}\right),
\end{split}
\end{equation}
\begin{equation}
\begin{split}
\label{powerty}
\quad \quad \quad \:\: (T_i - \xi) X_i^{-j} = 
&X_{i+1}^{-j}(T_i +\xi C_i C_{i+1})\\
&+ \xi \sum_{k=1}^{j-1} \left(
X_i^{-k} X_{i+1}^{k-j} + X_i^{j-k} X_{i+1}^{-k} C_i C_{i+1}\right).
\end{split}
\end{equation}

Finally, let $\fH_n$ denote the subalgebra of $\H_n$ generated
by all $C_i, T_j$ for $i=1,\dots,n,j=1,\dots,n-1$.
Alternatively, as follows easily from Theorem~\ref{base} below, 
$\fH_n$ can be defined as the superalgebra generated by elements
$C_i,T_j$ subject only to the relations (\ref{crel1}), (\ref{crel2}),
(\ref{hrel1}), (\ref{hrel2})
and (\ref{wrel0}), (\ref{wrel1}).
Hence, $\fH_n$ is the (finite) {\em Hecke-Clifford superalgebra}
first introduced by Olshanski \cite{Ol} as the $q$-analogue of the
Sergeev superalgebra of \cite{se}.

\Point{Basis theorem}
Now we proceed to study the algebra $\H_n$ in more detail.
The first goal is to construct a basis.
There are obvious homomorphisms $f:\A_n \rightarrow \H_n$
and $g:
\cH_n \rightarrow \H_n$ under which the $X_i,C_i$ or $T_j$ map
to the same elements of $\H_n$. 
We write $X^{\alpha}C^\beta$ also
for the image under $f$ of the basis element 
$X^{\alpha}C^\beta$ of $\A_n$,
and $T_w$ for the image under $g$ of $T_w \in \cH_n$.
This 
notation will be justified shortly, when we show that $f$ and $g$
are both algebra monomorphisms.
The following lemma is obvious from
the relations:

\vspace{1mm}
\begin{Lemma}\label{240800_0} Let $f \in \A_n$, $x \in \Sym_n$.
Then, in $\H_n$,
$$
T_x f = (x \cdot f) T_x + \sum_{y < x} f_y T_y,
\qquad\quad
fT_x = T_x(x^{-1} \cdot f)  + \sum_{y < x} T_y f_y'
$$
for some $f_y, f_y' \in \A_n$.
\end{Lemma}

\vspace{1mm}

It follows easily that $\H_n$ is at least
{\em spanned} by all $X^{\alpha} C^\beta T_w$ for all $\alpha \in \Z^n, 
\beta \in \Z_2^n$ and $w \in \Sym_n$.
We wish to prove that these elements are linearly independent too:

\vspace{1mm}
\begin{Theorem}\label{base} The
$\{X^{\alpha}C^\beta  T_w\:|\:\alpha \in \Z^n, \beta \in \Z_2^n, 
w \in \Sym_n\}$
form a basis for $\H_n$.
\end{Theorem}

\proof
Consider instead the algebra $\tilde \H_n$ on generators
$\tilde X_i^{\pm 1}, \tilde C_i, \tilde T_j$ for $1\leq i\leq n, 1\leq j <n$
subject to relations (\ref{prel1})--(\ref{cprel}), (\ref{hrel1}), 
(\ref{wrel0})--(\ref{affrel1}) and 
(\ref{wrel2})--(\ref{E290800_2}).
Thus we have all the relations of $\H_n$ {\em except for}
the braid relations (\ref{hrel2}).
Using precisely these relations as the reduction system,
it is a routine if tedious
exercise using Bergman's diamond lemma \cite[1.2]{Berg} to prove 
that $\tilde \H_n$ has a basis given by
all
$\tilde 
X^\alpha \tilde C^\beta \tilde T$ for all $\alpha \in \Z^n, \beta \in \Z_2^n$ and
all words $\tilde T$ in the $\tilde T_j$ which do not
involve a subword of the form $\tilde T_j^2$ for any $j$.
Hence, 
the subalgebra $\widetilde{A}_n$
of $\tilde \H_n$ generated by
the $\tilde X_i^{\pm 1}, \tilde C_i$ is isomorphic to $\A_n$.
Also let $\tilde \cH_n$ 
denote the subalgebra of $\tilde \H_n$ generated by the $\tilde T_j$,
so that $\tilde \cH_n$ is isomorphic to the algebra on generators
$\tilde T_1,\dots, \tilde T_{n-1}$ subject to relations 
$\tilde T_j^2 = \xi \tilde T_j +1$ for each $j$.

Now, by definition, $\H_n$ is the quotient of $\tilde \H_n$ by
the two-sided ideal $\mathcal I$ generated by the elements
$$
a_{i,j} = \tilde T_i \tilde T_j - \tilde T_j \tilde T_i,\quad
b_i = \tilde T_i \tilde T_{i+1} \tilde T_i - \tilde T_{i+1} \tilde T_i
\tilde T_{i+1}
$$
for $i,j$ as in (\ref{hrel2}). Let $\mathcal 
J$ be the two-sided ideal of $\tilde \cH_n$
generated by the same elements $a_{i,j}, b_i$ for all $i,j$.
Then, by the basis theorem for $\cH_n$, $\tilde \cH_n 
/ {\mathcal J} \cong \cH_n$, with basis
given by elements $T_w$ for $w \in \Sym_n$ defined in the usual way.
It follows immediately that to prove the theorem, it suffices to show that
${\mathcal I} ={ \widetilde{A}_n }{\mathcal J}$
in $\widetilde{\H}_n$. In turn, this follows if we can show that 
$rt = t'r$ for each 
$r \in \{a_{i,j}, b_i\}$ and 
each generator $t$ of $\widetilde{A}_n$,
where $t'$ is some other element of $\widetilde{A}_n$.
Now again this is a routine check: 
for example, the most complicated case involves verifying that
$(\tilde T_{i+1} \tilde T_i \tilde T_{i+1} - \tilde T_i \tilde T_{i+1} \tilde T_i) \tilde X_{i+2} = \tilde X_i
(\tilde T_{i+1} \tilde T_i \tilde T_{i+1} - \tilde T_i \tilde T_{i+1} \tilde T_i)$ 
using only the relations in $\tilde \H_n$.
\endproof

So we 
have a right from now on
to {\em identify} $\A_n$  and $\cH_n$
with the corresponding subalgebras of $\H_n$.
The theorem also shows that the superalgebra $\fH_n$ has
basis $\{C^\beta T_w \:|\:\beta\in\Z_2^n, w \in \Sym_n\}$.
So Theorem~\ref{base} can be restated as saying
that $\H_n$ is a free right $\fH_n$-module on basis
$\{X^\alpha \:|\:\alpha \in \Z^n\}$.

As another consequence, it makes sense to consider the
tower of superalgebras
$$
\H_0 \subset \H_1 \subset \H_2 \subset \dots \subset \H_n \subset \dots
$$
where for $i \leq n$,
$\H_i$ is identified with the subalgebra of $\H_{n}$
generated by $C_1,\dots,C_i$,
$X_1^{\pm 1},\dots,X_i^{\pm 1}$,
$T_1,\dots,T_{i-1}$.
Similarly, we can consider $\A_i \subseteq \A_n,
\P_i \subseteq \P_n$, etc....

Finally, we point out that
there are obvious variants of the basis of Theorem~\ref{base}, 
reordering the $C$'s, $X$'s and $T$'s. For instance, 
$\H_n$ also has $\{T_w X^\alpha C^\beta \:|\:w \in \Sym_n, \alpha \in \Z^n,
\beta \in \Z_2^n\}$ as a basis. 
This follows using Lemma~\ref{240800_0}.

\Point{\boldmath The center of $\H_n$}
The next theorem was first established in \cite[Prop. 3.2]{JNaz}
(for the case $F = \C$).

\vspace{1mm}
\begin{Theorem}\label{LCenter}
The (super)center of $\H_n$ consists of all symmetric polynomials
in $X_1+X_1^{-1}, \dots, X_n+X_n^{-1}$.
\end{Theorem}

\begin{proof}
Let $Z$ denote the center of $\H_n$, i.e. the $z$ such that
$zy = yz$ for all $y \in \H_n$. (The argument applies equally well
to the supercenter defined by $zy = (-1)^{\bar z \bar y} yz$.)
One first checks that
symmetric polynomials in $X_1+X_1^{-1},\dots,X_n+X_n^{-1}$ {\em are} central
following 
the argument in the proof of \cite[Prop. 3.2(b)]{JNaz}.

Conversely, take
$z = \sum_{w \in \Sym_n} f_w T_w \in Z$
where each $f_w \in \A_n$. 
Let $w$ be maximal with respect to the Bruhat
order such that $f_w \neq 0$. Assume for a contradiction that $w \neq 1$.
Then, there exists some $i \in\{1,\dots,n\}$ with $wi \neq i$.
Consider $(X_i+X_i^{-1})z - z (X_i + X_i^{-1})$. 
By Lemma~\ref{240800_0}, this looks like
$f_w(X_i+X_i^{-1} - X_{wi}-X_{wi}^{-1})T_w$ 
plus a linear combination of terms of the form
$f_x T_x$ for $f_x \in \A_n$ and $x \in \Sym_n$ with $x \not\geq w$
in the Bruhat order. So in view of Theorem~\ref{base}, 
$z$ is not central, a 
contradiction.

Hence, we must have that $z \in \A_n$.
Considering the form of the center of $\A_n$ one easily shows that
$z$ in fact lies in $F[X_1+X_1^{-1},\dots,X_n+X_n^{-1}]$. 
To see that $z$ is actually a symmetric
polynomial, write
$z = 
\sum_{i,j \geq 0} a_{i,j} (X_1+X_1^{-1})^i (X_2+X_2^{-1})^j$
where the coefficients $a_{i,j}$
lie in $F[X_3+X_3^{-1},\dots,X_n+X_n^{-1}]$.
Applying Lemma~\ref{240800_0} to $T_1 z = z T_1$ now gives that $a_{i,j} = a_{j,i}$
for each $i,j$, hence $z$ is symmetric in $X_1+X_1^{-1}$
and $X_2+X_2^{-1}$.
Similar argument shows that $z$ is symmetric in $X_i+X_i^{-1}$ 
and $X_{i+1}+X_{i+1}^{-1}$
for all $i =1,\dots,n-1$ to complete the proof.
\end{proof}

\Point{Parabolic subalgebras}
Suppose that $\mu = (\mu_1,\dots,\mu_u)$ is a composition of $n$, i.e. 
a sequence of positive integers summing to $n$.
Let $\Sym_\mu \cong \Sym_{\mu_1} \times \dots \times \Sym_{\mu_u}$
denote the corresponding Young subgroup of $\Sym_n$,
and $\cH_\mu \subseteq \cH_n$ denote its Hecke algebra.
So $\cH_\mu \cong \cH_{\mu_1} \otimes \dots\otimes \cH_{\mu_u}$
is the subalgebra of $\cH_n$ generated
by the $T_j$ for which $s_j \in \Sym_\mu$.

We define the parabolic subalgebra $\H_\mu$ of the
affine Hecke-Clifford superalgebra $\H_n$ in a similar way:
it is the subalgebra of $\H_n$ generated by $\A_n$
and all $T_j$ for which $s_j \in \Sym_\mu$.
It follows easily from Theorem~\ref{base} that the elements
$$
\{X^\alpha C^\beta T_w\:|\:\alpha \in \Z^n, \beta \in \Z_2^n, 
w \in \Sym_\mu\}
$$
form a basis for $\H_\mu$. In particular, 
$\H_\mu \cong \H_{\mu_1} \otimes \dots \otimes \H_{\mu_u}.$
Note that the parabolic subalgebra
$\H_{(1,1,\dots,1)}$ is precisely the subalgebra $\A_n$.

We will need the usual induction 
and restriction functors between
$\H_n$ and $\H_\mu$. These will be denoted simply
\begin{equation}\label{indres}
\ind^{n}_\mu: \mod{\H_\mu} \rightarrow \mod{\H_n},
\qquad
\res^{n}_\mu: \mod{\H_n} \rightarrow \mod{\H_\mu},
\end{equation}
the former being the tensor functor
$\H_n \otimes_{\H_\mu} ?$ which is left adjoint to $\res^{n}_\mu$.
More generally, 
we will consider induction and restriction between nested
parabolic subalgebras, with obvious notation.
We will also occasionally consider the restriction functor
\begin{equation}\label{blah}
\res^n_{n-1}:\mod{\H_n} \rightarrow \mod{\H_{n-1}}
\end{equation}
where $\H_{n-1}$ denotes the subalgebra of $\H_n$ generated
by $X_i^{\pm 1}, C_i$ and $T_j$ for $i = 1,\dots,n-1,
j = 1,\dots,n-2$.

\Point{Mackey theorem}\label{SMackey}
Let $\mu,\nu$ be compositions of $n$.
We let $D_\nu$ denote the set of minimal length
left $\Sym_\nu$-coset representatives in $\Sym_n$,
and $D_{\mu}^{-1}$ denote the set of minimal length right 
$\Sym_\mu$-coset representatives.
Then $D_{\mu,\nu} := D_\mu^{-1} \cap D_\nu$ is the set of
minimal length $(\Sym_\mu,\Sym_\nu)$-double coset representatives in $\Sym_n$.
We recall some well-known properties,
see e.g. \cite[\S 1]{DJ}.
First, for $x \in D_{\mu,\nu}$,
$\Sym_\mu \cap x \Sym_\nu x^{-1}$ and
$x^{-1} \Sym_\mu x \cap \Sym_\nu$ are Young subgroups of $\Sym_n$.
So we can define 
compositions $\mu\cap{x}\nu$ and ${{x^{-1}}}\mu\cap\nu$ of $n$ from 
$$
\Sym_\mu\cap x \Sym_\nu x^{-1}=\Sym_{\mu\cap{x}\nu} \quad \text{and} \quad x^{-1}\Sym_\mu x\cap \Sym_\nu=\Sym_{{{x^{-1}}}\mu\cap\nu}. 
$$
Then, for $x \in D_{\mu,\nu}$, every
$w \in \Sym_\mu x \Sym_\nu$ can be written as $w = u x v$ for unique
elements $u \in \Sym_\mu$ and $v \in \Sym_\nu \cap D^{-1}_{x^{-1}\mu\cap\nu}$.
Moreover, when this is done, 
$\ell(w) = \ell(u)+\ell(x)+\ell(v)$.
This fact implies the following well-known lemma which is essentially
equivalent to the Mackey theorem for $\cH_n$ (see e.g. \cite[Theorem 2.7]{DJ}):

\vspace{1mm}
\begin{Lemma}\label{clmackey} For $x \in D_{\mu,\nu}$, the subspace
$\cH_\mu T_x \cH_\nu$ of $\cH_n$ has basis
$\{T_w\:|\:w \in \Sym_\mu x \Sym_\nu\}.$
\end{Lemma}

\vspace{1mm}

This is our starting point for proving a version of the Mackey theorem
for $\H_n$.

\vspace{1mm}
\begin{Lemma}\label{soup}
For $x \in D_{\mu,\nu}$, the subspace
$\H_\mu T_x \cH_\nu$ of $\H_n$ has basis
$\{X^\alpha C^\beta T_w\:|\:
\al\in\Z^n, \beta\in \Z_2^n, w\in \Sym_\mu x \Sym_\nu\}$. 
Moreover,
$$
\H_n = \bigoplus_{x \in D_{\mu,\nu}} \H_\mu T_x \cH_\nu.
$$
\end{Lemma}

\begin{proof}
Since $\H_\mu = \A_n \cH_\mu$, Lemma~\ref{clmackey} implies at once that
the given $\{X^\alpha C^\beta T_w\}$ span
$\H_\mu T_x \cH_\nu$. But they are linearly independent too by
Theorem~\ref{base}, proving the first statement. The second follows immediately
using Theorem~\ref{base} once more.
\end{proof}

Now fix some total order $\prec$ refining the Bruhat order $<$
on $D_{\mu,\nu}$. For $x \in D_{\mu,\nu}$, set
\begin{align}\label{f1}
\B_{\preceq x} &= \bigoplus_{y \in D_{\mu,\nu},\ y \preceq x} \H_\mu T_y \cH_\nu,\\
\label{f2}
\B_{\prec x} &= \bigoplus_{y \in D_{\mu,\nu},\ y \prec x} \H_\mu T_y \cH_\nu,\\
\B_{x} &= \B_{\preceq x} / \B_{\prec x}.
\label{f3}
\end{align}
It follows immediately from Lemma~\ref{240800_0} that $\B_{\preceq x}$ (resp. 
$\B_{\prec x}$) is invariant under right multiplication by $\A_n$.
Hence, since $\H_\nu = \cH_\nu \A_n$, we have defined a filtration of
$\H_n$ as an $(\H_\mu,\H_\nu)$-bimodule. 
We want to describe the quotients $\B_x$ more explicitly
as $(\H_\mu,\H_\nu)$-bimodules. 

\vspace{1mm}
\begin{Lemma}\label{240800_4}
For each $x \in D_{\mu,\nu}$, 
there exists an algebra isomorphism
$$
\phi=\phi_{x^{-1}}:\H_{\mu\cap x\nu} \rightarrow \H_{{x^{-1}}\mu\cap\nu}
$$
with $\phi(T_w)=T_{x^{-1}wx}$,  $\phi(X_i)=X_{x^{-1}i}$, and 
$\phi(C_i)=C_{x^{-1}i}$ for $w\in \Sym_{\mu\cap x\nu}$, $1\leq i\leq n$. 
\end{Lemma}
\begin{proof}
The isomorphism $\psi:\Sym_{\mu\cap x\nu} 
\rightarrow \Sym_{{x^{-1}}\mu\cap\nu}, u\mapsto x^{-1}ux$ is length 
preserving. Equivalently,
$x^{-1}(i+1)=(x^{-1}i)+1$ for each $i$ with $s_i \in \Sym_{\mu\cap x\nu}$.
Using this, it is straightforward to check that the 
map $\phi$ defined as above respects the 
defining relations on generators. 
\end{proof}

Let $N$ be a left $\H_{{x^{-1}}\mu\cap\nu}$-module. 
By twisting the action with the isomorphism 
$\phi_{x^{-1}}:\H_{\mu\cap x\nu} \rightarrow \H_{{x^{-1}}\mu\cap\nu}$ from 
Lemma~\ref{240800_4}, we get a left $\H_{\mu\cap x\nu}$-module, 
which will be denoted ${^x N}$. 
Now we can identify the module $\B_x$ introduced above.

\vspace{1mm}
\begin{Lemma}\label{240800_5}
View 
$\H_\mu$ as an $(\H_\mu,\H_{\mu\cap x\nu})$-bimodule and 
$\H_\nu$ as an $(\H_{{x^{-1}}\mu\cap\nu},\H_\nu)$-bimodule in the natural ways. 
Then, ${^x}\H_\nu$ is an $(\H_{\mu\cap x\nu},\H_\nu)$-bimodule and
$$
\B_x\simeq \H_\mu\otimes_{\H_{\mu\cap x\nu}} {^x  \H_\nu}
$$
as an $(\H_\mu,\H_\nu)$-bimodule.
\end{Lemma}
\begin{proof}
We define a bilinear map $\H_\mu\times \H_\nu\rightarrow \B_x=\B_{\preceq x}/ 
\B_{\prec x}$ by 
$(u,v)\mapsto u T_x v +\B_{\prec x}$. 
For $y \in \Sym_{\mu} \cap x \Sym_\nu x^{-1}$,
$$
T_y T_x = T_{yx} = T_{xx^{-1}yx} = T_{x} T_{x^{-1} yx}
$$
which is all that is required to check that the map 
is $\H_{\mu\cap x\nu}$-balanced. Hence there is an induced
$(\H_\mu,\H_\nu)$-bimodule map
 $\Phi:\H_\mu\otimes_{\H_{\mu\cap x\nu}} {^x \H_\nu}
\rightarrow \B_x$. Finally, to prove that $\Phi$ is bijective, note that
$$
\{X^\al C^\beta T_u\otimes T_v  \:|\: \al\in \Z^n, \beta\in\Z_2^n,
 u\in \Sym_\mu, v \in \Sym_\nu \cap D^{-1}_{x^{-1}\mu\cap\nu}\}
$$
is a basis of the induced module  $\H_\mu\otimes_{\H_{\mu\cap x\nu}} {^x 
\H_\mu}$ as a vector space.
In view of Lemma~\ref{soup}, the image of these elements under $\Phi$ is a
basis of $\B_x$.
\end{proof}

Now we can prove the Mackey theorem.

\vspace{1mm}
\begin{Theorem} {\rm (``Mackey Theorem'')}
\label{TMackey}
Let $M$ be an $\H_\nu$-module. Then
$\res^{n}_{\mu}\ind^{n}_{\nu} M$
admits a filtration with subquotients $\:\simeq\:$ to
$\ind_{{\mu\cap x\nu}}^{\mu}{^x}
(\res_{{{x^{-1}}\mu\cap\nu}}^{\nu}M),$
one for each $x\in D_{\mu,\nu}$. 
Moreover, the subquotients can be taken in any order refining the 
Bruhat order on $D_{\mu,\nu}$,
in particular
$\ind_{{\mu\cap \nu}}^{\mu}
\res_{\mu\cap\nu}^{\nu}M$ appears as a submodule.
\end{Theorem}

\begin{proof} 
This follows from Lemma~\ref{240800_5} and the isomorphism 
$$
(\H_\mu\otimes_{\H_{\mu\cap x\nu}} {^x} \H_\nu)\otimes_{\H_\nu} M\simeq 
\ind_{\H_{\mu\cap x\nu}}^{\H_\mu}{^x}
(\res_{\H_{{x^{-1}}\mu\cap\nu}}^{\H_\nu}M),
$$
which is easy to check.
\end{proof} 

\Point{Some (anti)automorphisms}
A check of relations shows that $\H_n$ 
posesses an automorphism $\si$ and an 
antiautomorphism $\tau$ defined on the generators as follows:
\begin{align}\label{Esi}
\si&: T_i\mapsto -T_{n-i}+\xi,\quad \:\:\:\:\:\,C_j\mapsto C_{n+1-j},\quad
X_j\mapsto X_{n+1-j};\\
\label{Etau}
\tau&: T_i\mapsto T_i+\xi C_iC_{i+1},\quad 
C_j\mapsto C_j,\quad \quad \:\:\:\:\,X_j\mapsto X_j,
\end{align}
for all $i = 1,\dots,n-1, j = 1,\dots,n$.

If $M$ is a finite dimensional $\H_n$-module, 
we can use $\tau$ to make the dual space $M^*$ into an $\H_n$-module
denoted $M^\tau$, see
\ref{theory}.
Note $\tau$ leaves invariant every parabolic subalgebra of $\H_n$,
so also induces a duality on finite dimensional $\H_\mu$-modules
for each composition $\mu$ of $n$.

Instead, given any $\H_n$-module $M$, we can twist the action 
with $\si$ to get a new module denoted $M^\si$. 
More generally, for any composition $\nu = (\nu_1,\dots,\nu_u)$ 
of $n$ we denote by $\nu^*$ the 
composition with the same non-zero parts but taken in the opposite order. 
For example $(3,2,1)^*=(1,2,3)$. Then $\si$ induces an isomorphism of 
parabolic subalgebras $\H_{\nu^*}\rightarrow \H_{\nu}$. 
So if $M$ is an $\H_{\nu}$-module, we can inflate through 
$\si$ to get an $\H_{\nu^*}$-module denoted $M^\si$. 
If $M=M_1\boxtimes\dots\boxtimes 
M_u$ is an outer tensor 
product module over $\H_{\nu}$ then $M^\si\cong M_u^\si\boxtimes\dots\boxtimes 
M_1^\si$. The same holds if each $M_i$ is irreducible and
$\boxtimes$ is replaced with $\circledast$. 
These observations imply:

\vspace{1mm}
\begin{Lemma}
\label{L031000}
Let $M\in\mod{\H_m}$ and $N\in\mod{\H_n}$. Then 
\begin{equation*}
(\ind_{m,n}^{m+n} M\boxtimes N)^\si\cong
\ind_{n,m}^{m+n} N^\si\boxtimes M^\si.
\end{equation*}
Moreover, if $M$ and $N$ are irreducible, the same holds for $\circledast$ in place of $\boxtimes$. 
\end{Lemma}

\Point{Duality}
Thoughout this subsection, let $\mu$ be a composition of $n$
and set $\nu = \mu^*$.
Let $d \in D_{\mu,\nu}$ be the longest double coset representative.
Note that $\mu \cap d \nu = \mu$ and $d^{-1} \mu \cap \nu = \nu$,
so $\Sym_\mu d \Sym_\nu = \Sym_\mu d = d \Sym_\nu$.
There is an isomorphism
\begin{equation}\label{meme}
\phi = \phi_{d^{-1}}:\H_\mu \rightarrow \H_\nu,
\end{equation}
see Lemma~\ref{240800_4}.
As in \ref{SMackey}, for an $\H_\nu$-module $M$,
${^d M}$ denotes the $\H_\mu$-module obtained by pulling back the
action through $\phi$.
We begin by considering the classical situation,
adopting the obvious analogous 
notation for modules over $\cH_\mu,\cH_\nu$.

\vspace{1mm}
\begin{Lemma}\label{joe}
Define a linear map
$\theta^{\cl}:\cH_n \rightarrow {^d \cH_\nu}$
by
$$
\displaystyle
\theta^{\cl}(T_w) = \left\{
\begin{array}{ll}
T_{d^{-1}w}&\hbox{if $w \in d \Sym_\nu$,}\\
0&\hbox{otherwise,}
\end{array}\right.
$$
for each $w \in \Sym_n$.
Then:
\begin{enumerate}
\item[(i)] $\theta^{\cl}$ is an even homomorphism of 
$(\cH_\mu,\cH_\nu)$-bimodules;
\item[(ii)]  $\ker \theta^{\cl}$ contains no non-zero left ideals of $\cH_n$;
\item[(iii)] the map $$
f^{\cl}:
\cH_n \rightarrow \hom_{\cH_\mu}(\cH_n, {^d\cH_\nu}),\qquad
h \mapsto h \theta^{\cl}
$$ 
is an even isomorphism of $(\cH_n, \cH_\nu)$-bimodules.
\end{enumerate}
\end{Lemma}

\begin{proof}
(i) Since $d^{-1} \mu \cap \nu = \nu$,
$\cH_\mu T_d \cH_\nu = T_d \cH_\nu$ is isomorphic
as an $(\cH_\mu,\cH_\nu)$-bimodule to $^d\cH_\nu$, 
the isomorphism being simply the map
$T_w \mapsto T_{d^{-1}w}$ for $w \in d \Sym_\nu$, 
compare Lemma~\ref{240800_5}. Now (i) follows because $\theta^{\cl}$ is 
simply this isomorphism
composed with the projection from $\cH_n$
to $\cH_\mu T_d \cH_\nu$ along the bimodule decomposition 
$\cH_n = \bigoplus_{x \in D_{\mu,\nu}} \cH_\mu T_x \cH_\nu$.

(ii) We show by downward induction on $\ell(x)$ that
$\theta^{\cl}(\cH_n t) \neq 0$ whenever we are given
$x \in D_\nu$ and
$$
t = \sum_{y \in D_\nu\:\mathrm{with}\:
\ell(y) \leq \ell(x)} T_y h_y
$$
with each $h_y \in \cH_\nu$ and $h_x \neq 0$.
Since $d$ is the 
longest element of $D_\nu$, the induction starts with
$x = d$: in this case, the conclusion 
is clear as $\theta^{\cl}(t) = h_x \neq 0$.
So now suppose $x < d$ 
and that the claim has been proved for all higher $x \in D_\nu$.
Pick a basic transposition $s$ such that
$sx > x$ and $sx \in D_\nu$.
Then,
$$
T_s t = \sum_{y \in D_\nu\:\mathrm{with}\:\ell(y) \leq \ell(sx), } T_y h_y'
$$
for $h_y' \in \cH_\nu$ with $h_{sx}' = h_x \neq 0$.
But now the induction hypothesis shows that
$\theta^{\cl}(\cH_n t) = \theta^{\cl}(\cH_n T_s t) \neq 0$.

(iii) We remind the reader that $h \theta^{\cl}:\cH_n \rightarrow {^d \cH_\mu}$
denotes the map with $(h\theta^{\cl})(t) = 
(-1)^{\bar h \bar t}\theta^{\cl}(th)$ (the sign being $+$ always 
in this case).
Given this and (i), it is straightforward to check that $f^{\cl}$
is a homomorphism of $(\cH_n,\cH_\nu)$-bimodules.
To see that it is an isomorphism, it suffices by dimension to show that
it is injective. Suppose $h$ lies in the kernel.
Then, $(f^{\cl}(h))(t) = \theta^{\cl}(th) = 0$
for all $t \in \cH_n$. Hence $h = 0$ by (ii).
\end{proof}

Now we extend this result to $\H_n$, recalling the
definition of $\phi$ from (\ref{meme}).

\vspace{1mm}
\begin{Lemma}\label{joe3}
Define a linear map $\theta:\H_n \rightarrow {^d \H_\nu}$
by
$$
\displaystyle
\theta(f T_w) = \left\{
\begin{array}{ll}
\phi(f) T_{d^{-1}w}&\hbox{if $w \in d \Sym_\nu$,}\\
0&\hbox{otherwise,}
\end{array}\right.
$$
for each $f \in \A_n, w \in \Sym_n$.
Then:
\begin{enumerate}
\item[(i)] $\theta$ is an even homomorphism of $(\H_\mu,\H_\nu)$-bimodules;
\item[(ii)] the map $$
f:\H_n \rightarrow \hom_{\H_\mu}(\H_n,{^d \H_\nu}),\qquad
h \mapsto h \theta
$$ is an even isomorphism of $(\H_n,\H_\nu)$-bimodules.
\end{enumerate}
\end{Lemma}

\begin{proof}
(i) According to a special case of Lemma~\ref{240800_5},
the top factor $\B_d$ in the bimodule filtration of $\H_n$
defined in (\ref{f3})
is isomorphic to $^d \H_\nu$ as an $(\H_\mu,\H_\nu)$-bimodule.
The map $\theta$ is simply the composite of this isomorphism with the
quotient map $\H_n \rightarrow \B_d$.

(ii)
Recall that $\{T_w\:|\:w \in D_\mu^{-1}\}$
forms a basis for $\H_n$ as a free left $\H_\mu$-module,
and $^d \H_\nu$ is isomorphic to $\H_\mu$ as a left $\H_\mu$-module.
It follows that the maps
$$
\{\psi_w\:|\:w \in D_\mu^{-1}\}
$$
form a basis for $\hom_{\H_\mu}(\H_n,{^d \H_\nu})$ as a free right
$\H_\nu$-module,
where $\psi_w:\H_n \rightarrow {^d \H_\nu}$ is the unique left
$\H_\mu$-module homomorphism with 
$\psi_w(T_u) = \delta_{w,u}.1$ for all $u \in D_{\mu}^{-1}$.

The analogous maps $\psi_w^\cl \in \hom_{\cH_\mu}(\cH_n,{^d\cH_\nu})$
defined by $\psi_w^\cl(T_u) = \delta_{w,u}.1$ for $u \in D_\mu^{-1}$
form a basis for $\hom_{\cH_\mu}(\cH_n,{^d\cH_\nu})$
as a free right $\cH_\nu$-module.
So in view of Lemma~\ref{joe}(iii), we can find
a basis $\{a_w\:|\:w \in D_\mu^{-1}\}$
for $\cH_n$ viewed as a right $\cH_\nu$-module such that
$f^\cl(a_w) = \psi_w^\cl$ for each $w \in D_\mu^{-1}$,
i.e.
$$
\theta^{\cl}(T_u a_w) = \left\{
\begin{array}{ll}
1&\hbox{if $u = w$,}\\
0&\hbox{otherwise}
\end{array}\right.
$$
for every $u \in D_\mu^{-1}$.
But $\H_\nu = \cH_\nu \A_n$, so the elements $\{a_w\:|\:w \in D_\mu^{-1}\}$
also form a basis for $\H_n$ as a right $\H_\nu$-module, and
$f(a_w) = \psi_w$ since $\theta = \theta^\cl$ on $\cH_n$.
Thus $f$  maps a basis of $\H_n$
to a basis
of $\hom_{\H_\mu}(\H_n,{^d \H_\nu})$  (as free right $\H_\nu$-modules), 
hence $f$
is an isomorphism of $(\H_n,\H_\nu)$-bimodules.
\end{proof}

\begin{Corollary}\label{c1}
There is a natural isomorphism
$\hom_{\H_\mu}(\H_n,^d\!M)
\simeq
\H_n \otimes_{\H_\nu} M$ of $\H_n$-modules,
for every left $\H_\nu$-module $M$.
\end{Corollary}

\begin{proof}
Let $f:\H_n \rightarrow \hom_{\H_\mu}(\H_n,{^d}\H_\nu)$ be the bimodule
isomorphism constructed in Lemma~\ref{joe3}.
Then, there are natural isomorphisms
\ifbook@
\begin{align*}
\H_n \otimes_{\H_\nu} M \stackrel{f \otimes \id}{\longrightarrow} 
\hom_{\H_\mu}(\H_n,{^d} \H_\nu) \otimes_{\H_\nu} M
&\simeq \hom_{\H_{\mu}}(\H_n, {^d} \H_\nu \otimes_{\H_\nu} M)\\
&\simeq \hom_{\H_n}(\H_{n}, {^d} M),
\end{align*}
\else
$$
\H_n \otimes_{\H_\nu} M \stackrel{f \otimes \id}{\longrightarrow} 
\hom_{\H_\mu}(\H_n,{^d} \H_\nu) \otimes_{\H_\nu} M
\simeq \hom_{\H_{\mu}}(\H_n, {^d} \H_\nu \otimes_{\H_\nu} M)
\simeq \hom_{\H_n}(\H_{n}, {^d} M),
$$
\fi
the second isomorphism depending on the fact that $\H_n$ is a free
left $\H_\mu$-module,
see e.g. \cite[20.10]{AF}.
\end{proof}

Recall the duality $(\ref{Etau})$ 
on finite dimensional $\H_n$- (resp. $\H_\nu$-) modules.

\vspace{1mm}
\begin{Corollary}\label{c2}
There is a natural isomorphism $\ind_\nu^n (M^\tau) \simeq
(\ind_\mu^n ({^d}M))^\tau$
for every finite dimensional $\H_\nu$-module $M$.
\end{Corollary}

\begin{proof}
One routinely checks using (\ref{switch}) that the functor
$\tau \circ \ind_\mu^n \circ \tau$ (from the category of finite dimensional
$\H_\mu$-modules to finite dimensional $\H_n$-modules) is right adjoint to
$\res^n_\mu$. Hence,
it is isomorphic to $\hom_{\H_\mu}(\H_n, ?)$ by uniqueness of adjoint functors.
Now combine this natural isomorphism with the previous corollary
(with $\mu$ and $\nu$ swapped and $d$ replaced by $d^{-1}$).
\end{proof}

We finally record a special case, which is the analogue of
\cite[Proposition~5.8]{G}:

\vspace{1mm}
\begin{Theorem}\label{duals}
Given a finite dimensional $\H_m$-module $M$ and a
finite dimensional $\H_n$-module $N$,
$$
(\ind_{m,n}^{m+n} M \boxtimes N)^\tau
\cong \ind_{n,m}^{m+n} (N^\tau \boxtimes M^\tau).
$$
Moreover, if $M$ and $N$ are irreducible, the same is true
with $\boxtimes$ replaced by $\circledast$.
\end{Theorem}

\Point{Modifications in the degenerate case}\label{degmod}
Now we summarize the necessary changes 
in the degenerate case $q = 1$. So now $F$ is an algebraically
closed field of odd characteristic $h$ and $q = 1$.
The superalgebras $\P_n, \A_n, \H_n, \H_\mu, \fH_n$
need to be replaced with their degenerate analogues (but we keep the same
symbols).

First of all, $\P_n$ becomes the ordinary polynomial
algebra $F[x_1,\dots,x_n]$ concentrated in degree $\0$. Then, $\A_n$ is 
replaced by the algebra on even generators
$x_1,\dots,x_n$ and odd generators $c_1,\dots,c_n$, where the $x_i$ satisfy
the relations of a polynomial ring and the $c_i$ satisfy the same Clifford
relations as before (\ref{crel1}), (\ref{crel2}).
The relations (\ref{cprel}) become instead:
\begin{equation}\label{abcd}
c_i x_i = - x_i c_i,
\quad
c_i x_j = x_j c_i
\end{equation}
for $1 \leq i, j \leq n$ with $i \neq j$.
The algebra $\cH_n$ is replaced by the group algebra $F \Sym_n$ of the symmetric
group, always writing simply $w$ in place of $T_w$ before.

Now the affine Hecke-Clifford superalgebra becomes
the {\em affine Sergeev superalgebra} of \cite{Naz}. This is
the superalgebra $\H_n$ defined on
even generators $x_1,\dots,x_n, s_1,\dots,s_{n-1}$ and odd generators
$c_1,\dots,c_n$. The relations between the $x$'s and $c$'s are as in the
new $\A_n$, the relations between the $s_i$ are the usual relations of the
symmetric group, i.e. (\ref{hrel1}), (\ref{hrel2}) with $q$ there equal to $1$,
and there are new relations replacing (\ref{wrel0})--(\ref{affrel2}):
\begin{align}
s_i c_i = c_{i+1} s_i,
\qquad
&s_{i} c_{i+1} = c_{i} s_{i},
\quad
&s_i c_j = c_j s_i,\label{abcd1}\\
s_ix_{i}=x_{i+1}s_i-1-c_ic_{i+1},
\qquad &s_ix_{i+1}=x_is_i+1-c_ic_{i+1},
\quad &s_ix_j=x_js_i
\end{align}
for all admissible $i$, $j$ with $j\neq i,i+1$.
We record the useful formula, being the analogue of (\ref{powertx}):
\begin{equation}
\label{powertx2}
s_i x_i^j = x_{i+1}^j s_i - \sum_{k=0}^{j-1} (
x_i^{j-k-1}x_{i+1}^k + x_i^{j-k-1}(-x_{i+1})^k c_ic_{i+1})
\end{equation}
for all $j \geq 1, 1 \leq i < n$.
Finally, $\fH_n$ is replaced by the subalgebra of $\H_n$ generated by
the $c_i, s_j$: it is alternatively the twisted tensor product
of a Clifford algebra with the group algebra of the symmetric group, i.e.
the original {\em Sergeev superalgebra} considered in \cite{se},
also see \cite[$\S$3]{BK}.

The automorphism $\sigma:\H_n \rightarrow \H_n$ and
antiautomorphism $\tau:\H_n \rightarrow \H_n$ are now defined by:
\begin{align}
\si&: s_i\mapsto -s_{n-i},\quad c_j\mapsto c_{n+1-j},\quad
x_j\mapsto x_{n+1-j};\\
\tau&: s_i\mapsto s_i,\quad \quad\quad\!
c_j\mapsto c_j,\quad \quad \quad x_j\mapsto x_j,\label{tttdef}
\end{align}
for all $i = 1,\dots,n-1, j = 1,\dots,n$.

The basis theorem, proved in an analogous way to Theorem~\ref{base}, now says
that $\H_n$ has a basis given by $\{x^\alpha c^\beta w\:|\:
\alpha \in \Z_{\geq 0}^n,
\beta \in \Z_2^n, w \in \Sym_n\}$.
The center of $\H_n$, proved as in Theorem~\ref{LCenter}, 
is now the set of all symmetric
polynomials in $x_1^2, x_2^2, \dots, x_n^2$ (see also \cite[Prop. 3.1]{Naz}).
Parabolic subalgebras $\H_\mu$
of $\H_n$ are defined in the same way as before.
The Mackey Theorem and Theorem~\ref{duals} are proved in an entirely similar way.

\ifbook@\pagebreak\fi

\section{Cyclotomic Hecke-Clifford superalgebras}

\Point{Cyclotomic Hecke-Clifford superalgebras}
Keep all the notation from the previous section.
Suppose now that $f \in F[X_1] \subset \H_n$ 
is a polynomial of the form 
\begin{equation}\label{myform}
a_d X_1^d + a_{d-1} X_1^{d-1} + \dots + a_1 X_1 + a_0
\end{equation}
for coefficients $a_i \in F$ with $a_d = 1$
and $a_i = a_0 a_{d-i}$
for each $i = 0,1,\dots,d$.
This assumption implies that 
\begin{equation}\label{specia}
C_1 f = a_0 X_1^{-d} f C_1.
\end{equation}
Define ${\mathcal I}_f$ to be the two-sided ideal of $\H_n$ generated
by $f$ and let
\begin{equation}
\H_n^f := \H_n / {\mathcal I}_f.
\end{equation}
We call $\H_n^f$ the {\em cyclotomic Hecke-Clifford superalgebra}
corresponding to $f$. It is the analogue in our setting of the
cyclotomic Hecke algebra of \cite{AK}.

\Point{Basis theorem}
The goal in this subsection is to  
describe an explicit basis for $\H_n^f$, analogous to the 
Ariki-Koike basis \cite{AK} for cyclotomic Hecke algebras.
We introduce some further notation for the proof.
Set $f_1 = f$ and for $i = 2,\dots,n$, define inductively
$f_{i} = (T_{i-1} + \xi C_{i-1} C_{i}) f_{i-1} T_{i-1}$.
The first lemma follows easily by induction from 
(\ref{powertx}):

\vspace{1mm}
\begin{Lemma}\label{fi} For $i = 1,\dots,n$,
$$
f_i = X_i^d + (\hbox{terms 
lying in $ \P_{i-1}X_i^e \fH_i$ for $0 < e < d$}) + u_i
$$
where $u_i \in \fH_i$ is a unit.
\end{Lemma}

\vspace{1mm}

Given $Z = \{z_1 < \dots < z_u\} \subseteq \{1,\dots,n\}$, let
$
f_Z = f_{z_1} f_{z_2} \dots f_{z_u} \in \H_n.
$
Define
\begin{align}
\Pi_n &= \{(\alpha, Z)\:|\:Z \subseteq \{1,\dots,n\}, \alpha\in \Z^n
\hbox{ with } 0 \leq \alpha_i < d\hbox{ whenever $i \notin Z$}\},\\
\Pi_n^+ &= \{(\alpha,Z) \in \Pi_n\:|\:Z \neq \emptyset\}.
\end{align}

\vspace{1mm}
\begin{Lemma}\label{hardb}
$\H_n$ is a free right $\fH_n$-module on basis
$
\{X^\alpha f_Z\:|\:(\alpha,Z) \in \Pi_n\}.
$
\end{Lemma}

\begin{proof}
Define a total order $\prec$ on $\Z$ so that
$$
\textstyle\lfloor{\frac{d}{2}}\rfloor
 \prec \lfloor{\frac{d}{2}}\rfloor
-1 \prec \lfloor{\frac{d}{2}}\rfloor
+1 \prec \lfloor{\frac{d}{2}}\rfloor
-2 \prec \lfloor{\frac{d}{2}}\rfloor
+2 \prec \dots.
$$
We have a corresponding
reverse lexicographic ordering on $\Z^n$:
$\alpha \prec \alpha'$
if and only if
$
\alpha_n = \alpha_n', \dots, \alpha_{k+1} = \alpha_{k+1}',
\alpha_k \prec \alpha_k'
$
for some $k = 1,\dots,n$. It is important that $\Z^n$ has a smallest
element with respect to this total order.
Define a function
$\gamma:\Pi_n \rightarrow \Z^n$
by $\gamma(\alpha, Z) := (\gamma_1,\dots,\gamma_n)$ where
$$
\gamma_i = \left\{
\begin{array}{ll}
\alpha_i&\hbox{if $i \notin Z$ or $\alpha_i < 0$,}\\
\alpha_i + d&\hbox{if $i \in Z$ and $\alpha_i \geq 0$.}
\end{array}
\right.
$$
We claim that for $(\alpha,Z) \in \Pi_n$,
$$
X^\alpha f_Z = X^{\gamma(\alpha,Z)} u
+(\hbox{terms lying in
$X^{\beta} \fH_n$ for $\beta \prec \gamma(\alpha,Z)$}),
$$
where $u$ is some unit in $\fH_n$.
In other words, $X^\alpha f_Z = X^{\gamma(\alpha,Z)}u + $(lower terms).
Since $\gamma:\Pi_n \rightarrow \Z^n$ is a bijection and
we already know that 
the $\{X^\alpha\:|\:\alpha \in \Z^n\}$ form a basis for $\H_n$
viewed as a right $\fH_n$-module by Theorem~\ref{base}, the claim immediately
implies the lemma.

To prove the claim,
proceed by induction on $n$.
If $n = 1$, the statement is quite obvious.
Now assume $n > 1$ and the statement has been proved for $(n-1)$.
We need to consider $X^\alpha f_Z$ for $(\alpha, Z) \in \Pi_n$.
If $n \notin Z$, then the conclusion follows without difficulty from
the induction hypothesis, so assume $n \in Z$.
Let $\alpha' = (\alpha_1,\dots,\alpha_{n-1})$
and $Z' = Z - \{n\}$. Then by induction
$$
X^{\alpha'} f_{Z'} = X^{\gamma(\alpha',Z')} u'
+(\hbox{terms lying in
$X^{\beta'} \fH_{n-1}$ for $\beta' \prec \gamma(\alpha',Z')$}),
$$
where $u'$ is some unit in $\fH_{n-1}$.
Multiplying on the left by $X_n^{\alpha_n}$ and on the right by
$f_n$, using Lemma~\ref{fi}, one deduces that
\ifbook@
$$
X^\alpha f_Z = X_n^{\alpha_n+d} X^{\gamma(\alpha',Z')}u'+
X_n^{\alpha_n} X^{\gamma(\alpha',Z')} u'u_n 
+(\hbox{terms in $X^\beta \fH_n$ for $\beta \prec \gamma(\alpha,Z)$})
$$
\else
$$
X^\alpha f_Z = X_n^{\alpha_n+d} X^{\gamma(\alpha',Z')}u'+
X_n^{\alpha_n} X^{\gamma(\alpha',Z')} u'u_n 
+(\hbox{terms lying in $X^\beta \fH_n$ for $\beta \prec \gamma(\alpha,Z)$})
$$
\fi
where both $u'$ and $u'u_n$ are units in $\fH_n$.
Noting that
$$
X^{\gamma(\alpha,Z)} = 
\left\{
\begin{array}{ll}
X_n^{\alpha_n+d} X^{\gamma(\alpha',Z')}&\hbox{if $\alpha_n \geq 0$,}\\
X_n^{\alpha_n} X^{\gamma(\alpha',Z')}&\hbox{if $\alpha_n < 0$,}\\
\end{array}\right.
$$
this is of the desired form.
\end{proof}

\begin{Lemma}\label{little}
{\rm (i)} For $1 \leq i \leq n$,
$$
\displaystyle C_i f_i \fH_n \subseteq
X_i^{-d} f_i \fH_n + \sum_{k=1}^{i-1} \P_i f_k \fH_n.
$$
\begin{enumerate}
\item[(ii)] For $1 \leq i \leq j \leq n$,
$$\displaystyle \P_j (T_{j-1} - \xi) \dots (T_{i+1}-\xi)(T_i - \xi) X_i^{-d} f_i \fH_n
\subseteq
\sum_{k=1}^j \P_j f_k \fH_n.
$$
\item[(iii)] For $n > 1$, $\fH_{n-1} f_{n} \fH_n =  f_n \fH_n$.
\end{enumerate}
\end{Lemma}

\begin{proof}
(i) Proceed by induction on $i$, the case $i = 1$ being immediate from
(\ref{specia}). For $i > 1$, 
\begin{align*}
C_i f_i \fH_n &= C_i (T_{i-1}+\xi C_{i-1}C_i) f_{i-1} \fH_n
= (T_{i-1} - \xi) C_{i-1} f_{i-1} \fH_n\\
&\subseteq (T_{i-1} - \xi) X_{i-1}^{-d} f_{i-1} \fH_n
+ \sum_{k=1}^{i-2} (T_{i-1}-\xi) \P_{i-1} f_k \fH_n\\
&\subseteq X_i^{-d} (T_{i-1}+\xi C_{i-1}C_i) f_{i-1} \fH_n
+ \sum_{k=1}^{i-1} \P_i f_k \fH_n
= X_i^{-d} f_i \fH_n+ \sum_{k=1}^{i-1} \P_i f_k \fH_n,
\end{align*}
applying 
Lemma~\ref{240800_0}, the relations in $\H_n$ especially
(\ref{powerty}), and the induction hypothesis. 

(ii) Proceed by induction on $(j-i)$, the conclusion being immediate
in case $j = i$.
If $j > i$, an application of (\ref{powerty}), combined with (i)
to commute $C_i$ past $f_i$, gives that
\ifbook@
\begin{multline*}
\P_j (T_{j-1} - \xi) \dots (T_i - \xi) X_i^{-d} f_i \fH_n
\subseteq
\P_j (T_{j-1} - \xi) \dots (T_{i+1}-\xi) X_{i+1}^{-d}f_{i+1} \fH_n
+\\ 
\sum_{k=1}^j \P_j f_k \fH_n.\qquad
\end{multline*}
\else
$$
\P_j (T_{j-1} - \xi) \dots (T_i - \xi) X_i^{-d} f_i \fH_n
\subseteq
\P_j (T_{j-1} - \xi) \dots (T_{i+1}-\xi) X_{i+1}^{-d}f_{i+1} \fH_n
+ \sum_{k=1}^j \P_j f_k \fH_n.
$$
\fi
Now apply the induction hypothesis.

(iii) 
By considering the antiautomorphism $\tau$ of $\fH_{n-1}$,
one sees that $\fH_{n-1}$ is generated by the elements
$C_i$ for $1 \leq i \leq n-1$ and 
$(T_j  + \xi C_j C_{j+1})$ for $1 \leq j < n-1$, and the latter
satisfy the braid relations.
We show that each of these generators of $\fH_{n-1}$
leave $f_n \fH_n$ invariant.

First, consider $C_i f_n \fH_n$ for $1 \leq i \leq n-1$.
Expanding the definition of $f_n$ and commuting $C_i$ past the leading
terms, it equals
$$
(T_{n-1} + \xi C_{n-1}C_n) \dots C_i (T_i + \xi C_i C_{i+1}) f_{i} \fH_n.
$$
By the relations,
$C_i (T_i + \xi C_i C_{i+1}) = (T_i+ \xi C_i C_{i+1}) C_{i+1}$.
Now the conclusion in this case follows immediately since $C_{i+1} f_{i} \fH_n = f_{i} \fH_n$.

Next, consider $(T_j + \xi C_j C_{j+1}) f_n \fH_n$ for $1 \leq j < n-1$. 
Expanding and commuting again gives that it equals
$$
(T_{n-1} + \xi C_{n-1}C_n) \dots  (T_j + \xi C_j C_{j+1}) 
(T_{j+1} + \xi C_{j+1} C_{j+2})(T_j + \xi C_j C_{j+1}) 
f_j \fH_n.
$$
Applying the braid relation gives
$$
(T_{n-1} + \xi C_{n-1}C_n) \dots  (T_{j+1} + \xi C_{j+1} C_{j+2}) 
(T_{j} + \xi C_{j} C_{j+1})(T_{j+1} + \xi C_{j+1} C_{j+2}) 
f_j \fH_n
$$
which again equals $f_n \fH_n$ as required.
\end{proof}

\begin{Lemma}\label{before}
$\displaystyle{\mathcal I}_f = \sum_{i=1}^n \P_n f_i \fH_n$.
\end{Lemma}

\begin{proof}
Let $\fH_{2\dots n}$ denote the subalgebra of $\fH_n$
generated by $C_2,\dots,C_n, T_2, \dots, T_{n-1}$, so $\fH_{2\dots n} 
\cong \fH_{n-1}$.
We first claim that
$$
\fH_n \subseteq
\sum_{i = 1}^n (T_{i-1} + \xi C_{i-1}C_i) \dots (T_1 + \xi C_1 C_2)
\fH_{2\dots n}
+
\sum_{i = 1}^n (T_{i-1} - \xi) \dots (T_1 - \xi)
C_1 \fH_{2\dots n}.
$$
To prove this, 
it suffices to show that $T_w C^\beta$ lies in the right hand side
for each $w \in \Sym_n, \beta \in \Z_2^n$.
Proceed by induction on the Bruhat order on $w \in \Sym_n$,
the case $w = 1$ being trivial.
For the induction step,
we can find $1 \leq i \leq n$ and $u \in \Sym_{2\dots n}$ such that
$$
T_w C^\beta = T_{i-1} \dots T_1 T_u C^\beta = 
T_{i-1}\dots T_1 C_1^{\beta_1} T_u C^{\beta'}
$$
where $\beta' = (\0, \beta_2,\dots,\beta_n)$.
Now an argument using Lemma~\ref{240800_0} gives that 
$$
T_w C^\beta =
\left\{
\begin{array}{ll}
(T_{i-1} + \xi C_{i-1}C_i) \dots (T_1 + \xi C_1 C_2) T_u C^{\beta'}
&\hbox{if $\beta_1 = \0$,}\\
(T_{i-1} - \xi) \dots (T_1 - \xi) C_1 T_u C^{\beta'}
&\hbox{if $\beta_1 = \1$,}
\end{array}
\right.
$$
modulo lower terms of the form $T_{w'} C^{\beta''}$ with $w' < w$.
The claim follows.

Now let ${\mathcal J} = \sum_{i=1}^n \P_n f_i \fH_n$.
Clearly ${\mathcal J} \subseteq {\mathcal I}_f$.
So it suffices for the lemma 
to show that ${\mathcal I}_f \subseteq {\mathcal J}$.
Noting that $\H_n = \P_n \fH_n$ and using
the result in the previous paragraph,
we get that
\ifbook@
\begin{align*}
{\mathcal I}_f &= \H_n f_1 \H_n =
\H_n f_1 \P_n \fH_n = 
\H_n f_1 \fH_n = \P_n \fH_n f_1 \fH_n\\
& \subseteq
\sum_{i=1}^n
\P_n (T_{i-1} +\xi C_{i-1}C_i) \dots (T_1 + \xi C_1 C_2) f_1 \fH_n
+\\
&\qquad\qquad\qquad\qquad\qquad\qquad\sum_{i=1}^n
\P_n (T_{i-1} - \xi) \dots (T_1 - \xi) C_1 f_1 \fH_n.
\end{align*}
\else
\begin{align*}
{\mathcal I}_f &= \H_n f_1 \H_n =
\H_n f_1 \P_n \fH_n = 
\H_n f_1 \fH_n = \P_n \fH_n f_1 \fH_n\\
& \subseteq
\!\sum_{i=1}^n\!
\P_n (T_{i-1} + \xi C_{i-1}C_i) \dots (T_1 + \xi C_1 C_2) f_1 \fH_n
\!+
\!\sum_{i=1}^n\!
\P_n (T_{i-1} - \xi) \dots (T_1 - \xi) C_1 f_1 \fH_n.
\end{align*}
\fi
Now, each
$\P_n(T_{i-1} + \xi C_{i-1}C_i) \dots (T_1 + \xi C_1 C_2) f_1 \fH_n
=
\P_nf_i \fH_n \subseteq {\mathcal J}$.
Also each
$\P_n (T_{i-1} - \xi) \dots (T_1 - \xi) C_1 f_1 \fH_n$
is contained in ${\mathcal J}$ thanks to Lemma~\ref{little}(i),(ii).
\end{proof}

\begin{Lemma}\label{ke}
$\displaystyle
{\mathcal I}_f = \sum_{(\alpha,Z) \in \Pi_n^+} X^\alpha f_Z \fH_n.$
\end{Lemma}

\begin{proof}
We proceed by induction on $n$, the case $n = 1$ being almost obvious.
So now suppose that $n > 1$. Clearly we can assume that $d > 0$.
It will be convenient to write ${\mathcal I}_f'$ for the two-sided ideal
of $\H_{n-1}$ generated by $f$, so
\begin{equation}\label{hyp}
{\mathcal I}_f' = \sum_{(\alpha',Z') \in \Pi_{n-1}^+} X^{\alpha'} f_{Z'} 
\fH_{n-1}
\end{equation}
by the induction hypothesis.
Let ${\mathcal J} = \sum_{(\alpha,Z) \in \Pi_n^+} X^\alpha f_Z \fH_n.$
Obviously ${\mathcal J} \subseteq {\mathcal I}_f$.
So in view of Lemma~\ref{before}, it suffices to show that 
$X^\alpha f_i \fH_n \subseteq {\mathcal J}$ for each $\alpha \in \Z^n$
and each $i = 1,\dots,n$.

Consider first $X^\alpha f_n \fH_n$.
Write $X^\alpha = X_n^{\alpha_n} X^\beta$ for $\beta \in \Z^{n-1}$.
Expanding $X^\beta$ in terms of the basis of $\H_{n-1}$ from 
Lemma~\ref{hardb}, we see that
$$
X^\alpha f_n \fH_n
\subseteq \sum_{(\alpha',Z')\in\Pi_{n-1}}
X_n^{\alpha_n}X^{\alpha'} f_{Z'} \fH_{n-1} f_n \fH_n.
$$
This is contained in ${\mathcal J}$ thanks to Lemma~\ref{little}(iii).

Finally, consider $X^\alpha f_i \fH_n$ with $i < n$.
Write $X^\alpha = X_n^{\alpha_n} X^\beta$ for $\beta \in \Z^{n-1}$.
By the induction hypothesis,
$$
X^\alpha f_i \fH_n = X_n^{\alpha_n} X^\beta f_i \fH_n
\subseteq \sum_{(\alpha',Z') \in \Pi_{n-1}^+}
X_n^{\alpha_n} X^{\alpha'} f_{Z'} \fH_n.
$$
Now we need to consider the cases $\alpha_n \geq 0$ and $\alpha_n < d$
separately. The argument is entirely similar in each case, so suppose
in fact that $\alpha_n \geq 0$. Then we show by induction on
$\alpha_n$ that 
$X_n^{\alpha_n} X^{\alpha'} f_{Z'} \fH_n \subseteq {\mathcal J}$
for each $(\alpha',Z')\in\Pi_{n-1}^+$.
This is immediate if $\alpha_n < d$, so take $\alpha_n \geq d$ and
consider the induction step.
Expanding $f_n$ using Lemma~\ref{fi}, the set
$$
X_n^{\alpha_n-d} X^{\alpha'} f_{Z'} f_n \fH_n \subseteq {\mathcal J}
$$
looks like the desired $X_n^{\alpha_n} X^{\alpha'} f_{Z'} \fH_n$
plus a sum of terms lying in
$X_n^{\alpha_n - d + e} {\mathcal I}_f' \fH_n$
with $0 \leq e < d$. It now suffices to show that
each such $X_n^{\alpha_n - d+e} {\mathcal I}_f' \fH_n \subseteq {\mathcal J}$.
But by (\ref{hyp}), 
$$
X_n^{\alpha_n - d+e} {\mathcal I}_f' \fH_n
\subseteq
\sum_{(\alpha',Z') \in \Pi_{n-1}^+} 
X_n^{\alpha_n - d+e} X^{\alpha'} f_{Z'} \fH_n
$$
and each such term lies in ${\mathcal J}$ by induction,
since $0 \leq \alpha_n - d+ e < \alpha_n$.
\end{proof}

\begin{Theorem}\label{find}
The canonical images of the elements
$$
\{
X^\alpha C^\beta T_w\:|\:\alpha \in \Z^n\hbox{ with }0 \leq \alpha_1,\dots,\alpha_n < d, \beta \in \Z_2^n, w \in \Sym_n\}
$$
form a basis for $\H_n^f$.
\end{Theorem}

\begin{proof}
By Lemmas~\ref{hardb} and \ref{ke},
the elements $\{X^\alpha f_Z\:|\:(\alpha,Z) \in \Pi_n^+\}$
form a basis for ${\mathcal I}_f$ viewed as a right $\fH_n$-module.
Hence Lemma~\ref{hardb} implies that the elements
$$\{X^\alpha\:|\:\alpha \in \Z^n\hbox{ with }0 \leq \alpha_1,\dots,\alpha_n < d\}$$
form a basis for a complement to ${\mathcal I}_f$ in $\H_n$ 
viewed as a right $\fH_n$-module.
The theorem follows at once.
\end{proof}

\Point{Cyclotomic Mackey theorem}\label{cycloind}
We will need a special case of a Mackey theorem
for cyclotomic Hecke-Clifford  superalgebras.
Let $f \in \P_1$ be a polynomial of degree $d > 0$ of
the special form (\ref{myform}).
Let $\H_n^f$ denote the corresponding cyclotomic Hecke-Clifford superalgebra.
Given any $y \in \H_n$, we will write $\tilde y$ for its canonical image
in $\H_n^f$. Thus, Theorem~\ref{find} says that the elements
$$
\{\tilde X^\alpha \tilde C^\beta \tilde T_w\:|\:\alpha \in \Z^n\hbox{ with }
0 \leq \alpha_1,\dots,\alpha_n < d,
\beta \in \Z_2^n, w \in \Sym_n\}
$$ 
form a basis for $\H_n^f$.

It is obvious from Theorem~\ref{find} that the subalgebra of
$\H_{n+1}^f$ spanned by the
$\tilde X^\alpha \tilde C^\beta \tilde T_w$
for $\alpha \in \Z^{n}$ with $0 \leq \alpha_1,\dots,\alpha_{n} < d$,
$\beta \in \Z_2^{n}$ and $w \in \Sym_{n}$ is isomorphic to $\H_{n}^f$.
We will 
write $\ind_{\H_n^f}^{\H_{n+1}^f}$ and $\res_{\H_n^f}^{\H_{n+1}^f}$ 
for the induction and restriction
functors between $\H_{n}^f$ and $\H_{n+1}^f$, to avoid confusion with
the affine analogue from (\ref{blah}).
So,
$$
\ind_{\H_n^f}^{\H_{n+1}^f} M = \H_{n+1}^f \otimes_{\H_{n}^f} M.
$$

\vspace{1mm}
\begin{Lemma}\label{power}
For all $1 \leq i < n, j \geq 0, k \in \Z_2$,
$$
X_{i+1}^j C_{i+1}^k T_i - T_i X_i^j C_i^k \in \sum_{h=1}^j \left(X_{i+1}^h \A_i
+ X_{i+1}^{h-1} C_{i+1} \A_i\right).
$$
\end{Lemma}

\begin{proof}
Rearrange (\ref{powertx}) and (\ref{powertz}).
\end{proof}

\begin{Lemma}\label{list}
{\rm(i)} $\H_{n+1}^f$ is a free right $\H_{n}^f$-module on basis
$$
\{
\tilde X_j^a \tilde C_j^b
\tilde T_j \dots \tilde T_{n} \:|\:
0 \leq a < d, b \in \Z_2, 1 \leq j \leq n+1\}.
$$
\begin{enumerate}
\item[(ii)]
$\displaystyle
\H_{n+1}^f = \H_n^f \tilde T_n \H_n^f \oplus \bigoplus_{0 \leq a < d,
b \in \Z_2} \tilde X_{n+1}^a \tilde C_{n+1}^b \H_n^f $
as an $(\H_n^f,\H_n^f)$-bimodule.
\item[(iii)]
For any $0 \leq a < d$, there are isomorphisms
$$
\tilde X_{n+1}^a\H_n^f  \simeq \H_n^f,
\qquad
\tilde X_{n+1}^a \tilde C_{n+1} \H_n^f \simeq \Pi \H_n^f,
\qquad
\H_n^f \tilde T_n \H_n^f \simeq \H_n^f \otimes_{\H_{n-1}^f} \H_n^f,
$$
as $(\H_n^f,\H_n^f)$-bimodules.
\end{enumerate}
\end{Lemma}

\begin{proof}
For (i), by Theorem~\ref{find} and dimension
considerations, we just need to check that
$\H_{n+1}^f$ is generated as a right $\H_n^f$-module by the given
elements. For this, it suffices to show 
that all elements of the form
$$
\tilde X^\alpha \tilde C^\beta \tilde T_j \dots \tilde T_n
\qquad
(\alpha \in \Z^n, \beta \in \Z_2^n, 
0 \leq \alpha_1,\dots,\alpha_{n+1} < d)
$$
lie in the right $\H_n^f$-module generated by 
$\{
\tilde X_k^a \tilde C_k^b
\tilde T_k \dots \tilde T_{n} \:|\:
0 \leq a < d, b \in \Z_2, j \leq k \leq n+1\}.$
This involves considering terms of the form
$\tilde X_k^a \tilde C_k^b \tilde T_{k-1} \dots \tilde T_n$
for $j < k \leq n+1$ and $0 \leq a < d, b \in \Z_2$,
for which Lemma~\ref{power} is useful.

For (ii),(iii), 
define a map $\H_n^f \times \H_n^f 
\rightarrow \H_{n}^f \tilde T_n \H_n^f,
(u,v) \mapsto u \tilde T_n v$.
This is $\H_{n-1}^f$-balanced, so induces a well-defined
epimorphism
$$
\Phi:\H_n^f \otimes_{\H_{n-1}^f} \H_n^f \rightarrow \H_n^f \tilde T_n \H_n^f
$$
of $(\H_n^f,\H_n^f)$-bimodules.
We know from (i) that 
$\H_n^f \otimes_{\H_{n-1}^f} \H_n^f$
is a free right $\H_n^f$-module on basis
$\tilde X_j^a \tilde C_j^b
\tilde T_j \dots \tilde T_{n-1} \otimes 1$
for $1 \leq j \leq n,0 \leq a < d, b \in \Z_2$.
But $\Phi$ maps these elements to
$\tilde X_j^a \tilde C_j^b
\tilde T_j \dots \tilde T_{n-1} \tilde T_n$
which, again using (i), form a basis for $\H_n^f \tilde T_n \H_n^f$
as a free right $\H_n^f$-module.
This shows that
$\H_n^f \tilde T_n \H_n^f \simeq \H_n^f \otimes_{\H_{n-1}^f} \H_n^f.$
Now the remaining parts of (ii), (iii) are obvious consequences of (i).
\end{proof}

We have now decomposed $\H_{n+1}^f$ as an $(\H_n^f,\H_n^f)$-bimodule.
So the same argument as for Theorem~\ref{TMackey} easily gives:

\vspace{1mm}
\begin{Theorem}\label{cycloMackey}
Let $M$ be an $\H_{n}^f$-module.
Then, there is a natural isomorphism
$$
\res_{\H_n^f}^{\H_{n+1}^f} \ind_{\H_n^f}^{\H_{n+1}^f} M 
\simeq
(M\oplus\Pi M)^{\oplus d}\oplus \ind_{\H_{n-1}^f}^{\H_n^f} \res^{\H_n^f}_{
\H_{n-1}^f} M
$$
of $\H_{n}^f$-modules.
\end{Theorem}

\Point{Duality}
We wish next to prove that the induction functor $\ind_{\H_n^f}^{\H_{n+1}^f}$
commutes with the $\tau$-duality.
We need a little preliminary work.

\vspace{1mm}
\begin{Lemma}\label{newcleo}
For $1 \leq i \leq n$ and $a \geq 0$,
\begin{align}\label{ca}
(\tilde T_n+\xi \tilde C_n \tilde C_{n+1}) \dots (\tilde T_i+\xi \tilde C_i \tilde C_{i+1}) 
\tilde X_i^a \tilde
T_i \dots \tilde T_n
&=
\tilde X_{n+1}^a + (*)\\\intertext{where $(*)$ is a term lying in
$\displaystyle
\H_n^f \tilde T_n \H_n^f + 
\sum_{k = 1}^{a-1}\left(
\tilde X_{n+1}^{k} \H_n^f+\tilde X_{n+1}^{k-1} \tilde C_{n+1}\H_n^f\right)$,
and}
(\tilde T_n+\xi \tilde C_n \tilde C_{n+1}) \dots (\tilde T_i+\xi \tilde C_i \tilde C_{i+1}) 
\tilde X_i^a \tilde C_i \tilde
T_i \dots \tilde T_n
&=
\tilde X_{n+1}^a \tilde C_{n+1}+ (**)\label{cb}
\end{align}
where $(**)$ is a term lying in
$\displaystyle
\H_n^f \tilde T_n \H_n^f + 
\sum_{k = 1}^{a}\left(
\tilde X_{n+1}^{k} \H_n^f+\tilde X_{n+1}^{k-1} \tilde C_{n+1}\H_n^f\right)$.
\end{Lemma}

\begin{proof}
We prove (\ref{ca}) and (\ref{cb}) simultaneously by induction 
on $n = i,i+1,\dots$.
In case $n = i$, they follow from a calculation involving 
(\ref{powertx}) or (\ref{powertz}) respectively,
together with Lemma~\ref{power} to commute $\tilde X_{n+1}^a$ and
$\tilde X_{n+1}^a \tilde C_{n+1}$ past $\tilde T_n$.
The induction step is similar, 
noting that both $(\tilde T_{n+1} + \xi \tilde C_{n+1} \tilde C_{n+2})$
and $\tilde T_{n+1}$ centralize $\H_n^f$.
\end{proof}

\vspace{1mm}
\begin{Lemma} \label{cleo}
For any $s\in \H_n^f$,
$$
\tilde X_{n+1}^d s = -a_0  s
+
\left(\hbox{a term lying in }\displaystyle
\H_n^f \tilde T_n \H_n^f + 
\sum_{k = 1}^{d-1}
\left(
\tilde X_{n+1}^{k} \H_n^f+
\tilde X_{n+1}^{k-1} \tilde C_{n+1}\H_n^f\right)\right).
$$
\end{Lemma}

\begin{proof}
Recall the polynomial $f\in F[X_1]$ from (\ref{myform}).
Its image $\tilde f \in \H_n^f$ is zero.
Hence, 
$$
(\tilde T_n+\xi \tilde C_n \tilde C_{n+1})
\dots
(\tilde T_1+\xi \tilde C_1 \tilde C_{2}) \tilde f
\tilde T_1\dots\tilde T_n = 0.
$$
But a calculation using (\ref{ca}) shows that the left
hand side equals $\tilde X_{n+1}^d+a_0$ modulo
terms of the given form. The lemma follows easily on multiplying on the
right by $s$.
\end{proof}

\begin{Lemma}\label{grinch}
There exists an even $(\H_n^f,\H_n^f)$-bimodule homomorphism
$\theta:\H_{n+1}^f \rightarrow \H_n^f$ such that $\ker \theta$
contains no non-zero left ideals of $\H_{n+1}^f$.
\end{Lemma}

\begin{proof}
By Lemma~\ref{list}(ii), we know that
$$
\H_{n+1}^f = \H_n^f \oplus \tilde C_{n+1} \H_n^f 
\oplus \bigoplus_{a = 1}^{d-1} \left(\tilde X_{n+1}^a \H_n^f \oplus 
\tilde X_{n+1}^a \tilde C_{n+1} \H_n^f\right)
\oplus
\H_n^f \tilde T_n \H_n^f
$$
as an $(\H_n^f,\H_n^f)$-bimodule.
Let $\theta:\H_{n+1}^f \rightarrow \H_n^f$ be the projection onto the
first summand of this bimodule decomposition.
We just need to show that if $y \in \H_{n+1}^f$ has the property
that $\theta(hy) = 0$ for all $h \in \H_{n+1}^f$, then $y = 0$.
Using Lemma~\ref{list}(i), we may write
$$
y = \sum_{a=0}^{d-1} \left(\tilde X_{n+1}^a s_a + \tilde X_{n+1}^a \tilde C_{n+1} t_a\right)
+
\sum_{a=0}^{d-1} \sum_{j=1}^{n}
\left(\tilde X_j^a \tilde T_j \dots \tilde T_{n} 
u_{a,j}
+
\tilde X_j^a \tilde C_j \tilde T_j \dots \tilde T_n
v_{a,j}\right)
$$
for $s_a, t_a, u_{a,j},v_{a,j} \in \H_n^f$.
Consider $\theta(\tilde X_{n+1}^{d-1} \tilde C_{n+1} y)$.
An application of Lemma~\ref{power} reveals that this equals
$t_{d-1}$, hence $t_{d-1} = 0$.
Next consider $\theta(\tilde X_{n+1} y)$. Using 
Lemma~\ref{cleo} as well as Lemma~\ref{power} this time,
we get that $s_{d-1} = 0$.
Now consider similarly
$\theta(\tilde X_{n+1}^{d-2} \tilde C_{n+1} y)$,
$\theta(\tilde X_{n+1}^{2}y)$, $\theta(\tilde X_{n+1}^{d-3} \tilde C_{n+1} y),
\theta(\tilde X_{n+1}^{3}y), \dots$ in turn to deduce
$t_{d-2} = s_{d-2} = t_{d-3} = s_{d-3} = \dots = 0$.

We have now reduced to the case that
$$
y = 
\sum_{a=0}^{d-1} \sum_{j=1}^{n}
\left(\tilde X_j^a \tilde T_j \dots \tilde T_{n} 
u_{a,j}
+
\tilde X_j^a \tilde C_j \tilde T_j \dots \tilde T_n
v_{a,j}\right).
$$
Now consider $y' := 
(\tilde T_n + \xi \tilde C_n \tilde C_{n+1})y$.
Note
$$
(\tilde T_n + \xi \tilde C_n \tilde C_{n+1}) \tilde T_{n-1}
\tilde T_n 
= 
\tilde T_{n-1} \tilde T_n (\tilde T_{n-1}
+
\xi \tilde C_{n-1} \tilde C_{n})
$$
so the terms of $y$ with $j < n$ yield terms of $y'$ 
which lie
in $\H_n^f \tilde T_n \H_n^f$. Hence, by Lemma~\ref{newcleo} too,
$$
y' = \tilde X_{n+1}^{d-1} \tilde C_{n+1} v_{d-1,n}
+(*)
$$
where $(*)$ is a term lying in
$\displaystyle 
\H_n^f\tilde T_n \H_n^f+\sum_{k=1}^{d-1}
\left(\tilde X_{n+1}^k \H_n^f+
\tilde X_{n+1}^{k-1} \tilde C_{n+1} \H_n^f\right)$.
Now multiplying $y'$
by $\tilde X_{n+1}^{d-1} \tilde C_{n+1}$ and applying $\theta$,
as in the previous paragraph, gives that $v_{d-1,n} = 0$.
Hence in fact, by Lemma~\ref{newcleo} once more, we have that
$$
y' = 
\tilde X_{n+1}^{d-1} u_{d-1,n}
+(*)
$$
and now one gets $u_{d-1,n} =0$ on multiplying by $\tilde X_{n+1}$ and
applying $\theta$, again as in the previous paragraph. Continuing in this
way gives that all $u_{a,n} = v_{a,n} = 0$.

Now repeat the argument in the previous paragraph 
again, 
this time considering 
$y' := 
(\tilde T_n + \xi \tilde C_n \tilde C_{n+1})(\tilde T_{n-1} + 
\xi \tilde C_{n-1} \tilde C_{n})y$,
to get that all $u_{a,n-1}=v_{a,n-1} = 0$.
Continuing in this way eventually gives the desired conclusion:
$y = 0$.
\end{proof}

Now we are ready to prove the main result of the subsection:

\vspace{1mm}
\begin{Theorem}\label{indcoind}
There is a natural isomorphism
$\H_{n+1}^f\otimes_{\H_n^f} M \simeq 
\hom_{\H_n^f}(\H_{n+1}^f, M)$
for all $\H_{n}^f$-modules $M$.
\end{Theorem}

\begin{proof}
We show that there is an even isomorphism
$\phi:\H_{n+1}^f  \rightarrow \hom_{\H_n^f}(\H_{n+1}^f,\H_n^f)$
of $(\H_{n+1}^f, \H_n^f)$-bimodules.
The lemma then follows on applying the functor $? \otimes_{\H_n^f} M$:
one obtains natural isomorphisms
$$
\H_{n+1}^f \otimes_{\H_n^f} M \stackrel{\phi\otimes \id}{\longrightarrow} 
\hom_{\H_n^f}(\H_{n+1}^f,\H_n^f) \otimes_{\H_n^f} M 
\simeq \hom_{\H_{n}^f}(\H_{n+1}^f,M).
$$
Note the existence of the
second isomorphism here uses the fact that
$\H_{n+1}^f$ is a projective left $\H_n^f$-module,
see \cite[20.10]{AF}.

To construct $\phi$, let $\theta$ be as in Lemma~\ref{grinch}, and define
$\phi(h)$ to be the map $h \theta$, for each $h \in \H_{n+1}^f$.
One easily checks that $\phi:\H_{n+1}^f
\rightarrow \hom_{\H_n^f}(\H_{n+1}^f, \H_n^f)$ is then a well-defined
homomorphism
of $(\H_{n+1}^f, \H_n^f)$-bimodules.
To see that it is an isomorphism, it suffices by dimensions
to check it is injective.
Suppose $\phi (h) = 0$ for some $h \in \H_{n+1}^f$.
Then for every $x \in \H_{n+1}^f$, $\theta(xh) = 0$, 
i.e. the left ideal $\H_{n+1}^f h$
is contained in $\ker \theta$. So Lemma~\ref{grinch} implies $h = 0$.
\end{proof}

\begin{Corollary}
$\H_n^f$ is a Frobenius superalgebra, i.e. there is an even isomorphism
of left $\H_n^f$-modules
$\H_n^f \simeq \hom_F(\H_n^f, F)$ between the left regular module
and the $F$-linear dual of the right regular module.
\end{Corollary}

\begin{proof}
Proceed by induction on $n$.
For the induction step,
\begin{align*}
\H_{n}^f &\simeq 
\H_{n}^f \otimes_{\H_{n-1}^f} \H_{n-1}^f
\simeq
\H_{n}^f \otimes_{\H_{n-1}^f} \hom_F(\H_{n-1}^f, F)\\
&\simeq \hom_{\H_{n-1}^f}(\H_{n}^f, \hom_F(\H_{n-1}^f, F))
\simeq
\hom_F(\H_{n-1}^f \otimes_{\H_{n-1}^f} \H_{n}^f, F)\\
&\simeq
\hom_F(\H_{n}^f, F),
\end{align*}
applying Theorem~\ref{indcoind} and adjointness of tensor and Hom.
\end{proof}

For the next corollary, recall the duality induced by $\tau$ $(\ref{Etau})$ 
on finite dimensional $\H_n$-modules.
Since $\tau$ leaves the two-sided ideal ${\mathcal I}_f$
invariant, it induces a duality also denoted 
$\tau$ on finite dimensional $\H_n^f$-modules.

\vspace{1mm}
\begin{Corollary}\label{joey}
The exact functor $\ind_{\H_n^f}^{\H_{n+1}^f}$
is both left and right adjoint to $\res^{\H_{n+1}^f}_{\H_n^f}$.
Moreover, it commutes with duality in the sense that there is a natural 
isomorphism
$$
\ind_{\H_n^f}^{\H_{n+1}^f}(M^\tau) \simeq 
(\ind_{\H_n^f}^{\H_{n+1}^f} M)^\tau
$$
for all finite dimensional $\H_n^f$-modules $M$.
\end{Corollary}

\begin{proof}
The fact that $\ind_{\H_n^f}^{\H_{n+1}^f} = \H_{n+1}^f \otimes_{\H_n^f}?$ is right adjoint to
$\res^{\H_{n+1}^f}_{\H_n^f}$ is immediate from Theorem~\ref{indcoind},
since $\hom_{\H_n^f}(\H_{n+1}^f, ?)$ is right adjoint to restriction
by adjointness of tensor and Hom.
But on finite dimensional modules,
a standard check using (\ref{switch}) shows that
the functor
$\tau \circ \ind_{\H_n^f}^{\H_{n+1}^f} \circ \tau$
is also right adjoint to restriction.
Now the remaining part of the corollary follows by uniqueness of adjoint 
functors.
\end{proof}

\Point{Modifications in the degenerate case}
In the degenerate case, $\H_n^f$ becomes the {\em 
cyclotomic Sergeev superalgebra} defined for $f \in F[x_1] \subset \H_n$
a polynomial of the form
$$
x_1^d + a_{d-2} x_1^{d-2} + a_{d-4} x_1^{d-4}+\dots,
$$
i.e. the powers of $x_1$ appearing are either all even or all odd
and the leading coefficient
is $1$. By definition, $\H_n^f = \H_n / {\mathcal I}_f$ where ${\mathcal I}_f$
is the two-sided ideal generated by $f$.
The basis theorem says that
$\H_n^f$ has basis given by the images of
$$
\{
x^\alpha c^\beta w\:|\:\alpha \in \Z_{\geq 0}^n\hbox{ with }0 \leq \alpha_i < d, 
\beta \in \Z_2^n, w \in \Sym_n\},
$$
where ${\mathcal I}_f$ is the two-sided ideal of $\H_n$ generated by $f$.

The proof is entirely similar to that of Theorem~\ref{find}; actually
in this case it is much more straightforward. To give a little more detail,
one defines the element
$f_i = s_{i-1} \dots s_1 f s_1 \dots s_{i-1}$
for each $i = 1,\dots,n$ then defines 
$f_Z$ as before for $Z \subseteq \{1,\dots,n\}$.
By (\ref{powertx2}),
$$
f_i = x_i^d + (\hbox{a linear combination of terms
lying in $\P_{i-1} x_i^e \fH_i$ for $0 \leq e < d$}).
$$
Given this one easily proves as in Lemma~\ref{hardb} that the
$$
\{x^\alpha f_Z\:|\:Z \subseteq \{1,\dots,n\},
\alpha \in \Z_{\geq 0}^n\hbox{ with 
$0 \leq \alpha_i < d$ if $i \notin Z$}\}
$$
form a basis for $\H_n$ viewed as a right $\fH_n$-module.
Moreover, arguing as for Lemma~\ref{ke} the
$$
\{x^\alpha f_Z\:|\:\emptyset \neq Z \subseteq \{1,\dots,n\},
\alpha \in \Z_{\geq 0}^n\hbox{ with 
$0 \leq \alpha_i < d$ if $i \notin Z$}\}
$$
form a basis for ${\mathcal I}_f$ as a right $\fH_n$-module.
Finally, the proof is completed as in Theorem~\ref{find}.

Theorem~\ref{cycloMackey} 
goes through without significant alteration.
Note a suitable decomposition of $\H_{n+1}^f$ as an 
$(\H_n^f,\H_n^f)$-bimodule is
\begin{equation}\label{deed}
\H_{n+1}^f = \bigoplus_{0 \leq a < d, b \in \Z_2} x_{n+1}^a c_{n+1}^b 
\H_n^f \oplus
\H_n^f s_n \H_n^f.
\end{equation}
To prove that induction commutes with duality, i.e. Theorem~\ref{indcoind},
there is a slight twist in proving the analogue of Lemma~\ref{grinch}:
the map $\theta$ should be taken to be the projection
$$
\theta:\H_{n+1}^f \rightarrow x_{n+1}^{d-1} \H_n^f \simeq \H_n^f
$$
along the direct sum decomposition (\ref{deed}).

\ifbook@\pagebreak\fi

\section{The category of integral representations}

\Point{Affine Kac-Moody algebra}
Now we introduce some standard Lie theoretic notation. Let us treat the
case $h \neq \infty$ first,
when we let $\ell = (h-1) / 2$
and $\mathfrak g$ denote the twisted affine Kac-Moody algebra of
type $A_{2 \ell}^{(2)}$ (over $\mathbb C$), see 
\cite[ch. 4, table Aff 2]{Kac}.
In particular we label the Dynkin diagram
by the index set
$I = \{0,1,\dots,\ell\}$ as follows: 
$$
{\begin{picture}(340, 15)%
\put(6,5){\circle{4}}%
\put(101,2.3){$<$}%
\put(12,2.3){$<$}%
\put(256,2.3){$<$}%
\put(25, 5){\circle{4}}%
\put(44, 5){\circle{4}}%
\put(8, 4){\line(1, 0){15.5}}%
\put(8, 6){\line(1, 0){15.5}}%
\put(27, 5){\line(1, 0){15}}%
\put(46, 5){\line(1, 0){1}}%
\put(49, 5){\line(1, 0){1}}%
\put(52, 5){\line(1, 0){1}}%
\put(55, 5){\line(1, 0){1}}%
\put(58, 5){\line(1, 0){1}}%
\put(61, 5){\line(1, 0){1}}%
\put(64, 5){\line(1, 0){1}}%
\put(67, 5){\line(1, 0){1}}%
\put(70, 5){\line(1, 0){1}}%
\put(73, 5){\line(1, 0){1}}%
\put(76, 5){\circle{4}}%
\put(78, 5){\line(1, 0){15}}%
\put(95, 5){\circle{4}}%
\put(114,5){\circle{4}}%
\put(97, 4){\line(1, 0){15.5}}%
\put(97, 6){\line(1, 0){15.5}}%
\put(6, 11){\makebox(0, 0)[b]{$_{0}$}}%
\put(25, 11){\makebox(0, 0)[b]{$_{1}$}}%
\put(44, 11){\makebox(0, 0)[b]{$_{{2}}$}}%
\put(75, 11){\makebox(0, 0)[b]{$_{{\ell-2}}$}}%
\put(96, 11){\makebox(0, 0)[b]{$_{{\ell-1}}$}}%
\put(114, 11){\makebox(0, 0)[b]{$_{{\ell}}$}}%
\put(250,5){\circle{4}}%
\put(269,5){\circle{4}}%
\put(251.3,3.2){\line(1,0){16.6}}%
\put(252,4.4){\line(1,0){15.2}}%
\put(252,5.6){\line(1,0){15.2}}%
\put(251.3,6.8){\line(1,0){16.6}}%
\put(250, 11){\makebox(0, 0)[b]{$_{{0}}$}}%
\put(269, 11){\makebox(0, 0)[b]{$_{{1}}$}}%
\put(320,2){\makebox(0,0)[b]{if $\ell = 1$.}}%
\put(170,2){\makebox(0,0)[b]{if $\ell \geq 2$, and}}%
\end{picture}}
$$
The weight lattice is denoted $P$, the simple roots are 
$\{\alpha_i \:|\:i \in I\} \subset P$ and the corresponding
simple coroots are $\{h_i\:|\:i \in I\} \subset P^*$.
The Cartan matrix 
$\left(\langle h_i, \al_j\rangle\right)_{0\leq i,j\leq \ell}$ is 
$$
\left(
\begin{matrix}
2 & -2 & 0 & \cdots & 0 & 0 & 0 \\
-1 & 2 & -1 & \cdots & 0 & 0 & 0 \\
0 & -1 & 2 & \cdots & 0 & 0 & 0 \\
 & & & \ddots & & & \\
0 & 0 & 0 & \dots & 2 & -1& 0 \\
0 & 0 & 0 & \dots & -1 & 2& -2 \\
0 & 0 & 0 & \dots & 0 & -1& 2 \\
\end{matrix}
\right)
\quad
\text{if $\ell\geq 2$, and}\quad
\left(
\begin{matrix}
2 & -4 \\
-1 & 2
\end{matrix}
\right)
\quad\text{if $\ell=1$.}
$$
Let $\{\Lambda_i\:|\:i \in I\} \subset P$ 
denote fundamental dominant weights, so
that
$\langle h_i, \Lambda_j \rangle = \delta_{i,j}$, and let
$P_+\subset P$ denote the set of all dominant integral weights.
Set
\begin{equation}\label{cdelta}
c = h_0 + \sum_{i = 1}^\ell 2 h_i,\qquad
\delta = \sum_{i=0}^{\ell - 1} 2 \alpha_i + \alpha_\ell.
\end{equation}
Then the $\Lambda_0,\dots,\Lambda_\ell, \delta$ form a $\Z$-basis for $P$,
and $\langle c, \alpha_i \rangle = \langle h_i, \delta \rangle = 0$ for
all $i \in I$.

In the case $h = \infty$, 
we make the following changes to these definitions.
First, we let $\ell = \infty$, and $\mathfrak g$ denotes the
Kac-Moody algebra of type $B_\infty$, see \cite[$\S$7.11]{Kac}.
So $I = \{0,1,2,\dots\}$, corresponding to the nodes of the Dynkin diagram
$$
{\begin{picture}(100, 15)%
\put(6,5){\circle{4}}%
\put(12,2.3){$<$}%
\put(25, 5){\circle{4}}%
\put(44, 5){\circle{4}}%
\put(8, 4){\line(1, 0){15.5}}%
\put(8, 6){\line(1, 0){15.5}}%
\put(27, 5){\line(1, 0){15}}%
\put(46, 5){\line(1, 0){1}}%
\put(49, 5){\line(1, 0){1}}%
\put(52, 5){\line(1, 0){1}}%
\put(55, 5){\line(1, 0){1}}%
\put(58, 5){\line(1, 0){1}}%
\put(61, 5){\line(1, 0){1}}%
\put(64, 5){\line(1, 0){1}}%
\put(67, 5){\line(1, 0){1}}%
\put(70, 5){\line(1, 0){1}}%
\put(73, 5){\line(1, 0){1}}%
\put(6, 11){\makebox(0, 0)[b]{$_{0}$}}%
\put(25, 11){\makebox(0, 0)[b]{$_{1}$}}%
\put(44, 11){\makebox(0, 0)[b]{$_{{2}}$}}%
\end{picture}}
$$
\ifbook@
The Cartan matrix 
$\left(\langle h_i, \al_j\rangle\right)_{i,j\geq 0}$ is 
$$
\left(
\begin{matrix}
2 & -2 & 0 &&\\
-1 & 2 & -1 & 0&\\
0 & -1 & 2 & -1& \\
 &0 & -1& 2&\ddots\\
&&&\ddots&\ddots \\
\end{matrix}
\right)
$$
\fi
Note certain notions, for example the element
$c$ from (\ref{cdelta}), only make sense if one passes to the completed
algebra $b_{\infty}$, see \cite[$\S$7.12]{Kac}, though the intended
meaning whenever we make use of them should be obvious regardless.

\nopagebreak
Now, for either $h < \infty$ or $h = \infty$,
we let $U_\Q$ denote the $\Q$-subalgebra of the
universal enveloping algebra of $\mathfrak g$  
generated by the Chevalley generators $e_i,f_i,h_i\:(i \in I)$.
Recall these are subject only to the relations
\begin{align}\label{cr1}
[h_i, h_j] = 0,
\qquad &[e_i, f_j] = \delta_{i,j} h_i,\\\label{cr2}
[h_i, e_j] = \langle h_i, \alpha_j \rangle e_j,
\qquad& [h_i, f_j] = -\langle h_i, \alpha_j \rangle f_j,\\
\label{cr3}
(\operatorname{ad} e_i)^{1-\langle h_i,\alpha_k\rangle} e_k = 0,\qquad&
(\operatorname{ad} f_i)^{1-\langle h_i,\alpha_k \rangle} f_k = 0
\end{align}
for all $i,j,k \in I$ with $i \neq k$.
We let $U_\Z$ denote the $\Z$-form of $U_\Q$ generated by the divided powers
$e_i^{(n)} = e_i^n / n!$ and $f_i^{(n)} = f_i^n / n!$.
Then, $U_\Z$ has the usual triangular decomposition
$$
U_\Z = U_\Z^- U_\Z^0 U_\Z^+.
$$
We are particularly concerned here with the plus part $U_\Z^+$, generated by
all $e_i^{(n)}$.
It is a graded Hopf algebra over $\Z$ via 
the {\em principal grading}
$\deg(e_i^{(n)}) = n$ for all $i \in I, n \geq 0$.

\Point{Cyclotomic Hecke-Clifford superalgebras revisited}
Given $i \in I$, 
define
\begin{equation}\label{resdef}
q(i) := 2 \frac{q^{2i+1} + q^{-2i-1}}{q+q^{-1}} \in F.
\end{equation}
Note in particular that $q(0) = 2$.
For $\la \in P_+$, let ${\mathcal I}_\lambda$ denote the two-sided ideal
of $\H_n$ generated by the element
\begin{equation}\label{E210900_4}
(X_1 - 1)^{\langle h_0, \lambda\rangle} 
\prod_{i=1}^\ell (X_1 +X_1^{-1}- q(i))^{\langle h_i, \lambda \rangle}
\end{equation}
Up to a power of the unit $X_1$, this is an element of the
form (\ref{myform}), so the quotient superalgebra
$$
\H_n^\la := \H_n / {\mathcal I}_\la
$$
is a special case of the cyclotomic Hecke-Clifford superalgebras
introduced in the previous section.
This {quotient} of $\H_n$ defined for $\la \in P_+$
should not be confused with the
parabolic {subalgebra} $\H_\mu$ defined earlier for $\mu$ a composition
of $n$. 

Theorem~\ref{find} immediately gives the following basis
theorem for $\H_n^\la$:

\vspace{1mm}
\begin{Theorem}\label{cyclo}
For any $\la \in P_+$, the canonical images of the elements
$$
\{X^\alpha C^\beta T_w \:|\:\alpha \in \Z^n
\hbox{ with }0 \leq \alpha_1,\dots,\alpha_n < \langle c, \la \rangle, 
\beta \in \Z_2^n,
w \in \Sym_n\}
$$
form a basis for $\H_n^\la$. In particular, 
$\dim \H_n^\la = (2\langle c, \la \rangle)^n (n!).$
\end{Theorem}

\begin{Remark}\label{fhrem}\rm
In the special case $\lambda = \Lambda_0$ is the first fundamental
dominant weight, the cyclotomic Hecke-Clifford superalgebra $\H_n^\la$
can be identified with the superalgebra $\fH_n$.
This follows easily from Theorem~\ref{cyclo}: 
since $\langle c, \Lambda_0\rangle = 1$,
$\H_n^{\Lambda_0}$ has basis given by the images of the 
elements $\{C^\beta T_w\:|\:\beta \in \Z_2^n, w \in \Sym_n\}$ 
just as $\fH_n$, and the multiplications are the same by construction.
See also  \cite[Prop. 3.5]{JNaz}.
\end{Remark}

Introduce the functors
\begin{equation}\label{prind}
\pr^\la:\mod{\H_n} \rightarrow \mod{\H_n^\lambda},
\qquad
\inf^\la:\mod{\H_n^\lambda} \rightarrow \mod{\H_n}.
\end{equation}
Here, $\inf^\la$ is simply inflation along the canonical
epimorphism $\H_n \rightarrow \H_n^\la$, while
on a module $M$, $\pr^\la M = M / {\mathcal I}_\la M$ with the induced
action of $\H_n^\la$.
The functor $\inf^\la$ is right adjoint to $\pr^\la$, i.e.
there is a functorial isomorphism
\begin{equation}
\label{Eprinfl}
\hom_{\H_n^\la}(\pr^\la M,N)\simeq\hom_{\H_n}(M,\inf^\la N). 
\end{equation}
Note we will generally be sloppy and omit the functor $\inf^\la$ 
in our notation. In other words, we generally
identify $\mod{\H_n^\la}$ with the full subcategory of
$\mod{\H_n}$ consisting of all modules $M$ with ${\mathcal I}_\la M = 0$.

\Point{\boldmath Elements $\tilde \Phi_j$}\label{nazelts}
We will need certain elements of $\H_n$
defined originally by Jones and Nazarov.
Given $1 \leq j < n$, define
\begin{align}
z_j &
:= 
X_j+X_j^{-1}-X_{j+1}-X_{j+1}^{-1}
= X_j^{-1} (X_j X_{j+1}-1) (X_j X_{j+1}^{-1}-1),\\\label{nazelt}
\tilde\Phi_j
&:=
z_j^2 T_j + \xi \frac{z_j^2}{X_jX_{j+1}^{-1}-1}
- \xi \frac{z_j^2}{X_j X_{j+1}-1}C_j C_{j+1}.
\end{align}
Then $\tilde \Phi_j$ is equal to $z_j^2 \Phi_j$ where $\Phi_j$ is the element
defined by Jones and Nazarov in \cite[(3.6)]{JNaz}.
Note $\tilde \Phi_j$ really does make sense as an element of $\H_n$, unlike $\Phi_j$ which
belongs to a certain localization.
An easy calculation as in
\cite[(3.7)]{JNaz} gives that
\begin{align}\label{nazeltp}
\tilde \Phi_j X_j^{\pm 1} = X_{j+1}^{\pm 1} \tilde \Phi_j,\qquad
\tilde \Phi_{j} X_{j+1}^{\pm 1} &= X_{j}^{\pm 1} \tilde \Phi_j,\qquad
\tilde \Phi_j X_k^{\pm 1} = X_k^{\pm 1} \tilde \Phi_j,\\
\tilde \Phi_j C_j = C_{j+1} \tilde \Phi_j,
\qquad\;
\tilde \Phi_j C_{j+1} &= C_{j} \tilde \Phi_j,
\qquad\:\:\:\:\:\:\:
\tilde \Phi_j C_k = C_k \tilde \Phi_j\label{nazeltpp}
\end{align}
for $k \neq j,j+1$.
Moreover, 
\cite[Prop. 3.1]{JNaz} implies that
\begin{multline}\label{nazeltsq}
\tilde \Phi_j^2 
= 
z_j^2\big(X_j^{-2} (X_j X_{j+1}-1)^2 (X_j X_{j+1}^{-1}-1)^2
- \xi^2 X_j^{-1}X_{j+1}^{-1} (X_j X_{j+1}-1)^2\\
- \xi^2 X_j^{-1}
X_{j+1} (X_j X_{j+1}^{-1}-1)^2\big).\qquad
\end{multline}
In order to make use of this, we need the following
technical lemma:

\vspace{1mm}
\begin{Lemma}\label{nazelttech}
Suppose $a,b \in F^\times$ with $b+b^{-1} = q(i)$
for some $i \in I$.
If
\begin{multline*}
a^{-2}(ab-1)^2(ab^{-1}-1)^2\big(a^{-2} (ab-1)^2 (ab^{-1}-1)^2
- \xi^2 a^{-1}b^{-1}(ab-1)^2\\
- \xi^2 a^{-1} b(ab^{-1}-1)^2\big)=0.\qquad
\end{multline*}
then $a+a^{-1} = q(j)$ for $j \in I$ with $|i-j| \leq 1$.
\end{Lemma}

\begin{proof}
Follow \cite[(4.1)--(4.4)]{JNaz}.
\end{proof}

A more lengthy calculation also as in \cite[Prop. 3.1]{JNaz} shows that
the elements $\tilde \Phi_j$ satisfy the braid relations, i.e.
\begin{equation}\label{nazeltbraid}
\tilde \Phi_i \tilde \Phi_j = \tilde \Phi_j \tilde \Phi_i,\qquad
\tilde \Phi_i \tilde \Phi_{i+1} \tilde \Phi_i = 
\tilde \Phi_{i+1} \tilde \Phi_{i} \tilde \Phi_{i+1},
\end{equation}
for all admissible $i,j$ with $|i-j|> 1$.
This means that for any $w \in \Sym_n$, we obtain well-defined elements
$\tilde \Phi_w \in \H_n$, namely,
$\tilde \Phi_w := \tilde \Phi_{i_1} \dots \tilde \Phi_{i_m}$
where $w = s_{i_1} \dots s_{i_m}$ is any reduced expression for $w$.
According to (\ref{nazeltp}),(\ref{nazeltpp}), 
these elements have the property that
\begin{equation}\label{nazeltp2}
\tilde \Phi_w X_i^{\pm 1} = X_{wi}^{\pm 1}\tilde \Phi_w,\qquad
\tilde \Phi_w C_i = C_{wi}\tilde \Phi_w,\qquad
\end{equation}
for all $w \in S_n, 1 \leq i \leq n$.
Note we will not make essential use of (\ref{nazeltp2}) or 
the fact that the
$\tilde \Phi_i$ satisfy the braid relations in what follows.

\Point{Integral representations}\label{intrep}
Now call an $\A_n$-module $M$ {\em integral} if it is finite dimensional
and moreover all eigenvalues of 
$X_1+X_1^{-1},\dots,X_n+X_n^{-1}$ on $M$ are of the form $q(i)$ 
for $i \in I$, see (\ref{resdef}).
Call an $\H_n$-module, or more generally an $\H_\mu$-module
for $\mu$ a composition of $n$, 
{\em integral} if it is integral on restriction to $\A_n$.
In what follows we will restrict our attention to these modules,
and write
$\rep_I{\H_n}$ (resp. $\rep_I{\A_n}$, $\rep_I{\H_\mu}$)
for the full subcategory of $\mod{\H_n}$ (resp.
$\mod{\A_n}$, $\mod{\H_\mu}$) consisting of all integral modules.

\vspace{1mm}
\begin{Lemma}
\label{L210900}
Let $M$ be a finite dimensional $\H_n$-module, and 
$1\leq j\leq n$. 
Assume that the eigenvalues of $X_j+X_j^{-1}$ on $M$ are of the form $q(i)$, 
$i\in I$. Then the same is true for the eigenvalues of all 
other $X_k+X_k^{-1}$, 
$k=1,2,\dots,n$. 
\end{Lemma}
\begin{proof}
It suffices to show that the eigenvalues of $X_j+X_j^{-1}$ 
are of the form $q(i)$ if and only if the eigenvalues of 
$X_{j+1}+X_{j+1}^{-1}$ are of the same form, for $1\leq j<n$. 
Actually, by an argument involving conjugation with the
automorphism $\sigma$, it suffices just to prove the `if' part.
So assume that all eigenvalues of $X_{j+1}+X_{j+1}^{-1}$ on $M$ are
of the form $q(i)$ for various $i \in I$.
Let $a \neq 0$ be an eigenvalue for the action of $X_j$ on $M$.
We have to prove that $a+a^{-1}$ is also of the form $q(i)$.
Since $X_j$ and $X_{j+1}$ commute, we can pick $v$ lying in the
$a$-eigenspace of $X_j$ so that $v$ is also
an eigenvector for $X_{j+1}$, of eigenvalue $b$ say. By assumption,
$b + b^{-1} = q(i)$ for some $i \in I$. 
Now let $\tilde \Phi_j$ be the element (\ref{nazelt}).
By (\ref{nazeltp}),
$(X_{j+1}+X_{j+1}^{-1})\tilde\Phi_j=\tilde\Phi_j (X_j+X_j^{-1}).$
So if $\tilde\Phi_j v \neq 0$, we get that 
$$
(a+a^{-1}) \tilde\Phi_j v = \tilde\Phi_j (X_j + X_j^{-1}) v = 
(X_{j+1}+X_{j+1}^{-1}) \tilde\Phi_j v
$$
so that $a + a^{-1} = q(i')$ for some $i' \in I$ by assumption.
Else, $\tilde \Phi_j v = 0$ so $\tilde\Phi_j^2 v=0$. 
So applying (\ref{nazeltsq}) and Lemma~\ref{nazelttech},
we again get that $a+a^{-1} = q(i')$ for some $i' \in I$.
\end{proof}

\begin{Corollary}
\label{L210900_1}
Let $M$ be a finite dimensional $\H_n$-module. 
Then $M$ is integral if and only if ${\mathcal I}_\la M = 0$
for some $\la \in P_+$.
\end{Corollary}

\begin{proof}
If ${\mathcal I}_\la M = 0$, then
the eigenvalues of $X_1+X_1^{-1}$ on $M$ are of the form $q(i)$ for 
$i\in I$, by definition of ${\mathcal I}_\la$. 
Hence $M$ is integral in view of Lemma~\ref{L210900}. 
Conversely, suppose that $M$ is integral.
Then the minimal polynomial of $X_1+X_1^{-1}$ on $M$ is of the form 
$\prod_{i\in I}(t-q(i))^{\la_i}$ for some $\la_i \geq 0$. 
So if we set $\la = 2 \la_0 \Lambda_0 + \la_1 \Lambda_1 + \dots + 
\la_{\ell}\Lambda_{\ell} \in P_+$, we certainly have that 
the element (\ref{E210900_4}) acts as zero 
on $M$.
\end{proof}

Recall from \ref{ggp} that $\rep \H_n^\la$ denotes the category
of all finite dimensional $\H_n^\la$-modules.
Corollary~\ref{L210900_1} implies that the functors
$\pr^\la$ and $\inf^\la$ from (\ref{prind})
restrict to a well-defined adjoint pair of
functors at the level of integral representations:
\begin{equation}
\pr^\la:\rep_I{\H_n} \rightarrow \rep{\H_n^\lambda},
\qquad
\inf^\la:\rep{\H_n^\lambda} \rightarrow \rep_I{\H_n}.
\end{equation}
Let us also check at this point that induction from a parabolic subalgebra
of $\H_n$ preserves integral representations,
the analogous fact for restriction being obvious.

\vspace{1mm}
\begin{Lemma}\label{indind} Let $\mu$ be a composition of $n$ and $M$ be an integral
$\H_\mu$-module. Then, $\ind_\mu^n M$ is an integral $\H_n$-module.
\end{Lemma}

\begin{proof}
By Theorem~\ref{base}, $\ind_\mu^n M$ is spanned by elements $T_w \otimes m$
for $m \in M$, in particular it is finite dimensional.
Let 
$$
Y_j = \prod_{i \in I} (X_j+X_j^{-1} - q(i)).
$$
By Corollary~\ref{L210900_1}, it suffices to show that $Y_1^N$ annihilates
$\ind_\la^n M$ for sufficiently large $N$.
Consider $Y_1^N T_w \otimes m$ for $w \in \Sym_n, m \in M$.
We may write $T_w = T_u T_1 \dots T_k$ for $u \in 
\Sym_{2\dots n} \cong 
\Sym_{n-1}$ and $0 \leq k < n$. 
Then, $Y_1^N$ commutes with $T_u$, so we just need to consider
$Y_1^N T_1 \dots T_k \otimes m$. Now using the commutation relations,
one checks that
$Y_1^N T_1 \dots T_k \otimes m$ can be rewritten as an $\H_n$-linear
combination of elements of the form $1 \otimes Y_j^{N'} m$
for $1\leq j \leq n$ and $N - k \leq N' \leq N$.
Since $M$ is integral by assumption, 
we can choose $N$ sufficiently large
so that each such term is zero.
\end{proof}

It follows that the functors $\ind^n_\mu, \res^n_\mu$ restrict
to well-defined functors
\begin{equation}
\ind^n_\mu: \rep_I{\H_\mu} \rightarrow \rep_I{\H_n},
\qquad
\res^n_\mu: \rep_I{\H_n} \rightarrow \rep_I{\H_\mu}
\end{equation}
on integral representations. Similar remarks apply to more general
induction and restriction between nested parabolic subalgebras of
$\H_n$.

\Point{\boldmath Modules over $\A_n$} \label{sa1}
Let $i \in I$ and define
\begin{equation}
\label{Eb}
b_\pm(i)=\frac{q(i)}{2}\pm\sqrt{\frac{q(i)^2}{4}-1},
\end{equation}
i.e. the roots of the equation $x+x^{-1} = q(i)$.
Let $L(i)$ denote the vector superspace
on basis $w,w'$, where $w$ is even and $w'$ is odd,
made into an $\A_1$-module so that
$$
C_1w=w',\ C_1 w'=w,\ X_1^{\pm 1} w=b_\pm(i)w,\ X_1^{\pm 1} w'=b_\mp(i) w'.
$$
One easily checks:

\vspace{1mm}
\begin{Lemma}
\label{L290800_1} 
For each $i \in I$,
$L(i)$ is an irreducible $\A_1$-module, 
of type $\Mtype$ if $i \neq 0$ and of type $\Qtype$ if $i = 0$.
Moreover, the modules
$\{L(i)\:|\: i\in I\}$
form a complete set of pairwise non-isomorphic irreducibles 
in $\rep_I{\A_1}$. 
\end{Lemma}

\vspace{1mm}

Recall that
$\A_n\cong \A_1\otimes\dots\otimes \A_1$  ($n$ times)
as superalgebras. 
Hence, for ${\bi}=(i_1,\dots,i_n)\in I^n$, 
we can consider the irreducible $\A_n$-module
$L(i_1)\circledast\dots\circledast L(i_n)$.
By Lemma~\ref{L290800_1} and the general theory of outer tensor products
\ref{theory}, 
one obtains:

\vspace{1mm}
\begin{Lemma} \label{compset}
The $\A_n$-modules $\{L(i_1)\circledast\dots\circledast L(i_n)\:|\:\bi\in I^n\}$
form a complete set of pairwise non-isomorphic irreducible $\A_n$-modules.
Moreover, let $\gamma_0$ denote the number of $j = 1,\dots,n$
such that $i_j = 0$. Then,
$L(i_1)\circledast\dots\circledast L(i_n)$ 
is of type $\Mtype$ if $\gamma_0$ is even and
type $\Qtype$ if $\gamma_0$ is odd. Finally,
$\dim L(i_1)\circledast\dots\circledast L(i_n) 
= 2^{n - \lfloor \gamma_0/2\rfloor}.$
\end{Lemma}

\vspace{1mm}

Now let $M$ be any module in $\rep_I{\A_n}$.
For any ${\bi}\in I^n$,
let $M[{\bi}]$ be the largest submodule of $M$ all of whose 
composition factors are isomorphic to $L(i_1)\circledast\dots\circledast L(i_n)$. 
Alternatively, since each $X_k+X_k^{-1}$ acts on
$L(i_1)\circledast\dots\circledast L(i_n)$ 
by the scalar $q(i_k)$ and all the scalars $q(i)$ for $i \in I$
are distinct, we can describe $M[\bi]$ as the
simultaneous generalized eigenspace for the 
commuting operators $X_1+X_1^{-1},\dots, X_n+X_n^{-1}$ 
corresponding to eigenvalues $q(i_1),\dots,q(i_n)$,
respectively.
Hence:

\vspace{1mm}
\begin{Lemma}\label{decomp}
For any $M \in \rep_I{\A_n}$,
$\displaystyle M = \bigoplus_{\bi\in I^n} M[\bi]$ as an $\A_n$-module.
\end{Lemma}

\vspace{1mm}

We write $K(\rep_I{\A_n})$, $K(\rep_I{\H_n})$, \dots for the Grothendieck
groups of the categories $\rep_I{\A_n}, \rep_I{\H_n}$, \dots,
defined as in \ref{ggp}.
Note for an integral $\A_n$-module $M$,
knowledge of the dimensions of the spaces $M[\bi]$ for all $\bi$
is equivalent to knowing the coefficients $a_{\bi}$ when
$[M]\in K(\rep_I{\A_n})$ 
is expanded as
$$[M] = \sum_{\bi \in I^n} a_{\bi} 
[L(i_1)\circledast\dots\circledast L(i_n)]$$
in terms of the basis 
$\{[L(i_1)\circledast\dots\circledast L(i_n)]\:|\:\bi \in I^n\}$.

Now suppose instead that $M$ is an integral $\H_n$-module, so that
its restriction $\res^n_{1,\dots,1} M$ to $\A_n$ is in $\rep_I{\A_n}$.
We define the {\em formal character} of $M$ by:
\begin{equation}
\ch M = [\res^n_{1,\dots,1} M] \in K(\rep_I{\A_n}).
\end{equation}
Since the functor $\res^n_{1,\dots,1}$ is exact, 
$\ch$ induces a homomorphism
$$
\ch:K(\rep_I{\H_n}) \rightarrow K(\rep_I{\A_n})
$$
at the level of Grothendieck groups. 
We will later see that this map is
actually injective (Theorem~\ref{inj}), justifying the terminology.
Note we will occasionally consider characters of integral modules
over parabolic subalgebras $\H_\mu$ for $\mu$ a composition of $n$, 
or over the cyclotomic algebras $\H_n^\la$ for $\la \in P_+$.
The definitions are modified in these cases in obvious ways.

\vspace{1mm}
\begin{Lemma}
\label{LChari}
Let $\bi = (i_1,\dots,i_n)\in I^n$. Then
$$
\ch \ind_{1,\dots,1}^{n} L(i_1)\circledast\dots\circledast L(i_n)=\sum_{w\in \Sym_n} [L(i_{w^{-1}1})\circledast\dots\circledast L(i_{w^{-1}n})]
$$
\end{Lemma}
\begin{proof}
This follows from Theorem~\ref{TMackey} with $\mu = \nu = (1^n)$.
\end{proof}

\begin{Lemma}\label{LChar}
{\em (``Shuffle Lemma'') }
Let $n=m+k$, and let $M\in \rep_I{\H_m}$ and $K\in\rep_I{\H_k}$ 
be irreducible. Assume 
$$
\ch M=\sum_{{\bi}\in I^m}a_{\bi}[L(i_1)\circledast\dots\circledast L(i_m)],\qquad
\ch K=\sum_{{\bj}\in I^k}b_{\bj}[L(j_1)\circledast\dots\circledast L(j_k))].
$$
Then 
$$
\ch \ind_{m,k}^n M\circledast K = \sum_{{\bi}\in I^m}\sum_{{\bj}\in I^k}
a_{\bi} b_{\bj}(\sum_{\bh} L(h_1)\circledast \dots\circledast L(h_n)),
$$
where the last sum is over all $\bh=(h_1,\dots,h_n) \in I^n$ which are obtained
by shuffling $\bi$ and $\bj$,
i.e. there exist $1\leq u_1<\dots<u_m\leq n$ 
such that $(h_{u_1},\dots,h_{u_m})=(i_1,\dots,i_m)$, and $(h_1,\dots,\widehat {h}_{u_1},\dots,\widehat{h}_{u_m},\dots,h_n)=(j_1,\dots,j_k)$.
\end{Lemma}

\begin{proof}
This follows from Theorem~\ref{TMackey} with $\mu = (1^n)$ and
$\nu = (m,k)$.
\end{proof}

\Point{Central characters}\label{Defs}
Recall by Theorem~\ref{LCenter} that every element $z$ of
the center $Z(\H_n)$ of $\H_n$ 
can be 
written as a symmetric polynomial
$f(X_1+X_1^{-1},\dots,X_n+X_n^{-1})$ in the $X_k+X_k^{-1}$.
Given $\bi \in I^n$, we associate the
{\em central character}
$$
\chi_{\bi}:Z(\H_n) \rightarrow F,
\quad
f(X_1+X_1^{-1},\dots,X_n+X_n^{-1}) \mapsto f(q(i_1),\dots,q(i_n)).
$$
Consider the natural left action of $\Sym_n$ on $I^n$ by place permutation.
We write $\bi \sim \bj$ if $\bi,\bj$ lie in the same orbit.
The following lemma follows immediately from the fact that
the $q(i)$ are distinct as $i$ runs over the index set $I$.

\vspace{1mm}
\begin{Lemma}\label{ccc} For $\bi,\bj \in I^n$, 
$\chi_{\bi} = \chi_{\bj}$ if and only
if $\bi \sim \bj$.
\end{Lemma}

\vspace{1mm}

Given $\bi \in I^n$, we define its {\em weight}
$\wt(\bi) \in P$ by
\begin{equation}\label{wtdef}
\wt(\bi) = \sum_{i \in I} \gamma_i \alpha_i
\qquad\hbox{where}\qquad
\gamma_i = \sharp\{j=1,\dots,n\:|\:i_j = i\}.
\end{equation}
So $\wt(\bi)$ is an element of the set $\Gamma_n$ of 
non-negative integral linear combinations 
$\gamma = \sum_{i\in I} \gamma_i \alpha_i$ of the simple roots 
such that $\sum_{i \in I} \gamma_i = n$.
Obviously, the $\Sym_n$-orbit of $\bi$ is uniquely determined by its
weight, so we obtain a labelling of the orbits 
of $\Sym_n$ on $I^n$ by the elements of $\Gamma_n$.
We will also use the notation
$\chi_\gamma$ for the central character $\chi_{\bi}$
where $\bi$ is any element of $I^n$ with $\wt(\bi) = \gamma$. 
So $\chi_{\bi} = \chi_{\wt(\bi)}$.

Now let $M$ be an integral $\H_n$-module and $\gamma = \sum_{i \in I}
\gamma_i \alpha_i \in \Gamma_n$.
We let $M[\gamma]$ denote the generalized 
eigenspace of $M$ over $Z(\H_n)$ 
that corresponds to the central character $\chi_{\gamma}$, i.e.
$$
M[\gamma] = \{m \in M\:|\: (z - \chi_\gamma(z))^km = 0
\hbox{ for all $z \in Z(\H_n)$ and $k \gg 0$}\}.
$$
Observe this is an $\H_n$-submodule of $M$.
Now, for any $\bi \in I^n$ with $\wt(\bi) = \gamma$,
$Z(\H_n)$ acts on
$L(i_1)\circledast \dots\circledast L(i_n)$ via the central
character $\chi_\gamma$. So applying Lemma~\ref{ccc}, we see that
$$
M[\gamma] = \bigoplus_{\bi\:\rm{with}\:\wt(\bi)=\gamma} M[\bi],
$$
recalling the decomposition of $M$ as an $\A_n$-module from Lemma~\ref{decomp}.
Therefore:

\vspace{1mm}
\begin{Lemma}
\label{L070900}
Any integral $\H_n$-module
$M$ decomposes as
$$
M = \bigoplus_{\gamma \in \Gamma_n} M[\gamma]
$$
as an $\H_n$-module.
\end{Lemma}

\vspace{1mm}

Thus the $\{\chi_{\gamma}\:|\:\gamma \in \Gamma_n\}$
exhaust the possible central characters that can arise
in an integral $\H_n$-module, while
Lemma~\ref{LChari} shows that every such central character
does arise in some integral $\H_n$-module.

Let us write $\rep_\gamma \H_n$ for the full subcategory of
$\rep_I{\H_n}$ consisting of all modules $M$ with $M[\gamma] = M$.
Then, Lemma~\ref{L070900} implies that there is an equivalence of categories
\begin{equation}\label{block1}
\rep_I{\H_n} \cong
\bigoplus_{\gamma \in \Gamma_n} \rep_{\gamma} \H_n.
\end{equation}
We say that $\rep_{\gamma} \H_n$ is the {\em block} of $\rep_I{\H_n}$
corresponding to the central character $\chi_\gamma$.
In particular, if $M\neq 0$ is indecomposable then
$M$ belongs to $\rep_{\gamma} \H_n$, i.e. $M = M[\gamma]$, 
for a unique $\gamma \in \Gamma_n$.

We can extend some of these notions to $\H_n^\la$-modules,
for $\la \in P_+$.
In particular, if $M \in \rep{\H_n^\la}$, 
we also write
$M[\gamma]$ for the summand $M[\gamma]$ of $M$ defined by 
first viewing $M$ as an $\H_n$-module by inflation.
Also write $\rep_{\gamma} \H_n^\la$ for the full subcategory
of $\rep{\H_n^\la}$ consisting of the modules $M$ with $M = M[\gamma]$.
Thus we also have a decomposition
\begin{equation}\label{block2}
\rep{\H_n^\la} \cong \bigoplus_{\gamma \in \Gamma_n} \rep_{\gamma} \H_n^\la
\end{equation}
induced by (\ref{block1}).
Note though that we should not yet
refer to $\rep_{\gamma} \H_n^\la$ as a block of
$\rep{\H_n^\la}$: the center
of $\H_n^\la$ may be larger than the image of the
center of $\H_n$, so we cannot yet assert that
$Z(\H_n^\la)$ acts on $M[\gamma]$ by a single central character.
Also we no longer know precisely which
$\gamma\in\Gamma_n$ have the property that $\rep_{\gamma} 
\H_n^\la$ is non-trivial.
These questions will be settled in \ref{blockssect}.

\Point{Kato's theorem}
Let $i\in I$. 
Introduce the {\em principal series module}
\begin{equation}\label{mkato}
L(i^n):=\ind_{1,\dots,1}^{n} L(i)\circledast \dots\circledast L(i).
\end{equation}
By Lemma~\ref{LChari}, we know immediately that
$\ch L(i^n) = n! [L(i)\circledast \dots \circledast L(i)],$
hence $L(i^n)$ belongs to the block $\rep_{n \alpha_i} \H_n$.
In particular, for each $k = 1,\dots,n$,
the only eigenvalue of the element 
$X_k+X_k^{-1}$ on $L(i^n)$ is $q(i)$. 
In addition, if $i=0$, then the only eigenvalue of each $X_k$ is $1$. 

\vspace{1mm}
\begin{Lemma}
\label{L290800_2}
Let $n\geq 2$, $1\leq j<n$, $i\in I-\{0\}$, and $v\in 
L(i)\boxtimes\dots\boxtimes L(i)$ ($n$ copies). Then 
$X_j^{-1}(1-C_jC_{j+1})v\neq (1-C_jC_{j+1})X_{j+1}v.$
\end{Lemma}
\begin{proof}
The elements of $\A_n$ which are involved in the inequality act only on 
the positions $j$ and $j+1$ in the tensor product. 
So we may assume that $n=2$ and $j=1$. Let
$$
v=a w\otimes w + b w\otimes w'+c w'\otimes w+d w'\otimes w'
$$
for $a,b,c,d\in F$. Then 
\begin{multline*}
X_1^{-1}(1-C_1C_{2})v = (b_-(i)a+ b_-(i)d)w\otimes w + (b_-(i)b +b_-(i)c) w\otimes w' \\
+ (b_+(i)c- b_+(i)b)  w'\otimes w  
+(b_+(i)d - b_+(i)a)  w'\otimes w',
\end{multline*}
\begin{multline*}
\:\:\;(1-C_1C_{2})X_2v = (b_+(i)a+ b_-(i)d) w\otimes w + (b_-(i)b + b_+(i)c) w\otimes w' \\
+ (b_+(i)c - b_-(i)b)  w'\otimes w + 
(b_-(i)d - b_+(i)a)  w'\otimes  w'.
\end{multline*}
Now the lemma follows from the inequality $b_-(i)\neq b_+(i)$ for $i\neq 0$. 
\end{proof}

\begin{Lemma}
\label{L050900_1}
Let $i\in I$. Set $L=L(i)\circledast \dots\circledast L(i)$, so 
$L(i^n)=\H_n \otimes_{\A_n} L$.
\begin{enumerate}
\item[(i)] If $i\neq 0$, the common $q(i)$-eigenspace of the operators 
$X_1+X_1^{-1},\dots,X_{n-1}+X_{n-1}^{-1}$ on $L(i^n)$
is precisely $1\otimes L$, which is contained in the $q(i)$-eigenspace of
$X_n+X_n^{-1}$ too.
Moreover, all Jordan blocks of $X_1+X_1^{-1}$ on $L(i^n)$ are of size $n$. 
\item[(ii)] If $i=0$, the common $1$-eigenspace of the operators 
$X_1,\dots,X_{n-1}$ on $L(i^n)$ is precisely $1\otimes L$, which is contained
in the $1$-eigenspace of $X_n$ too.
Moreover, all Jordan blocks of $X_1$ on $L(i^n)$ are of size $n$.
\end{enumerate}
\end{Lemma}

\begin{proof}
We prove (i), (ii) being similar.
Note $L(i^n) = \bigoplus_{x \in \Sym_n} 
T_x \otimes L$, since by Theorem~\ref{base} we know that
$\H_n$ is a free right $\A_n$-module on basis $\{T_x\:|\:x \in \Sym_n\}$.

We first show that the eigenspace of $X_1+X_1^{-1}$ is a sum of the 
subpaces of the form $T_y\otimes L$, where $y\in \Sym_{2\dots n}\cong 
\Sym_{n-1}$ is the subgroup of $\Sym_n$ generated by $s_2,\dots,s_{n-1}$.
Well, any $T_x$ can be written as $T_y T_1T_2\dots T_j$ 
for some $y \in \Sym_{2\dots n}$ and $0\leq j< n$. 
Note
$$
(X_{j+1}+X_{j+1}^{-1}-q(i)) v = 0
$$
for any $v \in L$, by definition of $L$.
Now the defining relations of $\H_n$ especially (\ref{E290900_3}),  
(\ref{E290800_2}) imply
\begin{multline*}
(X_1+X_1^{-1}-q(i))T_y T_1T_2\dots T_j\otimes v
=\\\xi T_y T_1\dots T_{j-1}\otimes (X_j^{-1}(1-C_jC_{j+1})-(1-C_jC_{j+1})X_{j+1})v + (*),
\end{multline*}
where $(*)$ stands for a sum of terms which belong to subspaces of the form 
$T_{y'} T_1 \dots T_k\otimes L$ for 
$y' \in \Sym_{2\dots n}$ and $0 \leq k < j-1$.

Now assume that a linear combination 
$$
z:=\sum_{y\in \Sym_{2\dots n}}\sum_{0 \leq j < n} \sum_{v \in L} 
c_{y,j,v} T_y T_1T_2\dots T_j\otimes v
$$
is an eigenvector for $X_1+X_1^{-1}$. 
Then it must be annihilated by $X_1+X_1^{-1}-q(i)$. 
Choose the maximal $j$ for which the coefficient $c_{y,j,v}$ is non-zero, 
and for this $j$ choose the maximal (with respect to the Bruhat order) 
$y$ such that $c_{y,j,v}$ is non-zero. Then the calculation above and 
Lemma~\ref{L290800_2} show that $(X_1+X_1^{-1}-q(i))z\neq 0$ unless $j=0$. 
This proves our claim on the eigenspace of $X_1+X_1^{-1}$. 

Now apply the same 
argument to see that the common eigenspace of $X_1+X_1^{-1}$ and $X_2+X_2^{-1}$
is spanned by $T_y\otimes L$ for $y\in \Sym_{3\dots n}$, and so on, yielding 
the first claim in (i). 
Finally, define
$$
V(m):=\{z\in L(i^n)\mid (X_1+X_1^{-1} - q(i))^m z=0\}.
$$
It follows by induction from the calculation above and 
Lemma~\ref{L290800_2} that 
$$
V(m)=\operatorname{span}\{
T_y T_1T_2\dots T_j\otimes v\mid y\in \Sym_{2\dots n},\ j< m,\ v\in L\},
$$
giving the second claim.
\end{proof}

Now we are ready to prove the main theorem giving the structure of
the principal series module $L(i^n)$, compare
\cite{Kato}.

\vspace{1mm}
\begin{Theorem}\label{TKato} Let $i\in I$
and $\mu=(\mu_1,\dots,\mu_u)$ be a composition of $n$. 
\begin{enumerate}
\item[(i)]
$L(i^n)$ is irreducible of the same type as
$L(i)\circledast \dots\circledast L(i)$, i.e. type $\Mtype$
if either $i \neq 0$ or $i = 0$ and $n$ is even,
type $\Qtype$ otherwise, and it is the only irreducible module in its
block.
\item[(ii)]
All composition factors of 
$\res^n_\mu L(i^n)$ 
are isomorphic to $L(i^{\mu_1})\circledast \dots\circledast L(i^{\mu_u})$,
and $\soc \res^n_\mu L(i^n)$
is irreducible.
\item[(iii)]
$\soc \res^n_{n-1} L(i^n) \cong
\res^{n-1,1}_{n-1} L(i^{n-1}) \circledast L(i)$.
\end{enumerate}
\end{Theorem}

\begin{proof}
Denote $L(i)\circledast \dots\circledast L(i)$ by $L$.

(i) Let $M$ be a non-zero $\H_n$-submodule of $L(i^n)$.
Then, $\res^n_{1,\dots,1} M$ must contain an $\A_1$-submodule $N$ isomorphic
to $L$.
But the commuting operators $X_1+X_1^{-1},\dots, X_n+X_n^{-1}$ 
(or $X_1,\dots,X_n$ if $i=0$) act on $L$ as scalars, giving that
$N$ is contained in their common eigenspace on $L(i^n)$.
But by Lemma~\ref{L050900_1}, this implies that $N = 1 \otimes L$.
This shows that $M$ contains $1 \otimes L$, but this generates the whole of
$L(i^n)$ over $\H_n$. So $M=L(i^n)$. 

To see that the type of $L(i^n)$ is the same as the type of $L$,
the functor $\ind_{1,\dots,1}^n$ determines a map
$$
\End_{\A_n}(L) \rightarrow \End_{\H_n}(L(i^n)).
$$
We just need to see that this is an isomorphism, which we do by 
constructing the inverse map.
Let $f \in \End_{\H_n}(L(i^n))$.
Then $f$ leaves $1 \otimes L$ invariant by Lemma~\ref{L050900_1}, so $f$ restricts
to an $\A_n$-endomorphism $\bar f$ of $L$.

Finally, to see that $L(i^n)$ is the only irreducible in its block,
we have already observed using Lemma~\ref{LChari}
that $\ch L(i^n) = n! [L(i)\circledast\dots\circledast L(i)]$. 
Hence all composition factors of
$\res^n_{1,\dots,1} L(i^n)$ are isomorphic to 
$L(i)\circledast\dots\circledast L(i)$.
Now apply Frobenius reciprocity and the fact just proved
that $L(i^n)$ is irreducible. 

(ii) That all composition factors of $\res^n_\mu L(i^n)$
are isomorphic to $L(i^{\mu_1}) \circledast\dots\circledast L(i^{\mu_u})$
follows from the parabolic analogue of (i).
To see that $\soc \res^n_\mu L(i^n)$ is simple,
note that the submodule $\H_\mu\otimes L$ of 
$\res_{\H_\mu} L(i^n)$ is isomorphic to 
$L(i^{\mu_1})\circledast\dots\circledast L(i^{\mu_l})$. 
This module is irreducible, and so it is contained in the socle. 
Conversely, let $M$ be an irreducible $\H_\mu$-submodule of $L(i^n)$.
Then using Lemma~\ref{L050900_1} as in the proof of (i),
we see that $M$ must contain $1\otimes L$, hence $\H_\mu\otimes L$. 

(iii) By (ii), $L(i^n)$ has a unique $\H_{n-1,1}$-submodule
isomorphic to $L(i^{n-1}) \circledast L(i)$, namely, 
$\H_{n-1,1} \otimes L$.
Since this is completely reducible on restriction to $\H_{n-1}$,
it follows that $\H_{n-1,1} \otimes L \subseteq \soc \res^{n}_{n-1} L(i^n)$.
Conversely, take any irreducible $\H_{n-1}$-submodule $M$ of $L(i^n)$.
The common eigenspace of $X_1+X_1^{-1}, \dots, X_{n-1}+X_{n-1}^{-1}$
(resp. $X_1,\dots,X_{n-1}$ if $i = 0$) on $M$
must lie in $1 \otimes L$ by Lemma~\ref{L050900_1}.
Hence, $M \subseteq \H_{n-1,1} \otimes L$ which completes the proof.
\end{proof}

\Point{Covering modules}\label{cov}
Fix $i \in I$ and $n \geq 1$ throughout the subsection.
We will construct for each $m \geq 1$ 
an $\H_n$-module $L_m(i^n)$ with
irreducible cosocle isomorphic to $L(i^n)$.
Let ${\mathcal J}(i^n)$ denote the annihilator in $\H_n$
of $L(i^n)$. 
Introduce the quotient superalgebra
\begin{equation}\label{rmdef}
\R_m(i^n) := \H_n / {\mathcal J}(i^n)^m
\end{equation}
for each $m \geq 1$.
One checks that ${\mathcal J}(i^n)$ contains
$(X_k + X_k^{-1} - q(i))^{n!}$ for each $k = 1,\dots,n$. It follows
easily from this that each superalgebra 
$\R_m(i^n)$ is finite dimensional.
Moreover, by Theorem~\ref{TKato}, $L(i^n)$
is the unique irreducible $\R_m(i^n)$-module up to isomorphism.

Let $L_m(i^n)$ denote a projective cover of $L(i^n)$ in the category
$\mod{\R_m(i^n)}$.
For convenience, we also define $L_0(i^n) = \R_0(i^n) = 0$.
Note we know the dimension of $L(i^n)$ from Lemma~\ref{compset}, and moreover 
$L(i^n)$ is of type $\Qtype$ if $i = 0$ and $n$ is odd, type $\Mtype$
otherwise. Using this and the 
general theory of finite dimensional superalgebras, one shows:

\vspace{1mm}
\begin{Lemma}\label{lmir} 
For each $m \geq 1$,
$$
\R_m(i^n) \simeq
\left\{
\begin{array}{ll}
(L_m(i^n)\oplus \Pi L_m(i^n))^{\oplus 2^{n-1} (n!)}
&\hbox{if $i \neq 0$,}\\
(L_m(i^n) \oplus \Pi L_m(i^n))^{\oplus 2^{(n-2)/2} (n!)}
&\hbox{if $i = 0$ and $n$ is even,}\\
L_m(i^n)^{\oplus 2^{(n-1)/2} (n!)}
&\hbox{if $i = 0$ and $n$ is odd,}
\end{array}
\right.
$$
as left $\H_n$-modules.
Moreover, 
$L_m(i^n)$ admits an odd involution
if and only if $i = 0$ and $n$ is odd.
\end{Lemma}

\vspace{1mm}

There are obvious surjections
\begin{align}
&\R_1(i^n) \twoheadleftarrow\label{rseq}
\R_2(i^n) \twoheadleftarrow\dots.\\
\intertext{In the case $m = 1$, we certainly
have that $L_1(i^n) \cong L(i^n)$; let us assume by the choice of
$L_1(i^n)$ that in fact $L_1(i^n) = L(i^n)$.
Then we can choose the $L_m(i^n)$ for $m > 1$ 
so that the the maps (\ref{rseq}) induce even maps}
&L_1(i^n) \twoheadleftarrow
L_2(i^n) \twoheadleftarrow\dots.\label{lseq}
\end{align}
Moreover, in case $i = 0$ and $n$ is odd, we can choose 
the odd involutions 
\begin{equation}\label{thetaj}
\theta_m:L_m(i^n) \rightarrow L_m(i^n)
\end{equation} 
given by Lemma~\ref{lmir}
in such a way that they are compatible with the maps in (\ref{lseq}).

The significance of the $\H_n$-modules $\R_m(i^n)$ 
is explained by the following lemma:

\vspace{1mm}
\begin{Lemma}\label{nastystab}
Let $M$ be an $\H_n$-module annihilated by ${\mathcal J(i^n)}^k$ for some $k$.
Then, there is a natural isomorphism
$$
\hom_{\H_n}(\R_m(i^n), M) \stackrel{\sim}{\longrightarrow}
M
$$
for all $m \geq k$.
\end{Lemma}

\begin{proof}
The assumption implies that 
$M$ is the inflation of an $\R_m(i^n)$-module.
So
$$
\hom_{\H_n}(\R_m(i^n), M)
\simeq \hom_{\R_m(i^n)}(\R_m(i^n),M)
\simeq M,
$$
all isomorphisms being the natural ones.
\end{proof}

Let us finally consider the most important case $n = 1$
in more detail. In this case, one easily checks that
the ideal ${\mathcal J}(i)^m$
is generated by $(X_1 + X_1^{-1} - q(i))^m$
if $i \neq 0$ or $(X_1 - 1)^m$ if $i = 0$.
It follows that 
\begin{equation}\label{dimm}
\dim \R_m(i) = 
\left\{
\begin{array}{ll}
4m&\hbox{if $i \neq 0$,}\\
2m&\hbox{if $i = 0$.}
\end{array}
\right.
\end{equation}
Moreover, Lemma~\ref{lmir} shows in this case that
\begin{equation}\label{zoo}
\R_m(i) \simeq \left\{
\begin{array}{ll}
L_m(i) \oplus \Pi L_m(i)&\hbox{if $i \neq 0$,}\\
L_m(i)&\hbox{if $i = 0$.}\\
\end{array}\right.
\end{equation}
Hence, $\dim L_m(i) = 2m$ in either case.
Using this, it follows easily that $L_m(i)$ can be described alternatively
as the vector superspace on basis 
$w_1,\dots,w_m$, 
$w_1',\dots, w_m',$
where each $w_k$ is even and each $w_k'$ is odd,
with $\H_1$-module structure uniquely determined by
$$
X_1 w_k=b_+(i)w_k+w_{k+1},\qquad
C_1 w_k=w_k'
$$
for each $k = 1,\dots,m$,
interpreting $w_{m+1}$ as $0$.
Using this explicit description, one now
checks routinely that $L_m(i)$ is uniserial 
with $m$ composition factors all $\simeq L(i)$.

We can also describe the map $L_m(i) \twoheadleftarrow L_{m+1}(i)$
from (\ref{lseq}) explicitly: it is
the identity
on $w_1,\dots,w_m,w_1',\dots,w_m'$ but 
maps $w_{m+1}$ and $w_{m+1}'$
to zero.
Also, the map $\theta_m$ from (\ref{thetaj}) can be chosen so that
$$
w_k
\mapsto
\sqrt{-1}w_k', 
\qquad
w_k' \mapsto -\sqrt{-1}w_k
$$
for each $k = 1,\dots,m$.

\Point{Modifications in the degenerate case}
For $i \in I$, the definition (\ref{resdef})
of $q(i)$ becomes
\begin{equation}\label{newresdef}
q(i) = i(i+1) \in F
\end{equation}
for each $i \in I$.
For $\la \in P_+$,
the quotient superalgebra $\H_n^\la$, 
is defined to be 
the quotient $\H_n / {\mathcal I}_\la$ of the affine Sergeev
superalgebra $\H_n$ by the two-sided ideal ${\mathcal I}_\la$  generated
by
\begin{equation}
x_1^{\langle h_0,\la \rangle} \prod_{i = 1}^\ell (x_1^2 - q(i))^{\langle h_i,\la\rangle}.
\end{equation}
The basis theorem for $\H_n^\la$
gives that
$\H_n^\la$ has a basis given by the images of the elements
\begin{equation}
\{
x^\alpha c^\beta w\:|\:\alpha \in \Z_{\geq 0}^n\hbox{ with }0 \leq \alpha_i < \langle c,\la \rangle, \beta \in \Z_2^n, w \in \Sym_n\}.
\end{equation}
The definition of the category $\rep_I{\H_n}$ of
integral representations is modified in the appropriate way
to ensure that the integral representations of $\H_n$ are precisely the
inflations of finite dimensional representations of $\H_n^\la$ for 
$\la \in P_+$. To be precise, an integral $\A_n$-module 
now means
a finite dimensional $\A_n$-module
such that the eigenvalues of all $x_1^2, \dots, x_n^2$ are of the form
$q(i)$ for $i \in I$.
The appropriate analogue of Lemma~\ref{L210900}, 
involving elements $x_j^2$ now of course,
is proved using the elements 
$$
\Phi_j = s_j(x_j^2 - x_{j+1}^2) + (x_j+x_{j+1}) + c_j c_{j+1}(x_j - x_{j+1})
$$ from 
\cite[(3.4)]{Naz}.
Proofs of 
the basic properties of $\Phi_j$, analogous to those in \ref{nazelts},
can be found in
\cite[Prop. 3.2]{Naz} and at the end of \cite[$\S$4]{Naz}.

The module $L(i)$ in \ref{sa1}
is now defined to be the $\A_1$-module on basis $w,w'$ with action
$$
cw = w', \qquad cw' = w,\qquad xw = \sqrt{q(i)} w, \qquad xw' = -\sqrt{q(i)}w'.
$$
The remaining definitions go through more or less unchanged:
for example for an integral $\H_n$-module $M$,
$M[\bi]$ is the simultaneous generalized eigenspace of the
operators $x_1^2,\dots,x_n^2$ corresponding to eigenvalues
$q(i_1),\dots,q(i_n)$ respectively.

Kato's theorem (Theorem~\ref{TKato}) is the same, as is 
Lemma~\ref{L050900_1} 
when $X_i+X_i^{-1}$ is replaced by $x_i^2$ and eigenvalue $1$ for $X_i$
is replaced by
eigenvalue $0$ for $x_i$ in the usual way.
Note the main technical fact needed in the proof of
Lemma~\ref{L050900_1} is the following:
for $i \neq 0$, and 
every $0 \neq v \in L(i)\boxtimes\dots\boxtimes L(i)$ ($n$ copies),
$$
(x_j(1-c_jc_{j+1}) + (1-c_jc_{j+1})x_{j+1})v \neq 0
$$
for each $j=1,\dots,n-1$.

\ifbook@\pagebreak\fi

\section{Crystal operators}

\Point{\boldmath Multiplicity-free socles}\label{socs}
The arguments in this subsection are based on \cite{GV}.
Let $M \in \rep_I{\H_n}$ and $i \in I$.
Define
$\De_i M$ to be the generalized $q(i)$-eigenspace of $X_n+X_n^{-1}$ on $M$.
Alternatively,
$$
\De_i M = \bigoplus_{\bi \in I^n,\ i_n = i} M[\bi],
$$
recalling the decomposition from Lemma~\ref{decomp}.
Note since $X_n+X_n^{-1}$ is central in the parabolic subalgebra
$\H_{n-1,1}$ of $\H_n$, $\De_i M$ is invariant under this subalgebra.
So in fact, $\De_i$ can be viewed as an exact functor
\begin{equation}
\De_i:\rep_I{\H_n} \rightarrow \rep_I{\H_{n-1,1}},
\end{equation}
being defined on morphisms simply as restriction.
Clearly, there is an isomorphism of functors
$$
\res^{n}_{n-1,1} \simeq \De_0 \oplus \De_1 \oplus \dots \oplus \De_{\ell}.
$$
Slightly more generally, given $m \geq 0$,
define 
\begin{equation}\label{gende}
\De_{i^m}:\rep_I{\H_{n}} \rightarrow \rep_I{\H_{n-m,m}}
\end{equation}
so that $\De_{i^m} M$
is the simultaneous generalized $q(i)$-eigenspace of
of the commuting operators $X_k+X_k^{-1}$ for $k = n-m+1,\dots,n$.
In view of Theorem~\ref{TKato}(i), $\De_{i^m} M$ can also 
be characterized as the largest submodule of $\res^n_{n-m,m} M$ 
all of whose composition factors are of the form $N\circledast L(i^m)$
for irreducible $N \in \rep_I{\H_{n-m}}$. 

The definition of $\De_{i^m}$ implies functorial isomorphisms
\begin{equation*}
\hom_{\H_{n-m,m}}(N\boxtimes L(i^m),\De_{i^m} M)
\simeq \hom_{\H_n}(\ind_{n-m,m}^n N\boxtimes L(i^m), M)
\end{equation*}
for $N\in \rep_I{\H_{n-m}}$, $M\in\rep_I{\H_n}$. For irreducible $N$ this 
immediately imples:
\begin{equation}
\label{E090900_1}
\hom_{\H_{n-m,m}}(N\circledast L(i^m),\De_{i^m} M)
\cong \hom_{\H_n}(\ind_{n-m,m}^n N\circledast L(i^m), M).
\end{equation}
Also from definitions we get:

\vspace{1mm}
\begin{Lemma}\label{L310800_2}
Let $M\in\rep_I{\H_n}$ with
$\ch M=\sum_{{\bi}\in I^n}a_{\bi}[L(i_1)\circledast\dots\circledast L(i_n)].$
Then we have that
$\ch \De_{i^m} M=
\sum_{\bj} a_{\bj}[L(j_1)\circledast\dots\circledast L(j_n)],
$
summing over all $\bj \in I^n$ with
$j_{n-m+1} =\dots=j_n = i$.
\end{Lemma}

\vspace{1mm}

Now for $i \in I$ and $M \in \rep_I{\H_n}$, 
define
\begin{equation}
\label{E280900}
\eps_i(M)=\max\{m\geq 0\:|\: \De_{i^m} M\neq 0\}. 
\end{equation}
Note Lemma~\ref{L310800_2} shows that
$\eps_i(M)$ can be worked out just from knowledge
of the character $\ch M$.

\vspace{1mm}
\begin{Lemma}\label{sick}
Let $M\in \rep_I{\H_n}$ be irreducible, $i\in I$, $\eps=\eps_i(M)$.
If $N \circledast L(i^m)$ is an irreducible
submodule of $\De_{i^m} M$ for some $0 \leq m \leq \eps$,
then $\eps_i(N) =\eps-m$. 
\end{Lemma}

\begin{proof}
The definitions imply immediately that
$\eps_i(N) \leq \eps - m$. For the reverse inequality,
(\ref{E090900_1}) and the irreducibility of $M$ gives that
$M$ is a  quotient of
$\ind_{n-m,m}^n N \circledast L(i^m)$.
So applying the exact functor $\De_{i^\eps}$, 
$\De_{i^\eps} M\neq 0$ is a quotient of
$\De_{i^\eps} (\ind_{n-m,m}^n N \circledast L(i^m))$.
In particular, $\De_{i^\eps} (\ind_{n-m,m}^n N\circledast L(i^m))\neq 0$. 
Now one gets that $\eps_i(N) \geq \eps - m$ applying
the Shuffle Lemma (Lemma~\ref{LChar}) and Lemma~\ref{L310800_2}.
\end{proof}

\begin{Lemma}\label{L010900}
Let $m \geq 0$, $i \in I$ and $N$ be an irreducible module
in $\rep_I{\H_n}$ with $\eps_i(N) = 0$.
Set $M := \ind_{n,m}^{n+m} N \circledast L(i^m)$.
Then:
\begin{enumerate}
\item[(i)] $\De_{i^m} M \cong N \circledast L(i^m)$;
\item[(ii)] $\cosoc M$ is irreducible with $\eps_i(\cosoc M) = m$;
\item[(iii)] all other composition factors $L$ of $M$
have $\eps_i(L) < m$.
\end{enumerate}
\end{Lemma}

\begin{proof}
(i) Clearly a copy of 
$N \circledast L(i^m)$ appears in $\De_{i^m} M$.
But by the Shuffle Lemma and Lemma~\ref{L310800_2},
$\dim \De_{i^m} M = \dim N \circledast L(i^m)$, hence
$\De_{i^m} M \cong N \circledast L(i^m).$

(ii) By (\ref{E090900_1}), a copy of $N \circledast L(i^m)$
appears in $\De_{i^m} Q$ for any non-zero quotient $Q$ of $M$,
in particular for any constituent of $\cosoc M$.
But by (i), $N \circledast L(i^m)$ only appears once in
$\De_{i^m} M$, hence $\cosoc M$ must be irreducible.

(iii) We have shown that $\De_{i^m} M = \De_{i^m} (\cosoc M)$.
Hence, $\De_{i^m} L = 0$ for any other composition factor
of $M$ by exactness of $\De_{i^m}$.
\end{proof}

\begin{Lemma}
\label{L310800_1}
Let $M\in \rep_I{\H_n}$ be irreducible, $i\in I$, $\eps=\eps_i(M)$.
Then, $\De_{i^\eps} M$ is 
isomorphic to $N \circledast L(i^\eps)$ for some
irreducible $\H_{n-\eps}$-module $N$ with $\eps_i(N) = 0.$
\end{Lemma}

\begin{proof}
Pick an irreducible submodule of $\De_{i^\eps} M$. In
view of Theorem~\ref{TKato}(i), it
must be of the form $N\circledast L(i^\eps)$ for some irreducible
$\H_{n-\eps}$-module $N$.
Moreover, $\eps_i(N) = 0$ by Lemma~\ref{sick}. 
By (\ref{E090900_1}) and the irreducibility of $M$, 
$M$ is a quotient of  $\ind_{n-\eps,\eps}^n N\circledast L(i^\eps)$.
Hence, $\De_{i^\eps} M$ is a quotient of
$\De_{i^\eps} \ind_{n-\eps,\eps}^n N\circledast L(i^\eps)$.
But this is isomorphic to $N \circledast L(i^\eps)$
by Lemma~\ref{L010900}(i).
This shows that
$\De_{i^\eps} M \cong N \circledast L(i^\eps).$
\end{proof}

\begin{Lemma}
\label{L090900_1}
Let $m\geq 0$, 
$i\in I$ and $N$ be an irreducible module in $\rep_I{\H_n}$. 
Set $M:=\ind_{n,m}^{n+m}(N \circledast L(i^m))$.
Then, $\cosoc M$ is irreducible with $\eps_i(\cosoc M)=
\eps_i(N)+m$, and all other composition factors $L$ of 
$M$ have $\eps_i(L)<\eps_i(N)+m$. 
\end{Lemma}

\begin{proof}
Let $\eps = \eps_i(N)$.
By Lemma~\ref{L310800_1}, we have that
$\De_{i^\eps} N=K\circledast L(i^\eps)$
for an irreducible $K\in\rep_I{\H_{n-\eps}}$ with 
$\eps_i(K)=0$. 
By (\ref{E090900_1}) and the irreducibility of $N$, $N$ is a quotient of 
$\ind_{n-\eps,\eps}^n K\circledast L(i^\eps)$. 
So the transitivity of induction implies that 
$\ind_{n,m}^{n+m} N\circledast L(i^m)$ is a quotient of 
$\ind_{n-\eps,\eps+m}^{n+m}K\circledast L(i^{\eps+m})$. 
Now everything follows from 
Lemma~\ref{L010900}.
\end{proof}

\begin{Theorem}\label{T020900}
Let $M$ be an irreducible module in $\rep_I{\H_n}$ and $i \in I$.
Then, for any $0 \leq m \leq \eps_i(M)$,
$\soc \De_{i^m} M$ is an
irreducible $\H_{n-m,m}$-module of the same type as $M$,
and is isomorphic to $L \circledast L(i^m)$
for some irreducible $\H_{n-m}$-module $L$ with $\eps_i(L) =
\eps_i(M) - m$.
\end{Theorem}

\begin{proof}
Let $\eps = \eps_i(M)$.
Suppose $K_1 \circledast L(i^m) \oplus K_2 \circledast L(i^m)$ is a submodule
of $\De_{i^m} M$, for irreducibles $K_1, K_2$.
By Lemma~\ref{sick}, we have $\eps_i(K_j) = \eps - m$, hence
$\De_{i^{m-\eps}} K_j \cong K_j' \circledast L(i^{m-\eps})$ for some irreducible $K_j'$, for each $j = 1,2$.
This shows that $\res^{n - \eps, \eps}_{n - \eps, 1,\dots,1} \De_{i^\eps} M$ contains $K_1' \circledast L(i)^{\circledast \eps} \oplus K_2' \circledast
L(i)^{\circledast \eps}$ as a submodule. But according to Lemma~\ref{L310800_1}
and Theorem~\ref{TKato}, $\res^{n-\eps, \eps}_{n-\eps,1,\dots,1} \De_{i^\eps} M$ has irreducible socle, so this is a contradiction.
Hence, $\soc \De_{i^m} M$ is irreducible.
Now certainly 
$\soc \De_{i^m} M \cong L \circledast L(i^m)$ for some
irreducible $\H_{n-m}$-module $L$, by Theorem~\ref{TKato}(i).
It remains to show that $L \circledast L(i^m)$ has the same
type as $M$.
Note by Lemma~\ref{L090900_1},
$\ind_{n-m,m}^{n} L\circledast L(i^m)$ has irreducible cosocle, 
necessarily isomorphic to $M$ by Frobenius reciprocity.
So applying (\ref{E090900_1})
we have that
\begin{align*}
\End_{\H_{n-m,m}}(L\circledast L(i^m))
&\simeq
\hom_{\H_{n-m,m}}(L\circledast L(i^m), \De_{i^m} M)\\
&\simeq
\hom_{\H_{n}}(\ind_{n-m,m}^{n} L\circledast L(i^m), M)
\simeq
\End_{\H_{n}}(M),
\end{align*}
which implies the statement concerning types.
\end{proof}

The theorem has the following consequence:

\vspace{1mm}
\begin{Corollary}
For irreducible $M \in \rep_I \H_n$,
the socle of $\res^n_{n-1,1} M$ is multiplicity-free.
\end{Corollary}

\vspace{1mm}

We can also apply the theorem to study $\res^n_{n-1} M$, meaning the
restriction of $M$ to the subalgebra $\H_{n-1} \subset \H_n$,
see (\ref{blah}).

\vspace{1mm}
\begin{Corollary}\label{karg}
For an irreducible $M \in \rep_I \H_n$ with $\eps_i(M) > 0$, 
$$
\soc \res^{n-1,1}_{n-1} \circ \De_i (M)
\simeq
\left\{
\begin{array}{ll}
L\oplus \Pi L&\hbox{if $M$ is of type $\Qtype$ or $i \neq 0$,}\\
L&\hbox{if $M$ is of type $\Mtype$ and $i = 0$,}\\
\end{array}
\right.
$$
for some irreducible $\H_{n-1}$-module $L$ of the same type 
as $M$ if $i \neq 0$ and of the opposite type to $M$ if $i = 0$.
\end{Corollary}

\begin{proof}
Let
$\delta := 1$ if $M$ is of type $\Mtype$ and $i = 0$, 
$\delta := 2$ otherwise.
By Theorem~\ref{T020900}, the socle of $\De_i M$ is isomorphic to
$L \circledast L(i)$ for some irreducible $\H_{n-1}$-module $L$,
and $$
\res^{n-1,1}_{n-1} L \circledast L(i)
\cong L^{\oplus \delta},
$$
indeed it is exactly as in the statement of the corollary.
Now take any irreducible submodule $K$ of $\res^{n-1,1}_{n-1} \circ
\De_i (M)$.
Consider the $\H_{n-1,1}$-submodule $\H_1' K$, where
$\H_1'$ is the subalgebra generated by $C_n, X_n^{\pm 1}$.
All composition 
factors of $\H_1' K$ are isomorphic to $K \circledast L(i)$.
In particular, the socle of $\H_1' K$ is isomorphic to $K \circledast L(i)$
which implies $K \cong L$.
We have now shown that $\soc \res^{n-1,1}_{n-1} \De_i (M)
\cong L^{\oplus{\delta'}}$ for some $\delta' \geq \delta$.

By Lemma~\ref{sick}, $\eps_i (L) = \eps - 1$ where $\eps = \eps_i(M)$, 
and $\De_{i^{\eps - 1}} L$
is irreducible
by Lemma~\ref{L310800_1}.
So at least $\delta'$ copies of $\De_{i^{\eps -1}} L$ appear in
$\soc \res^{n-\eps,\eps}_{n-\eps,\eps-1} \De_{i^\eps} (M).$
But $\De_{i^\eps} M \cong N \circledast L(i^\eps)$ for some irreducible $N$,
so applying Theorem~\ref{TKato}(iii) and the facts about type
in Theorem~\ref{T020900},
\begin{align*}
\soc \res^{n-\eps,\eps}_{n-\eps,\eps-1}
\De_{i^\eps} (M)
&\cong \soc \res^{n-\eps,\eps}_{n-\eps,\eps-1}
N \circledast L(i^\eps)\\
&\cong \res^{n-\eps,\eps-1,1}_{n-\eps,\eps-1}
N \circledast L(i^{\eps-1}) \circledast L(i)
\cong 
(N \circledast L(i^{\eps-1}))^{\oplus \delta}.
\end{align*}
Hence $\delta' \leq \delta$ too.
\end{proof}

\Point{\boldmath Operators $\tilde e_i$ and $\tilde f_i$}\label{ctr}
Let $\Irr_{\!\!I} \H_n$ denote the set of isomorphism classes of irreducible
modules in $\rep_I \H_n$.
We define
\begin{equation}\label{binfty}
B(\infty) = \bigcup_{n \geq 0} \Irr_{\!\!I} \H_n,
\end{equation}
For each $i \in I$, we define the {\em affine crystal operators}
\begin{align}\label{affinecrystal}
\tilde e_i&:B(\infty)\cup\{0\} 
\rightarrow B(\infty)\cup\{0\},\\
\tilde f_i&:B(\infty)\cup\{0\} 
\rightarrow B(\infty)\cup\{0\}
\end{align}
as follows. First, we set $\tilde e_i(0) = \tilde f_i(0) = 0$.
Now let $M$ be an irreducible module in $\rep_I \H_n$.
Then, $\tilde f_i M$ is defined by 
\begin{equation}
\tilde f_i M:=\cosoc \ind_{n,1}^{n+1} M \circledast L(i),
\end{equation}
which is irreducible by 
Lemma~\ref{L090900_1}.
To define $\tilde e_i M$,
Theorem~\ref{T020900} shows that either 
$\De_i M=0$ or $\soc \De_i M\cong N\circledast L(i)$ for 
an irreducible module $N\in \rep_I{\H_{n-1}}$. 
In the former case, we define $\tilde e_i M = 0$;
in the latter case, we define $\tilde e_i M = N$.
Thus:
\begin{equation}
\soc \De_i M \cong (\tilde e_i M) \circledast L(i).
\end{equation}
Note right away from Lemma~\ref{sick} that
\begin{equation}
\label{E300900}
\eps_i(M)=\max\{m\geq 0\:|\: \tilde e_i^m M\neq 0\},
\end{equation}
while a special case of Lemma~\ref{L090900_1} shows that
\begin{equation}
\label{hi}
\eps_i(\tilde f_i M)=\eps_i (M) + 1.
\end{equation}

\vspace{1mm}
\begin{Lemma}\label{sunny}
Let $M$ be an irreducible in $\rep_I \H_n$, $i \in I$ and $m \geq 0$.
\begin{enumerate}
\item[(i)]
$\soc \De_{i^m} M \cong (\tilde e_i^m M)\circledast L(i^m).$
\item[(ii)]
$\cosoc \ind_{n,m}^{n+m} M \circledast L(i^m) \cong \tilde f_i^m M.$
\end{enumerate}
\end{Lemma}
\begin{proof}
(i)
If $m>\eps_i(M)$, 
then both parts in the equality above are zero. 
Let $m\leq \eps_i(M)$. 
Clearly, $(\tilde e_i^m M)$ is a submodule of $\res^{n-m,m}_{n-m} \De_{i^m} M$.
Hence, $(\tilde e_i^m M) \circledast L(i)^{\circledast m}$ is a submodule of
$\res^{n-m,m}_{n-m,1,\dots,1} \De_{i^m} M$. So Frobenius reciprocity gives that
$(\tilde e_i^m M) \circledast L(i^m)$ is a submodule of $\De_{i^m} M$. Now we
are done by Theorem~\ref{T020900}.

(ii) 
By exactness of induction, 
$\tilde f_i^m M$ is a quotient of $\ind_{n,m}^{n+m} M \circledast L(i^m)$. 
Now the result follows from the simplicity of the cosocle, 
see Lemma~\ref{L090900_1}.
\end{proof}

\begin{Lemma}
\label{L290900_2}
Let $M\in \rep_I \H_n$ and $N \in \rep_I \H_{n+1}$ be irreducible,
and $i \in I$. Then, 
$\tilde f_i M\cong N$ if and only if $\tilde e_i N\cong M$. 
\end{Lemma}

\begin{proof}
Suppose $\tilde f_i M \cong N$.
Then by (\ref{E090900_1}),
$\hom_{\H_{n,1}}(M \circledast L(i), \De_i N) 
\neq 0,$
so $M \circledast L(i)$ appears in the socle of $\De_i N$.
This means that $M \circledast L(i) \cong (\tilde e_i N) \circledast L(i)$,
whence $M \cong \tilde e_i N$.
The converse is similar.
\end{proof}

\begin{Corollary}\label{imm}
Let $M, N$ be irreducible modules in $\rep_I \H_n$.
Then, $\tilde f_i M \cong \tilde f_i N$ if and only if
$M \cong N$. Similarly, providing $\eps_i(M),\eps_i(N) > 0$,
$\tilde e_i M \cong \tilde e_i N$ if and only if $M \cong N$.
\end{Corollary}

\vspace{1mm}

We can also define cyclotomic analogues of the crystal operators.
So now suppose that $\la \in P_+$, and 
let $\Irr \H_n^\la$
denote the set of isomorphism classes of irreducible $\H_n^\la$-modules.
Define
\begin{equation}\label{bla}
B(\la) = \bigcup_{n \geq 0} \Irr \H_n^\la.
\end{equation}
The functors $\inf^\la$ and $\pr^\la$ induce
maps
\begin{equation}\label{inff}
\inf^\la:B(\la)\cup\{0\} \rightarrow B(\infty)\cup\{0\},
\qquad
\pr^\la:B(\infty)\cup\{0\} \rightarrow B(\la)\cup\{0\},
\end{equation}
with $\pr^\la \circ \inf^\la (L) = L$ for each $L \in B(\la)$.
In other words, we can view $B(\la)$ as a subset of $B(\infty)$ via
the embedding $\inf^\la$.
Now by restricting $\tilde e_i$ and $\tilde f_i$ to 
$B(\la) \subset B(\infty)$, we obtain the 
{\em cyclotomic crystal operators}, namely,
\begin{align}\label{cyclocrystal}
\tilde e_i^\la &= \pr^\la\circ \tilde e_i\circ \inf^\la:
B(\la)\cup\{0\} \rightarrow B(\la)\cup\{0\},\\
\tilde f_i^\la &= \pr^\la\circ \tilde f_i\circ \inf^\la:
B(\la)\cup\{0\} \rightarrow B(\la)\cup\{0\}
\end{align}
for each $i \in I$ and $\la \in P_+$.
Note that $\tilde e_i$ already maps $B(\la)$ into $B(\la) \cup\{0\}$, 
so we always have that $\tilde e_i M = \tilde e_i^\la M$ for $M \in B(\la)$.
This is certainly not the case for $\tilde f_i$:
for $M \in B(\la)$, it will often be the case that
$\tilde f_i^\la(M) = 0$ even though $\tilde f_i(M)$ is never zero.

\Point{Independence of irreducible characters}
We can now prove an important theorem, compare \cite[$\S$5.5]{Vt}:

\vspace{1mm}
\begin{Theorem}\label{inj}
The map
$\ch:K(\rep_I \H_n) \rightarrow K(\rep_I \A_n)$
is injective.
\end{Theorem}

\begin{proof}
We need to show that $\{\ch L\mid[L] \in \Irr_{\!\!I} \H_n\}$ 
is a linearly independent set in $K(\rep_I \A_n)$.
Proceed by induction on $n$, the case $n = 0$ being trivial.
Suppose $n > 0$ and 
$$
\sum_{L \in \Irr_{\!\!I} \H_n} a_L \ch L = 0
$$
for some $a_L \in \Z$.
Choose any $i \in I$. 
We will show by downward induction on 
$k = n,\dots,1$ that $a_L = 0$ for all $L$ with $\eps_i(L) = k$.
Since every irreducible $L$ has $\eps_i(L) > 0$ for at least
one $i \in I$, this is enough to complete the proof.

Consider first the case that $k = n$.
Then, $\De_{i^n} L = 0$ except if $L \cong L(i^n)$, by
Theorem~\ref{TKato}(i). Since $\ch \De_{i^n} L$ can be worked out
just from knowledge of $\ch L$ using Lemma~\ref{L310800_2},
we deduce on applying $\De_{i^n}$ to the equation that
the coefficient of $\ch L(i^n)$ is zero.
Thus the induction starts.
Now suppose $1 \leq k < n$ and that we have shown $a_L = 0$
for all $L$ with $\eps_i(L) > k$.
Apply $\De_{i^k}$ to the equation to deduce that
$$
\sum_{L\:\rm{with}\:\eps_i(L)=k} a_L \ch \De_{i^k} L = 0.
$$
Now each such $\De_{i^k} L$ is irreducible, hence
isomorphic to $(\tilde e_i^k L) 
\circledast L(i^k)$,
according to Lemmas~\ref{L310800_1} and \ref{sunny}(i). 
Moreover, for $L \not\cong L'$, $\tilde e_i^k L \not\cong \tilde e_i^k L'$
by Corollary~\ref{imm}.
So now the induction hypothesis on $n$ gives that all such coefficients
$a_L$ are zero, as required.
\end{proof}

\begin{Corollary}\label{selfdual}
If $L$ is an irreducible module in $\rep_I \H_n$,
then $L \cong L^\tau$.
\end{Corollary}

\begin{proof}
Since $\tau(X_i) = X_i$, $\tau$ leaves characters invariant.
Hence it leaves irreducibles invariant since they are determined
up to isomorphism by their character according to the theorem.
\end{proof}

We can also show at this point that
the type of an irreducible module $L$ is determined by the type
of its central character:

\vspace{1mm}
\begin{Lemma}\label{cccc}
Suppose $L \in \rep_I \H_n$ is irreducible with 
central character $\chi_\gamma$ where
$\gamma = \sum_{i \in I} \gamma_i \alpha_i \in \Gamma_n$.
Then, $L$ is of type $\Qtype$ is $\gamma_0$ is odd, type $\Mtype$ if
$\gamma_0$ is even.
\end{Lemma}

\begin{proof}
Proceed by induction on $n$, the case $n = 0$ being trivial.
If $n > 1$, let $L$ be an irreducible $\H_n$-module with central character
$\chi_\gamma$.
Choose $i \in I$ so that $\tilde e_i L \neq 0$.
By definition, $\tilde e_i L$ has central character $\gamma - \alpha_i$.
So by the induction hypothesis, $\tilde e_i L$
is of type $\Qtype$ if $\gamma_0 - \delta_{i,0}$ is odd,
type $\Mtype$ otherwise.
But by general theory \ref{theory} and
Lemma~\ref{L290800_1},
$(\tilde e_i L) \circledast L(i)$ is of the opposite type
to $\tilde e_i L$ if $i = 0$, of the same
type if $i \neq 0$.
Hence,
$(\tilde e_i L) \circledast L(i)$ is of
type $\Qtype$ if $\gamma_0$ is odd,
type $\Mtype$ otherwise.
Finally, the proof is completed by Theorem~\ref{T020900}, 
since this shows that
$L$ has the same type as $\soc \De_i L = (\tilde e_i L) \circledast L(i)$.
\end{proof}

\Point{Crystal graphs}\label{cg}
We can view the datum
$(B(\infty), \tilde e_i, \tilde f_i)$
as a combinatorial structure: the {\em crystal graph}.
This is the directed graph with vertices the set $B(\infty)$ and an
edge 
$$
[M] \stackrel{i}{\longrightarrow} [N]
$$
whenever $[M], [N] \in B(\infty)$ satisfy $\tilde f_i M \cong N$, or
equivalently, by Lemma~\ref{L290900_2}, $M \cong \tilde e_i N$.
Similarly, $(B(\la), \tilde e_i^\la, \tilde f_i^\la)$
can be viewed as a crystal graph.

Motivated by this, we introduce some notation to label the isomorphism classes
of irreducible representations, or equivalently the vertices of the
crystal graph.
Write ${\bf 1}$ for the (trivial) irreducible module
of $\H_0$.
If  $L$ is an irreducible module in $\rep_I \H_n$,
one easily shows using
Lemma~\ref{L290900_2} repeatedly that
$$
L \cong \tilde f_{i_n} \dots \tilde f_{i_2} \tilde f_{i_1} {\bf 1}
$$
for at least one tuple $\bi = (i_1,i_2,\dots,i_n) \in I^n$.
So if we define
$$
L(\bi) = L(i_1,\dots,i_n) := 
\tilde f_{i_n} \dots \tilde f_{i_2} \tilde f_{i_1} {\bf 1},
$$
we obtain a labelling of all irreducibles in $\rep_I \H_n$ by
tuples in $I^n$.
For example, $L(i,i,\dots,i)$ ($n$ times) is precisely the principal series
module
$L(i^n)$ introduced in (\ref{mkato}).
Similarly, any irreducible $\H_n^\la$-module can be represented
as 
$$
L^\la(\bi) = L^\la(i_1,\dots,i_n) :=
\tilde f_{i_n}^\la \tilde f_{i_{n-1}}^\la \dots \tilde f_{i_1}^\la
\bil,
$$
where $\bil$ is the irreducible $\H_0^\la$-module.
Of course, $L^\la(\bi) \cong \pr^\la L(\bi)$.

Thus, our labelling of the irreducible modules in $\rep_I \H_n$ 
is by {\em paths} in the crystal graph starting from $\mathbf 1$,
and similarly for $\rep \H_n^\la$.
Of course, the problem with this labelling
is that a given irreducible $L$ will in general
be parametrized by {several} {\em different} tuples $\bi \in I^n$, 
corresponding to different paths from $\mathbf 1$ to $L$. 
But basic properties of $L(\bi)$ are easy to
read off from the notation: for instance
the central character of $L(\bi)$ is 
$\chi_{\bi}$, so by Lemma~\ref{cccc}, $L(\bi)$ is 
of type $\Qtype$ if an odd number of the $i_j$ are zero,
type $\Mtype$ otherwise.

\Point{Boring central characters}
Given $\bi = (i_1,\dots,i_n)$, let
$$
\ind(\bi) = \ind(i_1,\dots,i_n) := \ind_{1,\dots,1}^{n} L(i_1)\circledast
\dots\circledast L(i_n).
$$
Note the character of $\ind(\bi)$ is
$\sum_{w \in \Sym_n} L(i_{w^{-1}1})\circledast\dots\circledast L(i_{w^{-1}n})$.
So every irreducible constituent of $\ind(\bi)$ belongs to the block
$\rep_\gamma \H_n$, where $\gamma = \wt(\bi)$.

\vspace{1mm}
\begin{Lemma}\label{boring1}
Let $\gamma \in \Gamma_n$ and pick any $\bi \in I^n$ with $\wt(\bi) = \gamma$.
Then:
\begin{enumerate}
\item[(i)] $\res^{n}_{1,\dots,1} L(\bi)$ has a submodule isomorphic to
$L(i_1) \circledast \dots \circledast L(i_n)$;
\item[(ii)] $\ind(\bi)$ contains a copy of $L(\bi)$ in its cosocle;
\item[(iii)] every irreducible module in the block $\rep_\gamma \H_n$
appears at least once as a constituent of $\ind(\bi)$.
\end{enumerate}
\end{Lemma}

\begin{proof}
(i) Proceed by induction on $n$.
For the induction step, let $\bj = (i_1,\dots,i_{n-1})$.
By Frobenius reciprocity, there
is a non-zero (hence necessarily injective) $\H_{n-1,1}$-module
homomorphism from
$L(\bj) \circledast L(i_n)$ to $\res^n_{n-1,1} L(\bi)$.
Hence by induction we get a copy of
$L(i_1)\circledast \dots \circledast L(i_{n-1}) \circledast L(i_n)$
in $\res^{\H_n}_{\A_n} L(\bi)$.

(ii) Use (i) and Frobenius reciprocity.

(iii) We have just shown that $L(\bi)$ appears in $\ind(\bi)$.
But for any other $\bj$ with the same weight as $\bi$,
$\ind(\bj)$ has the same character as $\ind(\bi)$, hence they
have the same set of composition factors  
thanks to Theorem~\ref{inj}.
Hence, $L(\bi)$ appears in $\ind(\bj)$. 
\end{proof}

We also need the following {criterion} for irreducibility,
compare \cite[Lemma~5.9]{G}:

\vspace{1mm}
\begin{Lemma}\label{useful}
Let $M,N$ be irreducibles in $\rep_I \H_m, \rep_I \H_n$ respectively.
Suppose
\begin{enumerate}
\item[(i)] $\ind_{m,n}^{m+n} M \circledast N
\cong \ind_{n,m}^{n+m} N \circledast M$;
\item[(ii)] $M \circledast N$ appears in
$\res^{m+n}_{m,n} \ind_{m,n}^{m+n} M\circledast N$
with multiplicity one.
\end{enumerate}
\noindent
Then, $\ind_{m,n}^{m+n} M\circledast N$
is irreducible.
\end{Lemma}

\begin{proof}
Suppose for a contradiction that
$K = \ind_{m,n}^{m+n} M \circledast N$ is reducible. 
Then we can find a proper irreducible submodule $S$, and set $Q = K/S$.
By Frobenius reciprocity, $M\circledast N$
appears in $\res^{m+n}_{m,n}Q$ 
with non-zero multiplicity. Hence, it cannot appear in $\res^{m+n}_{m,n}S$
by assumption (ii).
But assumption (i), Corollary~\ref{selfdual} 
and Theorem~\ref{duals} show that
$K \cong K^\tau$.
Hence, $K$ also has a quotient isomorphic to $S^\tau \cong S$,
and the Frobenius reciprocity argument implies that
$M \circledast N$ appears in $\res^{m+n}_{m,n}S$.
\end{proof}

\begin{Lemma}\label{boring2}
Let $\bi \in I^m, \bj \in I^n$ be tuples 
such that $|i_a - j_b| > 1$ for all $1 \leq a \leq m,
1 \leq b \leq n$.
Then,
$\ind_{m,n}^{m+n} L(\bi) \circledast L(\bj)
\cong
\ind_{n,m}^{m+n} L(\bj) \circledast L(\bi)$
is irreducible.
\end{Lemma}

\begin{proof}
By Lemma~\ref{useful} and the Shuffle Lemma,
it suffices to show that
$$
\ind_{m,n}^{m+n} L(\bi) \circledast L(\bj)
\cong
\ind_{n,m}^{m+n} L(\bj) \circledast L(\bi).
$$
By the Mackey Theorem,
$\res^{m+n}_{m,n} \ind_{n,m}^{m+n} L(\bj) \circledast L(\bi)$
contains $L(\bi) \circledast L(\bj)$ as a summand with multiplicity one,
all other constituents lying in different blocks.
Hence by Frobenius reciprocity, there exists a non-zero homomorphism
$$
f:\ind_{m,n}^{m+n} L(\bi) \circledast L(\bj) \rightarrow
\ind_{n,m}^{m+n} L(\bj) \circledast L(\bi).
$$
Every homomorphic image of
$\ind_{m,n}^{m+n} L(\bi) \circledast L(\bj)$
contains an $\H_{m,n}$-submodule isomorphic to $L(\bi) \circledast
L(\bj)$.
So, by Lemma~\ref{boring1}(i), we see that the image of $f$ contains
an $\A_{m+n}$-submodule $V$ isomorphic to 
$L(i_1) \circledast \dots \circledast L(i_m) \circledast
L(j_1) \circledast \dots \circledast L(j_n)$.

Now we claim that the image of $f$ also contains an
$\A_{m+n}$-submodule isomorphic to
$L(j_1) \circledast \dots \circledast L(j_n) \circledast
L(i_1) \circledast \dots \circledast L(i_m)$.
To see this, consider $\tilde \Phi_m V$.
Pick a common eigenvector $v \in V$ for the operators
$X_m^{\pm 1}$ and $X_{m+1}^{\pm 1}$;
then $(X_m + X_m^{-1})v = q(i_m) v$ and
$(X_{m+1}+X_{m+1}^{-1})v = q(j_1) v$.
So according to (\ref{nazeltsq}),
$\tilde \Phi_m^2$ acts on $v$ by a scalar, and the assumption that
$|i_m - j_1| > 1$ combined with Lemma~\ref{nazelttech} shows that this scalar
is necessarily non-zero.
Thus, $\tilde \Phi_m V \neq 0$, so by 
(\ref{nazeltp2}) it is an irreducible $\A_{m+n}$-module, namely,
$$
\tilde \Phi_m V \cong
L(i_1) \circledast \dots \circledast L(i_{m-1}) \circledast
L(j_1) \circledast L(i_m) \circledast L(j_2) 
\circledast\dots\circledast L(j_n).
$$
Next apply $\tilde \Phi_{m-1}, \dots, \tilde \Phi_1$ to move $L(j_1)$
to the first position, and continue in this way to complete the proof
of the claim.

We have now shown that the image of $f$ contains 
$L(j_1) \circledast \dots \circledast L(j_n)\circledast
L(i_1) \circledast \dots \circledast L(i_m).$
But by the Shuffle Lemma, all such
composition factors of 
$\res^{m+n}_{1,\dots,1}
\ind_{n,m}^{m+n} L(\bj) \circledast L(\bi)$ necessarily lie in the
irreducible 
$\H_{n,m}$-submodule $1 \otimes L(\bj) \circledast L(\bi)$
of the induced module.
Since this generates all of
$\ind_{n,m}^{m+n} L(\bj) \circledast L(\bi)$ as an $\H_{m+n}$-module,
this shows that $f$ is surjective.
Hence $f$ is an isomorphism by dimension,
which completes the proof.
\end{proof}

\begin{Theorem}
\label{TBoring}
Let $\bi \in I^m, \bj \in I^n$ be tuples 
such that $|i_a - j_b| \neq 1$ for all $1 \leq a \leq m,
1 \leq b \leq n$.
Then,
$\ind_{m,n}^{m+n} L(\bi) \circledast L(\bj)
\cong
\ind_{n,m}^{m+n} L(\bj) \circledast L(\bi)
$
is irreducible.
Moreover, every other irreducible module lying in the same block as
$\ind_{m,n}^{m+n} L(\bi) \circledast L(\bj)$ is of the form
$\ind_{m,n}^{m+n} L(\bi') \circledast L(\bj')$
for permutations $\bi'$ of $\bi$ and $\bj'$ of $\bj$.
\end{Theorem}

\begin{proof}
The second statement of the theorem is an easy consequence of
the first and Lemma~\ref{boring1}(iii).
For the first statement, 
proceed by induction on $m+n$, the case $m+n = 1$ being trivial.
For $m+n > 1$, we may assume by Lemma~\ref{boring2} that there exists
$k \in I$ that appears in both the tuples
$\bi$ and $\bj$. 
Note then that for every $a = 1,\dots,m$, either $i_a =k$
or $|i_a -k| > 1$, and similarly for every $b = 1,\dots,n$,
either $j_b = k$ or $|j_b - k| > 1$.
So by the induction hypothesis, we have that
$$
L(\bi)\cong \ind_{n-r,r}^n L(\bi')\circledast L(k^r),
\qquad
L(\bj) \cong \ind_{m-s,s}^m L(\bj') \circledast L(k^s)
$$
for some $r,s \geq 1$, where $\bi',\bj'$ are tuples with no entries equal to $k$.
By Theorem~\ref{duals}, Corollary~\ref{selfdual} and 
Lemma~\ref{boring2},
$$
\ind^{r+m-s}_{r,m-s} L(k^r) \circledast L(\bj')
\cong
\ind^{m-s+r}_{m-s,r} L(\bj') \circledast L(k^r).
$$
So using Theorem~\ref{TKato}(i) and
transitivity of induction,
\begin{align*}
\ind_{m,n}^{m+n} L(\bi)\circledast L(\bj)
&\cong
\ind_{n-r,r,m-s,s}^{m+n} 
L(\bi')\circledast L(k^r)
\circledast L(\bj') \circledast L(k^s)\\
&\cong
\ind_{n-r,m-s,r+s}^{m+n} 
L(\bi')\circledast L(\bj')
\circledast L(k^{r+s}).
\end{align*}
Finally this is irreducible by the induction hypothesis and 
Lemma~\ref{boring2}.
\end{proof}

\Point{Some character calculations}\label{calcsec}
At this point we need to compute the characters of certain very special
$\H_n$-modules explicitly.

\vspace{1mm}
\begin{Lemma}\label{calc1} Let $i,j \in I$ with $|i-j|=1$.
Then, for all $a,b \geq 0$ with $a+b < - \langle h_i,\alpha_j \rangle$,
there is a non-split short exact sequence
$$
0 \longrightarrow L(i^{a+1} j i^b) 
\longrightarrow \ind_{a+b+1,1}^{a+b+2} L(i^{a}j i^b) \circledast L(i)
\longrightarrow L(i^{a} j i^{b+1}) \longrightarrow 0.
$$
Moreover, for every $a,b \geq 0$ with $a+b \leq - \langle h_i,\alpha_j \rangle$,
$$
\ch L(i^a j i^b) = (a!)(b!) [L(i)^{\circledast a} \circledast L(j)
\circledast L(i)^{\circledast b}].
$$
\end{Lemma}

\begin{proof}
We proceed by induction on $n = 0,1,\dots,-\langle h_i,\alpha_j \rangle$
to show that
$$
\ch L(i^n j) = n! [L(i)^{\circledast n} \circledast L(j)],
$$
this being immediate in case $n = 0$.
For $n > 0$, 
let $M := \ind_{n,1}^{n+1} L(i^{n-1} j) \circledast L(i)$.
We know by the inductive hypothesis and the Shuffle Lemma
that
$$
\ch M
= 
n!
[L(i)^{\circledast n}\circledast  L(j)]
+ (n-1)!
[L(i)^{\circledast (n-1)}\circledast L(j)\circledast  L(i)].
$$
Now consider the $\H_{n,1}$-submodule
$$
N := (X_{n+1}+X_{n+1}^{-1} - q(i)) M \cong L(i^n) \circledast L(j)
$$
of $M$. The key point is that $N$ 
is stable under the action of $T_n$, hence all of $\H_{n+1}$.
Although this is in principle an elementary calculation, 
it turns out to be extremely lengthy. It was carried out by hand
for $n \leq 2$ but for the cases $n = 3,4$ (when $\ell$ necessarily
equals $1$),
we had to resort to a computer calculation using the
{\sc Gap} computer algebra package.
This proves the existence of an irreducible $\H_{n+1}$-module $N$ 
with character $n! [L(i)^{\circledast n} \circledast L(j)]$.
This must be $L(i^n j)$, by Lemma~\ref{boring1}(i), completing the proof
of the induction step.

Now we explain how to deduce the characters of the remaining irreducibles
in the block.
In the argument just given, the quotient module
$M / N$ has character
$(n-1)! [L(i)^{\circledast (n-1)} \circledast L(j)\circledast L(i)]$,
so must be $L(i^{n-1} j i)$.
Twisting with the automorphism $\sigma$ proves that there exist irreducibles
with characters
$n! [L(j) \circledast L(i)^{\circledast n}]$
and $(n-1)! [L(i)\circledast L(j)\circledast L(i)^{\circledast(n-1)}]$,
which must be $L(j i ^n)$ and $L(iji^{n-1})$ respectively by
Lemma~\ref{boring1}(i) once more.
This covers everything unless $n = 4$, when we necessarily have that
$i = 0, j = 1$ and $\ell = 1$. In this case, we have shown already that
there exist four irreducibles with characters
\begin{align*}
\ch L(00001) = 24[L(0)^{\circledast 4} \circledast L(1)],\quad
&\ch L(00010) = 6[L(0)^{\circledast 3} \circledast L(1)\circledast L(0)],\\
\ch L(10000) = 24[L(1)\circledast L(0)^{\circledast 4}],
\quad
&\ch L(01000) = 6[L(0)\circledast L(1)\circledast L(0)^{\circledast 3}].
\end{align*}
So by Lemma~\ref{boring1}, there must be exactly one more irreducible
module in the block, namely $L(00100)$,
since none of the above involve the character
$[L(0)^{\circledast 2} \circledast L(1)\circledast L(0)^{\circledast 2}]$.
Considering the character of 
$\ind_{n,1}^{n+1} L(0010)\circledast L(0)$
shows that $\ch L(00100)$ is either
$4 [L(0)^{\circledast 2} \circledast L(1)\circledast L(0)^{\circledast 2}]$
or
$4[L(0)^{\circledast 2} \circledast L(1)\circledast L(0)^{\circledast 2}]
+ 6 [L(0)^{\circledast 3} \circledast L(1) \circledast L(0)]$.
But the latter is not $\sigma$-invariant so cannot occur as there would
then be too many irreducibles.

Now that the characters are known, it is finally a routine matter using
the Shuffle Lemma and Lemma~\ref{L090900_1} to prove the existence
of the required non-split sequence.
\end{proof}

\begin{Lemma}\label{calc2}
Let $i, j \in I$ with $|i-j|=1$ and set $n = 1-\langle h_i, \alpha_j 
\rangle$.
Then $$
L(i^n j) \cong L(i^{n-1} j i).
$$
Moreover,
for every $a,b \geq 0$ with $a+b = -\langle h_i,\alpha_j \rangle$,
$$
L(i^a j i^{b+1}) \cong \ind_{n,1}^{n+1} L(i^a j i^b) \circledast L(i)
\cong \ind_{1,n}^{n+1} L(i)\circledast L(i^a j i^b)
$$
with character
$a! (b+1)! [L(i)^{\circledast a} \circledast L(j)
\circledast L(i)^{\circledast (b+1)}]+
(a+1)! b! [L(i)^{\circledast (a+1)} \circledast L(j)
\circledast L(i)^{\circledast b}].$
\end{Lemma}

\begin{proof}
Let $M = \ind_{n,1}^{n+1} L(i^{n-1}j) \circledast L(i)$.
We first claim that $M$ is irreducible.
To prove this, arguing in the same way as the proof of Lemma~\ref{calc1},
it suffices to show that the $\H_{n,1}$-submodule
$$
(X_{n+1}+X_{n+1}^{-1} - q(i)) M \cong L(i^n)\circledast L(j)
$$
of $M$ is {\em not} invariant under $T_n$.
Again this was checked by an explicit 
computer calculation.
Hence, there is an irreducible $\H_n$-module $M$ with character
$$
n!
[L(i)^{\circledast n}\circledast  L(j)]
+ (n-1)!
[L(i)^{\circledast (n-1)}\circledast L(j)\circledast  L(i)].
$$
Hence $\tilde e_i M \cong L(i^{n-1}j)$ and
$\tilde e_j M \cong L(i^n)$ by Theorem~\ref{T020900}.
So we deduce that $M \cong L(i^{n-1}j i) \cong L(i^n j)$ thanks to
Lemma~\ref{L290900_2}.

Now consider the remaining irreducibles in the block.
There are at most $(n-1)$ remaining, namely
$L(i^a j i^{b+1})$ for $a\geq 0,b\geq 1$ with $a+b=-\langle h_i,
\alpha_j\rangle$.
Considering the known characters
of $\ind_{n,1}^{n+1} L(i^a j i^b) \circledast L(i)$ and arguing in a similar
way to the second paragraph of the proof of the preceeding lemma, the remainder
of the lemma follows without further calculation.
\end{proof}

\Point{Higher crystal operators}
In this subsection we will introduce certain generalizations
of the crystal operators $\tilde f_j$, following the ideas of
\cite[$\S$10]{G}.
The results of this subsection are only needed 
in \ref{crel} below.
To simplify notation, we will write simply $\ind$
in place of $\ind_\mu^n$ throughout the subsection.

\vspace{1mm}
\begin{Lemma}\label{serre1}
Let $i,j \in I$ with $i \neq j$.
For any $a,b\geq 0$ with $a+b = -\langle h_i,\alpha_j \rangle$,
$$
\ind\: L(i^a j i^{b}) \circledast L(i^m)
\cong \ind\: L(i^m)\circledast L(i^a j i^{b})
$$
is irreducible.
\end{Lemma}

\begin{proof}
We first claim that
$$
\ind\: L(i^a j i^{b}) \circledast L(i^m)
\cong \ind\:
L(i^m)\circledast L(i^a j i^{b}).
$$
This is immediate from Lemma~\ref{boring2} in case $|i-j| > 1$.
If $|i-j| = 1$, then 
transitivity of induction and Lemma~\ref{calc2} give
that
\begin{align*}
\ind\: L(i^a j i^{b}) \circledast L(i^m)
&\cong
\ind\: L(i^a j i^{b}) 
\circledast L(i)\circledast L(i^{m-1})\\
&\cong
\ind\: L(i)\circledast L(i^a j i^{b}) \circledast  
L(i^{m-1}),
\end{align*}
and now repeating this argument $(m-1)$ more times gives the claim.

Hence, by
Corollary~\ref{selfdual} 
and Theorem~\ref{duals},
$K := \ind\: L(i^a j i^{b}) \circledast L(i^m)$ is self-dual.
Now suppose for a contradiction that $K$ is reducible. Then
we can pick  a  proper irreducible submodule
$S$ of $K$, and set $Q := K / S$.
Applying Lemmas~\ref{calc1} and \ref{LChar},
$$
\ch K = \sum_{k=0}^m \binom{m}{k}(a+k)!(b+m-k)! [L(i)^{\circledast (a+k)}
\circledast L(j) \circledast L(i)^{\circledast (b+m-k)}].
$$
By Frobenius reciprocity, $Q$ contains an $\H_{a+b+1,m}$-submodule
isomorphic to $L(i^aji^{b}) \circledast L(i^m)$.
So by Lemma~\ref{boring1}(i),
the irreducible $\A_{a+b+m+1}$-module
$L(i)^{\circledast a} \circledast L(j) \circledast L(i)^{\circledast (b+m)}$
appears in $Q$ with non-zero multiplicity, hence in fact by Theorem~\ref{TKato}(i)
it must appear with multiplicity $a!(b+m)!$ (viewing 
$Q$ as a module over $\H_{a,1,b+m}$).
It follows that 
$L(i)^{\circledast a}
\circledast L(j) \circledast L(i)^{\circledast (b+m)}$ is not a composition
factor of $S$. 
But this is a contradiction, since as $K$ is self-dual, $S \cong 
S^\tau$ is a quotient
module of $K$ hence must contain
$L(i)^{\circledast a} \circledast L(j) \circledast L(i)^{\circledast (b+m)}$
by the Frobenius reciprocity argument again.
\end{proof}

\begin{Lemma}\label{serre2}
Let $i,j \in I$ with $i \neq j$.
For any $a \geq 1$ and $b \geq 0$ with $a+b = -\langle h_i,\alpha_j\rangle$, 
any irreducible module $M$ in $\rep_I \H_n$ and any $m \geq 0$,
$$
\cosoc \ind\: 
M \circledast L(i^m) \circledast L(i^a j i^b)
$$
is irreducible.
\end{Lemma}

\begin{proof}
By the argument in the proof of Lemma~\ref{L090900_1},
it suffices to prove this in the special case that $\eps_i(M) = 0$.
Let $k = m+a+b+1$.
Recall from the previous lemma that
$$
N := \ind\: L(i^m) \circledast L(i^a j i^b)
$$ 
is an irreducible $\H_{k}$-module.
Moreover by Lemma~\ref{calc1}, 
$\ch L(i^a j i^b) = (a!)(b!)[L(i)^{\circledast a}\circledast L(j)
\circledast L(i)^{\circledast b}]$.
So since $\eps_i(M) = 0$ and $a > 0$, 
the Mackey Theorem and a block argument
shows that
$$
\res^{n+k}_{n,k} (\ind\:
M \circledast L(i^m) \circledast L(i^a j i^b))
\cong  (M \circledast N) \oplus U
$$
for some $\H_{n,k}$-module $U$
all of whose composition factors lie in different blocks to those of
$M \circledast N$.
Now let
$H := \cosoc
\ind\:
M \circledast L(i^m) \circledast L(i^a j i^b)$.
It follows from above that
$\res^{n+k}_{n,k} H \cong (M \circledast N)
\oplus \tilde U$
where $\tilde U$ is some quotient module of $U$.
Then:
\begin{align*}
\hom_{\H_{n+k}}(H,H)
&\cong
\hom_{\H_{n+k}}(\ind\: M \circledast L(i^m) \circledast L(i^a j i^b),H)\\&
\cong
\hom_{\H_{n,k}}(M \circledast N,
\res^{n+k}_{n,k} H)\\
&\cong
\hom_{\H_{n,k}}(M \circledast N,
M \circledast N
\oplus \tilde U
) \cong
\hom_{\H_{n,k}}(M \circledast N,M \circledast N).
\end{align*}
Since $H$ is completely reducible and $M \circledast N$ is irreducible,
this implies that $H$ is irreducible too, as required.
\end{proof}

Now we can define the higher crystal operators.
Let $i,j \in I$ with $i\neq j$, and $a \geq 1, b\geq 0$
with $a+b = -\langle h_i,\alpha_j\rangle$.
Then the special case $m = 0$ of the theorem shows that
$$
\tilde f_{i^a j i^b} M := \cosoc 
\ind\: M \circledast L(i^a j i^b)
$$
is irreducible for every irreducible $M$ in $\rep_I \H_n$.
Thus we have defined an operator
\begin{equation}\label{hop}
\tilde f_{i^a j i^b}:B(\infty)\rightarrow B(\infty).
\end{equation}

\vspace{1mm}
\begin{Lemma}\label{serre3}
Take $i,j \in I$ with $i \neq j$ and set $k = -\langle
h_i,\alpha_j\rangle$.
Let $M$ be an irreducible module in $\rep_I \H_n$.
\begin{enumerate}
\item[(i)] There exists a unique integer 
$a$ with $0\leq a \leq k$ such that for every $m\geq 0$ we have
$$
\eps_i (\tilde f_i^m 
\tilde f_j M) = m+\eps_i(M) - a.
$$
\item[(ii)] Assume $m \geq k$. Then
a copy of $\tilde f_i^m \tilde f_j M$ appears in
the cosocle of
$$
\ind\:
\tilde f_i^{m-k} M \circledast L(i^a j i^{k-a}).
$$
In particular, if $a \geq 1$, then
$\tilde f_i^m \tilde f_j M 
\cong \tilde f_{i^a j i^{k-a}} \tilde f_i^{m-k} M$.
\end{enumerate}
\end{Lemma}

\begin{proof}
Let $\eps=\eps_i(M)$ and write $M=\tilde f_i^\eps N$ for irreducible
$N\in\rep_I \H_{n-\eps}$ with $\eps_i(N)=0$.
It suffices to prove (i) 
for any fixed choice of $m$, the conclusion for all other $m \geq 0$
then following immediately by (\ref{hi}).
So take $m \geq k$. Note that $\tilde f_i^m\tilde f_j M=\tilde f_i^m\tilde
f_j\tilde f_i^\eps N$ is a quotient of 
$$
\ind\: N\circledast L(i^\eps)\circledast L(j)\circledast L(i)^{\circledast
k}\circledast L(i^{m-k}),
$$
which by Lemma~\ref{calc1} has a filtration with factors isomorphic to 
$$
F_a:=\ind\: N\circledast L(i^\eps)\circledast L(i^aji^{k-a})\circledast
L(i^{m-k}),\quad 0\leq a\leq k.
$$
So $\tilde f_i^m\tilde f_j M$ is a quotient of some such factor, and to
prove (i) it remains to show that $\eps_i(L)=\eps+m-a$ for any irreducible
quotient $L$ of $F_a$. 
The inequality $\eps_i(L)\leq \eps+m-a$ is clear from the
Shuffle Lemma.
On the other hand, by transitivity of induction and Lemma~\ref{serre1},
$F_a\cong\ind\, N\circledast(\ind\, L(i^aji^{k-a})\circledast
L(i^{\eps+m-k}))$. So by Frobenius Reciprocity, the irreducible module
$N\circledast(\ind\, L(i^aji^{k-a})\circledast L(i^{\eps+m-k}))$ is
contained in $\res_{n-\eps,m+1+\eps}L$. Hence $\eps_i(L)\geq \eps+m-a$.

For (ii), by Lemma~\ref{serre1}, we also have $F_a\cong\ind\, N\circledast
L(i^{m-k+\eps})\circledast L(i^aji^{k-a})$, and by the
Shuffle Lemma, the only
irreducible factors $K$ of $F_a$ with $\eps_i(K)=\eps+m-a$ come from its
quotient 
$$
\ind\: \tilde f_i^{m-k+\eps} N\circledast L(i^aji^{k-a})\cong\ind \tilde
f_i^{m-k} M\circledast L(i^aji^{k-a}).
$$
Finally, in case $a \geq 1$, the cosocle of the last module is precisely
$\tilde f_{i^a j i^{k-a}} \tilde f_i^{m-k} M$.
\end{proof}

\Point{Modifications in the degenerate case}\label{trivmods}
Other than replacing $X_i+X_i^{-1}$ by $x_i^2$ everywhere, everything in this
section goes through in the degenerate case in exactly the same way.
Note the computer calculations in \ref{calcsec} were checked separately
in the degenerate case.

\ifbook@\pagebreak\fi

\section{Induction and restriction}

\Point{\boldmath $i$-induction and $i$-restriction}
Fix $\la \in P_+$ throughout the subsection.
Recall the definition of the functors $\De_i$ 
for $i \in I$, see \ref{socs}.
Let us denote the composite functor $\res^{n-1,1}_{n-1} \circ \De_i$
instead by
\begin{equation}\label{relres}
\res_i:\rep_I{\H_n} \rightarrow \rep_I{\H_{n-1}},
\end{equation}
for any $n$.
Note if $M$ is an $\H_n^\la$-module, then $\res_i M$ is automatically
an $\H_{n-1}^\la$-module. So the restriction of the functor $\res_i$
gives a functor which we also denote
\begin{equation}\label{relres2}
\res_i:\rep{\H_n^\la} \rightarrow \rep \H_{n-1}^\la.
\end{equation}
We now focus on this cyclotomic case.

There is an alternative definition of $\res_i$:
if $M$ is a module in $\rep_\gamma \H_n^\la$ for some fixed
$\gamma = \sum_{j \in I} \gamma_j \alpha_j \in \Gamma_n$ then
\begin{equation}\label{sat1}
\res_i M = \left\{
\begin{array}{ll}
 (\res^{\H_n^\la}_{\H_{n-1}^\la} M)[\gamma - \alpha_i]&\hbox{if $\gamma_i > 0$,}\\
0&\hbox{if $\gamma_i = 0$.}
\end{array}
\right.
\end{equation}
This description makes it clear how to define an analogous (additive) functor
\begin{equation}\label{relind}
\ind_i:\rep{\H_n^\la} \rightarrow \rep\H_{n+1}^\la.
\end{equation}
Using (\ref{block2}) and 
additivity, it suffices to define this on an object $M$ belonging
to $\rep_\gamma \H_n^\la$ for fixed $\gamma = \sum_{j \in I}\gamma_j \alpha_j
\in \Gamma_n$.
Then, we set
\begin{equation}\label{sat2}
\ind_i M = 
(\ind_{\H_n^\la}^{\H_{n+1}^\la} M)[\gamma+\alpha_i].
\end{equation}
By the definitions (\ref{sat1}) and (\ref{sat2})
and Lemma~\ref{decomp}, we have that
\begin{equation}\label{indresdecomp}
\ind_{\H_n^\la}^{\H_{n+1}^\la} M = 
\bigoplus_{i \in I} \ind_i M,
\qquad
\res^{\H_n^\la}_{\H_{n-1}^\la} 
M = \bigoplus_{i \in I} \res_i M.
\end{equation}
To complete the definition of the functor $\ind_i$, it is defined
on a  morphism $f$ simply by restriction of the
corresponding morphism $\ind_{\H_n^\la}^{\H_{n+1}^\la} f$.
We stress that the functor $\ind_i$ depends fundamentally on the
fixed choice of $\la$, 
unlike $\res_i$ which is just the restriction of its affine counterpart.

\vspace{1mm}
\begin{Lemma}\label{sleepy} For $\la \in P_+$ and 
each $i \in I$,
\begin{enumerate}
\item[(i)] $\ind_i$ 
and $\res_i$ are
both left and right adjoint to each other, hence they are exact
and send projectives to projectives;
\item[(ii)] $\ind_i$ 
and $\res_i$ 
commute with duality, i.e.
there are natural isomorphisms
$$
\ind_i (M^\tau) \simeq (\ind_i M)^\tau,
\qquad
\res_i (M^\tau) \simeq (\res_i M)^\tau
$$
for each finite dimensional $\H_n^\la$-module $M$.
\end{enumerate}
\end{Lemma}

\begin{proof} We know that $\ind_i
M$ and $\res_i M$ 
are summands of $\ind_{\H_n^\la}^{\H_{n+1}^\la} M$
and $\res^{\H_n^\la}_{\H_{n-1}^\la} M$ respectively.
Moreover, $\tau$-duality leaves
central characters invariant because $\tau(X_j) = X_j$ for each $j$. 
Now everything follows easily applying Corollary~\ref{joey}.
\end{proof}

In order to refine the definitions of 
$\ind_i$ and 
$\res_i$ in the next subsection, 
we need to give an alternative definition, due to 
Grojnowski \cite[$\S$8]{G} in the untwisted case.
Recall the definition of the left $\H_1$-modules
$\R_m(i)$ for $i \in I, m \geq 0$ from (\ref{rmdef}).
The limits in the next lemma are taken with respect to the
systems induced by the maps (\ref{rseq}).

\vspace{1mm}
\begin{Lemma}\label{geq}
For every finite dimensional
$\H_n^\la$-module $M$ and $i \in I$, there are natural isomorphisms
\begin{align*}
\ind_i M&\simeq \varprojlim \pr^\la \ind_{n,1}^{n+1}
M \boxtimes \R_m(i),\\
\res_i M & \simeq \varinjlim \pr^\la \hom_{\H_1'}(\R_m(i), M)
\end{align*}
(in the second case, $\H_1'$ denotes the subalgebra
of $\H_{n-1,1}$ generated by $C_n, X_n^{\pm 1}$
and the $\H_{n-1}$-module structure is defined by
$(h f)(r) = h(f(r))$ for $f \in \hom_{\H_1'}(\R_m(i), M)$
and $r \in \R_m(i)$).
\end{Lemma}

\begin{proof}
For $\res_i$, it suffices to consider the effect
on $M \in \rep_\gamma \H_n^\la$ for $\gamma =\sum_{j\in I} \gamma_j \alpha_j
\in \Gamma_n$
with $\gamma_i > 0$, both sides of what we are trying to prove 
clearly being zero if $\gamma_i = 0$.
Then, for all sufficiently large $m$, Lemma~\ref{nastystab} (in the 
special case $n = 1$)
implies that
$$
\hom_{\H_1'}(\R_m(i), M) \simeq 
(\res^{\H_n^\la}_{\H_{n-1}^\la} M)[\gamma - \alpha_i].
$$
Hence,
$$
\varinjlim \pr^\la \hom_{\H_1'}(\R_m(i), M)
\simeq
(\res^{\H_n^\la}_{\H_{n-1}^\la}M)[\gamma - \alpha_i] = \res_i M.
$$
This proves the lemma for $\res_i$.

To deduce the statement about induction, it now suffices by uniqueness of
adjoint functors to show that
$\varprojlim \pr^\la \ind_{n,1}^{n+1}
 ? \boxtimes \R_m(i)$
is left adjoint to
$\varinjlim \pr^\la \hom_{\H_1'}(\R_m(i), ?).$
Let $N\in \rep{\H_{n-1}^\la}$ and $M\in \rep{\H_n^\la}$. 
First observe as explained in the previous paragraph that the
direct system
$$
\pr^\la\hom_{\H_1'}(\R_1(i), M)
\hookrightarrow
\pr^\la\hom_{\H_1'}(\R_2(i), M)
\hookrightarrow\dots
$$
stabilizes after finitely many terms.
We claim that the inverse system
$$
\pr^\la \ind_{n,1}^{n+1} N\boxtimes \R_1(i)
\twoheadleftarrow 
\pr^\la \ind_{n,1}^{n+1} N\boxtimes \R_2(i)
\twoheadleftarrow\dots
$$
also stabilizes after finitely many terms.
To see this,
it suffices to show that the dimension of 
$\pr^\la\ind_{n,1}^{n+1}  N \boxtimes \R_m(i)$
is bounded above independently of $m$.
Well, each $\R_m(i)$ is generated as an $\H_1'$-module
by a subspace $W$ isomorphic (as a vector space)
to the cosocle $\R_1(i)$ of $\R_m(i)$.
Then $\ind_{n,1}^{n+1}  N \boxtimes \R_m(i)$
is generated as an $\H_{n+1}$-module
by the subspace $W' = 1 \otimes (N \otimes W)$, also of dimension
independent of $m$.
Finally, $\pr^\la \ind_{n,1}^{n+1}
 N \boxtimes \R_m(i)$ is a quotient
of the vector space $\H_{n+1}^\la \otimes_F W'$,
whose dimension is independent of $m$.

Now we can complete the proof of adjointness.
Using the fact from the previous paragraph
that the direct and inverse systems stabilize after finitely many
terms,
we have natural isomorphisms
\begin{align*}
\hom_{\H_n^\la}(\varprojlim \pr^\la\ind_{n-1,1}^{n}  N\boxtimes \R_m(i), M)
&\simeq
\varinjlim \hom_{\H_n^\la}(\pr^\la\ind_{n-1,1}^{n}  N\boxtimes \R_m(i),  
M)\\
&\simeq \varinjlim \hom_{\H_n}(\ind_{n-1,1}^{n}
 N\boxtimes \R_m(i), M) \\
&\simeq\varinjlim \hom_{\H_{n-1,1}}(N\boxtimes \R_m(i), 
\res^{n}_{n-1,1}  M) \\
 &\simeq \varinjlim \hom_{\H_{n-1}}( N,
\hom_{\H_1'}(\R_m(i), M)) \\
 &\simeq \varinjlim \hom_{\H_{n-1}^\la}( N,
\pr^\la \hom_{\H_1'}(\R_m(i), M)) \\
 &\simeq \hom_{\H_{n-1}^\la}( N,\varinjlim \pr^\la\hom_{\H_1'}(\R_m(i),  M)).
\end{align*}
This completes the argument.
\end{proof}

\Point{\boldmath Operators $e_i$ and $f_i$}\label{eded}
Continue with $\la \in P_+$ being fixed.
We wish to refine the definitions of the functors
$\res_i$ and $\ind_i$ 
to give operators, denoted
$e_i$ and $f_i$ respectively, 
from {\em irreducible} $\H_n^\la$-modules
to {\em isomorphism classes} 
of $\H_{n-1}^\la$- (resp. $\H_{n+1}^\la$-) modules.

Actually, $e_i$ is simply the restriction to $\Irr \H_n^\la$ of an operator
also denoted $e_i$ on the irreducible modules in $\rep_I \H_n$.
We define this first;
recall the definition of the module $L_m(i)$
from \ref{cov}. Let $M$ be an irreducible in $\rep_I \H_n$.
Let $\H_1'$ denote the subalgebra of $\H_n$ generated by $C_n, X_n^{\pm 1}$.
For each $m \geq 1$, we define an $\H_{n-1}$-module
\begin{equation}\label{hbar}
\overline\hom_{\H_1'}(L_m(i), M)
\end{equation}
as follows. If $M$ is of type $\Mtype$ or $i \neq 0$, this is simply the
space
$\hom_{\H_1'}(L_m(i), M)$
viewed as an $\H_{n-1}$-module in the same way as in Lemma~\ref{geq}.
But if $M$ is of type $\Qtype$ and $i = 0$, we can pick an odd involution
$\theta_M:M \rightarrow M$ and also have the odd involutions
$\theta_m:L_m(i) \rightarrow L_m(i)$ from (\ref{thetaj}).
Let
$$
\theta_M \otimes \theta_m
:\hom_{\H_1'}(L_m(i),  M)\rightarrow \hom_{\H_1'}(L_m(i),  M)
$$
denote the map defined by $((\theta_M \otimes \theta_m) f)(v) = 
(-1)^{\bar f} \theta_M(f(\theta_m v)).$
One checks that $(\theta_M\otimes\theta_m)^2 = 1$, 
hence the $\pm 1$-eigenspaces
of $\theta_M\otimes \theta_m$ split
$\hom_{\H_1'}(L_m(i), M)$ into a direct sum of two isomorphic 
$\H_{n-1}$-modules (because there is an obvious odd automorphism
swapping the two eigenspaces). Now in this case, we define
$\overline{\hom}_{\H_1'}(L_m(i), M)$
to be the ${1}$-eigenspace (say).

In either case, we have a direct system
$$
\overline{\hom}_{\H_1'}(L_{1}(i), M)
\hookrightarrow
\overline{\hom}_{\H_1'}(L_2(i), M)
\hookrightarrow
\dots
$$
induced by the inverse system (\ref{lseq}).
Now define
\begin{equation}\label{eid}
e_i M = 
\varinjlim \overline{\hom}_{\H_1'}(L_m(i), M),
\end{equation}
giving us the affine version of the operator $e_i$.
Note if $M$ is an $\H_n^\la$-module then 
each
$\overline{\hom}_{\H_1'}(L_m(i), M)$ is an $\H_{n-1}^\la$-module, so
\begin{equation}\label{eida}
e_i M =
\varinjlim \pr^\la \overline{\hom}_{\H_1'}(L_m(i), M).
\end{equation}
We take (\ref{eida}) as our definition of the operator $e_i$
in the cyclotomic case.
Comparing (\ref{eida}) 
with Lemma~\ref{geq} and using (\ref{zoo}), 
one sees at once that:

\vspace{1mm}
\begin{Lemma}\label{relate1} Let $i \in I$ and 
$M$ be an irreducible module in $\rep_I \H_n$, or an irreducible
$\H_n^\la$-module.
Then,
$$
\res_i M \simeq
\left\{
\begin{array}{ll}
e_i M&\hbox{if $i = 0$ and $M$ is of type $\Mtype$,}\\
e_i M \oplus \Pi e_i M&\hbox{otherwise.}
\end{array}
\right.
$$
\end{Lemma}

\vspace{1mm}

Now we turn to the definition of $f_i M$ which, just like $\ind_i M$,
only makes sense in the cyclotomic case.
So, let $M$ an irreducible $\H_n^\la$-module.
We need to extend the definition of the operation $\circledast$
to give meaning to the notation
$M \circledast L_m(i)$,
for each $m \geq 1$.
If either $M$ is of type $\Mtype$ or $i \neq 0$, then
$M \circledast L_m(i) := M \boxtimes L_m(i)$.
But if $M$ is of type $\Qtype$ and $i = 0$,
pick an odd involution $\theta_M:M \rightarrow M$.
Then, the $\pm \sqrt{-1}$-eigenspaces of
$\theta_M \otimes \theta_m$ acting on the left 
on $M \boxtimes L_m(i)$ split
it into a direct sum of two isomorphic $\H_{n,1}$-modules.
Let $M \circledast  L_m(i)$ denote
the $\sqrt{-1}$-eigenspace (say) for each $m$.

We then have an inverse system
$M \circledast L_1(i) \twoheadleftarrow 
 M \circledast  L_{2}(i) \twoheadleftarrow \dots$ of $\H_{n,1}$-modules
induced by the maps from (\ref{lseq}). Now we can define
\begin{equation}\label{fidef}
f_i M = \varprojlim \pr^\la \ind_{n,1}^{n+1}
 M \circledast L_m(i).
\end{equation}
Comparing the definition with the proof of Lemma~\ref{geq}, one sees that the
inverse limit stabilizes after finitely many terms, hence that
$f_i M$ really is a well-defined finite dimensional 
$\H_{n+1}^\la$-module.
Indeed:

\vspace{1mm}
\begin{Lemma}\label{relate2} Let $i \in I$ and 
$M$ be an irreducible $\H_n^\la$-module.
Then,
$$
\ind_i M \simeq
\left\{
\begin{array}{ll}
f_i M&\hbox{if $i = 0$ and $M$ is of type $\Mtype$,}\\
f_i M \oplus \Pi f_i M&\hbox{otherwise.}
\end{array}
\right.
$$
\end{Lemma}

\vspace{1mm}
\begin{Lemma} \label{eprops}
Let $i \in I$ and $M$ be an irreducible module in $\rep_I \H_n$.
Then,
$e_i M$ is non-zero if and only if
$\tilde e_i M \neq 0$,
in which case 
it is a 
self-dual indecomposable module with irreducible socle and cosocle 
isomorphic to $\tilde e_i M$.
\end{Lemma}

\begin{proof}
To see that $e_i M$ has irreducible
socle $\tilde e_i M$ whenever it is non-zero, 
combine Lemma~\ref{relate1} with Corollary~\ref{karg}.
The remaining facts follow since $M$ is self-dual by 
Lemma~\ref{selfdual}, and $\res_i$ commutes with duality
by Lemma~\ref{sleepy}(ii).
\end{proof}

Before the next theorem, we recall again that the operator $f_i$ 
depends critically on the fixed choice of $\la$.

\begin{Theorem}\label{crystallize} 
Let $\la \in P_+$ and $i \in I$.
Then, for any irreducible $\H_n^\la$-module $M$,
\begin{enumerate}
\item[(i)] $e_i M$ is non-zero if and only if
$\tilde e_i^\la M \neq 0$,
in which case 
it is a 
self-dual indecomposable module with irreducible socle and cosocle 
isomorphic to $\tilde e_i^\la M$;
\item[(ii)] $f_i M$ is non-zero if and only if
$\tilde f_i^\la M \neq 0$,
in which case 
it is a 
self-dual indecomposable module with irreducible socle and cosocle 
isomorphic to $\tilde f_i^\la M$.
\end{enumerate}
\end{Theorem}

\begin{proof}
(i) This is immediate from Lemma~\ref{eprops}.

(ii) We deduce this from (i) by an adjointness argument.
Let $M$ be an irreducible $\H_n^\la$-module, and $N$ be an
irreducible $\H_{n+1}^\la$-module.
Let $\delta_M$ equal $1$ if $i =0$ and $M$ is of type $\Mtype$,
$2$ otherwise, and define $\delta_N$ similarly.
Then, by Lemmas~\ref{sleepy}(i), \ref{relate1} and \ref{relate2},
\begin{align*}
\dim \hom_{\H_{n+1}^\la} (f_i M, N)
&=\frac{1}{\delta_M} \dim \hom_{\H_{n+1}^\la} (\ind_i M, 
N)\\
&=\frac{1}{\delta_M} \dim \hom_{\H_{n}^\la} (M, 
\res_i N)=\frac{\delta_N}{\delta_M} 
\dim \hom_{\H_{n}^\la} (M, e_i N).
\end{align*}
By (i), the latter is zero unless $M = \tilde e_i N$, 
or equivalently $N = \tilde f_i M$ by Lemma~\ref{L290900_2}.
Taking into account the superalgebra analogue
of Schur's lemma using Lemma~\ref{cccc}, one deduces that
$\cosoc f_i M \cong \tilde f_i M$.
Finally, note $f_i M$ is self-dual by Lemma~\ref{sleepy}(ii)
so everything else follows.
\end{proof}

\begin{Remark}\label{feed}\rm
Let us also point out, as follows easily from the definitions,
that $e_i M$ and $f_i M$ admit odd involutions if either
$i \neq 0$ and $M$ is of type $\Qtype$, or $i = 0$ and $M$ is of
type $\Mtype$.
\end{Remark}

\Point{Divided powers}\label{divp}
Continue with $\la \in P_+$ and fix $i \in I$.
We can generalize the definitions of $e_i, f_i$
to define operators denoted $e_{i}^{(r)}, f_{i}^{(r)}$ on irreducible
$\H_n^\la$-modules, for each $r \geq 1$. It will be the case that
$e_{i}^{(1)} = e_i, f_i^{(1)} = f_i$.
For the definitions, we make use of the covering modules
$L_m(i^r)$ from \ref{cov}.

Let $M$ be an irreducible module in $\rep_I \H_n$.
If $r > n$, we set $e_i^{(r)} M = 0$. Otherwise, let
$\H_r'$ denote the subalgeba of $\H_n$ generated
by $X^{\pm 1}_{n-r+1}, \dots, X_n^{\pm 1}$,
$C_{n-r+1}, \dots, C_{n}$, $T_{n-r+1}, \dots, T_{n-1}$.
We have a direct system
$$
\overline{\hom}_{\H_r'}(L_{1}(i^r), M)
\hookrightarrow 
\overline{\hom}_{\H_r'}(L_2(i^r), M)
\hookrightarrow 
\dots
$$
induced by the inverse system (\ref{lseq}), where
$\overline{\hom}$ is interpreted in exactly the same way as in
\ref{eded} using the generalized maps $\theta_m$ from (\ref{thetaj})
in case $i = 0$.
Now define
\begin{equation}\label{eird}
e_i^{(r)} M = 
\varinjlim \overline{\hom}_{\H_r'}(L_m(i^r), M).
\end{equation}
As in \ref{eded}, if $M$ is 
an $\H_n^\la$-module then $e_i^{(r)} M$ is too, so that
\begin{equation}\label{eirda}
e_i^{(r)} M =
\varinjlim \pr^\la \overline{\hom}_{\H_r'}(L_m(i^r), M)
\end{equation}
in the cyclotomic case.

To define $f_i^{(r)}$, which as usual {only} makes sense in the cyclotomic
case, let $M$ be an irreducible $\H_n^\la$-module. 
We have an inverse system
$M \circledast L_1(i^r) \twoheadleftarrow 
 M \circledast  L_{2}(i^r) \twoheadleftarrow \dots$ of $\H_{n,r}$-modules
induced by the maps from (\ref{lseq}),
again interpreting $\circledast$ as in \ref{eded}.
Now we can define
\begin{equation}\label{firda}
f_i^{(r)} M = \varprojlim \pr^\la \ind_{n,r}^{n+r}
 M \circledast L_m(i^r).
\end{equation}

\vspace{1mm}
\begin{Lemma}\label{divpower}
Let $i \in I, r \geq 1$ and $M$ be an irreducible $\H_n^\la$-module.
Then:
\begin{align*}
(\res_i)^r M &\simeq \left\{\begin{array}{ll}
(e_i^{(r)} M \oplus \Pi e_i^{(r)} M)^{\oplus 2^{r-1} (r!)}&\hbox{if $i \neq 0$,}\\
e_i^{(r)} M^{\oplus 2^{(r-1)/2} (r!)}&\hbox{if $i = 0$, $r$ is odd,
$M$ is of type $\Mtype$,}\\
(e_i^{(r)} M \oplus \Pi e_i^{(r)} M)^{\oplus 2^{\lfloor (r-1)/2\rfloor} (r!)}&\hbox{otherwise;}
\end{array}\right.\\
(\ind_i)^r M &\simeq \left\{\begin{array}{ll}
(f_i^{(r)} M \oplus \Pi f_i^{(r)} M)^{\oplus 2^{r-1} (r!)}&\hbox{if $i \neq 0$,}\\
f_i^{(r)} M^{\oplus 2^{(r-1)/2} (r!)}&\hbox{if $i = 0$, $r$ is odd,
$M$ is of type $\Mtype$,}\\
(f_i^{(r)} M \oplus \Pi f_i^{(r)} M)^{\oplus 2^{\lfloor (r-1)/2\rfloor} (r!)}&\hbox{otherwise.}
\end{array}\right.
\end{align*}
\end{Lemma}

\begin{proof}
Using Lemma~\ref{lmir} and the definitions, it suffices to show that
\begin{align*}
(\res_i)^r M &\simeq
\varinjlim \pr^\la \hom_{\H_r'}(\R_m(i^r), M),\\
(\ind_i)^r M &\simeq 
\varprojlim \pr^\la \ind_{n,r}^{n+r} M \boxtimes \R_m(i^r).
\end{align*}
For $(\res_i)^r$, this follows from Lemma~\ref{nastystab}
in exactly the same way as in the proof of Lemma~\ref{geq}.
Now $(\ind_i)^r$ is left adjoint to $(\res_i)^r$, 
so the statement for induction
follows from uniqueness of adjoint functors on checking that
the functor
$\varprojlim \pr^\la \ind_{n,r}^{n+r} ? \boxtimes \R_m(i^r)$
is left adjoint to
$\varinjlim \pr^\la \hom_{\H_r'}(\R_m(i^r), ?).$
The latter follows as in the proof of Lemma~\ref{geq}.
\end{proof}

Since we have defined the operator $e_i^{(r)}$ on irreducible
modules we get induced operators also denoted $e_i^{(r)}$
at the level of Grothendieck groups, namely,
$$
e_i^{(r)}: K(\rep_I \H_n) \rightarrow K(\rep_I \H_{n-r}),
\qquad
e_i^{(r)}: K(\rep \H_n^\la) \rightarrow K(\rep \H_{n-r}^\la),
$$
in the affine and cyclotomic cases respectively.
Similarly $f_i^{(r)}$ induces an operator 
$$
f_i^{(r)}: K(\rep \H_n^\la) \rightarrow K(\rep \H_{n+r}^\la)
$$
on Grothendieck groups.
We record:

\vspace{1mm}
\begin{Lemma}\label{dp}
As operators on the Grothendieck group $K(\rep \H_n^\la)$ (or on
$K(\rep_I \H_n)$ in the case of $e_i$), we have that
$e_i^r = (r!) e_i^{(r)}$ and
$f_i^r = (r!) f_i^{(r)}$.
\end{Lemma}

\begin{proof}
Let us prove in the case $i \neq 0$.
The proof for $i = 0$ is the same idea, though the details are 
more delicate (one needs to use Lemma~\ref{cccc} too).
By Lemmas~\ref{relate1} and \ref{relate2}, we have that
$$
(\res_i)^r = 2^r e_i^r, \qquad
(\ind_i)^r = 2^r f_i^r
$$
as operators on the Grothendieck group.
By Lemma~\ref{divpower}, we have that
$$
(\res_i)^r = 2^r (r!) e_i^{(r)}, \qquad
(\ind_i)^r = 2^r (r!) f_i^{(r)}.
$$
This is enough.
\end{proof}

Let us finally note that we have only defined the operators
$e_i^{(r)}$ and $f_i^{(r)}$ on irreducible modules.
However, the definitions could be made more generally
on pairs $(M, \theta_M)$, where $M$ is an $\H_n^\la$-module
(or an integral $\H_n$-module in the case of $e_i$) and
$\theta_M:M \rightarrow M$ is either the identity map
or else an odd involution of $M$.
In case $\theta_M = \id_M$, the definitions of
$e_i^{(r)} M$ and $f_i^{(r)} M$ are exactly the same as the
case $M$ is irreducible of type $\Mtype$ above.
In case $\theta_M$ is an odd involution, the definitions of
$e_i^{(r)} M$ and $f_i^{(r)} M$ are exactly the same as the
case $M$ is irreducible of type $\Qtype$ above, substiting the
given map $\theta_M$ for the canonical odd involution of $M$
in the situation above.

This remark applies especially to give us
modules $e_i^{(r)} P_M, f_i^{(r)} P_M$, where $P_M$ is the
projective cover of an irreducible $\H_n^\la$-module:
in this case, if $M$ is of type $\Qtype$, the odd involution
$\theta_M$ of $M$ lifts to a unique odd involution also denoted
$\theta_M$ of the projective cover.
On doing this, we have that
\begin{equation}\label{okg}
[e_i^{(r)} P_M] = e_i^{(r)} [P_M],
\qquad
[f_i^{(r)} P_M] = f_i^{(r)} [P_M]
\end{equation}
where the equalities are written in 
$K(\rep \H_{n-r}^\la)$ and $K(\rep \H_{n+r}^\la)$
respectively.
To prove this, one needs to observe that all composition
factors of $P_M$ are of the same type as $M$
by Lemma~\ref{cccc}.

Note Lemma~\ref{divpower} is also true if $M$ is replaced by its projective
cover $P_M$, the proof being the same as above.
In particular, this shows that
$e_i^{(r)} P_M$ is a summand of $(\res_i)^r P_M$,
and similarly for $f_i$. So Lemma~\ref{sleepy}(i) gives that
$e_i^{(r)} P_M$ and $f_i^{(r)} P_M$ are also projective modules.
Hence $e_i^{(r)}$ and $f_i^{(r)}$ induce operators with the same
names on the Grothendieck groups of {\em projective} modules too:
$$
e_i^{(r)}:K(\proj \H_n^\la) \rightarrow K(\proj \H_{n-r}^\la),
\qquad
f_i^{(r)}:K(\proj \H_n^\la) \rightarrow K(\proj \H_{n+r}^\la).
$$
Moreover, by the same argument as in the proof of Lemma~\ref{dp},
we have:

\vspace{1mm}
\begin{Lemma}\label{dp2}
As operators on the Grothendieck group $K(\proj \H_n^\la)$,
we have that
$$
e_i^r [P_M] = (r!) [e_i^{(r)} P_M],
\qquad
f_i^r [P_M] = (r!) [f_i^{(r)} P_M],
$$
for all irreducible $\H_n^\la$-modules $M$.
\end{Lemma}

\Point{\boldmath Alternative descriptions of $\eps_i$}\label{ar}
In this subsection we give three new interpretations of the functions 
$\eps_i$, precisely as in \cite[Theorem~9.13]{G}.

\vspace{1mm}
\begin{Theorem}
\label{T280900}
Let $i\in I$ and $M$ be an irreducible module in $\rep_I{\H_n}$.
Then
\begin{enumerate}
\item[(i)] $[e_i M] = \eps_i(M) [\tilde e_i M] + \sum c_a [N_a]$
where the $N_a$ are irreducibles with 
$\eps_i(N_a) < \eps_i(\tilde e_i M)$;
\item[(ii)] $\eps_i(M)$ is the maximal size of a Jordan block of
$X_n+X_n^{-1}$ (resp. $X_n$ if $i = 0$) on $M$ with eigenvalue $q(i)$
(resp. eigenvalue $1$ if $i = 0$);
\item[(iii)] $\End_{\H_{n-1}}(e_i M)
\simeq
\End_{\H_{n-1}}(\tilde e_i M)^{\oplus \eps_i(M)}$
as vector superspaces.
\end{enumerate}
\end{Theorem}

\begin{proof}
Let $\eps = \eps_i(M)$ and $N = \tilde e_i^\eps M$.

(i) 
By Lemma~\ref{L310800_1} and Frobenius reciprocity, there is
a short exact sequence
$$
0 \longrightarrow R \longrightarrow \ind_{n-\eps,\eps}^n N \circledast 
L(i^\eps) \longrightarrow M \longrightarrow 0.
$$
Moreover, all composition factors $L$ of $R$ have
$\eps_i(L) < \eps$ by Lemma~\ref{L010900}(iii).
Applying the exact functor $\De_i$, we obtain an exact sequence
$$
0 \longrightarrow \De_i R \longrightarrow 
\De_i \ind_{n-\eps,\eps} N \circledast L(i^\eps) \longrightarrow \De_i M
\longrightarrow 0.
$$
By the Mackey Theorem,
$\De_i\ind_{n-\eps,\eps}^{n} N\circledast L(i^\eps) \cong  
\ind_{n-\eps,\eps-1,1}^{n-1,1} N\circledast \De_i L(i^\eps).$
By considering characters
$[\De_i L(i^\eps)] = \eps [L(i^{\eps - 1}) \circledast L(i)]$.
Hence,
\begin{equation}\label{E300900_3}
[\De_i\ind_{n-\eps,\eps}^{n} N\circledast L(i^\eps)] =  
\eps[\ind_{n-\eps,\eps-1,1}^{n-1,1} 
N\circledast L(i^{\eps-1})\circledast L(i)].
\end{equation}
Using Lemma~\ref{L010900} again, 
the cosocle of $\ind_{n-\eps, \eps-1,1}^{n-1,1} N \circledast L(i^{\eps-1})
\circledast L(i)$
is $(\tilde e_i M) \circledast L(i)$, and all other composition
factors of this module are of the form $L \circledast L(i)$
with $\eps_i(L) < \eps - 1$.
Moreover, all composition factors of $\De_i R$
are of the form $L \circledast L(i)$
with $\eps_i(L) < \eps - 1$.
So we have now shown that
$$
[\De_i M] = \eps [\tilde e_i M \circledast L(i)] + \sum c_a [N_a 
\circledast L(i)]
$$
for irreducibles $N_a$ with $\eps_i(N_a) < \eps_i(\tilde e_i M)$.
The conclusion follows on applying Lemma~\ref{relate1}.

(ii) We give the argument for the case $i\neq 0$,
the case $i= 0$ being similar but using Lemma~\ref{L050900_1}(ii) instead
of
Lemma~\ref{L050900_1}(i).
We know that $\De_{i^\eps} M\cong N
\circledast L(i^\eps)$. 
So, applying the automorphism $\sigma$ to Lemma~\ref{L050900_1}(i),
we deduce that the maximal size of a Jordan block of 
$X_n+X_n^{-1}$ on $\De_{i^\eps} M$ 
is $\eps$. 
Hence the maximal size of a Jordan block of $X_n+X_n^{-1}$ on $\De_i M$ 
is at least $\eps$. 

On the other hand, 
the argument given above in deriving
(\ref{E300900_3}) shows that the module
$\De_i\ind_{n-\eps,\eps}^{n} N\circledast L(i^\eps)$
has a filtration with $\eps$ factors, each of which is isomorphic to
$\ind_{n-\eps,\eps-1,1}^{n-1,1} 
N\circledast L(i^{\eps-1})\circledast L(i)$.
Since $(X_n +X_n^{-1}- q(i))$
annihilates
$\ind_{n-\eps,\eps-1,1}^{n-1,1} 
N\circledast L(i^{\eps-1})\circledast L(i)$, it follows that
$(X_n +X_n^{-1}- q(i))^\eps$ annihilates
$\De_i\ind_{n-\eps,\eps}^{n} N\circledast L(i^\eps)$.
So certainly
$(X_n +X_n^{-1}- q(i))^\eps$ annihilates its quotient
$\De_i M$.
So the maximal size of a Jordan block of $X_n+X_n^{-1}$ 
on $\De_i M$ is at most $\eps$.

(iii) 
Let $z = (X_n+X_n^{-1} - q(i))$ if $i \neq 0$
and $(X_n - 1)$ if $i = 0$.
Consider the effect of left multiplication by $z$ on the $\H_{n-1}$-module
$R = \res^{n-1,1}_{n-1} \circ \De_i (M).$
Note $R$ is equal to either $e_i M$ or $e_i M \oplus \Pi e_i M$,
by Lemma~\ref{relate1}.
In the latter case, $(X_n+X_n^{-1})$ (resp. $X_n$) acts as a scalar
on $\soc R \simeq \tilde e_i M \oplus \Pi \tilde e_i M$, hence it leaves
the two indecomposable summands invariant.
This shows that in any case 
left multiplication by $z$ (which centralizes
the subalgebra $\H_{n-1}$ of $\H_n$) induces an $\H_{n-1}$-endomorphism
$\theta:e_i M \rightarrow e_i M$.
But by (ii), $\theta^{\eps-1} \neq 0$ and $\theta^\eps = 0$.
Hence, $1, \theta, \dots, \theta^{\eps-1}$ give $\eps$ linearly independent
even $\H_{n-1}$-endomorphisms of $e_i M$.
In view of Remark~\ref{feed}, we automatically get from these
$\eps$ linearly independent odd endomorphisms in case
$\tilde e_i M$ is of type $\Qtype$, so we have now shown that
$$
\dim \End_{\H_{n-1}}(e_i M) \geq \eps 
\dim \End_{\H_{n-1}}(\tilde e_i M).
$$
On the other hand, $e_i M$ has irreducible cosocle $\tilde e_i M$, 
and this appears in $e_i M$ with multiplicity $\eps$ by (i),
so the reverse inequality also holds.
\end{proof}

\begin{Corollary}\label{epr}
Let $M, N$ be irreducible $\H_n$-modules with $M \not\cong N$.
Then, for every $i \in I$,
$\hom_{\H_{n-1}}(e_i M, e_i N) = 0.$
\end{Corollary}

\begin{proof}
Suppose there is a non-zero homomorphism
$\theta:e_i M \rightarrow e_i N$.
Then, since $e_i M$ has simple head $\tilde e_i M$,
we see that $e_i N$ has $\tilde e_i M$ as a composition factor.
Hence, by Theorem~\ref{T280900}(i),
$\eps_i(\tilde e_i N) \geq \eps_i(\tilde e_i M)$.
Dualizing and applying the same argument gives the inequality the other way
round, hence $\eps_i(\tilde e_i N) = \eps_i(\tilde e_i M)$.
But then, $\tilde e_i M$ is a composition factor of $e_i N$
with $\eps_i(\tilde e_i M) = \eps_i(\tilde e_i N)$, hence by
Theorem~\ref{T280900}(i) again, $\tilde e_i M \cong \tilde e_i N$.
But this contradicts Corollary~\ref{imm}.
\end{proof}

To state the next corollary, we first need to introduce
the $*$-operation on the crystal graph.
This will play a fundamental role later on.
Suppose $M$ is an irreducible module in 
$\rep_I{\H_n}$ and $0 \leq m \leq n$. 
Using Lemma~\ref{L031000} 
for the second equality in (\ref{staropse}), define
\begin{align}
\tilde e_i^* M &= (\tilde e_i (M^\sigma))^\sigma,\\\label{staropse}
\tilde f_i^* M &= (\tilde f_i (M^\sigma))^\sigma
= \cosoc \ind_{1,n}^{n+1} L(i) \circledast M,\\\label{staropsf}
\eps_i^*(M) &= \eps_i(M^\sigma) =
\max\{m\geq 0\:|\:(\tilde e_i^*)^m M\neq 0\}. 
\end{align}
Note $\eps_i^*(M)$ can be worked out just from knowledge of
the character $M$:
$\eps_i^*(M)$ is the maximum $k$ such that
$[L(i)^{\circledast k} \circledast \dots]$
appears in $\ch M$.

Recalling 
the definition of the
ideal ${\mathcal I}_\la$ generated by the element (\ref{E210900_4}),
Theorem~\ref{T280900}(ii) has the following important corollary:

\vspace{1mm}
\begin{Corollary}
\label{P300900}
Let $\la \in P_+$ and
$M$ be an irreducible in $\rep_I{\H_n}$.
Then ${\mathcal I}_\la M=0$ if and only if $\eps_i^*(M)\leq 
\langle h_i, \la\rangle$ for all $i\in I$. 
\end{Corollary}

\Point{\boldmath Functions $\phi_i$}\label{crel}
Fix $\la \in P_+$ throughout this subsection.
Let $M$ be an irreducible $\H_n^\la$-module.
Recall from (\ref{E300900}) that for $i \in I$,
\begin{align}
\eps_i(M) &= \max\{m\geq 0\:|\: (\tilde e_i^\la)^m M\neq 0\},\label{ep}\\\intertext{since $\tilde e_i^\la$ is simply the restriction of $\tilde e_i$.
Analogously, we define}
\phi_i(M) &= \max\{m\geq 0\:|\: (\tilde f_i^\la)^m M\neq 0\}.\label{ph}
\end{align}
We will see shortly (Corollary~\ref{stickingpoint}) that 
$\phi_i(M) < \infty$ always so that the definition makes sense.
Note unlike $\eps_i(M)$, the integer $\phi_i(M)$ depends on the
fixed choice of $\lambda$.

As in \ref{cg}, $\bil$ denotes the irreducible $\H_0^\la$-module.

\vspace{1mm}
\begin{Lemma}\label{cough}
$\eps_i(\bil) = 0$
and $\phi_i(\bil) = \langle h_i,\la \rangle$.
\end{Lemma}

\begin{proof}
The statement involving $\eps_i$ is obvious. For
$\phi_i(\bil)$, note that
$\tilde f_i^m \bil = L(i^m)$ and
$$
\eps_i^*(L(i^m)) = m,\qquad
\eps_j^*(L(i^m)) = 0
$$
for every $j \neq i$.
Hence by Corollary~\ref{P300900}, $\pr^\la L(i^m) \neq 0$
if and only if $m \leq \langle h_i,\la \rangle$.
This implies that $\phi_i(\bil) = \langle h_i,\la \rangle$.
\end{proof}

\begin{Lemma}\label{stete}
Let $i,j\in I$ with $i \neq j$ and $M$ be an irreducible module in $\rep_I
\H_n$. Then, 
$\eps_j^*(\tilde f_i^m M) \leq \eps_j^*(M)$ for every $m \geq 0$. 
\end{Lemma}

\begin{proof}
Follows from the Shuffle Lemma.
\end{proof}

\begin{Lemma}\label{cre2}
Let $i,j \in I$ with $i \neq j$.
Let $M$ be an irreducible module in $\rep \H_n^\la$ such that
$\phi_j(M) > 0$.
Then,
$\phi_i(\tilde f_j M)-\eps_i(\tilde f_j M)
\leq \phi_i(M) - \eps_i(M)
-\langle h_i,\alpha_j \rangle.$
\end{Lemma}

\begin{proof}
Let $\eps = \eps_i(M), \phi = \phi_i(M)$ and $k=-\langle
h_i,\alpha_j\rangle$.
By Lemma~\ref{serre3}, there exist unique $a,b \geq 0$ with
$a+b = k$ such that
$\eps_i (\tilde f_j M) 
= \eps - a.$
We need to show that
$\phi_i(\tilde f_j M) \leq \phi +b$, 
which follows if we can show that
$\pr^\la \tilde f_i^m \tilde f_j M = 0$
for all $m > \phi+b$. We claim that
\begin{equation*}
\eps_i^*(\tilde f_i^{m+b} \tilde f_j M) \geq 
\eps_i^*(\tilde f_i^m M)
\end{equation*}
for all $m \geq 0$.
Given the claim, we know by the definition of $\phi$,
Corollary~\ref{P300900}
and Lemma~\ref{stete} that
$\eps_i^*(\tilde f_i^{m} M) > \langle h_i,\la \rangle$
for all $m > \phi$.
So the claim implies that
$\eps_i^*(\tilde f_i^{m} \tilde f_j M) >\langle h_i,\la \rangle$
for all $m > \phi+b$, hence by Corollary~\ref{P300900} once more,
$\pr^\la \tilde f_i^m \tilde f_j M = 0$
as required.

To prove the claim, note that $a \leq \eps$, so
$b+\eps \geq k$. Hence, Lemma~\ref{serre3}(ii) shows that there is a
surjection
$$
\ind_{n,m-a,a+b+1}^{n+m+b+1} 
M \circledast L(i^{m-a})
\circledast L(i^a j i^b)
\twoheadrightarrow
\tilde f_i^{m+b} \tilde f_j M.
$$
By Lemma~\ref{calc1},
$\res^{a+b+1}_{a,b+1} L(i^a j i^b) \cong L(i^a) \circledast L(j i^b).$
Hence by Frobenius reciprocity, there is a surjection
$$
\ind_{a,b+1}^{a+b+1} L(i^a) \circledast L(j i^b) \twoheadrightarrow
L(i^a j i^b).
$$
Combining, we have proved existence of a surjection
$$
\ind_{n,m,b+1}^{n+m+b+1} 
M \circledast L(i^{m})
\circledast L(j i^b)
\twoheadrightarrow
\tilde f_i^{m+b} \tilde f_j M.
$$
Hence by Frobenius reciprocity there is a non-zero map
$$
(\ind_{n,m}^{n+m} M \circledast L(i^{m}) )
\circledast L(j i^b)
\rightarrow
\res^{n+m+b+1}_{n+m,b+1} \tilde f_i^{m+b} \tilde f_j M.
$$
Since the left hand module has irreducible cosocle
$\tilde f_i^{m} M \circledast L(j i ^b)$, 
we deduce that
$\tilde f_i^{m+b} \tilde f_j M$ 
has a constituent isomorphic to
$\tilde f_i^m M$ on restriction to the subalgebra 
$\H_{n+m} \subseteq \H_{n+m+b+1}$.
This implies the claim.
\end{proof}

\begin{Corollary}\label{stickingpoint}
Let $\la \in P_+$ and $M$ be an irreducible $\H_n^\la$-module
with central character $\chi_\gamma$ for some $\gamma \in \Gamma_n$.
Then,
$\phi_i(M) - \eps_i(M) \leq \langle
h_i, \la - \gamma \rangle.$
\end{Corollary}

\begin{proof}
Proceed by induction on $n$, the case $n = 0$ being immediate
by Lemma~\ref{cough}.
For $n > 0$, we may write $M = \tilde f_j N$
for some irreducible $\H_{n-1}^\la$-module $N$ with $\phi_j(N) > 0$.
By induction,
$\phi_i(N) - \eps_i(N) \leq \langle h_i, \la - \gamma + \alpha_j\rangle.$
The conclusion follows from Lemma~\ref{cre2}.
\end{proof}

\Point{\boldmath Alternative descriptions of $\phi_i$}
Keep $\la \in P_+$ fixed.
Now we wish to prove the analogue of Theorem~\ref{T280900}
for the function $\phi_i$.
This is considerably more difficult to do.
Let $M$ be an irreducible $\H_n^\la$-module.
Recall that
$$
f_i M = \varprojlim \pr^\la \ind_{n,1}^{n+1}
 M \circledast L_m(i),
\qquad
\ind_i M = \varprojlim \pr^\la \ind_{n,1}^{n+1}
 M \boxtimes \R_m(i),
$$
and that the inverse limits stabilize after finitely many terms.
Define $\tilde \phi_i(M)$ to be the stabilization point of the limit, i.e. 
the least $m \geq 0$ such that
$f_i M = \pr^\la \ind_{n,1}^{n+1} M \circledast L_m(i)$, or equivalently
$\ind_i M = \pr^\la \ind_{n,1}^{n+1} M \boxtimes \R_m(i).$
The first lemma follows \cite[Theorem 9.15]{G}.

\vspace{1mm}
\begin{Lemma}\label{phiprops}
Let $M$ be an irreducible $\H_n^\la$-module and $i \in I$.
Then:
\begin{enumerate}
\item[(i)] $[f_i M] = \tilde\phi_i(M) [\tilde f_i M] + \sum c_a [N_a]$
where the $N_a$ are irreducible $\H_{n+1}^\la$-modules with
$\eps_i(N_a) < \eps_i(\tilde f_i M)$;
\item[(ii)] $\End_{\H_{n+1}^\la}(f_i M)
\simeq \End_{\H_{n+1}^\la}(\tilde f_i M)^{\oplus \tilde\phi_i(M)}$ as
vector superspaces.
\end{enumerate}
\end{Lemma}

\begin{proof}
(i) 
Take any $m \geq 1$.
Since $\pr^\la$ is right exact, the
natural surjection
$L_m(i) \twoheadrightarrow L_{m-1}(i)$ and
the obvious embedding $L_m(i) \hookrightarrow L_{m+1}(i)$ 
(see \ref{cov}) induce a commutative diagram
\ifbook@
$$
\begin{CD}
\pr^\la \ind_{n,1}^{n+1} M \circledast L(i)
&@>\alpha_m>> \pr^\la \ind_{n,1}^{n+1} M \circledast L_m(i)
&@>\beta_m>> \pr^\la \ind_{n,1}^{n+1} M \circledast L_{m-1}(i)\\
@|&@VVV\\
\pr^\la \ind_{n,1}^{n+1} M \circledast L(i)
&@>\alpha_{m+1}>> \pr^\la \ind_{n,1}^{n+1} M \circledast L_{m+1}(i)
&@>\beta_{m+1}>> \pr^\la \ind_{n,1}^{n+1} M \circledast L_{m}(i)
\end{CD}
$$
where the rows are exact and $\beta_m,\beta_{m+1}$ are epimorphisms.
\else
$$
\begin{CD}
\pr^\la \ind_{n,1}^{n+1} M \circledast L(i)
&@>\alpha_m>> \pr^\la \ind_{n,1}^{n+1} M \circledast L_m(i)
&@>\beta_m>> \pr^\la \ind_{n,1}^{n+1} M \circledast L_{m-1}(i)
\rightarrow \:0\\
@|&@VVV\\
\pr^\la \ind_{n,1}^{n+1} M \circledast L(i)
&@>\alpha_{m+1}>> \pr^\la \ind_{n,1}^{n+1} M \circledast L_{m+1}(i)
&@>\beta_{m+1}>> \pr^\la \ind_{n,1}^{n+1} M \circledast L_{m}(i)
\quad\rightarrow \:0
\end{CD}
$$
where the rows are exact.
\fi
Note if $\alpha_m = 0$ then $\alpha_{m+1} = 0$.
It follows that if $\beta_m$ is
an isomorphism so is $\beta_{m'}$ for every $m' \geq m$.
So by definition of $\tilde \phi_i(M)$,
the maps 
$\beta_1,\beta_2,\dots, \beta_{\tilde \phi_i(M)}$ are not isomorphisms
but all other $\beta_{m'}, m' > \tilde \phi_i(M)$ are isomorphisms.

Now to prove (i), we show by induction on
$m = 0,1,\dots,\tilde \phi_i(M)$ that
$$
[\pr^\la \ind_{n,1}^{n+1} M \circledast L_m(i)] = m[\tilde f_i M] + 
\hbox{lower terms},
$$
where the lower terms are irreducible $\H_{n+1}^\la$-modules $N$ with
$\eps_i(N) < \eps_i(\tilde f_i M)$.
This is vacuous if $m = 0$.
For $m > 0$, $\beta_m$ is not an isomorphism, so $\alpha_m \neq 0$.
Hence, by Lemma~\ref{L090900_1}, 
the image of $\alpha_m$ contains a copy of
$\tilde f_i M$ plus lower terms. Now the induction step is immediate.

(ii)
Take $m = \tilde \phi_i(M)$.
One easily shows 
using the explicit construction of $L_m(i)$ in \ref{cov} that
there is an even endomorphism
$\theta:L_m(i) \rightarrow L_m(i)$
of $\H_{1}$-modules,
such that the image of $\theta^k$ is $\simeq L_{m-k}(i)$ for each 
$0 \leq k \leq m$.
Frobenius reciprocity induces
superalgebra homomorphisms
$$
\End_{\H_1}(L_m(i))
\hookrightarrow \End_{\H_{n,1}}(M \circledast L_m(i))
\hookrightarrow \End_{\H_{n+1}}(\ind_{n,1}^{n+1} M \circledast L_m(i)).
$$
So $\theta$ induces an even $\H_{n+1}$-endomorphism
$\tilde\theta$ of $\ind_{n,1}^{n+1} M \circledast L_m(i)$,
such that the image of
$\tilde\theta^k$
is $\simeq\ind_{n,1}^{n+1} M \circledast L_{m-k}(i)$
for $0 \leq k \leq m$.
Now apply the right exact
functor $\pr^\la$ to get 
an even $\H_{n+1}^\la$-endomorphism
$$
\hat\theta:\pr^\la \ind_{n,1}^{n+1} M \circledast L_m(i)
\rightarrow \pr^\la \ind_{n,1}^{n+1} M \circledast L_m(i)
$$
induced by $\tilde \theta$.
Note $\hat\theta^m = 0$ and $\hat\theta^{m-1} \neq 0$
because its image coincides with the image of the non-zero map
$\alpha_m$ in the proof of (i).
Hence, $1, \hat\theta,\dots,\hat\theta^{m-1}$ are linearly independent,
even endomorphisms of $f_i M$.
Now the proof of (ii) is completed in the same way as in the proof of
Theorem~\ref{T280900}(iii).
\end{proof}

\begin{Corollary}\label{epr2}
Let $M, N$ be irreducible $\H_n^\la$-modules with $M \not\cong N$.
Then, for every $i \in I$,
$\hom_{\H_{n+1}^\la}(f_i M, f_i N)
= 0.$
\end{Corollary}

\begin{proof}
Repeat the argument in the proof of Corollary~\ref{epr},  but using
$\tilde \phi_i$ and Lemma~\ref{phiprops}(i) in place of
$\eps_i$ and Theorem~\ref{T280900}(i).
\end{proof}

The main thing now
is to prove that $\tilde \phi_i(M) = \phi_i(M)$.
Note right away from the definitions
that 
$\tilde\phi_i(M) = 0$ if and only if $\phi_i(M)= 0$.

\vspace{1mm}
\begin{Lemma}\label{phe}
If $M$ is an irreducible $\H_n^\la$-module then
$$
(\tilde \phi_0(M) - \eps_0(M)) + 
2 \sum_{i=1}^{\ell} (\tilde \phi_i(M) - \eps_i(M))
= \langle c, \la \rangle.
$$
\end{Lemma}

\begin{proof}
By Frobenius reciprocity and Theorem~\ref{cycloMackey},
we have that
\begin{align*}
\End_{\H_{n+1}^\la} (\ind_{\H_n^\la}^{\H_{n+1}^\la} M)
&\simeq
\hom_{\H_{n}^\la} (M, \res_{\H_n^\la}^{\H_{n+1}^\la}\ind_{\H_n^\la}^{\H_{n+1}^\la} M)\\
&\simeq
\hom_{\H_n^\la}(M, M \oplus \Pi M)^{\oplus \langle c,\lambda\rangle}
\oplus
\hom_{\H_{n}^\la} (M, \ind_{\H_{n-1}^\la}^{\H_{n}^\la}\res_{\H_{n-1}^\la}^{\H_{n}^\la} M)\\
&\simeq
\hom_{\H_n^\la}(M, M \oplus \Pi M)^{\oplus \langle c,\lambda\rangle}
\oplus
\End_{\H_{n-1}^\la} (\res_{\H_{n-1}^\la}^{\H_n^\la}M).
\end{align*}
Hence, by Schur's lemma,
$$
\dim \End_{\H_{n+1}^\la}(\ind_{\H_n^\la}^{\H_{n+1}^\la} M)
-\dim \End_{\H_{n-1}^\la}(\res_{\H_{n-1}^\la}^{\H_{n}^\la} M)
=
 \left\{
\begin{array}{ll}
2\langle c,\lambda\rangle
&\hbox{if $M$ is of type $\Mtype$,}\\
4 \langle c,\lambda\rangle
&\hbox{if $M$ is of type $\Qtype$.}
\end{array}
\right.
$$
Now if $M$ is of type $\Mtype$,
$$
\ind_{\H_n^\la}^{\H_{n+1}^\la} M
\simeq f_0 M \oplus \bigoplus_{i=1}^{\ell} (f_i M \oplus \Pi f_i M),
\qquad
\res_{\H_{n-1}^\la}^{\H_{n}^\la} M
\simeq e_0 M \oplus \bigoplus_{i=1}^{\ell} (e_i M \oplus \Pi e_i M).
$$
by Lemmas~\ref{relate1}, \ref{relate2} and (\ref{indresdecomp}).
Hence by Lemma~\ref{phiprops}(ii) and Theorem~\ref{T280900}(iii),
\begin{align*}
\dim \End_{\H_{n+1}^\la}(\ind_{\H_n^\la}^{\H_{n+1}^\la} M)
&=2\tilde \phi_0(M) + 4 \sum_{i=1}^\ell \tilde \phi_i(M),\\
\dim \End_{\H_{n-1}^\la}(\res_{\H_{n-1}^\la}^{\H_{n}^\la} M)
&=2\eps_0(M) + 4 \sum_{i=1}^\ell \eps_i(M),
\end{align*}
and the conclusion follows in this case. The argument for
$M$ of type $\Qtype$ is similar.
\end{proof}

\begin{Lemma}\label{cnt}
Let $M$ be an irreducible $\H_n^\la$-module and $i \in I$.
Let $c_i = 1$ if $i = 0$, $c_i = 2$ otherwise.
Then, 
\begin{align*}
[\res_i \ind_i M:M] = 2 c_i \eps_i(\tilde f_i M)\tilde \phi_i(M),
\qquad
&[\ind_i \res_i M:M] = 2 c_i \eps_i(M)\tilde \phi_i(\tilde e_i M),\\
\soc \res_i \ind_i M
 \simeq (M\oplus \Pi M)^{\oplus c_i \tilde \phi_i(M)},\qquad
&\soc \ind_i \res_i M  \simeq
(M\oplus \Pi M)^{\oplus c_i \eps_i(M)}.
\end{align*}
\end{Lemma}

\begin{proof}
The statement about composition multiplicities follows from
Theorem~\ref{T280900}(i) and Lemma~\ref{phiprops}(i), taking into account how
$\res_i$ and $\ind_i$ are related to $e_i$ and $f_i$ as explained in
Lemmas~\ref{relate1} and \ref{relate2}.
Now consider the statement about socles. We consider only
$\res_i \ind_i M$, the other case being entirely similar but using results
from \ref{ar} instead.
By adjointness, it suffices to be able to compute
$\hom_{\H_{n+1}^\la}(\ind_i N, \ind_i M)$
for any irreducible $\H_n^\la$-module $N$.
But in view of Lemma~\ref{relate2}, this can be computed from knowledge
of $\hom_{\H_{n+1}^\la}(f_i N, f_i M)$
which is known by Corollary~\ref{epr2} (if $M \not\cong N$) and
Lemma~\ref{phiprops}(ii) (if $M \cong N$).
The details are similar to those in the proof of Lemma~\ref{phe}, so we
omit them.
\end{proof}

The proof of the next lemma is based on \cite[Lemma 6.1]{V}.

\vspace{1mm}
\begin{Lemma}\label{vaz}
Let $M$ be an irreducible $\H_n^\la$-module and $i \in I$.
There are maps
$$
\ind_i\res_i M
\stackrel{\psi}{\longrightarrow}
\res_i\ind_i M
\stackrel{\operatorname{can}}{\longrightarrow}
\res_i \ind_i M / 
\soc \res_i \ind_i M,
$$
whose composite is surjective.
\end{Lemma}

\begin{proof}
Let $k = \tilde\phi_i(M)$ and
$$
\pi:\ind_{n,1}^{n+1} M \boxtimes \R_{k}(i) \twoheadrightarrow 
\pr^\la \ind_{n,1}^{n+1} M \boxtimes \R_{k}(i) = \ind_i M
$$
be the quotient map.
Let $\H_1'$ denote the subalgebra of $\H_n$ generated by $X_n^{\pm 1},
C_n$,
and set $z = X_n+X_n^{-1} - q(i)$ (resp. $X_n - 1$ if $i = 0$).
Recall from \ref{cov} that viewed as an $\H_1'$-module, we have that
$\R_k(i) \simeq \H_1' / (z^k)$.
In particular, $\R_k(i)$ is a cyclic module generated by the image
$\tilde 1$ of $1 \in \H_1'$.

We first observe that for any $m \geq \eps_i(M) + \tilde \phi_i(M)$,
$z^m$ annihilates the vector
$$\pi[T_n \otimes (u \otimes v)]
\in \ind_i M$$
for any $u \in M, v \in \R_{k}(i)$.
This follows from the relations in $\H_{n+1}$,
e.g. in case $i \neq 0$ one ultimately appeals to the facts that
$(X_n+X_n^{-1} - q(i))^{\eps_i(M)}$ 
annihilates $u$ (see Theorem~\ref{T280900}(ii)) and $(X_{n+1} +
X_{n+1}^{-1} - q(i))^{\tilde \phi_i(M)}$
annihilates $v$.
It follows that the unique $\H_{n-1,1}$-homomorphism
$(\res_i M) \boxtimes \H_1' \rightarrow \res^n_{n-1,1} \res_i \ind_i M$
such that
$u \otimes 1 \mapsto \pi[T_n \otimes (u\otimes \tilde 1)]$
for each $u \in \res_i M \subseteq M$ factors 
to induce a well-defined $\H_{n-1,1}$-module homorphism
$(\res_i M) \boxtimes \R_m(i) \rightarrow \res^n_{n-1,1} \res_i \ind_i M.$
We then get from Frobenius reciprocity an induced map
\begin{equation}\label{zxzx}
\psi_m:\ind_{n-1,1}^n (\res_i M) \boxtimes \R_m(i)
\rightarrow \res_i \ind_i M
\end{equation}
for each $m \geq \eps_i(M) + \tilde \phi_i(M)$.
Each $\psi_m$ factors through the quotient
\ifbook@
$$\pr^\la \ind_{n-1,1}^n (\res_i M) \boxtimes \R_m(i),
$$
\else
$\pr^\la \ind_{n-1,1}^n (\res_i M) \boxtimes \R_m(i)$,
\fi
so we get an induced map
$$
\psi:\ind_i \res_i M = \varprojlim \pr^\la \ind_{n-1,1}^n (\res_i M)
\boxtimes \R_m(i)
\rightarrow \res_i \ind_i M.
$$
It remains to show that the composite of $\psi$ with the 
canonical epimorphism from $\res_i \ind_i M$
to 
$\res_i \ind_i M / \soc \res_i \ind_i M$ is surjective.

Let $x=(n,n+1)\in S_{n+1}$. By Mackey Theorem there exists an exact
sequence
$$
0 \rightarrow M \boxtimes \R_k(i)
\rightarrow \res^{n+1}_{n,1} (\ind_{n,1}^{n+1} M \boxtimes \R_k(i))
\rightarrow \ind_{n-1,1,1}^{n,1} {^x}((\res^{n}_{n-1,1} M) \boxtimes
\R_k(i))
\rightarrow 0.
$$
In other words, there is an $\H_{n,1}$-isomorphism
\begin{eqnarray*}
\ind_{n-1,1,1}^{n,1} {^x}((\res^{n}_{n-1,1} M) \boxtimes \R_k(i))
 & \overset{\text{$\sim$}}{\longrightarrow} & 
\res^{n+1}_{n,1} (\ind_{n,1}^{n+1} M \boxtimes \R_k(i))
/(M \boxtimes \R_k(i)), \\
h\otimes (u\otimes v) & \mapsto & hT_n\otimes u\otimes v + M \boxtimes
\R_k(i)
\end{eqnarray*}
for $h\in\H_n,u\in M, v\in \R_k(i)$, where $M \boxtimes \R_k(i)$ is
embedded into $\res^{n+1}_{n,1} (\ind_{n,1}^{n+1} M \boxtimes \R_k(i))$ as
$1\otimes M \otimes \R_k(i)$. 
Recall from (\ref{dimm}) that
$\dim \R_k(i) = 2 k c_i$ where $c_i$ is as in Lemma~\ref{cnt}.
Hence, applying Lemma~\ref{cnt},
$$
\res^{n,1}_{n} M \boxtimes \R_k(i) \simeq (M \oplus \Pi M)^{k c_i}
\simeq \soc \res_i \ind_i M.
$$
So applying the exact functor $\res_i = \res^{n,1}_n \circ \De_i$
to the isomorphism above we get an isomorphism 
\begin{eqnarray*}
\ind_{n-1,1}^{n} (\res_i M) \boxtimes \R_k(i)
& \overset{\text{$\sim$}}{\longrightarrow} &
\res_i (\ind_{n,1}^{n+1} M \boxtimes \R_k(i))/\soc \res_i \ind_i M, \\
h\otimes u\otimes v 
& \mapsto & 
h T_n\otimes u\otimes v + \soc \res_i \ind_i M.
\end{eqnarray*}
It follows that there is a surjection
$$
\theta:\ind_{n-1,1}^n (\res_i M) \boxtimes \R_k(i)
\twoheadrightarrow \res_i \ind_i M / \soc \res_i \ind_i M
$$
such that the diagram
$$
\begin{CD}
\ind_{n-1,1}^n (\res_i M) \boxtimes \R_m(i)
&@>\psi_m>>&
\res_i \ind_i M\\
@VVV&&@VV\operatorname{can}V\\
\ind_{n-1,1}^n (\res_i M) \boxtimes \R_k(i)
&@>\theta>>&
\res_i \ind_i M / \soc \res_i \ind_i M
\end{CD}
$$
commutes for all $m \geq \eps_i(M) + \tilde \phi_i(M)$,
where $\psi_m$ is the map from (\ref{zxzx})
and the left hand arrow is the natural surjection.
Now surjectivity of $\theta$ immediately
implies the desired surjectivity of
the composite.
\end{proof}

\begin{Lemma}\label{hooray}
Let $M$ be an irreducible $\H_n^\la$-module
with $\eps_i(M) > 0$.
Then,
$$
\tilde \phi_i(\tilde e_i M) = \tilde \phi_i(M)+1.
$$
\end{Lemma}

\begin{proof}
Let us first show that
\begin{equation}\label{step1}
\tilde \phi_i(\tilde e_i M) \geq \tilde \phi_i(M)+1.
\end{equation}
Recall $\phi_i(M) = 0$ if and only if $\tilde \phi_i(M) = 0$.
Suppose first that  $\phi_i(M) = 0$, when
$\phi_i(\tilde e_i M) \neq 0$.
Then,
$\tilde \phi_i(M) = 0$ and
$\tilde \phi_i(\tilde e_i M) \neq 0
$, so the conclusion certainly holds in this case.
So we may assume that $\phi_i(M) > 0$, hence $\tilde \phi_i(M) \neq 0$.
Note by Lemma~\ref{cnt}, 
$$
[\res_i \ind_i M / \soc \res_i \ind_i M:M] = 
2 c_i \eps_i(\tilde f_i M)\tilde \phi_i(M)
-
2 c_i \tilde \phi_i(M) = 2 c_i \eps_i(M) \tilde \phi_i(M) \neq 0.
$$
In particular,
the map $\psi$ in Lemma~\ref{vaz} is non-zero.
Now Lemma~\ref{vaz} implies that the multiplicity of $M$ in $\im\:\psi$
is {\em strictly} greater than $2 c_i \eps_i (M) \tilde \phi_i(M)$,
since at least one composition factor
of $\soc \im\:\psi \subseteq \soc \res_i \ind_i M$ must be sent to zero
on composing with the second map $\operatorname{can}$.
Using another part of Lemma~\ref{cnt}, this shows that
$$
2 c_i \eps_i(M) \tilde \phi_i(\tilde e_i M)
> 
2 c_i \eps_i (M) \tilde \phi_i(M)
$$
and (\ref{step1}) follows.

Now using (\ref{step1}) and Lemma~\ref{cnt}, we see that in 
the Grothendieck group,
$$
[\res_i \ind_i M - \ind_i \res_i M:M]
\leq 2 c_i (\tilde \phi_i(M) - \eps_i(M)),
$$
with equality if and only if equality holds in (\ref{step1}).
By central character considerations, for $i \neq j$,
$[\res_i \ind_j M :M] = [\ind_j \res_i M:M] = 0$.
So using (\ref{indresdecomp}) we deduce that
\ifbook@
\begin{multline*}
\qquad
[\res^{\H_{n+1}^\la}_{\H_n^\la} \ind^{\H_{n+1}^\la}_{\H_n^\la} 
M-\ind_{\H_{n-1}^\la}^{\H_{n}^\la} \res^{\H_n^\la}_{\H_{n-1}^\la} M:M] \leq 
\\2 (\tilde \phi_0(M) - \eps_0(M))
+ 
4 \sum_{i=1}^{\ell} (\tilde \phi_i(M) - \eps_i (M))
\qquad
\end{multline*}
\else
$$
[\res^{\H_{n+1}^\la}_{\H_n^\la} \ind^{\H_{n+1}^\la}_{\H_n^\la} 
M-\ind_{\H_{n-1}^\la}^{\H_{n}^\la} \res^{\H_n^\la}_{\H_{n-1}^\la} M:M] \leq 2 (\tilde \phi_0(M) - \eps_0(M))
+ 4 \sum_{i=1}^{\ell} (\tilde \phi_i(M) - \eps_i (M))
$$
\fi
with equality if and only if equality holds in (\ref{step1}) for all 
$i \in I$.
Now Lemma~\ref{phe} shows that the right hand side equals
$2 \langle c,\lambda \rangle$, which does indeed equal the left hand side
thanks to Theorem~\ref{cycloMackey}.
\end{proof}

\begin{Corollary}\label{phiephi}
For any irreducible $\H_n^\la$-module $M$,
$\phi_i(M) = \tilde \phi_i(M)$.
\end{Corollary}

\begin{proof}
We proceed by induction on $\phi_i(M)$, the conclusion being known already
in case $\phi_i(M) = 0$.
For the induction step, take an irreducible $\H_n^\la$-module $N$ with
$\phi_i(N) > 0$, so $N = \tilde e_i M$ where $M = \tilde f_i N$ is
an irreducible $\H_{n+1}^\la$-module with $\eps_i(M) > 0$, 
$\phi_i(M) < \phi_i(N)$.
Then by Lemma~\ref{hooray} and the induction hypothesis,
$$
\tilde \phi_i(N)
= \tilde \phi_i(\tilde e_i M) = 
\tilde \phi_i(M)+1 = \phi_i(M) + 1 = \phi_i(\tilde e_i M)
= \phi_i(N).
$$
This completes the induction step.
\end{proof}

As a first consequence, we can improve Corollary~\ref{stickingpoint}:

\vspace{1mm}
\begin{Lemma}\label{unstuck}
Let $M$ be an irreducible $\H_n^\la$-module of central character 
$\chi_\gamma$ for $\gamma \in \Gamma_n$.
Then,
$\phi_i(M) - \eps_i(M) = \langle
h_i, \la - \gamma \rangle.$
\end{Lemma}

\begin{proof}
In view of Corollary~\ref{stickingpoint}, it suffices to show that
$$
(\phi_0(M) - \eps_0(M))
+ 2 \sum_{i=1}^{\ell}
(\phi_i(M) - \eps_i(M))
= 
\langle c, \la \rangle.
$$
But this is immediate from Lemma~\ref{phe} and Corollary~\ref{phiephi}.
\end{proof}

We are finally ready to assemble the results of the subsection
to obtain the full analogue of Theorem~\ref{T280900}
for $\phi_i$:

\vspace{1mm}
\begin{Theorem}\label{weredone}
Let $i\in I$ and $M$ be an irreducible $\H_n^\la$-module.
Then,
\begin{enumerate}
\item[(i)] $[f_i M] = \phi_i(M) [\tilde f_i M] + \sum c_a [N_a]$
where the $N_a$ are irreducibles with 
$\phi_i(N_a) < \phi_i(\tilde f_i M)$;
\item[(ii)] $\phi_i(M)$ is the least
$m \geq 0$ such that $f_i M =
\pr^\la \ind_{n,1}^{n+1} M \circledast L_m(i)$;
\item[(iii)] $\End_{\H_{n+1}}(f_i M)
\simeq
\End_{\H_{n+1}}(\tilde f_i M)^{\oplus \phi_i(M)}$
as vector superspaces.
\end{enumerate}
\end{Theorem}

\begin{proof}
(i) Since $\phi_i(M) = \tilde \phi_i(M)$ by Corollary~\ref{phiephi},
we know by Lemma~\ref{phiprops}(i) that
$$
[f_i M] = \phi_i(M) [\tilde f_i M] + \sum c_a [N_a]
$$
where the $N_a$ are irreducibles with 
$\eps_i(N_a) < \eps_i(\tilde f_i M)$.
Suppose that $M \in \rep_\gamma \H_n^\la$, for $\gamma \in \Gamma_n$.
Then $\tilde f_i M$ and each $N_a$ have central character
$\chi_{_{\gamma + \alpha_i}}$, since they are all composition factors
of $\ind_i M$.
So by Lemma~\ref{unstuck},
$$
\phi_i(N_a) = \langle h_i, \lambda - \gamma - \alpha_i \rangle + \eps_i(N_a),
\qquad
\phi_i(\tilde f_i M) = \langle h_i, \lambda - \gamma - \alpha_i \rangle + \eps_i(\tilde f_i M).
$$
It follows that $\phi_i(N_a) < \phi_i(\tilde f_i N_a)$ too.

(ii) This
is just the definition of $\tilde \phi_i(M)$
combined with Corollary~\ref{phiephi}.

(iii) This follows from Lemma \ref{phiprops}(ii).
\end{proof}

\ifbook@\pagebreak\fi

\section{Construction of $U_\Z^+$}

\Point{Grothendieck groups revisited}
Let us now write
\begin{align}\label{kinfty}
K(\infty) &= \bigoplus_{n \geq 0} K(\rep_I{\H_n})\\\intertext{for 
the sum over all $n$ of the Grothendieck groups of the categories
$\rep_I{\H_n}$. Also, write}
K(\infty)_{\Q} &= \Q \otimes_{\Z} K(\infty),
\end{align}
extending scalars.
Thus $K(\infty)$ is a free $\Z$-module with canonical
basis given by $B(\infty)$,
the isomorphism classes of irreducible modules,
and $K(\infty)_\Q$ is the $\Q$-vector space on basis $B(\infty)$.
We will always view $K(\infty)$ as a lattice in $K(\infty)_{\Q}$.

We let $K(\infty)^*$ 
denote the {\em restricted dual} of $K(\infty)$,
namely, the set of functions $f:K(\infty) \rightarrow \Z$
such that $f$ vanishes on all but finitely many elements of
$B(\infty)$.
Thus $K(\infty)^*$ is also a free $\Z$-module, with canonical basis 
$$
\{\delta_M\:|\:[M] \in B(\infty)\}
$$
dual to the basis $B(\infty)$ of $K(\infty)$,
i.e. $\delta_M([M]) = 1$, $\delta_M([N]) = 0$ for $[N] \in B(\infty)$
with $N \not\cong M$.
Note for an arbitrary $N \in \rep_I \H_n$,
$\delta_M([N])$ simply computes the {\em composition multiplicity}
$[N:M]$ of the irreducible module $M$ as a composition factor of $N$.
Finally, we write
$B(\infty)^*_\Q := \Q \otimes_{\Z} B(\infty)^*,$
which can be identified with the restricted 
dual of $B(\infty)_{\Q}$.

Entirely similar definitions can be made for each $\la \in P_+$:
\begin{equation}\label{kla}
K(\la) = \bigoplus_{n \geq 0} K(\rep_I{\H_n^\la})
\end{equation}
denotes
the Grothendieck groups of the categories $\rep{\H_n^\la}$ for all $n$.
Again, $K(\la)$ is a free $\Z$-module on the basis $B(\la)$
of isomorphism classes of irreducible modules.
Moreover, $\inf^\la$ induces a canonical embedding
$$
\inf^\la:K(\la) \hookrightarrow K(\infty).
$$
We will generally identify $K(\la)$ with its image under this embedding.
We also define
$K(\la)^*$ and $K(\la)_{\Q} = \Q \otimes_{\Z} K(\la)$ 
as above.

Recall
the operators $e_i$ and more generally the divided power
operators $e_i^{(r)}$ for $r \geq 1$,
defined on irreducible modules in $\res_I \H_n$
in (\ref{eid}) and (\ref{eird}) respectively. These induce
linear maps
\begin{equation}
e_i^{(r)}:K(\infty) \rightarrow K(\infty)
\end{equation}
for each $r \geq 1$.
Similarly, the operators $e_i^{(r)}$ and $f_i^{(r)}$ 
from (\ref{eirda}) and (\ref{firda})
respectively induce maps
\begin{equation}
e_i^{(r)}, f_i^{(r)} :K(\lambda) \rightarrow K(\lambda).
\end{equation}
Recall by Lemma~\ref{dp} that
\begin{equation}\label{diviv}
e_i^r = (r!) e_i^{(r)}, \qquad
f_i^r = (r!) f_i^{(r)}.
\end{equation}
Extending scalars, the maps $e_i^{(r)}, f_i^{(r)}$
induce linear maps on 
$K(\infty)_{\Q}$ and $K(\la)_{\Q}$ too.

\Point{Hopf algebra structure}
Now we wish to give $K(\infty)$ the structure of a graded Hopf algebra
over $\Z$.
To do this, recall 
the canonical isomorphism
\begin{equation}\label{canid}
K(\rep_I{\H_m}) \otimes_{\Z} K(\rep_I{\H_n})
\rightarrow
K(\rep_I{\H_{m,n}})
\end{equation}
from (\ref{stara}),
for each $m,n \geq 0$.
The exact functor $\ind_{m,n}^{m+n}$ induces a well-defined map
$$
\ind_{m,n}^{m+n}:K(\rep_I{\H_{m,n}}) \rightarrow K(\rep_I{\H_{m+n}}).
$$
Composing with the isomorphism (\ref{canid}) and taking the direct
sum over all
$m,n \geq 0$, we obtain a homogeneous map
\begin{equation}
\diamond:K(\infty) \otimes_{\Z} K(\infty) \rightarrow K(\infty).
\end{equation}
By transitivity of induction, 
this makes $K(\infty)$ into an associative graded $\Z$-algebra.
Moreover, there is a unit
\begin{equation}
\iota:\Z \rightarrow K(\infty)
\end{equation}
mapping $1$ to the identity module $\mathbf 1$ of $\H_0$.
Finally, since the duality $\tau$ induces the identity map
at the level of Grothendieck groups,
Theorem~\ref{duals} implies that the multiplication
$\diamond$ is commutative (in the usual unsigned sense).

Now consider how to define the comultiplication.
The exact functor $\res^n_{n_1,n_2}$
induces a map
$$
\res^n_{n_1,n_2}:
K(\rep_I \H_n) \rightarrow
K(\rep_I \H_{n_1,n_2}).
$$
On composing with the isomorphism
(\ref{canid}), we obtain maps
\begin{align}
\Delta^n_{n_1,n_2}:
&K(\rep_I \H_n) \rightarrow
K(\rep_I \H_{n_1}) \otimes_{\Z} K(\rep_I H_{n_2}),\\
\Delta^n = \sum_{n_1+n_2 = n} \Delta^n_{n_1,n_2}:
&K(\rep_I \H_n) \rightarrow
\bigoplus_{n_1+n_2 = n}
K(\rep_I \H_{n_1}) \otimes_{\Z} K(\rep_I H_{n_2}).
\end{align}
Now taking the direct sum over all $n \geq 0$
gives a homomogeneous map
\begin{equation}
\Delta: K(\infty)
\rightarrow K(\infty) \otimes_{\Z} K(\infty).
\end{equation}
Transitivity of restriction implies that $\Delta$ is coassociative,
while the homogeneous projection onto $K(\rep_I\H_0) \cong \Z$
gives a counit
\begin{equation}
\eps:
K(\infty) \rightarrow \Z.
\end{equation}
Thus $K(\infty)$ is also a graded coalgebra over $\Z$.
Now finally:

\vspace{1mm}
\begin{Theorem}
$(K(\infty), \diamond, \Delta, \iota, \eps)$
is a commutative, graded Hopf algebra
over $\Z$.
\end{Theorem}

\begin{proof}
It just remains to check that $\Delta$ is an algebra homomorphism,
which follows using the Mackey Theorem.
Note in checking the details, 
one needs to use Lemma~\ref{cccc} to take the definition of $\circledast$
into account correctly.
\end{proof}

We record 
the following lemma explaining how to compute the action of $e_i$
on $K(\infty)$ explicitly in terms of $\Delta$:

\vspace{1mm}
\begin{Lemma}\label{howto}
Let $M$ be a module in $\rep_I \H_n$. Write
$$
\Delta^n_{n-1,1} [M] = \sum_a [M_a] \otimes [N_a]
$$
for irreducible 
$\H_{n-1}$-modules $M_a$ and irreducible $\H_{1}$-modules $N_a$.
Then,
$$
e_i [M] = \sum_{a\:with\:N_a \cong L(i)} [M_a].
$$
\end{Lemma}

\begin{proof}
This is immediate from Lemma~\ref{relate1}.
\end{proof}

\begin{Lemma}\label{serre}
The operators $e_i:K(\infty) \rightarrow K(\infty)$ satisfy the Serre
relations, i.e.
\begin{align*}
e_i e_j &= e_j e_i 
&\hbox{if $|i-j| > 1$};\\
e_i^2 e_j + e_j e_i^2 &= 2 e_i e_j e_i 
&\hbox{if $|i-j|=1, i \neq 0, j \neq \ell$};\\
e_0^3 e_1 + 3 e_0 e_1 e_0^2 &= 3 e_0^2 e_1 e_0 + e_1 e_0^3
&\hbox{if $\ell \neq 1$};\\
e_{\ell-1}^3 e_\ell + 3 e_{\ell-1} e_\ell e_{\ell-1}^2 &= 
3 e_{\ell-1}^2 e_\ell e_{\ell-1} + e_{\ell} e_{\ell-1}^3
&\hbox{if $\ell \neq 1$};\\
e_0^5 e_1 + 5 e_0 e_1 e_0^4 + 10 e_0^3 e_1 e_0^2 &=
10 e_0^2 e_1 e_0^3 + 5 e_0^4 e_1 e_0 + e_1 e_0^5
&\hbox{if $\ell = 1$},
\end{align*}
for $i,j \in I$.
\end{Lemma}

\begin{proof}
In view of Lemma~\ref{howto} and coassociativity of $\De$, this reduces
to checking it on irreducible 
$\H_n$-modules for $n = 2,3,4,4,6$ respectively.
For this, the character information in Lemmas~\ref{calc1} and \ref{calc2}
is sufficient.
\end{proof}

Now consider $K(\la)$ for $\la \in P_+$.
This has a natural structure as $K(\infty)$-comodule:
viewing $K(\la)$ as a subset of $K(\infty)$, the
comodule structure map
is the restriction
$$
\Delta^\la:K(\la) \rightarrow K(\la) \otimes_{\Z} K(\infty)
$$
of $\Delta$. In 
other words, each $K(\la)$ is a subcomodule of the right regular 
$K(\infty)$-comodule.

One checks from the definitions that
the dual maps to $\diamond,\Delta, \iota, \eps$
induce on
$K(\infty)^*$ the structure of a cocommutative graded Hopf algebra.
From this point of view, each $K(\la)$
is a left $K(\infty)^*$-module in the natural way: 
$f \in K(\infty)^*$ acts on the left on $K(\la)$
as the map $(\id \bar\otimes f) \circ \Delta^\la$.
Similarly, $K(\infty)$ is itself a left $K(\infty)^*$-module, indeed in
this case the action is even {\em faithful}.

\vspace{1mm}
\begin{Lemma}\label{id}
The operator $e_i^{(r)}$
acts on $K(\infty)$ (resp. $K(\la)$ for any $\la \in P_+$) 
in the same way as
the basis element $\delta_{L(i^r)}$ of $K(\infty)^*$.
\end{Lemma}

\begin{proof}
Let $M$ be an irreducible module in $\rep_I \H_n$ or $\rep \H_n^\la$.
Let
$$
\Delta_{n-r,r} [M] = \sum_a [M_a] \otimes [N_a]
$$
for irreducible $\H_{n-r}^\la$-modules $M_a$
and irreducible $\H_{r}$-modules $N_a$.
By the definition of the action of $K(\infty)^*$ on $K(\la)$, it follows that
$$
\delta_{L(i^r)} [M] = \sum_{a\:with\:N_a \cong L(i^r)} [M_a].
$$
Hence, since $[\res^r_{1,\dots,1} L(i^r)] = (r!) [L(i)^{\circledast r}]$,
we get that $\delta_{L(i)}^r$ acts in the same way as
$(r!) \delta_{L(i^r)}$.
So in view of (\ref{diviv}), it just remains to check that
$\delta_{L(i)}$ acts in the same way as $e_i$, which follows by
Lemma~\ref{howto}.
\end{proof}

\begin{Lemma}\label{pidef} There is a unique homomorphism
$\pi:U_\Z^+ \rightarrow K(\infty)^*$
of graded Hopf algebras such that
$e_i^{(r)} \mapsto \delta_{L(i^r)}$ for each $i \in I$ and $r \geq 1$.
\end{Lemma}

\begin{proof}
Since $K(\infty)$ is a faithful $K(\infty)^*$-module, 
(\ref{diviv}) and Lemmas~\ref{serre} and \ref{id} imply
that the operators $\delta_{L(i^r)}$ satisfy the same relations
as the generators $e_i^{(r)}$ of $U_{\Z}^+$.
This implies existence of a unique such algebra homomorphism.
The fact that $\pi$ is a coalgebra map follows because
$$
\Delta(\delta_{L(i)}) 
= \delta_{L(i)} \otimes 1 + 1 \otimes \delta_{L(i)}
$$
by the definition of the comultiplication on $K(\infty)^*$.
\end{proof}

\Point{Shapovalov form}
Now focus on a fixed $\la \in P_+$.
For an $\H_n^\la$-module $M$, 
we let $P_M$ denote its projective cover in the category
$\mod{\H_n^\la}$.
Since $\H_n^\la$ is a finite dimensional superalgebra,
we can identify
\begin{equation}
K(\la)^* = \bigoplus_{n \geq 0} K(\proj \H_n^\la)
\end{equation}
so that the basis element $\delta_M$ corresponds to the 
isomorphism class $[P_M]$ for each irreducible $\H_n^\la$-module $M$
and each $n \geq 0$.
Moreover, under this identification,
the canonical pairing
\begin{equation}\label{cpa}
(.,.):K(\la)^* \times K(\la) \rightarrow \Z
\end{equation}
satisfies
\begin{equation}\label{asasa}
([P_M], [N]) = \left\{
\begin{array}{ll}
\dim \hom_{\H_n^\la}(P_M, N)&\hbox{if $M$ is of type $\Mtype$,}\\
\frac{1}{2}
\dim \hom_{\H_n^\la}(P_M, N)&\hbox{if $M$ is of type $\Qtype$,}
\end{array}
\right.
\end{equation}
for $\H_n^\la$-modules $M, N$ with $M$ irreducible
(since the right hand side clearly computes the composition multiplicity
$[N:M]$).

There is a homogeneous map
\begin{equation}\label{scd}
\omega:K(\la)^* \rightarrow K(\la)
\end{equation}
induced by the natural maps
$K(\proj \H_n^\la) \rightarrow K(\rep \H_n^\la)$ for each $n$.
We warn the reader that we do not yet know that $\omega$ is injective.

As explained at the end of \ref{divp}, we can define an action of
$e_i^{(r)}$ and $f_i^{(r)}$ on the projective indecomposable modules,
hence on $K(\la)^*$.
We know by Lemma~\ref{dp2} that (\ref{diviv}) holds for the operations
on $K(\la)^*$ as well as on $K(\la)$.
Also, the actions of $e_i^{(r)}$ and $f_i^{(r)}$ commute
with $\omega$ by (\ref{okg}).

\vspace{1mm}
\begin{Lemma}\label{Adad}
The operators $e_i, f_i$ on $K(\la)^*$ and $K(\la)$
satisfy $(e_i x, y) = (x, f_i y)$, $(f_i x, y) = (x, e_i y)$
for each $x \in K(\la)^*$ and $y \in K(\la)$.
\end{Lemma} 

\begin{proof}
Let $M$ be an irreducible $\H_n^\la$-module and $N$ be an irreducible
$\H_{n+1}^\la$-module. We check that
$(f_i [P_M], [N])  = ([P_M], e_i [N])$
in the special case that $i = 0$, $M$ is of type $\Qtype$ and
$N$ is of type $\Mtype$.
In this case,
Lemmas~\ref{relate2}, \ref{sleepy}(i) and \ref{relate1},
\begin{align*}
(f_i [P_M], [N])
&= \frac{1}{2}(\ind_i [P_M], [N]) 
= \frac{1}{2} \dim \hom_{\H_{n+1}^\la}(\ind_i P_M, N)\\
&= \frac{1}{2} \dim \hom_{\H_n^\la}(P_M, \res_i N)
= ([P_M], \res_i [N])
= ([P_M], e_i [N]).
\end{align*}
All the other situations that need to be considered follow similarly.
\end{proof}

\begin{Corollary}\label{pnew}
Suppose 
$$
e_i^{(r)} [M] = \sum_{[N] \in B(\la)} a_{M,N} [N],
\quad
f_i^{(r)} [M] = \sum_{[N] \in B(\la)} b_{M,N} [N].
$$
for $[M] \in B(\la)$.
Then,
$$
e_i^{(r)} [P_N] = \sum_{[M] \in B(\la)} b_{M,N} [P_M],
\quad
f_i^{(r)} [P_N] = \sum_{[M] \in B(\la)} a_{M,N} [P_M]
$$
for $[N] \in B(\la)$.
\end{Corollary}

\vspace{1mm}

Using this, the next lemma
is an easy consequence of
Theorem~\ref{weredone}(i), see
\cite[Lemma~11.2]{G} for the detailed proof.

\vspace{1mm}
\begin{Lemma}\label{g1}
Let $M$ be an irreducible $\H_n^\la$-module, set
$\eps = \eps_i(M), \phi = \phi_i(M)$.
Then, for any $m \geq 0$,
$$
e_i^{(m)} [P_M] = 
\sum_{N\:with\:\eps_i(N) \geq m} a_N [P_{\tilde e_i^m N}]
$$
for coefficients $a_N \in \Z_{\geq 0}$.
Moreover, in case $m = \eps$,
$$
e_i^{(\eps)} [P_M] = \binom{\eps + \phi}{\eps} [P_{\tilde e_i^\eps M}]
+
\sum_{N\:with\:\eps_i(N) > \eps} a_N [P_{\tilde e_i^\eps N}].
$$
\end{Lemma}

We also need:

\vspace{1mm}
\begin{Theorem}\label{gg}
Given an irreducible $\H_n^\la$-module $M$, 
$[P_M] \in K(\proj \H_n^\la)$ can be written as an integral
linear combination of terms of the form
$f_{i_1}^{(r_1)} \dots f_{i_a}^{(r_a)} [{\bil}]$.
\end{Theorem}

\begin{proof}
Proceed by induction on $n$, the conclusion being trivial for $n = 0$.
So suppose $n > 0$ and that the result is true for all smaller $n$.
Suppose for a contradiction that we can find an irreducible $\H_n^\la$-module
$M$ for which the result does not hold. Pick $i$ with $\eps := \eps_i(M) > 0$.
Since there are only finitely many irreducible $\H_n^\la$-modules,
we may assume by the choice of $M$ that the result holds for all irreducible
$\H_n^\la$-modules $L$ with $\eps_i(L) > \eps$.
Write 
$$
f_i^{(\eps)} [P_{\tilde e_i^\eps M}]
=
\sum_{[L] \in \Irr \H_n^\la}
a_L [P_L]
$$
for coefficients $a_L \in \Z$.
By Corollary~\ref{pnew}, 
$a_L = [e_i^{(\eps)} L : \tilde e_i^\eps M].$
In particular, $a_L = 0$ unless $\eps_i(L) \geq \eps$.
Moreover, if $a_L \neq 0$ for $L$ with $\eps_i(L) = \eps$,
then by Theorem~\ref{T280900}(i),
we have that $e_i^{(\eps)} L \cong \tilde e_i^{\eps} L
\cong \tilde e_i^\eps M$, whence $L \cong M$ and $a_M = 1$.
This shows:
\begin{equation}
[P_M] = 
f_i^{(\eps)} [P_{\tilde e_i^\eps M}]
-
\sum_{L\:with\:\eps_i(L) > \eps}
a_L [P_L]
\end{equation}
for some $a_L \in \Z$.
But the inductive hypothesis and choice of $M$ ensure that
all terms on the right hand side 
are integral linear combinations of terms
$f_{i_1}^{(r_1)} \dots f_{i_a}^{(r_a)} [\bil]$,
hence the same is true for $[P_M]$, a contradiction.
\end{proof}

The next two theorems are fundamental, so we include the proofs
even though they are 
identical to the argument in \cite[Theorem 11.1]{G}.

\vspace{1mm}
\begin{Theorem}\label{injectivity}
The map $\omega:K(\la)^* \rightarrow K(\la)$ from (\ref{scd}) 
is injective.
\end{Theorem}

\begin{proof}
We show by induction on $n$ that
the map $\omega:K(\proj \H_n^\la) \rightarrow K(\rep \H_n^\la)$
is injective. This is clear if $n = 0$, so assume $n > 0$ and the
result has been proved for all smaller $n$.
Suppose we have a relation
$$
\sum_{[M] \in \Irr \H_n^\la} a_M \omega ([P_M]) = 0
$$
with not all coefficients $a_M$ being zero.
We may choose $i \in I$ and an irreducible module $M$
such that $a_M \neq 0, \eps := \eps_i(M) > 0$ 
and $a_N = 0$ for all $N$ with $\eps_i(N) > \eps$.

Apply $e_i^{(\eps)}$ to the sum. By Lemma~\ref{g1}, we get
$$
\sum_{N\:with\:\eps_i(N) = \eps}
\binom{\eps + \phi_i(N)}{\eps} a_N \omega ([P_{\tilde e_i^\eps N}])
+ X = 0
$$
where $X$ is a sum of terms of the form
$\omega([P_{\tilde e_i^\eps L}])$ with $\eps_i(L) > \eps$.
Now the inductive hypothesis shows that $X = 0$ and that
$a_N = 0$ for each $N$ with $\eps_i(N) = \eps$.
In particular, $a_M = 0$, a contradiction.
\end{proof}

In view of Theorem~\ref{injectivity}, we may identify
$K(\la)^*$ with its image under $\omega$,
so $K(\la)^* \subseteq K(\la)$ are two different lattices in 
$K(\la)_{\Q}$: on tensoring with $\Q$ they become equal.
Extending scalars, the pairing (\ref{cpa}) induces a bilinear form
\begin{equation}\label{shap}
(.,.):K(\la)_{\Q} \times K(\la)_{\Q} \rightarrow \Q
\end{equation}
with respect to which the operators $e_i$ and $f_i$ are adjoint.

\vspace{1mm}
\begin{Theorem}\label{shapthm}
The form
$(.,.):K(\la)_\Q \times K(\la)_\Q \rightarrow \Q$
is symmetric and non-degenerate.
\end{Theorem}

\begin{proof}
It is non-degenerate by construction, being induced by the pairing
$(.,.)$ from (\ref{cpa}).
So we just need to check that it is symmetric.
Proceed by induction on $n$ to show that
$(.,.):K(\rep \H_n^\la)_\Q \times K(\rep \H_n^\la)_\Q \rightarrow \Q$
is symmetric. In view of Theorem~\ref{gg}, any element of
$K(\rep \H_n^\la)_{\Q}$ 
can be written as
$f_i x$ for $x \in K(\rep \H_{n-1}^\la)_{\Q}$.
Then for any other $y \in K(\rep \H_n^\la)_{\Q}$, we have that
$$
(f_i x, y) = (x, e_i y) = (e_i y, x) = (y, f_i x)
$$
by the induction hypothesis.
\end{proof}

\Point{\boldmath Chevalley relations}
Continue working with a fixed $\la \in P_+$.
We turn now to considering the relations satisfied by the operators
$e_i, f_i$ on $K(\la)$. 

\vspace{1mm}
\begin{Lemma}\label{serre0}
The operators $e_i, f_i:K(\lambda) \rightarrow K(\lambda)$
satisfy the Serre relations (\ref{cr3}).
\end{Lemma}

\begin{proof}
We know the $e_i$ satisfy the Serre
relations on all of $K(\infty)$ by Lemma~\ref{serre},
so they certainly satisfy the Serre relations on restriction to
$K(\lambda)$. Moreover, $e_i$ and $f_i$ are adjoint
operators for the bilinear form $(.,.)$ according to
Lemma~\ref{Adad}, and this form is non-degenerate
by Theorem~\ref{shapthm}. The lemma follows.
\end{proof}

Now we consider relations between the $e_i$ and $f_j$.
For $i \in I$ and an irreducible $\H_n^\la$-module $M$
with central character $\chi_\gamma$ for $\gamma \in \Gamma_n$,
define
\begin{equation}\label{hidef}
h_i [M] = \langle h_i, \lambda - \gamma \rangle [M].
\end{equation}
Recall according to Lemma~\ref{unstuck} that we have equivalently
that 
\begin{equation}
h_i [M] = (\phi_i(M) - \eps_i(M)) [M].
\end{equation}
More generally, define
\begin{equation}
\binom{h_i}{r}:K(\la) \rightarrow K(\la),\qquad
[M] \mapsto \binom{\phi_i(M) - \eps_i(M)}{r} [M]
\end{equation}
where $\binom{m}{r}$ denotes $m(m-1)\dots (m-r+1) / (r!)$.
Extending linearly, each $\binom{h_i}{r}$ can be viewed as a diagonal
linear operator $K(\la) \rightarrow K(\la)$.
The definition (\ref{hidef}) implies immediately that:

\vspace{1mm}
\begin{Lemma}\label{imm1}
As operators on $K(\lambda)$,
$[h_i, e_j] = \langle h_i,\alpha_j \rangle e_j$ and
$[h_i, f_j] = -\langle h_i, \alpha_j \rangle f_j$
for all $i,j \in I$.
\end{Lemma}

\vspace{1mm}

For the next lemma,  we follow \cite[Proposition~12.5]{G}.

\vspace{1mm}
\begin{Lemma}\label{imm2}
As operators on $K(\la)$, the relation
$[e_i, f_j] = \delta_{i,j} h_i$
holds for each $i,j \in I$.
\end{Lemma}

\begin{proof}
Let $M$ be an irreducible $\H_n^\la$-module.
It follows immediately from Theorems~\ref{T280900}(i)
and \ref{weredone}(i) (together with central character considerations in case
$i \neq j$)
that $[M]$ appears in 
$e_i f_j [M] - f_j e_i [M]$ with multiplicity
$\delta_{i,j} (\phi_i(M) - \eps_i(M)).$
Therefore, it suffices simply to show that
$e_i f_j [M] - f_j e_i [M]$
is a multiple of $[M]$.
Let us show equivalently that
$$
[\res_i \ind_j M - \ind_j \res_i M]
$$
is a multiple of $[M]$.

For $m \gg 0$, we have a surjection
$\ind_{n,1}^{n+1} M \boxtimes \R_m(j) \twoheadrightarrow \ind_j M.$
Apply $\pr^\la \circ \res_i$ to get a surjection
\begin{equation}\label{s1}
\pr^\la \res_i \ind_{n,1}^{n+1} M \boxtimes \R_m(j) \twoheadrightarrow \res_i \ind_j M.
\end{equation}
By the Mackey Theorem and (\ref{zoo}), there is an exact sequence
$$
0 \rightarrow (M \oplus \Pi M)^{\oplus \delta_{i,j}. m c_j}
\longrightarrow \res_i \ind_{n,1}^{n+1} M \boxtimes \R_m(j)
\longrightarrow \ind_{n-1,1}^n (\res_i M) \boxtimes \R_m(j)
\rightarrow 0,
$$
where $c_j$ is as in Lemma~\ref{cnt}.
For sufficiently large $m$,
$\pr^\la \ind_{n-1,1}^n (\res_i M) \boxtimes \R_m(j) = \ind_j \res_i M$.
So on applying the right exact functor $\pr^\la$ and using the irreducibility
of $M$, this implies that there is an exact sequence
\begin{equation}\label{s2}
0 \longrightarrow M^{\oplus m_1} \oplus \Pi M^{\oplus m_2}
\longrightarrow \pr^\la \res_i \ind_{n,1}^{n+1} M \boxtimes \R_m(j)
\longrightarrow \ind_j \res_i M \longrightarrow 0,
\end{equation}
for some $m_1,m_2$.
Now let $N$ be any irreducible $\H_n^\la$-module with $N \not\cong M$.
Combining (\ref{s1}) and (\ref{s2}) shows that 
\begin{equation}\label{all}
[\ind_j \res_i M - \res_i \ind_j M :N] \geq 0
\end{equation}
Now summing over all $i,j$ and using (\ref{indresdecomp})
gives that $[\ind \:\res \:M-\res \:\ind \:M:N] \geq 0.$
But Theorem~\ref{cycloMackey} shows that equality holds here,
hence it must hold in (\ref{all}) for all $i,j \in I$.
This completes the proof.
\end{proof}

To summarize, we have shown in (\ref{diviv}), 
Lemmas~\ref{serre0}, \ref{imm1} and \ref{imm2} that:

\vspace{1mm}
\begin{Theorem}\label{chevrel}
The action of the operators 
$e_i, f_i, h_i$ on $K(\la)$ satisfy the Chevalley relations
(\ref{cr1}), (\ref{cr2}) and (\ref{cr3}).
Hence, the actions of $e_i^{(r)}, f_i^{(r)}$ and $\binom{h_i}{r}$
for all $i \in I, r \geq 1$ make $K(\la)_{\Q}$
into a $U_{\Q}$-module so that
$K(\la)^*, K(\la)$ are $U_\Z$-submodules.
\end{Theorem}

\Point{\boldmath Identification of $K(\infty)^*$, $K(\la)^*$ and $K(\la)$}
Now we can prove Theorems A and B stated in the introduction.
Compare \cite[4.3, 4.4]{A1} and \cite[14.1, 14.2]{G}.

\vspace{1mm}
\begin{Theorem}\label{thmb}
For any $\la \in P_+$, 
\begin{enumerate}
\item[(i)] $K(\la)_{\Q}$ is precisely the integrable highest weight 
$U_\Q$-module
of highest weight $\la$, with highest weight vector $[\bil]$;
\item[(ii)] the bilinear form $(.,.)$ from (\ref{shap})
on the highest weight module $K(\la)_{\Q}$ coincides with
the usual Shapovalov form satisfying $([\bil],[\bil]) = 1$;
\item[(iii)] $K(\la)^* \subset K(\la)$ are integral forms
of $K(\la)_{\Q}$ containing $[\bil]$,
with $K(\la)^*$ being the minimal lattice $U_\Z^- [\bil]$
and $K(\la)$ being its dual under the Shapovalov form.
\end{enumerate}
\end{Theorem}

\begin{proof}
It makes sense to think of $K(\la)_{\Q}$ as a $U_{\Q}$-module according to
Theorem~\ref{chevrel}.
The actions of $e_i$ and $f_i$ are locally nilpotent by
Theorems~\ref{T280900}(i) and \ref{weredone}(i).
The action of $h_i$ is diagonal by definition.
Hence, $K(\la)_{\Q}$ is an integrable module.
Clearly $[\bil]$ is a highest weight vector of highest weight $\la$.
Moreover, $K(\la)_{\Q} = U_\Q^- [\bil]$ by Theorem~\ref{gg}.
This completes the proof of (i), and (ii) follows immediately from
Lemma~\ref{Adad}.
For (iii), we know already that $K(\la)^* \subset K(\la)$ are dual lattices
of $K(\la)_{\Q}$ which are invariant under $U_\Z$.
Moreover, Theorem~\ref{gg} again shows $K(\la)^* = U_{\Z}^- [\bil]$.
\end{proof}

\begin{Theorem}\label{thma}
The map $\pi:U_\Z^+ \rightarrow K(\infty)^*$ constructed
in Lemma~\ref{pidef} is an isomorphism.
\end{Theorem}

\begin{proof}
Note by Lemma~\ref{id} that the action of the $e_i^{(r)} \in U_{\Z}^+$
on $K(\infty)$, hence on each $K(\la)$, factors through the
map $\pi$.
So if $x \in \ker \pi$, we have by the previous theorem
that $x$ acts as
zero on all integrable highest weight modules $K(\la), \la \in P_+$.
Hence $x = 0$ and $\pi$ is injective.
To prove surjectivity, take $x \in K(\infty)$.
It suffices to show that
$(\pi(u),x) = 0$ for all $u \in U_\Z^+$ implies that
$x = 0$.
Note
$$
(\pi(u), x) = (\pi(1)\pi(u), x) 
=
(\pi(1), \pi(u) x)
= (\pi(1), ux)
$$
where the second equality follows because the right regular action
of $K(\infty)^*$ on itself is precisely the dual action to the left
action of $K(\infty)^*$ on $K(\infty)$, and the third equality
follows from Lemma~\ref{id}.
Hence, if $(\pi(u), x) = 0$ for all $u \in U_{\Z}^+$, we have that
$(\pi(1), ux) = 0$ for all $u \in U_{\Z}^+$.
Now choose $\la \in P_+$ sufficiently large so that in fact
$x \in K(\la) \subset K(\infty)$. Then, it follows that
$([\bil], ux) = 0$ for all $u \in U_{\Z}^+$,
where $(.,.)$ now is canonical 
pairing between $K(\la)^*$ and $K(\la)$.
Hence by Lemma~\ref{Adad}, $(v [\bil], x) = 0$ for all $v \in U_\Z^-$.
But then Theorem~\ref{gg} implies that $x = 0$.
\end{proof}

\ifbook@\pagebreak\fi

\section{Identification of the crystal}

\Point{\boldmath Final properties of $B(\infty)$}\label{ended}
Now we follow the ideas of \cite[$\S$13]{G}.

\vspace{1mm}
\begin{Lemma}\label{sunday}
Let $M \in \rep_I \H_m$ be irreducible.
\begin{enumerate}
\item[(i)] For any $i \in I$,
 either $\eps_i(\tilde f_i^* M) = \eps_i(M)$ or $\eps_i(M)+1$.

\item[(ii)] For any $i,j \in I$ with $i \neq j$, 
$\eps_i(\tilde f_j^* M) = \eps_i(M)$.


\end{enumerate}
\end{Lemma}

\begin{proof}
We prove (i), the proof of (ii) being similar.
By the Shuffle Lemma, we certainly have that
$\eps_i(\tilde f_i^* M) \leq \eps_i(M) + 1$.
Now let $N = \tilde f_i^* M$.
Then obviously, $\eps_i(\tilde e_i^* N) \leq \eps_i(N)$.
Hence, $\eps_i(M) \leq \eps_i(\tilde f_i^* M)$.
\end{proof}

\begin{Lemma}
\label{oldway}
Let $M\in\rep_I{\H_m}$ be irreducible and $i,j \in I$. 
If $\eps_i(\tilde f_j^* M) = \eps_i(M)$ 
then, writing $\eps  := \eps_i(M)$, we have 
$\tilde e_i^\eps  \tilde f_j^* M \cong \tilde f_j^* \tilde e_i^\eps  M$.
\end{Lemma}

\begin{proof}
Set $n = m-\eps $. Let $N = \tilde e_i^\eps  M$, so $N$ is an 
irreducible $\H_{n}$-module
with $\eps_i(N) = 0$
and $M = \tilde f_i^\eps  N$.
For $0 \leq b \leq \eps $, let
$Q_b = \res_i^{\eps -b} \tilde f_j^* M$.
Theorem~\ref{T280900}(i) and Lemma~\ref{relate1} 
imply that in the Grothendieck group,
$Q_b$ is some number of copies of
$\tilde e_i^{\eps -b} \tilde f_j^* M$ plus terms with strictly smaller
$\eps_i$. In particular, $\eps_i(L) \leq b$ for all composition factors
$L$ of $Q_b$, while $Q_0$ is consists only of 
copies of $\tilde e_i^\eps  \tilde f_j^* M$.

We will show by decreasing induction on $b = \eps ,\eps -1,\dots,0$ that
there is a non-zero $\H_{n+b+1}$-module homorphism
$$
\gamma_b:\ind_{1,n,b}^{n+b+1} 
L(j) \boxtimes N \boxtimes L(i^b)
\rightarrow Q_b.
$$
In case $b = \eps $, $Q_\eps  = \tilde f_j^* M$ is a quotient of
$\ind_{1,m}^{m+1} L(j) \boxtimes M$ and $M$ is a quotient
of $\ind_{n,\eps }^{m} N \boxtimes L(i^\eps )$, so the induction starts.
Now we suppose by induction that we have proved $\gamma_b\neq 0$ exists
for
$b \geq 1$ and construct $\gamma_{b-1}$.

Consider $\res^{n+b+1}_{n+b,1} \ind_{1,n,b}^{n+b+1}
L(j) \boxtimes N \boxtimes L(i^b)$.
By the Mackey Theorem, this has a 
filtration
$0 \subset F_1 \subset F_2 \subset F_3$ with successive quotients
\begin{align*}
F_1
& \simeq
\ind_{1,n,b-1,1}^{n+b,1} \res^{1,n,b}_{1,n,b-1,1}
L(j) \boxtimes N \boxtimes L(i^b),\\
F_2/F_1 & \simeq \ind_{1,n-1,b,1}^{n+b,1} {^w}\res^{1,n,b}_{1,n-1,1,b}
L(j) \boxtimes N \boxtimes L(i^b),\\
F_3 / F_2 &\simeq
\ind_{n,b,1}^{n+b,1} N \boxtimes L(i^b) \boxtimes L(j),
\end{align*}
where $w$ is the obvious permutation.
As $\gamma_b \neq 0$, Frobenius reciprocity implies that there is
a copy of the $\H_{1,n,b}$-module
$L(j) \circledast N \circledast L(i^b)$ in the image of $\gamma_b$.
Now $b > 0$, so the $q(i)$-eigenspace of
$X_{n+b+1}+X_{n+b+1}^{-1}$ acting on 
$L(j) \circledast N \circledast L(i^b)$ is non-trivial.
We conclude that the map
$$
\tilde \gamma_b = \res_i (\gamma_b):
\res_i \ind_{1,n,b}^{n+b+1} L(j)\boxtimes N \boxtimes L(i^b)
\rightarrow \res_i Q_b = Q_{b-1}
$$
is non-zero.

If $i \neq j$, then it follows from the description of $F_3/F_2$ and
$F_2/F_1$ above that $\res_i (F_3 / F_1) = 0$.
So in this case, we necessarily have that $\tilde \gamma_b(F_1) \neq 0$.
Similarly if $i = j$, $\res_i (F_2 / F_1) = 0$, so
if $\tilde \gamma_b(F_1) = 0$ we see that
$\tilde \gamma_b$ factors to a non-zero homomorphism
$$
\res^{n+b,1}_{n+b} \ind_{n,b,1}^{n+b,1} N \boxtimes L(i^b) \boxtimes L(i)
\rightarrow Q_{b-1}.
$$
But this implies that $Q_{b-1}$ has a constituent $L$ with $\eps_i(L) =
b$,
which we know is not the case.
Hence we have that $\tilde \gamma_b(F_1) \neq 0$ in the case $i = j$ too.

Hence, the restriction of $\tilde \gamma_b$ to $F_1$
gives us a non-zero homomorphism
$$
\res^{n+b,1}_{n+b}
\ind_{1,n,b-1,1}^{b+n,1} \res^{1,n,b}_{1,n,b-1,1}
L(j) \boxtimes N \boxtimes L(i^b)
\rightarrow Q_{b-1}.
$$
Now finally as all composition factors of
$\res^{b}_{b-1} L(i^b)$ are isomorphic to
$L(i^{b-1})$, 
this implies the existence of a non-zero homomorphism
$\gamma_{b-1}:
\ind_{1,n,b-1}^{b+n} 
L(j) \boxtimes N \boxtimes L(i^{b-1})
\rightarrow Q_{b-1}$
completing the induction.

Now taking $b=0$ we have a non-zero map 
$\ga_0: \ind_{1,n}^{n+1}L(j)\boxtimes N \rightarrow Q_0.$
But the left hand side has irreducible cosocle $\tilde f_j^* N
= \tilde f_j^* \tilde e_i^\eps  M$
while all composition factors of the right hand side are
isomorphic to $\tilde e_i^\eps  \tilde f_j^* M$.
This completes the proof.
\end{proof}

\begin{Corollary}
\label{middleway}
Let $M \in \rep_I(\H_m)$ be irreducible and $i \in I$.
Let $M_1 = \tilde e_i^{\eps_i(M)} M$ and
$M_2 = (\tilde e_i^*)^{\eps_i^*(M)} M$.
Then, $\eps_i^*(M) = \eps_i^*(M_1)$ if and only if 
$\eps_i(M) = \eps_i(M_2)$.
\end{Corollary}

\begin{proof}
Since we can apply the automorphism $\sigma$,
it suffices to check just one of the implications. So suppose
$\eps_i(M) = \eps_i(M_2)$.
Clearly, $\eps_i^*(M) \geq \eps_i^*(M_1)$.
For the reverse inequality, the preceeding two lemmas show that
$$
M_1 \cong
\tilde e_i^{\eps_i(M)} M \cong
\tilde e_i^{\eps_i(M)} (\tilde f_i^*)^{\eps_i^*(M)} M_2 \cong
(\tilde f_i^*)^{\eps_i^*(M)} \tilde e_i^{\eps_i(M)} 
M_2.
$$
This shows that $\eps_i^*(M_1) \geq \eps_i^*(M)$.
\end{proof}

\begin{Lemma}
\label{newway}
Let $M\in\rep_I{\H_m}$ be irreducible and $i \in I$ satisfy
$\eps_i(\tilde f_i^* M) = \eps_i(M) + 1$.
Then
$\tilde e_i \tilde f_i^* M = M$.
\end{Lemma}

\begin{proof}
Set $\eps :=\eps_i(M)$ and $N = \tilde e_i^\eps  M$. By the Shuffle Lemma
and Theorem~\ref{inj},
$$
[\ind_{1,m-\eps ,\eps }^{m+1} L(i) \circledast N\circledast L(i^\eps )]
=
[\ind_{m-\eps ,\eps +1}^{m+1} N \circledast L(i^{\eps +1})].
$$
Hence by Theorem~\ref{T280900}(i),
$[\ind_{1,m-\eps ,\eps }^{m+1}L(i)\circledast N\circledast L(i^\eps )]$ 
equals $[\tilde f_i^{\eps +1} N] = [\tilde f_i M]$ 
plus terms $[L]$ for irreducible $L$
with  $\eps_i(L)\leq \eps $. 
On the other hand, $\ind_{1,m-\eps ,\eps }^{m+1}L(i)\circledast
N\circledast L(i^\eps )$ 
surjects onto $\tilde f_i^* M$. So the assumption
$\eps_i(\tilde f_i^* M)=\eps +1$ 
implies $\tilde f_i M\cong\tilde f_i^* M$. 
\end{proof}

\Point{Crystals}\label{crs}
Let us now recall some definitions from \cite{Kas}.
A {\em crystal} is a set $B$ endowed with maps
\begin{align*}
\phi_i, \eps_i&:B \rightarrow \Z \cup \{-\infty\}\quad(i\in I),\\
\tilde e_i, \tilde f_i&:B \rightarrow B \cup \{0\}\quad(i\in I),\\
\wt&:B \rightarrow P
\end{align*}
such that
\begin{itemize}
\item[(C1)]
$\phi_i(b) = \eps_i(b) + \langle h_i, \wt(b)\rangle$ for any $i \in I$;
\item[(C2)]
if $b \in B$ satisfies $\tilde e_i b \neq 0$, then $\eps_i(\tilde e_i b)
= \eps_i(b) - 1$, $\phi_i(\tilde e_i b) = \phi_i(b)+1$ and
$
\wt(\tilde e_i b) = \wt(b) + \alpha_i$;
\item[(C3)]
if $b \in B$ satisfies $\tilde f_i b \neq 0$, then $\eps_i(\tilde f_i b)
= \eps_i(b) + 1$, $\phi_i(\tilde f_i b) = \phi_i(b)-1$ and
$
\wt(\tilde e_i b) = \wt(b) - \alpha_i$;
\item[(C4)] for $b_1,b_2 \in B$, $b_2 = \tilde f_i b_1$ if and only if
$b_1 = \tilde e_i b_2$;
\item[(C5)] if $\phi_i(b) = -\infty$, then $\tilde e_i b = \tilde f_i b =
0$.
\end{itemize}
For example, for each $i \in I$, we have the crystal $B_i$ defined as a
set
to be $\{b_i(n)\:|\:n \in \Z\}$ with
\begin{align*}
\eps_j(b_i(n)) = \left\{
\begin{array}{ll}
-n&\hbox{if $j = i$,}\\
-\infty&\hbox{if $j \neq i$;}
\end{array}\right.
\qquad\qquad&\phi_j(b_i(n)) = \left\{
\begin{array}{ll}
n&\hbox{if $j = i$,}\\
-\infty&\hbox{if $j \neq i$;}
\end{array}\right.\\
\tilde e_j(b_i(n)) = \left\{
\begin{array}{ll}
b_i(n+1)&\hbox{if $j = i$,}\\
0&\hbox{if $j \neq i$;}
\end{array}\right.
\qquad\qquad&\tilde f_j(b_i(n)) = \left\{
\begin{array}{ll}
b_i(n-1)&\hbox{if $j = i$,}\\
0&\hbox{if $j \neq i$;}
\end{array}\right.
\end{align*}
and $\wt(b_i(n)) = n \alpha_i$.
We abbreviate $b_i(0)$ by $b_i$.
Also for $\la \in P$, we have the crystal $T_\lambda$
equal as a set to $\{t_\lambda\}$, with 
$\eps_i(t_\lambda) = \phi_i(t_\lambda) = -\infty$,
$\tilde e_i t_\lambda = \tilde f_i t_\lambda = 0$
and $\wt(t_\lambda) = \lambda$.

A morphism $\psi:B \rightarrow B'$ of crystals
is a map $\psi:B \cup \{0\} \rightarrow B' \cup \{0\}$ such that
\begin{itemize}
\item[(H1)] $\psi(0) = 0$;
\item[(H2)] if $\psi(b) \neq 0$ for $b \in B$, then 
$\wt(\psi(b)) = \wt(b)$, $\eps_i(\psi(b)) = \eps_i(b)$ and
$\phi_i(\psi(b)) = \psi(b)$;
\item[(H3)] for $b \in B$ such that $\psi(b) \neq 0$ and $\psi(\tilde e_i
b) \neq 0$, we have that $\psi(\tilde e_i b) = \tilde e_i \psi(b)$;
\item[(H4)] for $b \in B$ such that $\psi(b) \neq 0$ and $\psi(\tilde f_i
b) 
\neq 0$, we have that $\psi(\tilde f_i b) = \tilde f_i \psi(b)$.
\end{itemize}
A morphism of crystals is called {\em strict}
if $\psi$ commutes with the $\tilde e_i$'s and $\tilde f_i$'s, and an
{\em embedding} if $\psi$ is injective.

We also need the notion of a tensor product of two crystals $B, B'$.
As a set, $B \otimes B'$ is equal to $\{b \otimes b'\:|\:b \in B, b' \in
B'\}$.
This is made into a crystal by
\ifbook@
\begin{align*}
\eps_i(b \otimes b') &= \max(\eps_i(b), \eps_i(b') - \langle h_i, \wt(b)
\rangle),\\
\phi_i(b \otimes b') &= \max(\phi_i(b) + \langle h_i, \wt(b') \rangle,
\phi_i(b')),\\
\tilde e_i (b \otimes b') &= \left\{\begin{array}{ll}
\tilde e_i b \otimes b'&\hbox{if $\phi_i(b) \geq \eps_i(b')$}\\
b \otimes \tilde e_i b'&\hbox{if $\phi_i(b) < \eps_i(b')$}
\end{array}\right.,\\
\tilde f_i (b \otimes b') &= \left\{\begin{array}{ll}
\tilde f_i b \otimes b'&\hbox{if $\phi_i(b) > \eps_i(b')$}\\
b \otimes \tilde f_i b'&\hbox{if $\phi_i(b) \leq \eps_i(b')$}
\end{array}\right.,\\
\wt(b \otimes b') &= \wt(b) + \wt(b').
\end{align*}
\else
\begin{align*}
\eps_i(b \otimes b') = \max(\eps_i(b), \eps_i(b') - \langle h_i, \wt(b)
\rangle),\quad
&\phi_i(b \otimes b') = \max(\phi_i(b) + \langle h_i, \wt(b') \rangle,
\phi_i(b')),\\
\tilde e_i (b \otimes b') = \left\{\begin{array}{ll}
\tilde e_i b \otimes b'&\hbox{if $\phi_i(b) \geq \eps_i(b')$}\\
b \otimes \tilde e_i b'&\hbox{if $\phi_i(b) < \eps_i(b')$}
\end{array}\right.,
\quad
&\tilde f_i (b \otimes b') = \left\{\begin{array}{ll}
\tilde f_i b \otimes b'&\hbox{if $\phi_i(b) > \eps_i(b')$}\\
b \otimes \tilde f_i b'&\hbox{if $\phi_i(b) \leq \eps_i(b')$}
\end{array}\right.,\\
\wt(b \otimes b') &= \wt(b) + \wt(b').
\end{align*}
\fi
Here, we understand $b \otimes 0 = 0 = 0 \otimes b$.

\vspace{2mm}
Having recalled these definitions, we now explain how to
make our sets $B(\infty)$
and $B(\la)$ from (\ref{binfty}) and (\ref{bla})
into crystals in the above sense.
To do this, it just remains to define the weight
functions 
on both $B(\infty)$ and $B(\la)$, as well as the function
$\phi_i$ on $B(\infty)$:
set 
\begin{align}
\wt(M) &= - \gamma,\\\intertext{for an irreducible $M \in \rep_\gamma
\H_n$, and}
\wt^\la(N) &= \la - \gamma,
\end{align}
for an irreducible $N \in \rep_\gamma \H_n^\la$.
Also for an irreducible $[M] \in B(\infty)$, define
\begin{equation}
\phi_i(M) = \eps_i(M) + \langle h_i, \wt(M) \rangle,
\end{equation}
so that property (C1) is automatic.
Thus we have defined $(B(\infty), \eps_i, \phi_i, \tilde e_i,
\tilde f_i, \wt)$ and $(B(\lambda), \eps_i, \phi_i, \tilde e_i^\la, 
\tilde f_i^\la, \wt^\la)$ 
purely in terms of the representation theory of the Hecke-Clifford
superalgebras.

\vspace{1mm}
\begin{Lemma}
$(B(\infty), \eps_i, \phi_i, \tilde e_i,
\tilde f_i, \wt)$ and each
$(B(\lambda), \eps_i, \phi_i, \tilde e_i^\la, \tilde f_i^\la,
\wt^\la)$ for $\la \in P_+$ are crystals in the sense of Kashiwara.
\end{Lemma}

\begin{proof}
Property (C1) is Lemma~\ref{unstuck} or the definition in the affine case.
Property (C4) is 
Lemma~\ref{L290900_2}.
The remaining properties are immediate.
\end{proof}

Recall the embedding $\inf^\la:
B(\la) \cup \{0\} \rightarrow B(\infty)\cup\{0\}$ from 
(\ref{inff}).

\vspace{1mm}
\begin{Lemma}\label{lastly}
The map $B(\la) \hookrightarrow B(\infty) \otimes T_\lambda
$,
$[N] \mapsto \inf^\la [N] \otimes t_\lambda$
is an embedding of crystals with image
$$
\left\{[M] \otimes t_\lambda \in B(\infty) \otimes T_\lambda
\:\big|\:
\eps_i^*(M) \leq \langle h_i, \lambda\rangle\hbox{ for each }i \in
I\right\}.
$$
\end{Lemma}

\begin{proof}
Since $\tilde e_i^\la$ and $\tilde f_i^\la$ are restrictions
of $\tilde e_i, \tilde f_i$ from $B(\infty)$ to $B(\la)$, respectively,
the first statement is immediate.
The second is a restatement of Corollary~\ref{P300900}.
\end{proof}

\Point{\boldmath Identification of $B(\infty)$ and $B(\la)$}
The first lemma follows directly from Lemmas~\ref{sunday}
and \ref{oldway}, 
applying the automorphism $\sigma$ to get (ii).

\vspace{1mm}
\begin{Lemma}\label{soc}
Let $M\in\rep_I \H_m$ be irreducible and $i,j \in I$ with $i\neq j$. 
Set $a = \eps_i^*(M)$. 
\begin{enumerate}

\item[(i)] $\eps_j(M) = \eps_j((\tilde e_i^*)^a M)$.

\item[(ii)] If $\eps_j(M) > 0$, then
$\eps_i^*(\tilde e_j M) = \eps_i^*(M)$
and 
$(\tilde e_i^*)^a \tilde e_j M \cong \tilde e_j (\tilde e_i^*)^a M.$
\end{enumerate}
\end{Lemma}

\vspace{1mm}


The proof of the next result is taken from \cite[Proposition 10.2]{G}.

\vspace{1mm}
\begin{Lemma}\label{sec}
Let $M\in\rep_I \H_m$ be irreducible and $i \in I$. 
Set $a = \eps_i^*(M)$ and $\bar M = (\tilde e_i^*)^a M$.
\begin{enumerate}
\item[(i)] $\eps_i(M) = \max(\eps_i(\bar M), a - \langle h_i, \wt(\bar M)
\rangle).$

\item[(ii)] If $\eps_i(M) > 0$,
$$
\eps_i^*(\tilde e_i M) = \left\{
\begin{array}{ll}
a&\hbox{if $\eps_i(\bar M) \geq a - \langle h_i, \wt(\bar M) \rangle$,}\\
a-1&\hbox{if $\eps_i(\bar M) < a - \langle h_i, \wt(\bar M) \rangle$.}
\end{array}
\right.
$$

\item[(iii)] If $\eps_i(M) > 0$,
$$
(\tilde e_i^*)^b \tilde e_i M \cong \left\{
\begin{array}{ll}
\tilde e_i(\bar M)&\hbox{if $\eps_i(\bar M) \geq a - \langle h_i, \wt(\bar
M) \rangle$,}\\
\bar M&\hbox{if $\eps_i(\bar M) < a - \langle h_i, \wt(\bar M) \rangle$,}
\end{array}
\right.
$$
where $b = \eps_i^*(\tilde e_i M)$.
\end{enumerate}
\end{Lemma}
\begin{proof}
Let $\eps=\eps_i(M)$, $n=m-\eps$, and $N=(\tilde e_i)^\eps M$. 

(i) By twisting with $\sigma$, it suffices to prove that (for arbitrary
$M$)
\begin{equation}
\label{Esec}
\eps_i^*(M) = \max(\eps_i^*(N), \eps - \langle h_i, \wt(N)
\rangle).
\end{equation}
Define the weights $\la(0),\la(1),\dots\in P_+$ by taking the
$\Lambda_i$-coefficient of $\la(r)$ to be $\eps_i^*(N)+r$ and the
$\Lambda_j$-coefficients of $\la(r)$ for $j\neq i$ to be $\gg 0$. 
Then ${\mathcal I}_{\la(r)} N=0$ for any $r\geq 0$, see
Corollory~\ref{P300900}. Moreover, the same corollary implies that
for
$k=\phi_i^{\la(r)}(N)$,
$$
\eps_i^*(\tilde f_i^k N)=\langle h_i, \la(r)\rangle
\quad \text{and}\quad \eps_i^*(\tilde
f_i^{k+1} N)=\langle h_i, \la(r) \rangle +1.
$$
Now,  
$
\phi_i^\la(N)-\eps_i^\la(N)=\langle h_i,\la+\wt(N)\rangle
$
and $\eps_i(N)=\eps_i^\la(N)$. Moreover, by Lemma~\ref{sunday}(i) 
(twisted with $\sigma$) we have $\eps_i^*(\tilde f^k N)\geq \eps_i^*(N)$
for any $k$. All of these applied consecutively to $\la(0),\la(1),\dots$
imply: 
$$
\eps_i^*(\tilde f_i^s N) = \left\{
\begin{array}{ll}
\eps_i^*(N)
&\hbox{if $s\leq\eps_i^*(N)+\eps_i(N)+\langle h_i,\wt(N)\rangle$,}\\
s-\eps_i(N)-\langle h_i,\wt(N)\rangle
&\hbox{if $s\geq\eps_i^*(N)+\eps_i(N)+\langle h_i,\wt(N)\rangle$}
\end{array}
\right.
$$
for all $s\geq 0$. For $s=\eps$, taking into account $\eps_i(N)=0$, this
gives (\ref{Esec}). 

(ii) Observe by (\ref{Esec}) that $\eps_i^*(\tilde e_i M)=\eps_i^*(M)-1$
if and only if $\eps>\langle h_i,\wt(N)\rangle+\eps_i^*(N)$, and that
otherwise $\eps_i^*(\tilde e_i M)=\eps_i^*(M)$. But $\langle
h_i,\wt(N)\rangle= \langle h_i,\wt(M)\rangle+2\eps$ and 
$\langle h_i,\wt(\bar M)\rangle= \langle h_i,\wt(M)\rangle+2a$
so (ii) follows if we show that
$\langle h_i, \wt(M)\rangle+\eps_i^*(N)+\eps<0$ if and only if 
$\langle h_i, \wt(M)\rangle+\eps_i(\bar M)+a<0$. 
But by (i) and (\ref{Esec}), 
\begin{align*}
\langle h_i,\wt(M)\rangle+\eps_i^*(N)+\eps
& = 
\max(\langle h_i,\wt(M)\rangle+\eps_i^*(N)+\eps_i(\bar M),\eps_i^*(N)-a),
\\
\langle h_i, \wt(M)\rangle+\eps_i(\bar M)+a
& = \max(\langle h_i, \wt(M)\rangle+\eps_i(\bar M)+\eps_i^*(N),
\eps_i(\bar M)-\eps).
\end{align*}
Moreover, obviously $\eps_i^*(N)-a\leq 0$ and $\eps_i(\bar M)-\eps\leq 0$,
and it remains to observe that $\eps_i^*(N)-a= 0$ if and only if
$\eps_i(\bar M)-\eps= 0$, thanks to Corollary~\ref{middleway}.

(iii) This follows from (ii) and Lemmas~\ref{oldway} and \ref{newway} (twisted
with $\sigma$).
\end{proof}

\vspace{1mm}

Now for each $i \in I$, define a map
\begin{equation}
\Psi_i:B(\infty) \rightarrow B(\infty) \otimes B_i
\end{equation}
mapping each $[M] \in B(\infty)$
to $[(\tilde e_i^*)^a M] \otimes \tilde f_i^a b_i$, where
$a = \eps_i^*(M)$.

\vspace{1mm}
\begin{Lemma}\label{prop}
The following properties hold:
\begin{enumerate}
\item[(i)] for every $[M] \in B(\infty)$, $\wt(M)$ is a negative
sum of simple roots;

\item[(ii)] $[{\bf 1}]$ is the unique element of $B(\infty)$ with weight
$0$;

\item[(iii)] $\eps_i({\bf 1}) = 0$ for every $i \in I$;

\item[(iv)] $\eps_i(M) \in \Z$ for every $[M] \in B(\infty)$ and
every $i \in I$;

\item[(v)] for every $i$, the map $\Psi_i:B(\infty) \rightarrow B(\infty)
\otimes B_i$ defined above is a strict embedding of crystals;

\item[(vi)] $\Psi_i(B(\infty)) \subseteq B(\infty) \times \{\tilde f_i^n
b_i\:|\:
n\geq 0\}$;

\item[(vii)] for any $[M] \in B(\infty)$ other than $[{\bf 1}]$,
there exists $i \in I$ such that $\Psi_i([M]) = [N] \otimes \tilde f_i^n
b_i$ for some $[N] \in B(\infty)$ and $n > 0$.
\end{enumerate}
\end{Lemma}

\begin{proof}
Properties (i)--(iv) are immediate from our construction of $B(\infty)$.
The information required to verify (v) is 
exactly contained in
Lemmas~\ref{soc} and \ref{sec}.
Finally, (vi) is immediate from the definition of $\Psi_i$,
and (vii) holds because every such $M$ has $\eps_i^*(M) > 0$
for at least one $i \in I$.
\end{proof}

The properties in Lemma~\ref{prop} exactly characterize the crystal
$B(\infty)$ by \cite[Proposition~3.2.3]{KS}.
Hence, we have proved:

\vspace{1mm}
\begin{Theorem}\label{ident1}
The crystal $B(\infty)$ is isomorphic to Kashiwara's
crystal $B(\infty)$ associated to the crystal base of $U_\Q^-$.
\end{Theorem}

\vspace{1mm}

In view of \cite[Theorem~8.2]{Kas}, we can also identify our maps
$\Psi_i$ with
those of \cite{Kas}. Taking into account \cite[Proposition 8.1]{Kas} we
can then identify
our functions $\eps_i^*$ on $B(\infty)$ with those in \cite{Kas}. 
It follows from this, Lemma~\ref{lastly} and \cite[Proposition~8.2]{Kas}
that:

\vspace{1mm}
\begin{Theorem}\label{ident2}
For each $\la \in P_+$,
the crystal $B(\lambda)$ is isomorphic to Kashiwara's crystal
$B(\la)$ associated to the integrable highest weight $U_\Q$-module
of highest weight $\la$.
\end{Theorem}

\Point{Blocks}\label{blockssect}
Let $\la \in P_+$.
As an application of the theory, we mention here another delightful
argument of Grojnowski from \cite{G2}, which in our setting
classifies the blocks
of the cyclotomic Hecke-Clifford superalgebras $\H_n^\la$.

Let us recall the definition in the case of a finite dimensional
superalgebra $A$.
Let $\sim$ be the equivalence relation on the set of isomorphism
classes of irreducible $A$-modules such that
$[L] \sim [M]$ if and only if there exists a chain 
$L \cong L_0, L_1, \dots, L_m \cong M$ of irreducible $A$-modules
with
either $\operatorname{Ext}^1_A(L_i, L_{i+1}) \neq 0$
or $\operatorname{Ext}^1_A(L_{i+1}, L_{i}) \neq 0$ for each $i$.
Given a $\sim$-equivalence class $b$, the corresponding
{\em block} $\rep_b A$ is the full subcategory
of $\rep A$ consisting of the $A$-modules all of whose composition factors
belong to the equivalence class $\sim$.
Thus, there is a decomposition
$$
\rep A = \bigoplus_{b} \rep_b A
$$
as $b$ runs over all $\sim$-equivalence classes.
As usual, there are various equivalent points of view: for instance,
$[L] \sim [M]$ if and only if the even central characters
$\chi_L, \chi_M:Z(A)_{\0} \rightarrow F$
arising from the actions 
on the irreducible modules $L, M$ are equal.
Alternatively, the blocks of $A$ can be defined in terms of 
a decomposition
of the identity $1_A$ into a sum of centrally primitive
{\em even} idempotents.

\vspace{1mm}
\begin{Theorem}
Let $M, N$ be irreducible $\H_n^\la$-modules with $M \not\cong N$.
Let
$$
0 \longrightarrow M \longrightarrow X \longrightarrow N \longrightarrow 0
$$
be an exact sequence of $\H_n$-modules. Then, $\pr^\la X = X$ and
the sequence is also an exact sequence of $\H_n^\la$-modules.
\end{Theorem}

\begin{proof}
We note that
for irreducible modules $M, N$ in $\rep_I \H_n$ with $M \not\cong N$,
we have that
$$
\hom_{\H_{n-1}}(\res^{n}_{n-1} M, 
\res^n_{n-1} N) = 0.
$$ 
This follows immediately from Corollary~\ref{epr}, Lemma~\ref{relate1}
and (\ref{indresdecomp}).
Using this in place of \cite[Lemma 2]{G2}, the proof of the theorem is now
completed by exactly the same argument as in \cite{G2}.
\end{proof}

Note the theorem can be reinterpreted as saying that
\begin{equation}
\ext^1_{\H_n}(M, N) \simeq \ext^1_{\H_n^\la}(M,N)
\end{equation}
for irreducible $\H_n^\la$-modules $M, N$ with $M \not\cong N$.
This is certainly not the case if $M \cong N$!
Recalling the definitions from \ref{Defs},
we immediately deduce the following
corollary which determines the blocks of $\H_n^\la$:

\vspace{1mm}
\begin{Corollary}\label{blockclass}
The 
blocks of $\H_n^\la$ are precisely the subcategories
$\rep_\gamma \H_n^\la$ for  $\gamma \in \Gamma_n$.
Moreover, $\rep_\gamma \H_n^\la$ is non-trivial if and only
if the $(\la - \gamma)$-weight space of the highest weight
module $K(\la)_{\Q}$ is non-zero.
\end{Corollary}

\ifbook@\pagebreak\fi

\section{Branching rules}\label{last}

In this final section, we deduce some important consequences
for modular representations of the 
double cover $\widehat S_n$ of the symmetric group.
In particular we obtain the classification of the irreducible
modules, describe the blocks of the group algebra 
$F\widehat S_n$ and prove analogues of the modular branching rules of
\cite{K1,K2,K3,K4,JWB:branching,BK1}.
Actually the results on branching are somewhat weaker here: there is no
representation theoretic interpretation at present to the
notion  of ``normal node'' introduced below.

\Point{Kang's description of the crystal graph}
Now we focus on the fundamental 
highest weight $\Lambda_0$.
In this case, Kang \cite{Kang} has given a
convenient combinatorial description
of the crystal $B(\Lambda_0)$ in terms of Young diagrams,
which we now describe.

For any $n \geq 0$, let $\la=(\la_1,\la_2,\dots)$ be a partition of $n$, 
i.e. a non-increasing sequence
of non-negative integers summing to $n$.
Recall that $\ell \in \Z_{>0} \cup \{\infty\}$ and $h = 2 \ell + 1$.
We call $\la$ an {\em $h$-strict partition} if $h$ divides 
$\la_r$ whenever $\la_r=\la_{r+1}$ for $r \geq 1$
(to interpret correctly in case $h = \infty$, we adopt the convention that
$\infty$ divides $0$ and nothing else).
Let ${\mathscr P}_h(n)$ denote the set of all $h$-strict partitions of $n$,
and ${\mathscr P}_h := \bigcup_{n \geq 0} {\mathscr P}_h(n)$.
We say that  $\la \in {\mathscr P}_h(n)$
is {\em restricted} if in addition 
$$
\left\{
\begin{array}{ll}
\la_r-\la_{r+1}< h &\hbox{if $h|\la_r$},\\
\la_r-\la_{r+1}\leq h &\hbox{if $h \nmid \la_r$}
\end{array}
\right.
$$
for each $r \geq 1$.
Let ${\mathscr{RP}}_h(n)$ 
denote the set of all restricted $h$-strict partitions
of $n$, and ${\mathscr{RP}}_h := \bigcup_{n \geq 0} {\mathscr{RP}}_h(n)$.

Let $\la \in {\mathscr P}_h$ 
be an $h$-strict partition. 
We identify $\la$ with its {\em Young diagram} 
$$
\la=\{(r,s)\in\Z_{>0}\times\Z_{>0}\mid s\leq\la_r\}.
$$
Elements $(r,s)\in\Z_{>0}\times\Z_{>0}$ are called {\em nodes}. 
We label the nodes of $\la$ with {\em residues},
which are the elements of the set
$I = \{0,1,\dots,\ell\}$.
The labelling depends only on the column and follows the repeating pattern
$$
0,1,\dots,\ell-1,\ell,\ell-1,\dots,1,0,
$$
starting fom the first column and going to the right, 
see Example~\ref{E1} below. The residue of the node $A$ is denoted $\Res A$. 
Define the {\em residue content} of $\la$ to be the tuple
\begin{equation}
\cont(\la)=(\gamma_i)_{i \in I}
\end{equation}
where for each $i \in I$,
$\gamma_i$ is the number of nodes of residue $i$ contained in the
diagram $\la$. 

Let $i \in I$ be some fixed residue.
A node $A = (r,s)
\in \la$ is called {\em $i$-removable} (for $\la$) if one of the following
holds:
\begin{itemize}
\item[(R1)] $\Res A = i$ and
$\la_A:=\la-\{A\}$ is again an $h$-strict partition;
\item[(R2)] the node $B = (r,s+1)$ immediately to the right of $A$
belongs to $\la$,
$\Res A = \Res B = i$,
and both $\la_B = \la - \{B\}$ and
$\la_{A,B} := \la - \{A,B\}$ are $h$-strict partitions.
\end{itemize}
Similarly, a node $B = (r,s)\notin\la$ is called 
{\em $i$-addable} (for $\la$) if one of the following holds:
\begin{itemize}
\item[(A1)] $\Res B = i$ and
$\la^B:=\la\cup\{B\}$ is again an $h$-strict partition;
\item[(A2)] 
the node $A = (r,s-1)$
immediately to the left of $B$ does not belong to $\la$,
$\Res A = \Res B = i$, and both 
$\la^A = \la \cup \{A\}$ and 
$\la^{A, B} := \la \cup\{A,B\}$ are $h$-strict partitions.
\end{itemize}
We note that (R2) and (A2) above are only possible in case $i = 0$.

Now label all $i$-addable
nodes of the diagram $\la$ by $+$ and all $i$-removable nodes by $-$.
Then, the {\em $i$-signature} of 
$\la$ is the sequence of pluses and minuses obtained by going along the rim of the Young diagram from bottom left to top right and reading off
all the signs.
The {\em reduced $i$-signature} of $\la$ is obtained 
from the $i$-signature
by successively erasing all neighbouring 
pairs of the form $+-$. 

Note the reduced $i$-signature always looks like a sequence
of $-$'s followed by $+$'s.
Nodes corresponding to a $-$ in the reduced $i$-signature are
called {\em $i$-normal}, nodes corresponding to a $+$ are
called {\em $i$-conormal}.
The rightmost $i$-normal node (corresponding to the rightmost $-$
in the reduced $i$-signature) is called {\em $i$-good}, 
and the leftmost $i$-conormal node (corresponding to the leftmost $+$
in the reduced $i$-signature) is called {\em $i$-cogood}.

\vspace{1mm}
\begin{Example}\label{E1}
{\rm 
Let $h=5$, so $\ell=2$.
The partition $\la=(16, 11,10,10,9,5,1)$ belongs to $\mathscr{RP}_h$,
and its residues are as follows:
$$
\diagram{
$0$ & $1$ & $2$ & $1$ & $0$ & $0$& $1$ & $2$ & $1$ & $0$ & $0$ & $1$ & $2$ & $1$ & $0$ & $0$\cr 
$0$ & $1$ & $2$ & $1$ & $0$ & $0$ & $1$ & $2$ & $1$ & $0$ & $0$ \cr
$0$ & $1$ & $2$ & $1$ & $0$ & $0$ & $1$ & $2$ & $1$ & $0$  \cr
$0$ & $1$ & $2$ & $1$ & $0$ & $0$ & $1$ & $2$ & $1$ & $0$  \cr
$0$ & $1$ & $2$ & $1$ & $0$ & $0$ & $1$ & $2$ & $1$ \cr
$0$ & $1$ & $2$ & $1$ & $0$\cr
$0$ \cr
}
$$
The $0$-addable and $0$-removable nodes are as labelled in the diagram:
$$
\begin{picture}(340,95)
\put(63,40)
{$
\diagram{
 &  &  &  &  & &  &  &  &  &  &  &  &  & $-$ & $-$ \cr 
 &  &  &  &  &  &  &  &  &  & $-$ \cr
 &  &  &  &  &  &  &  &  &   \cr 
 &  &  &  &  &  &  &  &  &  \cr
 &  &  &  &  &  &  &  & \cr 
 &  &  &  & $-$\cr
$-$ \cr
}
$}
\put(123.4, 14.4){\circle{9}}
\put(69.4, .8){\circle{9}}
\put(271.9, 82.3){\circle{9}}
\put(137,13.5){\makebox(0,0)[b]{$+$}}
\put(191,26){\makebox(0,0)[b]{$+$}}
\end{picture}
$$
Hence, the $0$-signature of $\la$ is
$-,-,+,+,-,-,-$ 
and the reduced $0$-signature is
$-,-,-$.
Note the nodes corresponding to the $-$'s in the reduced $0$-signature
have been circled in the above diagram.
So, there are three $0$-normal nodes, the rightmost of which is
$0$-good; there are no $0$-conormal or $0$-cogood nodes.
}\end{Example}

In general, we define
\begin{align}\label{neweps}
\eps_i(\la) & = \sharp\{\text{$i$-normal nodes in $\la$}\}=\sharp\{\text{$-$'s in the reduced $i$-signature of $\la$}\},\\
\phi_i(\la) & = \sharp\{\text{$i$-conormal nodes in $\la$}\}
=\sharp\{\text{$+$'s in the reduced $i$-signature of $\la$}\}.\label{newf}
\end{align}
Also set
\begin{align}\label{newet}
\tilde e_i(\la) &=
\left\{
\begin{array}{ll}
\la_A&\hbox{if $\eps_i(\la) > 0$ and $A$ is the (unique) $i$-good node,}\\
0&\hbox{if $\eps_i(\la) = 0$,}
\end{array}
\right.\\
\tilde f_i(\la) &=\label{newft}
\left\{
\begin{array}{ll}
\la^B&\hbox{if $\phi_i(\la) > 0$ and $B$ is the (unique) $i$-cogood 
node,}\\
0&\hbox{if $\phi_i(\la) = 0$.}
\end{array}
\right.
\end{align}
Finally define
\begin{equation}
\wt(\la) = \Lambda_0 - \sum_{i \in I} \gamma_i \alpha_i
\end{equation}
where $\cont(\la) = (\gamma_i)_{i \in I}$.
We have now defined a datum
$({\mathscr P}_h, \eps_i, \phi_i, \tilde e_i, \tilde f_i, \wt)$
which makes the set
${\mathscr P}_h$ of all $h$-strict partitions
into a crystal in the sense of \ref{crs}.
The definitions imply that $\tilde e_i(\la), \tilde f_i(\la)$ are restricted (or zero)
in case $\la$ is itself restricted.
Hence, 
$\mathscr{RP}_h$
is a sub-crystal of $\mathscr{P}_h$.
We can now state the main result of
Kang \cite[7.1]{Kang}
for type $A_{2\ell}^{(2)}$:

\vspace{1mm}
\begin{Theorem}
\label{TKang}
The set $\mathscr{RP}_h$ 
equipped with
\,$\eps_i,\varphi_i,\tilde e_i,\tilde f_i, \wt$ as above is 
isomorphic (in the unique way) to the crystal
$B(\Lambda_0)$ 
associated to the integrable highest weight $U_\Q$-module
of fundamental highest weight $\Lambda_0$.
\end{Theorem}

\vspace{1mm}
\begin{Example}
{\rm
The crystal graph of $\mathscr{RP}_h = B(\Lambda_0)$ 
in case $h = 3$, up to degree $10$, is as follows:

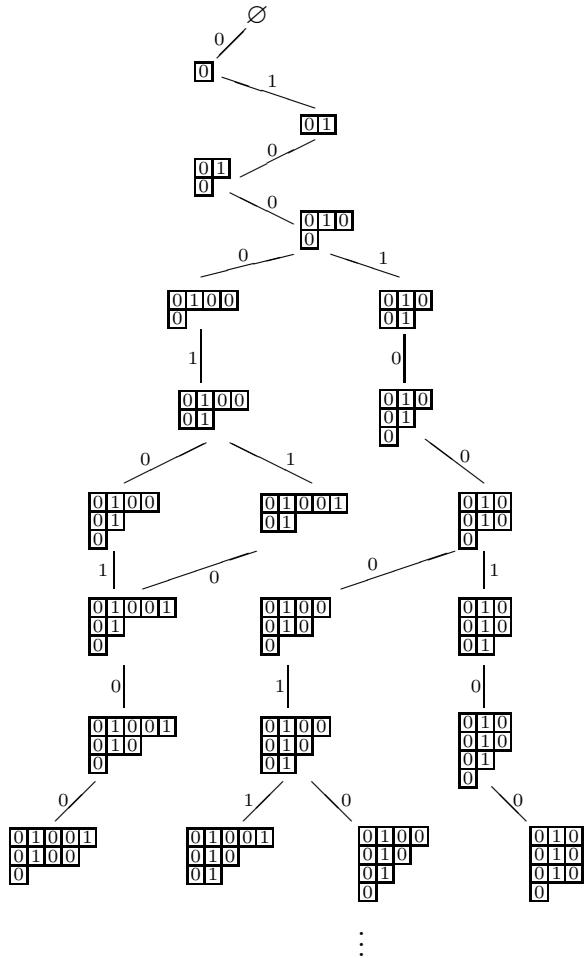
\vspace{-5mm}
\begin{picture}(340,325)

\put(162,301)
{$\varnothing$}

\put(143,280)
{$
\hoogte=6pt    \newdimen\breedte   \breedte=7pt  
\newdimen\dikte     \dikte=0.5pt 
\diagram{
\tiny{0}\cr
}
$}

\put(183,260)
{$
\hoogte=6pt    \newdimen\breedte   \breedte=7pt  
\newdimen\dikte     \dikte=0.5pt 
\diagram{
\tiny{0}& \tiny{1} \cr
}
$}

\put(143,240)
{$
\hoogte=6pt    \newdimen\breedte   \breedte=7pt  
\newdimen\dikte     \dikte=0.5pt 
\diagram{
\tiny{0}& \tiny{1} \cr
\tiny{0}\cr
}
$}

\put(183,220)
{$
\hoogte=6pt    \newdimen\breedte   \breedte=7pt  
\newdimen\dikte     \dikte=0.5pt 
\diagram{
\tiny{0}& \tiny{1} & \tiny{0}\cr
\tiny{0}\cr
}
$}

\put(133,190)
{$
\hoogte=6pt    \newdimen\breedte   \breedte=7pt  
\newdimen\dikte     \dikte=0.5pt 
\diagram{
\tiny{0}& \tiny{1} & \tiny{0} & \tiny{0}\cr
\tiny{0}\cr
}
$}
\put(213,190)
{$
\hoogte=6pt    \newdimen\breedte   \breedte=7pt  
\newdimen\dikte     \dikte=0.5pt 
\diagram{
\tiny{0}& \tiny{1} & \tiny{0}\cr
\tiny{0} & \tiny{1}\cr
}
$}

\put(137,152)
{$
\hoogte=6pt    \newdimen\breedte   \breedte=7pt  
\newdimen\dikte     \dikte=0.5pt 
\diagram{
\tiny{0}& \tiny{1} & \tiny{0} & \tiny{0} \cr
\tiny{0} & \tiny{1}\cr
}
$}
\put(213,149)
{$
\hoogte=6pt    \newdimen\breedte   \breedte=7pt  
\newdimen\dikte     \dikte=0.5pt 
\diagram{
\tiny{0}& \tiny{1} & \tiny{0}\cr
\tiny{0} & \tiny{1}\cr
\tiny{0}\cr
}
$}

\put(103,110)
{$
\hoogte=6pt    \newdimen\breedte   \breedte=7pt  
\newdimen\dikte     \dikte=0.5pt 
\diagram{
\tiny{0}& \tiny{1} & \tiny{0} & \tiny{0} \cr
\tiny{0} & \tiny{1}\cr
\tiny{0}\cr
}
$}
\put(243,110)
{$
\hoogte=6pt    \newdimen\breedte   \breedte=7pt  
\newdimen\dikte     \dikte=0.5pt 
\diagram{
\tiny{0}& \tiny{1} & \tiny{0}  \cr
\tiny{0} & \tiny{1} & \tiny{0} \cr
\tiny{0}\cr
}
$}
\put(168,113)
{$
\hoogte=6pt    \newdimen\breedte   \breedte=7pt  
\newdimen\dikte     \dikte=0.5pt 
\diagram{
\tiny{0}& \tiny{1} & \tiny{0} & \tiny{0} & \tiny{1}\cr
\tiny{0} & \tiny{1}\cr
}
$}

\put(168,70)
{$
\hoogte=6pt    \newdimen\breedte   \breedte=7pt  
\newdimen\dikte     \dikte=0.5pt 
\diagram{
\tiny{0}& \tiny{1} & \tiny{0} & \tiny{0} \cr
\tiny{0} & \tiny{1} & \tiny{0} \cr
\tiny{0}\cr
}
$}
\put(103,70)
{$
\hoogte=6pt    \newdimen\breedte   \breedte=7pt  
\newdimen\dikte     \dikte=0.5pt 
\diagram{
\tiny{0}& \tiny{1} & \tiny{0} & \tiny{0} & \tiny{1} \cr
\tiny{0} & \tiny{1} \cr
\tiny{0}\cr
}
$}
\put(243,70)
{$
\hoogte=6pt    \newdimen\breedte   \breedte=7pt  
\newdimen\dikte     \dikte=0.5pt 
\diagram{
\tiny{0}& \tiny{1} & \tiny{0} \cr
\tiny{0} & \tiny{1} & \tiny{0}\cr
\tiny{0} & \tiny{1}\cr
}
$}

\put(168,25)
{$
\hoogte=6pt    \newdimen\breedte   \breedte=7pt  
\newdimen\dikte     \dikte=0.5pt 
\diagram{
\tiny{0}& \tiny{1} & \tiny{0} & \tiny{0} \cr
\tiny{0} & \tiny{1} & \tiny{0} \cr
\tiny{0} & \tiny{1} \cr
}
$}
\put(103,25)
{$
\hoogte=6pt    \newdimen\breedte   \breedte=7pt  
\newdimen\dikte     \dikte=0.5pt 
\diagram{
\tiny{0}& \tiny{1} & \tiny{0} & \tiny{0} & \tiny{1} \cr
\tiny{0} & \tiny{1} & \tiny{0} \cr
\tiny{0}\cr
}
$}
\put(243,23)
{$
\hoogte=6pt    \newdimen\breedte   \breedte=7pt  
\newdimen\dikte     \dikte=0.5pt 
\diagram{
\tiny{0}& \tiny{1} & \tiny{0} \cr
\tiny{0} & \tiny{1} & \tiny{0}\cr
\tiny{0} & \tiny{1}\cr
\tiny{0}\cr
}
$}

\put(73,-17)
{$
\hoogte=6pt    \newdimen\breedte   \breedte=7pt  
\newdimen\dikte     \dikte=0.5pt 
\diagram{
\tiny{0}& \tiny{1} & \tiny{0} & \tiny{0} &\tiny{1}\cr
\tiny{0} & \tiny{1} & \tiny{0}&\tiny{0} \cr
\tiny{0}  \cr
}
$}
\put(140,-17)
{$
\hoogte=6pt    \newdimen\breedte   \breedte=7pt  
\newdimen\dikte     \dikte=0.5pt 
\diagram{
\tiny{0}& \tiny{1} & \tiny{0} & \tiny{0} & \tiny{1} \cr
\tiny{0} & \tiny{1} & \tiny{0} \cr
\tiny{0}&\tiny{1}\cr
}
$}
\put(205,-20)
{$
\hoogte=6pt    \newdimen\breedte   \breedte=7pt  
\newdimen\dikte     \dikte=0.5pt 
\diagram{
\tiny{0}& \tiny{1} & \tiny{0}&\tiny{0} \cr
\tiny{0} & \tiny{1} & \tiny{0}\cr
\tiny{0} & \tiny{1}\cr
\tiny{0}\cr
}
$}
\put(270,-20)
{$
\hoogte=6pt    \newdimen\breedte   \breedte=7pt  
\newdimen\dikte     \dikte=0.5pt 
\diagram{
\tiny{0}& \tiny{1} & \tiny{0} \cr
\tiny{0} & \tiny{1} & \tiny{0}\cr
\tiny{0} & \tiny{1}&\tiny{0}\cr
\tiny{0}\cr
}
$}

\put(90,-1){\line(1,1){15}}\put(91,5){\tiny{0}}
\put(161,-1){\line(1,1){15}}\put(161,5){\tiny{1}}
\put(202,-1){\line(-1,1){15}}\put(198,5){\tiny{0}}
\put(268,-1){\line(-1,1){13}}\put(263,5){\tiny{0}}

\put(252,42){\line(0,1){16}}\put(247,48){\tiny{0}}
\put(178,42){\line(0,1){16}}\put(173,48){\tiny{1}}
\put(116,42){\line(0,1){16}}\put(111,48){\tiny{0}}

\put(112,87){\line(0,1){14}}\put(106,92){\tiny{1}}
\put(123,87){\line(3,1){43}}\put(148,89){\tiny{0}}
\put(198,87){\line(3,1){43}}\put(208,94){\tiny{0}}
\put(252,87){\line(0,1){14}}\put(254,92){\tiny{1}}

\put(187,127){\line(-2,1){31}}\put(177,134){\tiny{1}}
\put(116,127){\line(2,1){31}}\put(243,135){\tiny{0}}
\put(252,127){\line(-4,3){22}}\put(122,134){\tiny{0}}

\put(222,165){\line(0,1){18}}\put(217,172){\tiny{0}}
\put(145,165){\line(0,1){20}}\put(140,172){\tiny{1}}

\put(220,205){\line(-3,1){26}}\put(212,210){\tiny{1}}
\put(145,205){\line(4,1){35}}\put(159,211){\tiny{0}}

\put(180,225){\line(-2,1){24}}\put(170,231){\tiny{0}}

\put(160,243){\line(2,1){28}}\put(170,251){\tiny{0}}

\put(188,269){\line(-3,1){35}}\put(170,277){\tiny{1}}

\put(151,288){\line(1,1){11}}\put(150,293){\tiny{0}}
\end{picture}
\vspace{11mm}
$$
\vdots
$$
}
\end{Example}

\vspace{2mm}

Let us also mention here the extension of
Morris' notion of $h$-bar core \cite{Morris}
to an
arbitrary $h$-strict partition $\la \in \mathscr{P}_h(n)$.
By an {\em $h$-bar} of $\la$, we mean one of the following:
\begin{itemize}
\item[(B1)] 
the rightmost $h$ nodes of row $i$ of $\la$ if $\la_i \geq h$ and
either  $h | \la_i$ or $\la$ has no row of length $(\la_i - h)$;
\item[(B2)] 
the set of nodes in rows $i$ and $j$ of $\la$
if $\la_i+\la_j = h$.
\end{itemize}
If $\la$ has no $h$-bars, it is called an {\em $h$-bar core}.
In general, the $h$-bar 
core $\tilde\la$ of $\la$ is obtained by successively removing $h$-bars,
reordering the rows each time so that the result still
lies in $\mathscr{P}_h$, until it is reduced to a core.
The {\em $h$-bar weight} of $\lambda$, denoted $w(\lambda)$,
is then the total number of $h$-bars that get removed.
For a Lie theoretic explanation of these notions,
we refer the reader to \cite[$\S$12.6]{Kac}.
In particular, as observed in \cite[$\S$4]{LT},
for $\mu,\la \in \mathscr{P}_h(n)$ we have that
\begin{equation}\label{coreblock}
\cont(\mu) = \cont(\la)\hbox{ if and only if }
\tilde\mu = \tilde \la.
\end{equation}
Also, bearing in mind Theorem~\ref{TKang},
we can state Kac' formula \cite[(12.13.5)]{Kac} 
for the character of
the highest weight
$U_\Q$-module of highest weight $\Lambda_0$ as follows:
for $\la \in \mathscr{RP}_h(n)$,
\begin{equation}\label{kacformula}
\sharp\{\mu \in \mathscr{RP}_h(n)\:|\:
\cont(\mu) =\cont(\la)\}
= \operatorname{Par}_\ell(w(\la)),
\end{equation}
where $\operatorname{Par}_\ell(N)$ 
denotes the number of partitions of $N$
as a sum of positive integers of $\ell$ different colors.

\Point{Representations of finite Hecke-Clifford superalgebras}\label{klk}
Now that we have an explicit description of the crystal $B(\Lambda_0)$,
we formulate a more combinatorial description of our main results for 
the representation theory of the finite Hecke-Clifford
superalgebras $\fH_n$.
Recall from Remark~\ref{fhrem} that this is precisely the
cyclotomic Hecke-Clifford superalgebra $\H_n^{\Lambda_0}$.
The results of this subsection also hold in the degenerate case,
when $\fH_n$ 
is the finite Sergeev superalgebra, see \ref{degmod}.

As explained in \ref{cg},
the isomorphism classes of irreducible $\fH_n$-modules
are parametrized by the nodes of the crystal graph $B(\Lambda_0)$.
By Theorems~\ref{TKang} and \ref{ident2}, we can identify
$B(\Lambda_0)$ with $\mathscr{RP}_h$.
In other words, we can use the set $\mathscr{RP}_h(n)$ 
of restricted $h$-strict partitions of $n$
to parametrize the irreducible $\fH_n$-modules for each $n \geq 0$.
Let us write $M(\lambda)$ for the irreducible
$\fH_n$-module corresponding to $\lambda \in \mathscr{RP}_h(n)$.
To be precise, 
$$
M(\lambda) := L(i_1,\dots,i_n)
$$
if $\lambda = \tilde f_{i_n} \dots \tilde f_{i_1} \varnothing$.
Here the operator $\tilde f_i$ is as defined in (\ref{newft}), 
corresponding under the identification
$\mathscr{RP}_h(n) = B(\Lambda_0)$
to the crystal operator denoted $\tilde f_i^{\Lambda_0}$ in earlier sections,
and $\varnothing$ denotes the empty partition,
corresponding to $\bil \in B(\Lambda_0)$.

For $\la \in \mathscr{RP}_h(n)$,
we also define
\begin{equation}\label{bbdef}
b(\la) := \sharp\{r \geq 1\:|\:h \nmid \la_r\},
\end{equation}
the number of parts of $\la$ that are not divisible by $h$.
The definition of residue content immediately gives that
\begin{equation}\label{resds}
b(\la) \equiv \gamma_0 \pmod{2},
\end{equation}
where $\gamma_0$ denotes the number of $0$'s in the residue content of $\la$.

\vspace{1mm}
\begin{Theorem}
\label{TBr} 
The modules
$\{M(\lambda)\:|\:\lambda \in \mathscr{RP}_h(n)\}$
form a complete set of pairwise non-isomorphic irreducible 
$\fH_n$-modules.
Moreover, for $\la,\mu \in \mathscr{RP}_h(n)$,
\begin{enumerate}
\item[(i)] $M(\la) \cong M(\la)^\tau$;

\item[(ii)] $M(\la)$ is of type $\Mtype$ if $b(\la)$ is even,
type $\Qtype$ if $b(\la)$ is odd;

\item[(iii)] $M(\mu)$ and $M(\la)$ belong to the same block if and only
if $\cont(\mu) = \cont(\la)$;

\item[(iv)] $M(\la)$ is a projective module if and only if $\la$
is an $h$-bar core.
\end{enumerate}
\end{Theorem}

\begin{proof}
We have already discussed the first statement of the theorem,
being a consequence of our main results combined with 
Kang's Theorem~\ref{TKang}.
For the rest, 
(i) follows from Corollary~\ref{selfdual},
(ii) is a special case of 
Lemma~\ref{cccc} combined with (\ref{resds}), and
(iii) is a special case of Corollary~\ref{blockclass}.
For (iv), note
that if $M(\la)$ is projective then it is the only irreducible in its
block, hence by (\ref{kacformula}),
$\operatorname{Par}_\ell(w(\la)) = 1$.
So either $w(\la) = 0$, or $\ell = 1$ and $w(\la) = 1$.
Now if $w(\la) = 0$ then $\la$ is an $h$-bar core
so the Shapovalov form on the ($1$-dimensional) 
$\wt(\la)$-weight 
space of $K(\Lambda_0)_{\Z}$ is $1$ (since $\wt(\la)$ is conjugate
to $\Lambda_0$ under the action of the affine Weyl group).
Hence, $M(\la)$ is projective by Theorem~\ref{thmb}(ii).
To rule out the remaining
possibilty $\ell = 1$ and $w(\la) = 1$, one checks
in that case that
the Shapovalov form on the $\wt(\la)$-weight space 
of $K(\Lambda_0)_{\Z}$ is $3$.
\end{proof}

The next two theorems summarize earlier results concerning
restriction and induction.

\vspace{1mm}
\begin{Theorem}\label{Te}
Let $\la \in {\mathscr{RP}}_h(n)$.
There exist
$\fH_{n-1}$-modules $e_i M(\la)$ for each $i \in I$,
unique up to isomorphism, such that
\begin{enumerate}
\item[(i)]
$\displaystyle
\res_{\fH_{n-1}}^{\fH_{n}} M(\la)
\cong
\left\{
\begin{array}{ll}
\displaystyle
2 e_0 M(\la) \oplus 2 e_1 M(\la) \oplus \dots \oplus 2 e_\ell M(\la)
&\hbox{if $b(\la)$ is odd,}\\
e_0 M(\la) \oplus 2 e_1 M(\la) \oplus \dots \oplus 2 e_\ell M(\la)
&\hbox{if $b(\la)$ is even;}\\
\end{array}
\right.$
\item[(ii)] 
for each $i \in I$,
$e_i M(\la)\neq 0$ if and only if $\la$ has an $i$-good node $A$,
in which case $e_i M(\la)$ is a self-dual indecomposable 
module with irreducible socle and cosocle isomorphic to $M(\la_A)$.
\end{enumerate}
Moreover, if $i \in I$ and $\la$ has an $i$-good node $A$, then
\begin{enumerate}
\item[(iii)] 
the multiplicity of $M(\la_A)$ in $e_i M(\la)$ is 
$\eps_i(\la)$, $\eps_i(\la_A)=\eps_i(\la)-1$, and 
$\eps_i(\mu)<\eps_i(\la)-1$ for all other composition factors 
$M(\mu)$ of $e_i M(\la)$;
\item[(iv)] 
$\End_{\fH_{n-1}}(e_i M(\la))
\simeq \End_{\fH_{n-1}}(M(\la_A))^{\oplus \eps_i(\la)}$
as a vector superspace;
\item[(v)]
$\hom_{\fH_{n-1}}(e_i M(\la), e_i M(\mu)) = 0$
for all $\mu \in \mathscr{RP}_h(n)$ with $\mu \neq \la$;
\item[(vi)] 
$e_i M(\la)$ is irreducible if and only if $\eps_i(\la)=1$. 
\end{enumerate}
Hence,
$\res^{\fH_n}_{\fH_{n-1}} M(\la)$ 
is completely reducible if and only if $\eps_i(\la)\leq 1$ for every $i\in I$.
\end{Theorem}

\begin{proof}
The existence of such modules $e_i M(\lambda)$
follows from
(\ref{indresdecomp}), Lemma~\ref{relate1} and Theorem~\ref{crystallize}(i),
combined as usual with Kang's Theorem~\ref{TKang}.
Uniqueness follows from Krull-Schmidt and the block classification
from Theorem~\ref{TBr}(iii).
For the remaining properties,
(iii),(iv) and (v) follow from  Theorem~\ref{T280900}
and Corollary~\ref{epr}.
Finally, (vi) follows from (iii) as $e_i M(\la)$ is a module with
simple socle and cosocle both isomorphic to $M(\la_A)$.
\end{proof}

\begin{Theorem}\label{Tf}
Let $\la \in {\mathscr{RP}}_h(n)$.
There exist
$\fH_{n+1}$-modules $f_i M(\la)$ for each $i \in I$,
unique up to isomorphism, such that
\begin{enumerate}
\item[(i)]
$\displaystyle
\ind_{\fH_{n}}^{\fH_{n+1}} M(\la)
\cong
\left\{
\begin{array}{ll}
\displaystyle
2 f_0 M(\la) \oplus 2 f_1 M(\la) \oplus \dots \oplus 2 f_\ell M(\la)
&\hbox{if $b(\la)$ is odd,}\\
f_0 M(\la) \oplus 2 f_1 M(\la) \oplus \dots \oplus 2 f_\ell M(\la)
&\hbox{if $b(\la)$ is even;}\\
\end{array}
\right.$
\item[(ii)] 
for each $i \in I$,
$f_i M(\la)\neq 0$ if and only if $\la$ has an $i$-cogood node $B$,
in which case $f_i M(\la)$ is a self-dual indecomposable 
module with irreducible socle and cosocle isomorphic to $M(\la^B)$.
\end{enumerate}
Moreover, if $i \in I$ and $\la$ has an $i$-cogood node $B$, then
\begin{enumerate}
\item[(iii)] 
the multiplicity of $M(\la^B)$ in $f_i M(\la)$ is 
$\phi_i(\la)$, $\phi_i(\la^B)=\phi_i(\la)-1$, and 
$\phi_i(\mu)<\phi_i(\la)-1$ for all other composition factors 
$M(\mu)$ of $f_i M(\la)$;
\item[(iv)] 
$\End_{\fH_{n+1}}(f_i M(\la))
\simeq \End_{\fH_{n+1}}(M(\la^B))^{\oplus \phi_i(\la)}$
as a vector superspace;
\item[(v)]
$\hom_{\fH_{n+1}}(f_i M(\la), f_i M(\mu)) = 0$
for all $\mu \in \mathscr{RP}_h(n)$ with $\mu \neq \la$;
\item[(vi)] 
$f_i M(\la)$ is irreducible if and only if $\phi_i(\la)=1$. 
\end{enumerate}
Hence,
$\ind^{\fH_{n+1}}_{\fH_{n}} M(\la)$ 
is completely reducible if and only if $\phi_i(\la)\leq 1$ for every $i\in I$.
\end{Theorem}

\begin{proof}
The argument is the 
same as Theorem~\ref{Te}, but using
(\ref{indresdecomp}), Lemma~\ref{relate2}, Theorem~\ref{crystallize}(ii),
Corollary~\ref{epr2}
and Theorem~\ref{weredone}.
\end{proof}

There is one $\fH_n$-module that deserves special mention,
the so-called {\em basic spin module}.
Recall from \ref{supalgs}
that the subalgebra $\cH_n$ of $\fH_n$ generated by
$T_1,\dots,T_{n-1}$ is the classical Hecke algebra associated
to the symmetric group. It has a one dimensional module
denoted $\bid$, on which each 
$T_i$ acts as multiplication by $q$.
For $n \geq 1$, we define
\begin{equation}\label{inn}
I(n) := \ind_{\cH_n}^{\fH_n} {\bid},
\end{equation}
giving an $\fH_n$-module of dimension $2^n$.
Also introduce the restricted $h$-strict partition
\begin{equation}\label{on}
\omega_n := \left\{
\begin{array}{ll}
(h^a, b)&\hbox{if $b \neq 0$,}\\
(h^{a-1}, h-1,1)&\hbox{if $b = 0$,}
\end{array}\right.
\end{equation}
where $n = ah + b$ with $0 \leq b < h$.

\begin{Lemma}\label{bs}
If $p \nmid n$ then 
$I(n) \cong M(\omega_n)$;
if $p | n$ then $I(n)$
is an indecomposable module with two composition factors
both isomorphic to $M(\omega_n)$.
In particular,
$$
\dim M(\omega_n) = \left\{
\begin{array}{ll}
2^n&\hbox{if $p \nmid n$,}\\
2^{n-1}&\hbox{if $p \mid n$.}
\end{array}\right.
$$
\end{Lemma}

\begin{proof}
This is obvious if $n=1,2$ and easy to check directly
if $n = 3$.
Now for $n > 3$ we proceed by induction using
Theorem~\ref{Te} together
with the observation that
$\res^{\fH_n}_{\fH_{n-1}} I(n) \simeq I(n-1) \oplus \Pi I(n-1)$.
We consider the four cases
$n \equiv 0,1$ or $2 \pmod{h}$ and $n \not\equiv 0,1,2\pmod{h}$
separately.

Suppose first that $n\not\equiv 0,1,2\pmod{h}$.
Considering the crystal graph
shows that $\tilde f_i \omega_{n-1} \neq 0$
only for $i = 0$ and for one other $i \in I$, for which
$\tilde f_i \omega_{n-1} = \omega_n$.
By the induction hypothesis, $\res^{\fH_n}_{\fH_{n-1}} I(n)
\cong 2 M(\omega_{n-1})$.
Hence by Theorem~\ref{Te}, $I(n)$ can only contain
$M(\omega_n)$ and $M(\tilde f_0 \omega_{n-1})$ as composition factors.
But the latter case cannot hold since by Theorem~\ref{Te} again,
$\res^{\fH_n}_{\fH_{n-1}} M(\tilde f_0 \omega_{n-1})$
is not isotypic.
Hence all composition factors of 
$I(n)$ are $\cong M(\omega_n)$,
and one easily gets that in fact 
$I(n)\cong M(\omega_n)$ by a dimension argument.

Next suppose that $n \equiv 0 \pmod{h}$.
This time, $\tilde f_0 \omega_{n-1} = \omega_n$
and all other $\tilde f_i \omega_{n-1}$ are zero.
Hence, by the induction hypothesis
and the branching rules,
$I(n)$ only involves $M(\omega_n)$ as a constituent.
But we have that $\res^{\fH_n}_{\fH_{n-1}} M(\omega_n)
=  e_0 M(\omega_n) \cong M(\omega_{n-1})$
so in fact that $I(n)$ must have $M(\omega_n)$ as a constituent
with multiplicty two.
Further consideration of the endomorphism ring of $I(n)$
shows moreover that it is an indecomposable module.

The argument in the remaining two cases
$n \equiv 1 \pmod{h}$ and $n \equiv 2 \pmod{h}$ is entirely similar.
\end{proof}

\Point{\boldmath Projective representations of $S_n$}
We specialize for the final applications to the degenerate case.
So now $F$ is an algebraically closed field of characteristic $p \neq 2$,
$h = p$ (if $p \neq 0$) or $h = \infty$ (if $p = 0$), and
$\ell = (h-1)/2$.
We are interested in the projective representations of the symmetric
group $S_n$ over the field $F$.
Equivalently, see for example \cite[$\S$3]{BK},
we consider the representations of the
{\em twisted group algebra} $S(n)$,
defined by generators $t_1,\dots,t_{n-1}$ subject to the relations
$$
t_i^2 = 1,\qquad
t_i t_{i+1} t_i = t_{i+1}t_i t_{i+1},
\qquad t_i t_j = -t_j t_i
$$
for all $1 \leq i \leq n-1$ and all $1 \leq j \leq n-1$ with $|i-j| > 1$.
We view $S(n)$ as a {superalgebra},
defining the grading by declaring the generators $t_1,\dots,t_{n-1}$
to be of degree $\1$. All modules in this subsection
will be $\Z_2$-graded as usual,
see \ref{theory}.

Let $W(n)$ denote the finite Sergeev superalgebra, 
replacing the notation $\fH_n$ used previously.
Recall from \ref{degmod} that
$W(n)$ is a twisted tensor product of the 
group algebra $F S_n$ of the symmetric group
and the Clifford superalgebra $C(n)$.
So $W(n)$ has even generators $s_1,\dots,s_{n-1}$ subject to the
usual relations of the basic transpositions in $S_n$,
odd generators $c_1,\dots,c_n$ subject to the Clifford relations
as in (\ref{crel1}), (\ref{crel2}), and the additional relations (\ref{abcd1}).
The connection between the superalgebras $S(n)$ and $W(n)$ 
is explained by the following 
observation, due originally to Sergeev
\cite{sergeev2}, see also \cite[3.3]{BK}:

\vspace{1mm}
\begin{Lemma}\label{TW}
There is an isomorphism of superalgebras 
$$
\varphi: S(n)\otimes C(n)\stackrel{\sim}{\longrightarrow} W(n)
$$
such that $\varphi(1 \otimes c_i)=c_i$ and 
$\varphi(t_j\otimes 1)=\frac{1}{\sqrt{-2}} s_j (c_j - c_{j+1})$, 
$i=1,\dots,n$, $j=1,\dots,n-1$. 
\end{Lemma}

We will from now on {\em identify} $W(n)$ with $S(n) \otimes C(n)$
according to the lemma, and view
$S(n)$ (resp. $C(n)$) as the subalgebra $S(n) \otimes 1$
(resp. $1 \otimes C(n)$) of $W(n)$.
Recall the antiautomorphism $\tau$ of $W(n)$ defined
as in (\ref{tttdef}) to be the unique antiautomorphism
which is the identity on the generators $s_1,\dots,s_{n-1}, c_1,\dots,c_n$.
This leaves invariant the subalgebra $S(n)$ of $W(n)$, so
induces an antiautomorphism
\begin{equation}\label{ttdef}
\tau:S(n) \rightarrow S(n).
\end{equation}
More explicitly, $\tau$ on $S(n)$ is the unique antiautomorphism
which maps the generator $t_i$ to $-t_i$ for each $i = 1,\dots,n-1$.
Given a finite dimensional $S(n)$- (resp. $W(n)$-)  module $M$
we write $M^\tau$ for the dual vector superspace viewed as a module
by twisting the natural right action into a left action via the
antiautomorphism $\tau$, see \ref{theory}.

We also need the {\em Clifford module} $U(n)$, which is
the unique irreducible $C(n)$-module up to isomorphism, 
see for instance \cite[2.10]{BK}.
We recall that $U(n)$ is of type $\Mtype$ if $n$ is even, 
type $\Qtype$ if $n$ is odd, and
$\dim U(n) = 2^{\lfloor (n+1)/2 \rfloor}$.

Now consider the exact functors
\begin{align*}
{\mathcal F}_n&:\rep{S(n)} \rightarrow \rep{W(n)},
\qquad
{\mathcal F}_n := ? \boxtimes U(n),\\
{\mathcal G}_n&:\rep{W(n)} \rightarrow \rep{S(n)},
\qquad
{\mathcal G}_n := \hom_{C(n)}(U(n), ?).
\end{align*}
Also let
\begin{align*}
\res^{W(n)}_{W(n-1)}:\rep W(n) \rightarrow \rep W(n-1),
\quad
&
\ind^{W(n)}_{W(n-1)}:\rep W(n-1) \rightarrow \rep W(n),\\
\res^{S(n)}_{S(n-1)}:\rep S(n) \rightarrow \rep S(n-1),
\quad
&
\ind^{S(n)}_{S(n-1)}:\rep S(n-1) \rightarrow \rep S(n)
\end{align*}
denote the (exact) induction and restriction functors, where 
$S(n-1) \subset S(n)$ and $W(n-1) \subset W(n)$ 
are the natural subalgebras generated by all but the last generators.
Recalling that $\Pi$ denotes the parity change functor
(\ref{pcf}),
the following lemma lists some basic properties:

\vspace{1mm}
\begin{Lemma}\label{ds}
The functors ${\mathcal F}_n$ and ${\mathcal G}_n$ are
left and right
adjoint to one another,
and both commute with $\tau$-duality.
Moreover:
\begin{enumerate}
\item[(i)] Suppose that $n$ is even.
Then ${\mathcal F}_n$ and ${\mathcal G}_n$ are inverse
equivalences of categories, 
so induce a 
type-preserving bijection between the isomorphism classes of 
irreducible $S(n)$-modules and of irreducible $W(n)$-modules.
Also,
\begin{align}\label{x1}
{\mathcal F}_{n-1}\circ\res^{S(n)}_{S(n-1)} & \simeq
\res^{W(n)}_{W(n-1)}\circ{\mathcal F}_n,\\\label{x2}
{\mathcal G}_{n-1}\circ\res^{W(n)}_{W(n-1)} & \simeq
\res^{S(n)}_{S(n-1)}\circ{\mathcal G}_n \oplus \Pi\circ\res^{S(n)}_{S(n-1)}\circ{\mathcal G}_n,\\
{\mathcal F}_{n+1}\circ\ind^{S(n+1)}_{S(n)} & \simeq
\ind^{W(n+1)}_{W(n)}\circ{\mathcal F}_n,\label{x3}\\
{\mathcal G}_{n+1}\circ\ind^{W(n+1)}_{W(n)}
 & \simeq
\ind^{S(n+1)}_{S(n)}\circ{\mathcal G}_n\oplus\Pi
\circ
\ind^{S(n+1)}_{S(n)}\circ{\mathcal G}_n.\label{x4}
\end{align}
\item[(ii)]
Suppose that $n$ is odd.
Then ${\mathcal F}_n \circ {\mathcal G}_n \simeq \Id \oplus \Pi$
and ${\mathcal G}_n \circ {\mathcal F}_n \simeq \Id \oplus \Pi$.
Furthermore, the functor $\mathcal{F}_n$ induces a bijection
between isomorphism classes of
irreducible $S(n)$-modules of type $\Mtype$
and irreducible
$W(n)$-modules of type $\Qtype$, while
the functor $\mathcal{G}_n$ induces a bijection
between isomorphism classes of
irreducible $W(n)$-modules of type $\Mtype$ 
and irreducible $S(n)$-modules of type $\Qtype$.
Finally,
\begin{align}
\res^{W(n)}_{W(n-1)}\circ \mathcal{F}_n
&\simeq\mathcal{F}_{n-1} \circ \res^{S(n)}_{S(n-1)}
\oplus\Pi\circ\mathcal{F}_{n-1} \circ \res^{S(n)}_{S(n-1)},\label{x5}\\
\label{x6}\res^{S(n)}_{S(n-1)}\circ{\mathcal G}_n & \simeq
{\mathcal G}_{n-1}\circ\res^{W(n)}_{W(n-1)},\\
\ind^{W(n+1)}_{W(n)}\circ{\mathcal F}_n & \simeq
{\mathcal F}_{n+1}\circ\ind^{S(n+1)}_{S(n)} \oplus \Pi\circ{\mathcal F}_{n+1}\circ\ind^{S(n+1)}_{S(n)},\label{x7}\\
\ind^{S(n+1)}_{S(n)}\circ{\mathcal G}_n & \simeq
{\mathcal G}_{n+1}\circ\ind^{W(n+1)}_{W(n)}.\label{x8}
\end{align} 
\end{enumerate}
\end{Lemma}

\begin{proof}
Most of these facts are proved in \cite[3.4,3.5]{BK}, but we recall some of the
details since we need to go slightly further.
Let us consider the proof that
${\mathcal F}_n \circ {\mathcal G}_n \simeq \operatorname{Id} \oplus \Pi$
and ${\mathcal G}_n \circ {\mathcal F}_n \simeq \operatorname{Id} \oplus \Pi$,
assuming that $n$ is odd.
Let $I, J$ be a basis for $\End_{C(n)}(U(n))$ with $I$ being the identity
and $J$ being an odd involution.
Then, there are natural isomorphisms
$$
\eta:{\mathcal F}_n \circ {\mathcal G}_n \stackrel{\sim}{\longrightarrow} 
\Id \oplus \Pi,
\qquad
\xi:\Id \oplus \Pi
\stackrel{\sim}{\longrightarrow} {\mathcal G}_n \circ {\mathcal F}_n.
$$
The first is defined for each $W(n)$-module $M$ by
$\eta_M:\hom_{C(n)}(U(n), M) \boxtimes U(n)
\rightarrow M \oplus \Pi M$,
$\theta \otimes u \mapsto (\theta(u),(-1)^{\bar \theta} \theta(Ju))$.
The second is defined for each $S(n)$-module $N$
by $\xi_N:N \oplus \Pi N \rightarrow
\hom_{C(n)}(U(n), N \boxtimes U(n))$,
$(n, n') \mapsto \theta_{n,n'}$, where
$\theta_{n,n'}
(u) = n \otimes u +(-1)^{\bar{n}'} n' \otimes J u$
for each $u \in U(n)$.
The proof that $\eta$ and $\xi$ really are isomorphisms
is similar to the argument in \cite[3.4]{BK}.
Now consider the composite natural transformations
\begin{align}
\Id \longrightarrow \Id \oplus \Pi \stackrel{\xi}{\longrightarrow}
\mathcal{G}_n \circ \mathcal{F}_n,\qquad&
\mathcal{F}_n \circ \mathcal{G}_n \stackrel{\eta}{\longrightarrow} 
\Id \oplus \Pi \longrightarrow \Id,\label{fir}\\
\Id \longrightarrow \Id \oplus \Pi \stackrel{\eta^{-1}}{\longrightarrow}
\mathcal{F}_n \circ \mathcal{G}_n,\qquad&
\mathcal{G}_n \circ \mathcal{F}_n \stackrel{\xi^{-1}}{\longrightarrow} 
\Id \oplus \Pi \longrightarrow \Id,\label{secs}
\end{align}
where the unmarked arrows are the obvious even ones.
Then, (\ref{fir}) (resp. \ref{secs})
gives the unit and counit of the adjunction 
needed to prove 
that $\mathcal{F}_n$ is left (resp. right) adoint to $\mathcal{G}_n$.
We leave the details to be checked to the reader,
see e.g. \cite[IV.1, Theorem 2(v)]{Macl}.

In case $n$ is even, a similar but easier argument shows that
${\mathcal F}_n \circ {\mathcal G}_n \simeq \Id,
{\mathcal G}_n \circ {\mathcal F}_n \simeq \Id$ so that
${\mathcal F}_n$ and ${\mathcal G}_n$ are inverse equivalences,
hence left and right adjoint to each other.
Now the statements in (i) and (ii) about isomorphism classes of
irreducible modules follow easily as in \cite[3.5]{BK}.

Let us next prove that ${\mathcal F}_n$ commutes with duality. 
The antiautomorphism $\tau$ of $W(n)$ induces the antiautomorphism
$\tau$ of the subalgebra $C(n)$ with $\tau(c_i) = c_i$ for each
$i = 1,\dots,n$.
Let $\phi:U(n)\rightarrow U(n)^\tau, u \mapsto \phi_u$ 
be a homogeneous isomorphism;
note that it is not always possible to choose $\phi$ to be even,
since $U(n)\not\simeq U(n)^\tau$ in case $n \equiv 2 \pmod{4}$.
We get a natural isomorphism
$$
\Phi_M: M^\tau \boxtimes U(n) \rightarrow (M \boxtimes U(n))^\tau,
\quad
f \otimes u \mapsto \theta_{f,u}
$$
for each finite dimensional $S(n)$-module $M$,
where
$\theta_{f,u}
(m \otimes v) =(-1)^{\bar m \bar \phi}
f(m)\phi_u(v)$ for each $m \in M, v \in U(n)$.
This shows that ${\mathcal F}_n$ commutes with duality,
i.e. $\tau \circ {\mathcal F}_n \circ \tau \cong {\mathcal F}_n$.
Hence, using (\ref{switch}) and the fact that
$\mathcal{G}_n$ is left adjoint to ${\mathcal F}_n$,
the composite functor
$\tau \circ {\mathcal G}_n \circ \tau$ is right adjoint to ${\mathcal F}_n$.
But we already know that $\mathcal{G}_n$ is right adjoint to 
${\mathcal F}_n$, so uniqueness of adjoints gives that
$\tau \circ {\mathcal G}_n \circ \tau \cong {\mathcal G}_n$.
This shows that ${\mathcal G}_n$ commutes with duality too.

It just remains to check the isomorphisms
(\ref{x1})--(\ref{x8}).
Well, (\ref{x1}) follows from the definition on noting that
$\res^{C(n)}_{C(n-1)} U(n) \simeq U(n-1)$ if $n$ is even.
Then (\ref{x2}) follows from (\ref{x1}) on composing on the left with
$\mathcal{G}_{n-1}$ and on the right with $\mathcal{G}_n$.
Next, (\ref{x6}) follows from the definition and an application of
Frobenius reciprocity, using the observation that
$U(n) \simeq \ind_{C(n-1)}^{C(n)} U(n-1)$ if $n$ is odd.
As before (\ref{x5}) then follows, composing with $\mathcal{F}_{n-1}$
and $\mathcal{F}_n$.
Finally, 
(\ref{x3}), (\ref{x4}), (\ref{x7}) and (\ref{x8}) follow from
(\ref{x6}), (\ref{x5}), (\ref{x2}) and (\ref{x1}) respectively
by uniqueness of adjoints.
\end{proof}

We are ready to derive the consequences for $S(n)$ of
the results of \ref{klk}, or rather, of the analogous results for $W(n)$ in the 
degenerate case.
By Theorem~\ref{TBr}, we have a parametrization
$\{M(\la)\:|\:\la \in \mathscr{RP}_p(n)\}$
of the irreducible $W(n)$-modules.
Lemma~\ref{ds} shows that the functors $\mathcal{F}_n$ and
$\mathcal{G}_n$ set up a natural correspondence between
classes of irreducible $S(n)$ and $W(n)$-modules, type-preserving
if $n$ is even and type-reversing if $n$ is odd.
Hence we have a parametrization $\{D(\la)\:|\:\la \in \mathscr{RP}_p(n)\}$
of the irreducible $S(n)$-modules, letting $D(\la)$ be an irreducible
$S(n)$-module corresponding to $M(\la)$ under the correspondence.
Also, recalling the definition (\ref{bbdef}), define
\begin{equation}
a(\la) := n - b(\la)
\end{equation}
for $\la \in \mathscr{RP}_p(n)$.
We observe by (\ref{resds}) that
\begin{equation}
a(\la) \equiv \gamma_1+\dots+\gamma_\ell \pmod{2},
\end{equation}
where $\gamma_1+\dots+\gamma_\ell$ counts the number of nodes in
the Young diagram $\la$ of residue {\em different from} $0$.
Then:

\vspace{1mm}
\begin{Theorem}
\label{SBr} 
The modules
$\{D(\lambda)\:|\:\lambda \in \mathscr{RP}_p(n)\}$
form a complete set of pairwise non-isomorphic irreducible 
$S(n)$-modules.
Moreover, for $\la,\mu \in \mathscr{RP}_p(n)$,
\begin{enumerate}
\item[(i)] $D(\la) \cong D(\la)^\tau$;

\item[(ii)] $D(\la)$ is of type $\Mtype$ if $a(\la)$ is even,
type $\Qtype$ if $a(\la)$ is odd;

\item[(iii)] $D(\mu)$ and $D(\la)$ belong to the same block if and only
if $\cont(\mu) = \cont(\la)$;

\item[(iv)] $D(\la)$ is a projective module if and only if $\la$
is a $p$-bar core.
\end{enumerate}
\end{Theorem}

\begin{proof}
It just remains to observe that (i)--(iv) follow
follow directly from Theorem~\ref{TBr}(i)--(iv) using
Lemma~\ref{ds}.
\end{proof}

\begin{Remark}{\rm
The $p$-blocks of the ordinary irreducible projective
representations of $S_n$ were described by Humphreys \cite{H}, in
terms of the notion of $p$-bar core.
However, unlike the case of $S_n$, Humphreys' result does not imply
Theorem~\ref{SBr}(iii) 
because of the lack of information on decomposition numbers.
}
\end{Remark}

Define $\omega_n \in \mathscr{RP}_p(n)$ as in
(\ref{on}). Then the irreducible $S(n)$-module
$D(\omega_n)$ is the {\em basic spin module}:

\vspace{1mm}
\begin{Lemma}
$D(\omega_n)$ 
is of
dimension $2^{\lfloor n/2\rfloor}$,
unless $p | n$ when its dimension is
$2^{\lfloor (n-1)/2\rfloor}$.
Moreover, $D(\omega_n)$ is equal to the reduction modulo $p$
of the basic spin module $D((n))_\C$ 
of $S(n)_{\C}$ over $\C$,
except if $p|n$ and $n$ is even
when the reduction modulo $p$ of $D((n))_\C$ has two composition factors
both isomorphic to $D(\omega_n)$.
\end{Lemma}

\begin{proof}
The statement about dimension is immediate from
Lemmas~\ref{bs} and \ref{ds}.
The final statement is easily proved
by working in terms of $W(n)$ and using the explicit construction 
given in (\ref{inn}).
\end{proof}

To motivate the next two theorems, note that
the map $[D(\la)] \mapsto [M(\la)]$ for each 
$\la \in \mathscr{RP}_p(n)$ 
extends linearly to an isomorphism
$$
K(\rep S(n)) \stackrel{\sim}{\longrightarrow} K(\rep W(n))
$$
at the level of Grothendieck groups.
Using this identification, we can lift the operators $e_i$ and $f_i$
on $K(\Lambda_0) = \bigoplus_{n \geq 0} K(\rep W(n))$ defined earlier
to define similar operators on
$\bigoplus_{n \geq 0} K(\rep S(n))$.
Then all our earlier results about $K(\Lambda_0)$, for instance 
Theorems~\ref{chevrel} and \ref{thmb},  
could be restated purely in terms of the representations of $S(n)$
instead of $W(n)$.
In fact, we can do slightly better and define the operators
$e_i$ and $f_i$ on irreducible $S(n)$-modules, 
not just on the Grothendieck group.
Theorems~\ref{Se} and \ref{Sf} below should be compared with
the parallel results 
\cite[Theorems E and E$'$]{BK1} for the symmetric group.

\vspace{1mm}
\begin{Theorem}\label{Se}
Let $\la \in {\mathscr{RP}}_p(n)$.
There exist
$S(n-1)$-modules $e_i D(\la)$ for each $i \in I$,
unique up to isomorphism, such that
\begin{enumerate}
\item[(i)]
$\displaystyle
\res_{S(n-1)}^{S(n)} D(\la)
\cong
\left\{
\begin{array}{ll}
\displaystyle
e_0 D(\la) \oplus 2 e_1 D(\la) \oplus \dots \oplus 2 e_\ell D(\la)
&\hbox{if $a(\la)$ is odd,}\\
e_0 D(\la) \oplus  e_1 D(\la) \oplus \dots \oplus  e_\ell D(\la)
&\hbox{if $a(\la)$ is even;}\\
\end{array}
\right.$
\item[(ii)] 
for each $i \in I$,
$e_i D(\la)\neq 0$ if and only if $\la$ has an $i$-good node $A$,
in which case $e_i D(\la)$ is a self-dual indecomposable 
module with irreducible socle and cosocle isomorphic to $D(\la_A)$.
\end{enumerate}
Moreover, if $i \in I$ and $\la$ has an $i$-good node $A$, then
\begin{enumerate}
\item[(iii)] 
the multiplicity of $D(\la_A)$ in $e_i D(\la)$ is 
$\eps_i(\la)$, $\eps_i(\la_A)=\eps_i(\la)-1$, and 
$\eps_i(\mu)<\eps_i(\la)-1$ for all other composition factors 
$D(\mu)$ of $e_i D(\la)$;
\item[(iv)] 
$\End_{S(n-1)}(e_i D(\la))
\simeq \End_{S(n-1)}(D(\la_A))^{\oplus \eps_i(\la)}$
as a vector superspace;
\item[(v)]
$\hom_{S(n-1)}(e_i D(\la), e_i D(\mu)) = 0$
for all $\mu \in \mathscr{RP}_p(n)$ with $\mu \neq \la$;
\item[(vi)] 
$e_i D(\la)$ is irreducible if and only if $\eps_i(\la)=1$. 
\end{enumerate}
Hence,
$\res^{S(n)}_{S(n-1)} D(\la)$ 
is completely reducible if and only if $\eps_i(\la)\leq 1$ for every $i\in I$.
\end{Theorem}

\begin{proof}
If $n$ is odd, we simply define $e_i D(\la) := \mathcal{G}_{n-1} (e_i M(\la))$
for each $i \in I, \la \in \mathscr{RP}_p(n)$.
If $n$ is even, take
$$
e_i D(\la) := \left\{
\begin{array}{ll}
\mathcal{G}_{n-1} (e_i M(\la))&\hbox{if 
$a(\la)$ is even and $i \neq 0$, or
$a(\la)$ is odd and $i = 0$,}\\
\overline{\mathcal{G}}_{n-1} (e_i M(\la))&\hbox{if 
$a(\la)$ is even and $i = 0$, or
$a(\la)$ is odd and $i \neq 0$.}
\end{array}
\right.
$$
We need to explain the notation $\overline{\mathcal{G}}_{n-1}$
used in the last two cases: here,
$e_i M(\la)$ admits an odd involution by
Remark~\ref{feed} and Theorem~\ref{TBr}(ii),
and also $U(n-1)$ has an odd involution since $n$ is even.
So in exactly the same way as in 
the definition of (\ref{hbar}), we can 
introduce the space
$$
\overline{\mathcal{G}}_{n-1} (e_i M(\la))
:= 
\overline{\hom}_{C(n-1)}(U(n-1), e_i M(\la)).
$$
It is then the case that
$\mathcal{G}_{n-1} (e_i M(\la)) \simeq
\overline{\mathcal{G}}_{n-1} (e_i M(\la)) \oplus
\Pi \overline{\mathcal{G}}_{n-1} (e_i M(\la))$.
Equivalently, 
by Lemma~\ref{ds}(ii)
$e_i D(\la)$ can be characterized by
$$
e_i M(\la) \cong \mathcal{F}_{n-1}(e_i D(\la))
$$
if 
$a(\la)$ is even and $i = 0$, or
$a(\la)$ is odd and $i \neq 0$.

With these definitions, it is now a straightforward matter to 
prove (i)--(vi) using Theorem~\ref{Te} and Lemma~\ref{ds}.
Finally, the uniqueness statement is immediate from 
Krull-Schmidt and the description of blocks from Theorem~\ref{SBr}(iii).
\end{proof}

\begin{Theorem}\label{Sf}
Let $\la \in {\mathscr{RP}}_p(n)$.
There exist
$S(n+1)$-modules $f_i D(\la)$ for each $i \in I$,
unique up to isomorphism, such that
\begin{enumerate}
\item[(i)]
$\displaystyle
\ind_{S(n)}^{S(n+1)} D(\la)
\cong
\left\{
\begin{array}{ll}
\displaystyle
 f_0 D(\la) \oplus 2 f_1 D(\la) \oplus \dots \oplus 2 f_\ell D(\la)
&\hbox{if $a(\la)$ is odd,}\\
f_0 D(\la) \oplus  f_1 D(\la) \oplus \dots \oplus  f_\ell D(\la)
&\hbox{if $a(\la)$ is even;}\\
\end{array}
\right.$
\item[(ii)] 
for each $i \in I$,
$f_i D(\la)\neq 0$ if and only if $\la$ has an $i$-cogood node $B$,
in which case $f_i D(\la)$ is a self-dual indecomposable 
module with irreducible socle and cosocle isomorphic to $D(\la^B)$.
\end{enumerate}
Moreover, if $i \in I$ and $\la$ has an $i$-cogood node $B$, then
\begin{enumerate}
\item[(iii)] 
the multiplicity of $D(\la^B)$ in $f_i D(\la)$ is 
$\phi_i(\la)$, $\phi_i(\la^B)=\phi_i(\la)-1$, and 
$\phi_i(\mu)<\phi_i(\la)-1$ for all other composition factors 
$D(\mu)$ of $f_i D(\la)$;
\item[(iv)] 
$\End_{S(n+1)}(f_i D(\la))
\simeq \End_{S(n+1)}(D(\la^B))^{\oplus \phi_i(\la)}$
as a vector superspace;
\item[(v)]
$\hom_{S(n+1)}(f_i D(\la), f_i D(\mu)) = 0$
for all $\mu \in \mathscr{RP}_p(n)$ with $\mu \neq \la$;
\item[(vi)] 
$f_i D(\la)$ is irreducible if and only if $\phi_i(\la)=1$. 
\end{enumerate}
Hence,
$\ind^{S(n+1)}_{S(n)} M(\la)$ 
is completely reducible if and only if $\phi_i(\la)\leq 1$ for every $i\in I$.
\end{Theorem}

\begin{proof}
This is deduced from Theorem~\ref{Tf} by similar argument to the proof
of Theorem~\ref{Se}.
\end{proof}

\begin{Remark}{\rm
(i) 
In \cite[10.3]{BK}, we gave an entirely different 
construction of the irreducible $S(n)$-modules,
which we also denoted by $D(\la)$ for
$\la \in \mathscr{RP}_p(n)$.
We warn the reader that we have not yet proved that
the modules denoted $D(\la)$ here are isomorphic to those
in \cite{BK}, though we expect this to be the case.

(ii) Over $\C$, the branching rules in the preceeding two theorems
are the same as Morris' branching rules, see \cite{Morris}.
In particular using this observation, one easily shows that
our labelling of irreducibles over $\C$ agrees with the standard
labelling. Hence over $\C$ the labelling here agrees with the labelling
in \cite{BK}, compare \cite{se}.

(iii) 
There is one other case where it is easy to see right away that
the labelling here agrees with \cite{BK}:
if $\la$ is a $p$-bar core then
consideration of central characters shows that
$D(\la)$ is equal to a
reduction modulo $p$ of the irreducible representation
of $S(n)_\C$ over $\C$ with the same label.
So it coincides with
the module $D(\la)$ of \cite{BK} thanks to \cite[10.8]{BK}.
}
\end{Remark}

\Point{The Jantzen-Seitz problem}
The results of the previous subsection give a solution to the
{\em Jantzen-Seitz problem} for projective representations of the
symmetric and alternating groups.
This problem originated in \cite{JS}, and is of interest in the study
of maximal subgroups of the finite classical groups.
To consider the Jantzen-Seitz problem, 
we first need to switch to studying {\em ungraded} representations
of the twisted group algebra $S(n)$.
The goal is to describe all ungraded irreducible $S(n)$-modules
which remain irreducible on restriction to the subalgebra $S(n-1)$.

As in \cite[$\S$10]{BK}, 
it is straightforward to obtain a parametrization of the ungraded
irreducible $S(n)$-modules from Theorem~\ref{SBr}:
if $a(\la)$ is odd then $D(\la)$ decomposes as an ungraded module
as $D(\la) = D(\la, +) \oplus D(\la, -)$ for two non-isomorphic
irreducible $S(n)$-modules $D(\la, +)$ and $D(\la, -)$.
If $a(\la)$ is even, then $D(\la)$ is irreducible viewed as an
ungraded $S(n)$-module, but we denote it instead by $D(\la,0)$
to make it clear that we are no longer considering a $\Z_2$-grading.
Then
$$
\{D(\la,0)\:|\:\la \in \mathscr{RP}_p(n), a(\la)\hbox{ even}\}
\sqcup
\{D(\la, +), D(\la, -)\:|\:\la \in \mathscr{RP}_p(n), a(\la)\hbox{ odd}\}
$$
gives a complete set of pairwise non-isomorphic 
ungraded irreducible $S(n)$-modules.

\vspace{1mm}
\begin{Remark}{\rm
Using Theorem~\ref{SBr} and a counting argument
involving Humphreys' block classification \cite{H},
it is not hard to obtain the following description of
the ungraded blocks of the algebra $S(n)$.
Let $D(\la,\eps)$ and $D(\mu,\delta)$ be ungraded
irreducible $S(n)$-modules. 
Then, with one exception,
$D(\la,\eps)$ and $D(\mu,\delta)$ 
lie in the same block if and only if $\la$ and $\mu$ have the same
$p$-bar core.
The exception is if $\la = \mu$ is a $p$-bar core,
$a(\la)$ is odd and $\eps = - \delta$, when
$D(\la,\eps)$ and $D(\mu,\delta)$ are in different blocks.
}\end{Remark}

\vspace{1mm}
Now we state the solution to
the Jantzen-Seitz problem for projective representations
of the symmetric group.
The proof is a straightforward consequence of Theorem~\ref{Se}.

\vspace{1mm}
\begin{Theorem}
Let $\la \in \mathscr{RP}_p(n)$.
Then:
\begin{enumerate}
\item[(i)] If $a(\la)$ is even, 
$\res^{S(n)}_{S(n-1)}
D(\la,0)$ is irreducible 
if and only if 
$\eps_0(\la) = \sum_{i \in I} \eps_i(\la) = 1$.
\item[(ii)] If $a(\la)$ is odd, 
$\res^{S(n)}_{S(n-1)}
D(\la, \pm)$ 
is irreducible 
if and only if $\sum_{i \in I} \eps_i(\la) = 1$.
\end{enumerate}
\end{Theorem}

Finally let us discuss the analogous problem for the projective representations
of the alternating group.
As explained in \cite[$\S$10]{BK}, it suffices for this to consider the
representation theory of the algebra $A(n) := S(n)_{\0}$, i.e. the
twisted group algebra of the alternating group.
The irreducible $A(n)$-modules (there being no ambiguity between graded
and ungraded modules since $A(n)$ is purely even) are constructed
out of those for $S(n)$
as in \cite[Theorem 10.4]{BK}.
More precisely, if $\la \in \mathscr{RP}_p(n)$ has
$a(\la)$ even, then $\res^{S(n)}_{A(n)} D(\la)$ decomposes as
a direct sum $E(\la, +)\oplus E(\la, -)$ of two non-isomorphic irreducible
$A(n)$-modules.
If $a(\la)$ is odd, then $\res^{S(n)}_{A(n)} D(\la)$ decomposes
as a direct sum of two copies of a single irreducible $A(n)$-module
denoted 
$E(\la,0)$. Then,
$$
\{E(\la,0)\:|\:\la \in \mathscr{RP}_p(n), a(\la)\hbox{ odd}\}
\sqcup
\{E(\la, +), E(\la, -)\:|\:\la \in \mathscr{RP}_p(n), a(\la)\hbox{ even}\}
$$
gives a complete set of pairwise non-isomorphic 
irreducible $A(n)$-modules.
Now the solution to the Jantzen-Seitz problem in this case, again
an easy consequence of Theorem~\ref{Se}, is as follows:

\vspace{1mm}
\begin{Theorem}
Let $\la \in \mathscr{RP}_p(n)$.
Then:
\begin{enumerate}
\item[(i)] If $a(\la)$ is even, 
$\res^{A(n)}_{A(n-1)}
E(\la, \pm)$ 
is irreducible 
if and only if 
$\sum_{i \in I} \eps_i(\la) = 1$;

\item[(ii)] If $a(\la)$ is odd, 
$\res^{A(n)}_{A(n-1)}
E(\la, 0)$ 
is irreducible 
if and only if 
$\eps_0(\la) = \sum_{i \in I} \eps_i(\la) = 1$.
\end{enumerate}
\end{Theorem}

\ifbook@\pagebreak\fi

\end{document}